\documentstyle[12pt]{article}
\input{amssym.def}
\input amssym.tex

\def\newtheorems{\newtheorem{theorem}{Theorem}[section]
\newtheorem{cor}[theorem]{Corollary}
\newtheorem{prop}[theorem]{Proposition}
\newtheorem{lemma}[theorem]{Lemma}

\newtheorem{defn}[theorem]{Definition}

\newtheorem{example}[theorem]{Example}

\newtheorem{question}[theorem]{Question}

}

\newtheorems

\def\fs#1{\mbox{\it #1\kern 1.3pt}}
\def\fss#1{\mbox{\scriptsize\it #1\kern 1.3pt}}
\def\fst#1{\hbox{\tiny$ #1$\kern 1.3pt}}
\newcommand{\ff}[1]{\hbox{\it #1} \kern 0.7pt}
\newcommand{\nff}[2]{\mbox{\it #1}\kern 2.2pt(\kern -0.5pt #2
\kern -0.3pt)}
\newcommand{\nnff}[2]{\mbox{\it #1}\kern 1.5pt(\kern -0.5pt #2
\kern -0.3pt)}
\newcommand{\sff}[1]{\hbox{\scriptsize\it #1} \kern 0.4pt}

\newcommand{\rcomp}{\hbox{$\sim \kern 1.0pt$}}
\newcommand{\comp}{\mathrel{\hbox{\kern2pt$<\kern -3pt
>$\kern2pt}}}
\newcommand{\ncomp}
{\;\hbox{\kern2pt\hbox{/}\kern -10.8pt \hbox{$<\kern -2.5pt
>$\kern2pt}}}
\newcommand{\scomp}{\mathrel{\hbox{\scriptsize$<\kern -2pt >$}}}
\newcommand{\nscomp}{\mathrel{\hbox{\scriptsize$<\kern -2.0pt
>\kern-10.1pt/\kern4.5pt$}}}

\newcommand{\lcparenth}
{\vbox{\hbox{$($ \kern-9pt \raise1.5pt \hbox{\tiny
$\circ$}}}\kern1.5pt}

\newcommand{\rcparenth}
{\vbox{\hbox{$)$ \kern-7.5pt \raise1.1pt \hbox{\tiny $\circ$}}}}

\newcommand{\cmpl}
{\raise 1.5pt \hbox{\tiny $\sim$\kern 0.5pt}}

\newcommand{\meet}
{\hbox{$\wedge \kern -5.75pt \raise 1.5pt \hbox{$.$}\,$}}
\newcommand{\Meet}
{\hbox{$\bigwedge \kern -8pt \raise 0.75pt \hbox{$.$}\:$}}
\newcommand{\ld}
{\hbox{$< \kern -6pt \raise 2pt \hbox{$.$}\,$}}
\newcommand{\sss}{\: \hbox{$
\underline{\hbox{$\subset$}}\kern -4pt\raise -2pt \hbox{$\tiny |$}
$}\: }
\newcommand{\almostcontained}{\: \hbox{$
\raise 1.5pt \hbox{\scriptsize $\subset$}\kern -6.3pt\raise -3.5pt
\hbox{\scriptsize $\sim$}
$}\: }
\newcommand{\rraro}[2]{\hbox{$\kern 3pt\raise 2pt
\hbox{$\raro$}
\kern -14pt \raise
-3.5pt\hbox{\tiny{$#1\raro #2$}}$}}

\newcommand{\ct}{\centerline}

\newcommand{\frc}{\hbox{$\parallel \kern -5.7pt \hbox{$-$}$}}

\newcommand{\nfrc}{\not \kern -5pt \frc}
\newcommand{\rest}{\vbox{\hbox{$\:\kern -
2pt\mathbin{\vert\kern-3.1pt\lower-1pt
\hbox{$\mathsurround=0pt\mathchar"0012$}\kern-4pt}\:$}}}

\newcommand{\rests}{\vbox{\hbox{\scriptsize$\:\kern
-1.4pt\mathbin{\vert\kern-2.4pt\lower-1pt
\hbox{\scriptsize$\mathsurround=0pt\mathchar"0012$}\kern-
3.6pt}\:$}}}

\newcommand{\cat}{\kern 1.8pt{\mathaccent 94 {\,\ }}\kern 0.7pt}

\newcommand{\cntd}{\subseteq}

\newcommand{\pcntda}{\lower5pt\hbox{$\stackrel{\subset}{\neq}$
}}
\newcommand{\pcntdb}{\lower5pt\hbox{$\stackrel{\supset}{\neq}$
}}
\newcommand{\pcntdc}{\lower5pt\hbox{$\stackrel{\subseteq}{\neq
}$}}
\newcommand{\wreath}{{\bf W\kern -2pt r}\kern 3pt}

\newcommand{\m}[1]{\hbox{$ #1 $}}

\newcommand{\gf}{\varphi}

\newcommand{\CC}{{\cal C}}

\newcommand{\CI}{{\cal I}}

\newcommand{\CL}{{\cal L}}

\newcommand{\CO}{{\cal O}}

\newcommand{\CR}{{\cal R}}
\newcommand{\CS}{{\cal S}}
\newcommand{\CT}{{\cal T}}

\newcommand{\smallr}{\mbox{\scriptsize rel}}

\newcommand{\barg}{\bar{g}}

\newcommand{\barB}{\bar{B}}

\newcommand{\barL}{\bar{L}}

\newcommand{\barN}{\bar{N}}

\newcommand{\veca}{\vec{a}}
\newcommand{\vecb}{\vec{b}}

\newcommand{\vecg}{\vec{g}}

\newcommand{\vecu}{\vec{u}}
\newcommand{\vecU}{\vec{U}}
\newcommand{\vecv}{\vec{v}}
\newcommand{\vecV}{\vec{V}}
\newcommand{\vecw}{\vec{w}}
\newcommand{\vecW}{\vec{W}}
\newcommand{\lvec}{\hbox{\raise 3.3 pt \hbox{\raise 0.584pt
\hbox{$\bar{\
}$}\hbox{\kern -1pt $\vec{\ }$}}}}
\newcommand{\lvecW}{W \kern -11pt \lvec \kern 5pt}
\newcommand{\vecx}{\vec{x}}
\newcommand{\vecy}{\vec{y}}
\newcommand{\vecz}{\vec{z}}

\newcommand{\tldgt}{\tilde{\tau}}

\newcommand{\tldF}{\tilde{F}}

\newcommand{\tldR}{\tilde{R}}

\newcommand{\raro}{\rightarrow}

\newcommand{\itmb}[1]{\item[[\kern 1.5pt #1\kern -4pt]]}
\newcommand{\itms}[1]{\item[[#1\kern -5pt]]}

\newcommand{\II}{{\bf I\kern -1pt I}}

\newcommand{\singcolb}[2]{\left(\begin{array}{c}#1\\#2\end{array
}\right)}

\newcommand{\doubcol}[4]{\left[\begin{array}{cc}#1&#2\\#3&#4\e
nd{array}\right]}
\newcommand{\doubcolb}[4]{\left(\begin{array}{cc}#1&#2\\#3&#4\
end{array}\right)}

\newcommand{\eqqdf}{\stackrel{def}{\equiv}}
\newcommand{\eqdf}{\stackrel{def}{=}}

\newcommand{\surj}{\vbox{\hbox{$\longrightarrow $
\kern -22pt \hbox{\lower 2.5pt  \hbox{\tiny onto}}
\kern -16pt \hbox{\raise 5pt  \hbox{\tiny 1-1}}
\kern 3pt}}}
\newcommand{\uarrow}[2]{\vbox{\hbox{$\longrightarrow $
\kern -16pt \hbox{\raise 5pt  \hbox{\tiny $#1$}}
\kern 10pt }}}

\newcommand{\simm}{\kern 4pt \sim \kern 4pt}
\newcommand{\spri}{\vbox{\hbox{\raise 2pt \hbox{\tiny $\|$}}}}
\newcommand{\tlr}{\vbox{\hbox{\raise 2pt \hbox{\tiny
$\leftarrow$}}}}
\newcommand{\trr}{\vbox{\hbox{\raise 2pt \hbox{\tiny
$\rightarrow$}}}}
\newcommand{\spr}{\mathrel{\vbox{\hbox{\tlr \kern -1.7pt \spri
\kern -1.7pt
\trr}}}}
\newcommand{\iseg}{\mathrel{\vbox{\hbox{$< \kern-5.4pt \lower
2.5pt\hbox{\tiny $|$} \kern -2.3pt \raise 5.0pt\hbox{\tiny
$|$}$}}\kern
0.5pt}}
\newcommand{\isegeq}{\mathrel{\vbox{\hbox{$< \kern-5.4pt
\lower
2.5pt\hbox{\tiny $|$} \kern -2.3pt \raise 5.0pt\hbox{\tiny
$|$}$}}\kern
3.5pt}}

\newcommand{\overrelation}{ \kern1pt / \kern-2.5pt}

\newcommand{\Nm}{\hbox{\kern -1.3pt \em I\kern-.2300em N\,}}
\newcommand{\Na}{\hbox{\kern -1.3pt \it I\kern-.2300em N\,}}
\newcommand{\Aa}{\hbox{\it A\kern -6.8pt\lower 1.0pt\hbox{-
}\,}}
\newcommand{\Ba}{\hbox{\kern -1.3pt \it I\kern-.2300em B\,}}
\newcommand{\Da}{\hbox{\kern -1.3pt \it I\kern-.2300em D\,}}
\newcommand{\Ka}{\hbox{\kern -1.3pt \it I\kern-.2300em K\,}}
\newcommand{\La}{\hbox{\kern -1.3pt \it I\kern-.2300em L\,}}
\newcommand{\Ta}{\hbox{\kern 0.5pt \it T\kern-.5550em T\,}}
\newcommand{\Td}
{\hbox{\kern 0.5pt \it T \kern-.9300em {\small\it I} \kern -.590em
{\small\it I}\kern 2.4pt}}
\newcommand{\jTa}{\hbox{\kern 0.5pt \it T\kern-.5550em T\,}}
\newcommand{\nTa}{\hbox{\kern 0.5pt \it T\kern-.6300em T\,}}
\newcommand{\Nmt}{\hbox{\kern -1.3pt I\kern-.2300em N\,}}
\newcommand{\N}{\hbox{I\kern-.1500em \hbox{\sf N}}}
\newcommand{\As}{\hbox{\scriptsize\it A\kern -5.3pt\lower
0.7pt\hbox{\it -}\,}}
\newcommand{\Ns}{\hbox{{\scriptsize\it I\kern-.1500em N}}}
\newcommand{\Ls}{\hbox{{\scriptsize\it I\kern-.1700em L}}}
\newcommand{\Rm}{\hbox{\kern -1.3pt \em I\kern-.1950em R\,}}
\newcommand{\Ra}{\hbox{\kern -1.3pt \it I\kern-.1950em R\,}}

\newcommand{\qed}{\kern 5pt\vrule height8pt width6.5pt
depth2pt}
\newcommand{\text}[1]{\mbox{#1}}

\newcommand{\bd}{\begin{description}}
\newcommand{\ed}{\end{description}}

\newcommand{\ben}{\begin{enumerate}}
\newcommand{\een}{\end{enumerate}}

\newcommand{\sngltn}[1]{\{ #1 \}}
\newcommand{\dbltn}[2]{\{ #1, #2 \}}

\newcommand{\setm}[2]{\{#1\:|\:#2\}}

\newcommand{\fsetn}[2]{\{\,#1,\ldots ,#2\}}

\newcommand{\sing}[1]{\langle\  #1\ \rangle}
\newcommand{\pair}[2]{\langle#1 ,#2\rangle}

\newcommand{\trpl}[3]{\langle#1 ,#2 ,#3 \rangle}
\newcommand{\qdrpl}[4]{\langle#1 ,#2 ,#3 ,#4 \rangle}

\newcommand{\change}[1]{\lower 0.7pt \hbox{\mbox{\large $#1$}}}

\newcommand{\fnn}[3]{#1:#2 \raro #3}
\newcommand{\tendin}[3]
{\lim_{#1\kern 1pt \ni\kern 1.5pt #2\kern 1pt \rightarrow\kern
1.5pt #3}}
\newcommand{\tend}[2]
{\lim_{#1\kern 1pt \rightarrow\kern 1.5pt #2}}

\newcommand{\iso}[3]{\m{ #1 : #2 \cong #3 }}

\newcommand{\titles}[1]{\bigskip\ct {\Large\bf #1}\bigskip
\bigskip}

\newcommand{\authorss}[4]{\bigskip\ct {#1} \ct{#2}
\ct{#3} \ct{#4} }

\newcommand{\inverse}{^{-1}}

\newcommand{\thinrd}[1]{\kern -1pt #1}
\newcommand{\thinthinrd}[1]{\kern -1.7pt #1}
\newcommand{\thinld}[1]{#1  \kern -1pt}
\newcommand{\thinthinld}[1]{#1  \kern -1.7pt}

\newcommand{\rcheck}{{\lower 1.7pt\hbox{\large $\check{\
}\kern1pt$}}}
\newcommand{\zspccheck}{
\kern3pt\lower1.7pt\hbox{\large $\check{}$\kern-4pt}}
\newcommand{\rhat}{{\lower 1.7pt\hbox{\large $\hat{\
}\kern1pt$}}}
\newcommand{\zspchat}{
\kern3pt\lower1.7pt\hbox{\large $\hat{}$\kern-4pt}}

\newcommand{\barC}{\bar{C}}
\newcommand{\barI}{\bar{I}}
\newcommand{\barJ}{\bar{J}}
\newcommand{\barK}{\bar{K}}
\newtheorem{note}[theorem]{Note}
\newcommand{\tldc}{\tilde{c}}
\newcommand{\R}{\bf R}
\newcommand{\Q}{\bf Q}
\newcommand{\Z}{\bf Z}
\newcommand{\C}{\bf C}
\newcommand{\LL}{\bf L}
\newcommand{\irr}{\mbox{\scriptsize irr}}
\newcommand{\rrestriction}{\kern-3pt\restriction\kern-3pt}
\def\endproof{{\kern 6pt\penalty 500
   \raise -0.5pt\hbox{\vrule height6pt width1.1pt \vbox to6pt
   {\hrule height1.1pt width 5pt
  \vfill\hrule height1.1pt}\vrule height6pt width1.1pt}}\ \,}
\def\eendproof{\hfill {\kern 6pt\penalty 500
   \raise -0.5pt\hbox{\vrule height6pt width1.1pt \vbox to6pt
   {\hrule height1.1pt width 5pt
  \vfill\hrule height1.1pt}\vrule height6pt width1.1pt}}}
\begin{document}
\baselineskip 18pt

\addtolength{\abovedisplayskip}{-3pt}
\addtolength{\belowdisplayskip}{-3pt}
\addtolength{\abovedisplayshortskip}{-3pt}
\addtolength{\belowdisplayshortskip}{-3pt}


%
\titles{Locally Moving Groups}
\vspace{-1cm}
\titles{and the Reconstruction Problem}
\vspace{-1cm}
\titles{for Chains and Circles\footnote{The authors would like to thank
Matthew Brin for inspiring this paper and making us write it.}}

    \bigskip

    \authorss
   {\bf Stephen McCleary}
{Department of Mathematics}{Bowling Green State University}
{Bowling Green,  Ohio 43403, USA}
    \bigskip

    \ct{ and}

    \authorss
   {\bf Matatyahu Rubin
\footnote{The results in this paper were found in 1995,
when the second author visited in Bowling Green State University.
The second author thanks Bowling Green State University.}
}
{Department of Mathematics}{Ben Gurion University of the
    Negev}{Beer Sheva, Israel}

\vfill


\newpage

%

    \begin{abstract}
    \begin{small}
For groups of automorphisms of chains (linear orders) and circles
which satisfy certain weak transitivity hypotheses,
the Dedekind completion of the chain or circle
can be reconstructed up to reversal of order or orientation from
the group alone.
 
Let $\pair{L}{<}$ be a dense chain. Throughout, an interval means
a nonempty open interval. A subgroup $G \leq \fs{Aut}(\pair{L}{<})$
is {\it locally moving}
if for every interval $I$ there exists
$g \in G - \sngltn{\fs{Id}}$
such that\break
$g \rrestriction (L - I) = \fs{Id}$.
$G$ is {\it $n$-interval-transitive} if for every pair of sequences of
intervals
$I_1 < \ldots < I_n$ and $J_1 < \ldots < J_n$
there exists $g \in G$ such that
$g(I_1) \cap J_1 \neq \emptyset$,\ $\ldots$\ ,
$g(I_n) \cap J_n \neq \emptyset$.
The Dedekind completion of $L$ is denoted by $\barL$.
It is assumed that $\barL$ does not add endpoints to $L$.
 
{\bf Theorem 1}
Assume $L_i$ has no endpoints
and let $G_i \leq \fs{Aut}(L_i)$ be locally moving and
$2$-interval-transitive, $i = 1,2$.
Suppose that $\iso{\alpha}{G_1}{G_2}$.
Then there is a monotonic
bijection
$\fnn{\tau}{\barL_1}{\barL_2}$ which
induces~$\alpha$, that is,
$\alpha(g) = \tau \circ g \circ \tau\inverse$
for every $g \in G$; and $\tau$ is unique.
(In particular, if $G \leq \fs{Aut}(L)$ is locally moving and\break
$2$-interval-transitive,
then eve\-ry $\alpha \in \fs{Aut}(G)$ is conjugation by a
unique monotonic bijection $\tau$ of $\barL$.)
 
The next theorem shows that the hypothesis of
2-interval-transi\-tivity in Theorem 1 is rather weak.
 
{\bf Theorem 2} There exist $L$ without endpoints
and $G \leq \fs{Aut}(L)$ such that $G$ is locally moving
and $2$-interval-transitive,
but (i) $G$ is not $3$-interval-transitive,
and (ii) $G$ does not act 2-order-transitively on any of its orbits in
$\barL$
(that is, there is no orbit $L_0 \subseteq \barL$
such that for every $x < y$ and $u < v$ in $L_0$
there exists $g \in G$ such that $g(x) = u$ and $g(y) = v$).
 
The ``circularization'' of $L$ is the structure
$\fs{CR}(L) = \pair{L}{\fs{Cr}^L}$,
where
$\fs{Cr}^L(x,y,z) \equiv
(x < y < z) \vee (y < z < x) \vee (z < x < y)$.
For pairwise disjoint circular intervals $I_1,\ldots,I_n$,
$\fs{Cr}(I_1,\ldots,I_n)$ means that
for every $1 \leq i < j < k \leq n$
and $x,y,z$ in
$I_i$, $I_j$, $I_k$
respectively, $\fs{Cr}(x,y,z)$ holds.
$n$-interval-transitivity is defined here as for chains,
but with ``$I_1 < \ldots < I_n$''
replaced by ``$\fs{Cr}(I_1,\ldots,I_n)$''.
 
{\bf Theorem 3}
Let $G_i \leq \fs{Aut}(\fs{CR}(L_i))$ be
locally moving and\break $3$-interval-transitive, \,$i = 1,2$.
Suppose that $\iso{\alpha}{G_1}{G_2}$. Then there is a unique
orientation-preserving or -reversing
bijection
\newline
$\fnn{\tau}{\fs{CR}(\barL_1)}{\fs{CR}(\barL_2)}$
which induces~$\alpha$.
 
{\bf Theorem~4}
In Theorems 1 and 3 there is a first order interpretation of the
action of $G$ on a dense subset of $\barL$
(respectively, $\fs{CR}(\barL)$).
 
{\thickmuskip=2.5mu \medmuskip=1.5mu \thinmuskip = 0.5mu
The ``equal direction'' structure of $L$
is $\fs{ED}(L) = \pair{L}{\fs{Ed}^L}$,
}
where
\newline
$\fs{Ed}^L(x_1,x_2,y_1,y_2) \equiv
 ((x_1 < x_2) \wedge (y_1 < y_2)) \vee
((x_1 > x_2) \wedge (y_1 > y_2))$.
The ``equal orientation'' structure $\fs{EO}(L)$ is defined from
$\fs{CR}(L)$ in an analogous way,
 
The analogues of the above theorems hold for $\fs{ED}(L)$'s given
$3$-interval-transitivity,
and for $\fs{EO}(L)$'s given $4$-interval-transitivity.
 
{\bf Theorem 5}
(a) If $N = \pair{L}{<}$ or $N = \fs{ED}(L)$,
and $G \leq \fs{Aut}(N)$ is locally moving and $3$-interval-transitive,
then $G$ is $n$-interval-transitive for every $n$.
 
{\thickmuskip=2.5mu \medmuskip=1.5mu \thinmuskip = 0.5mu
(b) If $N = \fs{CR}(L)$ or $N = \fs{EO}(L)$,
and $G \leq \fs{Aut}(N)$ is locally moving and $4$-interval-transitive,
then $G$ is $n$-interval-transitive for every $n$.
}
 
Other weak transitivity properties of the structures
$\pair{L}{<}$, $\fs{CR}(L)$, $\fs{ED}(L)$ and $\fs{EO}(L)$
are investigated.
 
In some cases, with different transitivity hypotheses,
reconstruction is obtained via second order interpretations,
that is, interpretations which use formulas of second order logic.
 
    \end{small}
    \end{abstract}
\newpage
\tableofcontents

\baselineskip 18pt

\newpage

%

\section{Introduction}
\label{INTRO}

This work deals with groups of automorphisms of four closely related
types of structures:
linearly ordered sets (chains), equal direction structures (essentially
the same as betweenness structures), circularly ordered sets, and
equal orientation structures.
We call such structures nearly ordered structures.

Our primary goal is to show that if $N$ is a nearly ordered structure,
$G \leq \fs{Aut}(N)$, and $\pair{N}{G}$ satisfies some rather mild
additional hypotheses, then up to possible reversal of its order
(orientation),
the Dedekind completion $\barN$ of $N$ can be reconstructed from
$G$;
and further that every group automorphism of $G$ is induced by a
unique order(orientation)-preserving or -reversing bijection of
$\barN$.  Foremost among these hypotheses is that $\pair{N}{G}$ is
``locally moving'' (see Definition~\ref{FIRSTDEFLM}).

There is a long history of such reconstruction problems, involving a
gradual falling away of the additional hypotheses.  Before
recapitulating some of this history, we define the four types of
nearly ordered structures and an array of additional hypotheses.
Throughout this paper $\pair{L}{<}$ will denote a dense linearly
ordered set (with more than one element), often abbreviated by $L$
alone.

\begin{defn} \label{}
\begin{rm}
A {\em nearly ordered structure} $N$ is a structure of any one of the
following types:

(1) A dense linearly ordered set $L = \pair{L}{<}$ without endpoints.

(2) The {\em equal direction structure} based on a chain $L$ as in
(1), defined as
$\fs{ED}(L) \eqdf \pair{L}{\fs{Ed}}$,
where
\begin{eqnarray*}
& \fs{Ed} = \fs{Ed}^L \eqdf \\*
& \setm{\vecx \cat \vecy}{\vecx, \vecy \in L^2 \mbox{ and }
((x_1 < x_2) \wedge (y_1 < y_2))
\: \vee \:
((x_1 > x_2) \wedge (y_1 > y_2))}.
\end{eqnarray*}
\fs{ED}$(L)$ is essentially the same as the usual betweenness
structure on $L$ (see Proposition~\ref{BETW}), and its automorphism
group consists of the order-preserving and order-reversing
permutations of $L$.  $\fs{Aut}(L) \leq \fs{Aut}(\fs{ED}(L))$, a
subgroup of index $\leq 2$.

(3) The {\em circularly ordered set}, or {\em circle}, based on a
dense chain $\pair{L}{<}$ having at most one endpoint, defined as
$\fs{CR}(L) \eqdf
\pair{L}{\fs{Cr}}$,
where
\[
\fs{Cr} = \fs{Cr}^L \eqdf
\setm{\trpl{x}{y}{z} \in L^3}
{(x < y < z) \vee (y < z < x) \vee (z < x < y)}.
\]
When $L$ is clear from context, we sometimes denote $\fs{CR}(L)$
by $C$.  The reader is encouraged to visualize the triples
$\trpl{x}{y}{z}$ in \fs{Cr} as ``oriented counterclockwise''.

(4) The {\em equal orientation structure} based on a chain $L$ as in
(3), defined as
$\fs{EO}(L) \eqdf
\pair{L}{\fs{Eo}}$, where
\begin{eqnarray*}
& \fs{Eo} = \fs{$Eo^L$} \eqdf
\setm{\vecx \cat \vecy}{\vecx, \vecy \in L^3 \mbox{ and }
\\*
&(\fs{Cr}(x_1,x_2,x_3) \wedge \fs{Cr}(y_1,y_2,y_3))
\: \vee \:
(\fs{Cr}(x_3,x_2,x_1) \wedge \fs{Cr}(y_3,y_2,y_1))}.
\end{eqnarray*}
Sometimes we denote $\fs{EO}(L)$ by $\fs{EO}(C)$ and speak of
$\fs{EO}(C)$ as being based on $C$.  $\fs{EO}(C)$ is essentially the
same as the usual separation structure on $C$ (see
Proposition~\ref{SEP}), and its automorphism group consists of the
orientation-preserving and -reversing permutations of $C$.
$\fs{Aut}(C) \leq \fs{Aut}(\fs{EO}(C))$, a subgroup of index $\leq 2$.

\end{rm}
\end{defn}

When $N$ is a nearly ordered structure with universe $|N|$, we
sometimes write $x \in N$ to mean $x \in |N|$, and similarly for $I
\subseteq N$.

\begin{example}
\label{REALCIR}
Let $\C$
denote the usual real circle in
${\R}^2$. Let
\[
\fs{Ccw} = \setm{\vecx \in {\C}^3}
{x_1,x_2,x_3 \mbox{ are oriented counterclockwise}}.
\]
Then $\pair
{\C}{\fs{Ccw}}$ is isomorphic to
$\fs{CR}(\pair{[0,1)}{<})$, but not to
$\fs{CR}(\pair{\R}{<})$.  In contrast, when $[0,1)$ denotes the {\em
rational} unit interval, $\fs{CR}(\pair{[0,1)}{<}) \cong
\fs{CR}(\pair{\Q}{<})$.
\end{example}

For $\pair{L}{<}$ a linear order, we denote the
reverse order of $<$ by $<^*$, and denote
$\pair{L}{<^*}$ by $L^*$.
Since \fs{ED}$(L)$ coincides with \fs{ED}$(L^*)$, $\pair{L}{<}$ can be
recovered from \fs{ED}$(L)$ only up to possible order reversal, so
$\fs{ED}(L)$ can equally well be based on $L^*$ (and these are the
only two linear orders on which \fs{ED}$(L)$ can be based).

There are many linear orders on which $C = \fs{CR}(L)$ can be based,
one for each place at which $C$ could be ``cut'' to give a chain.
$\fs{CR}(L^*) = C^*$, which denotes $C$ with the reverse orientation.
$\fs{EO}(C) = \fs{EO}(C^*)$, and these are the only two circular orders
on which $\fs{EO}(C)$ can be based.

For each type of nearly ordered structure $N$, all automorphisms of
$N$ are homeomorphisms with respect to the order topology (or
analogue thereof) on $N$.  Because of the Dedekind completeness of
$\pair{\R}{<}$, the automorphisms of $\fs{ED}(\pair{\R}{<})$
constitute {\em all} homeomorphisms; and similarly for the real
circle $\pair{\C}{\fs{Ccw}}$.

\begin{defn} \label{}
\begin{rm}
Let $N$ be a nearly ordered structure, and let $G \leq \fs{Aut}(N)$.
We call $\pair{N}{G}$ a {\em nearly ordered permutation group}
(though it is actually the structure being permuted that is nearly
ordered, not the group).  These come in four types.  We call
$\pair{N}{G}$ a

(1) {\em Linear permutation group} if $N=\pair{L}{<}$,

(2) {\em Monotonic permutation group} if $N=\fs{ED}(L)$,

(3) {\em Circular permutation group} if $N=\fs{CR}(L)$,

(4) {\em Monocircular permutation group} if $N=\fs{EO}(L)$.

\end{rm}
\end{defn}

$\pair{L}{G}$ is often called an {\em ordered permutation group}
when $G$ is partially ordered via the {\em pointwise order}, wherein
$f \leq g$ iff $f(x) \leq g(x)$ for all $x \in L$.  Here, however, $G$ is
{\em not} ordered and we use the term ``linear''.
If $G \leq \fs{Aut}(\pair{L}{<})$ then certainly $G \leq
\fs{Aut}(\fs{ED}(L))$, but (1) is not a special case of (2) because the
structures are different.

In case (2), in which the elements of $G$ are {\em monotonic}
permutations (either preserving or reversing order) of the chain
$\pair{L}{<}$, let
\[
\fs{Opp}(G) \eqdf \setm{g \in G}{g \mbox{ preserves order}} = G \cap
\fs{Aut}(\pair{L}{<}) \leq G,
\]
the subgroup of order-preserving permutations, a subgroup of index
$\leq 2$.  $\pair{\fs{ED}(L)}{\fs{Opp}(G)}$ is a monotonic
permutation group and is unaffected by whether $\fs{ED}(L)$ is
based on $L$ or $L^*$.  $\pair{L}{\fs{Opp}(G)}$ and
$\pair{L^*}{\fs{Opp}(G)}$ are linear permutation groups.

Similar comments apply to (3) and (4).  In (4), in which the elements
of $G$ are {\em monocircular} permutations of $C$ (either preserving
or reversing orientation), we use the same notation $\fs{Opp}(G)$ as
before, this time for the subgroup of {\em orientation}-preserving
permutations.

In cases (1) and (3) we let $\fs{Opp}(G) = G$.  In each of the four
cases, $\pair{N}{\fs{Opp}(G)}$ is a nearly ordered permutation group
of the same type as $\pair{N}{G}$.

We next set the conventions for the Dedekind completion $\barN$ of
a nearly ordered structure.  Recall that $L$ has no endpoints in cases
(1) and (2), and at most one endpoint in (3) and (4).

\begin{defn} \label{}
\begin{rm}
Let $N$ be a nearly ordered structure.

(1) If $N=L$ then $\barN \eqdf \barL$, the Dedekind completion of
$L$ without endpoints.

(2) If $N=\fs{ED}(L)$ then $\barN \eqdf \fs{ED}(\barL)$.

(3) If $N=\fs{CR}(L)$ then $\barN \eqdf \fs{CR}(\barL^{Cr})$, where
$\barL^{Cr}$ is $\barL$ with an endpoint adjoined at one end if
$\barL$ has no endpoints.

(4) If $N = \fs{EO}(L)$ then $\barN \eqdf \fs{EO}(\barL^{Cr})$.

\end{rm}
\end{defn}

Note that if $N$ is nearly ordered, then $\barN$ is nearly ordered
and of the same type.
It is obvious that if $g \in \fs{Aut}(N)$, then there is a
unique $\barg \in \fs{Aut}(\barN)$ which extends $g$.
Further, the function $g \mapsto \barg$ is an embedding of
$\fs{Aut}(N)$ in $\fs{Aut}(\barN)$.
We shall not distinguish between $g$ and $\barg$.

\begin{defn} \label{DEFINT}
\begin{rm}
The $N$-intervals of a nearly ordered structure $N$, all of which are
open and nonempty, are defined as follows:

(1) For distinct $x < y \in L$, let $(x,y) = \setm{z \in L}{x < z < y}$,
$(-\infty,x) = \setm{z \in L}{z < x}$,
$(x,\infty) = \setm{z \in L}{x < z}$,
and $(-\infty,\infty) = L$.
For $N = \pair{L}{<}$ or $N = \fs{ED}(L)$, the {\em bounded
$N$-intervals} are those of the first type, the {\em $N$-rays} are
those of
the second and third types, and the {\em $N$-intervals} include all
four types.   The terminology for $\fs{ED}(L)$-intervals is unaffected
by whether $\fs{ED}(L)$ is based on $L$ or $L^*$ (though these {\em
notations} apply only when the base is specified).

(2) For $N=C=\fs{CR}(L)$ or $N=\fs{EO}(C)$, the {\em bounded
$N$-intervals} are the sets $(x,y)$, where $x$ and $y$ are distinct
elements of $C$ and $(x,y)$ is the {\em circular interval} $\setm{z
\in C}{\fs{Cr}(x,z,y)}$.  The {\em unbounded $N$-intervals} are $C$
itself, and for each $x \in C$ the set $C_x \eqdf \setm{z \in C}{z \neq
x}$.
The terminology for $\fs{EO}(C)$-intervals is unaffected by whether
$\fs{EO}(C)$ is based on $C$ or $C^*$ (though the notation $(x,y)$
applies only when the base is specified).

(3) Sometimes we deal with $N$-intervals and sometimes with
$\barN$-intervals.  For an $N$-interval $I \neq N$, $\barI$ will
denote the unique $\barN$-interval such that $\barI \cap N = I$.

 (4) A subset of $\barN$ is {\em bounded} if it is contained in a
bounded  $\barN$-interval.

\end{rm}
\end{defn}

\begin{defn}\label{FIRSTDEFLM}
\begin{rm}
Let $N$ be a nearly ordered structure, and let $G \leq \fs{Aut}(N)$.

(1) For $g \in G$, the {\em support} of $g$ means $\setm{x \in
\barN}{g(x) \neq x}$, an open subset of $\barN$, and is denoted by
\fs{supp}$(g)$.  When $I$ is an $N$-interval, we sometimes write
$\fs{supp}(g) \subseteq I$ as short for
$\fs{supp}(g) \subseteq \barI$, and say
that $I$ {\em supports} $g$.  We refer to an element having bounded
support as a {\em bounded element}.  When $N = C$ or $N =
\fs{EO}(C)$, this is equivalent to the existence of a $C$-interval $I$
such that $g\rest I = \fs{Id}$.

(2) $G$ is a {\em locally moving subgroup} of $\fs{Aut}(N)$
if for every $N$-interval $I$ there exists
$g \in G - \sngltn{\fs{Id}}$
such that $g \rest (|N| - I) = \fs{Id}$; that is, if every $N$-interval
supports some $g \in G - \sngltn{\fs{Id}}$.  Here we also refer to the
nearly ordered permutation group $\pair{N}{G}$ as being locally
moving.
This is our most crucial hypothesis---in all our reconstruction results
it will be assumed that $\pair{N}{G}$ is locally moving.

\end{rm}
\end{defn}

\begin{defn} \label{}\begin{rm}
Let $N$ be a nearly ordered structure.

(1) Let $N=\pair{L}{<}$, and let $n \geq 2$.  Let $I_1,\ldots,I_n
\subseteq L$.   We write  $I_1 < \cdots < I_n$ if $x_1 < \cdots < x_n$
whenever $x_i \in I_i$ for $i = 1,\ldots,n$, and if each $I_i \neq
\emptyset$.

(2) Let $N = C = \fs{CR}(L)$, and let $n \geq 3$.

$\,\,\,\,\,\,\,\,\,\,$(a) Let
$x_1,\ldots,x_n \in C$.  We write $\fs{Cr}(x_1,\ldots,x_n)$ if
$\fs{Cr}(x_i,x_j,x_k)$ for all $i,j,k$ such that $1 \leq i < j < k \leq
n$.
(Informally, $\fs{Cr}(x_1,\ldots,x_n)$ iff $x_1,\ldots,x_n$ are distinct
and the sequence $x_1,x_2,\ldots,x_n,x_1$ wraps once
counterclockwise around the circle.)

$\,\,\,\,\,\,\,\,\,\,$(b) Let $I_1,\ldots,I_n \subseteq C$.  We write
$\fs{Crs}(I_1,\ldots,I_n)$ if
$\fs{Cr}(x_1,\ldots,x_n)$ whenever $x_i \in I_i$ for $i = 1,\ldots,n$,
and if each $I_i \neq \emptyset$.
$\fs{Crs}(I_1,\ldots,I_n)$ implies that $I_1,\ldots,I_n$ are pairwise
disjoint.

\end{rm}
\end{defn}

We turn now to transitivity assumptions, the exact nature of which
will vary from one result to another.  These transitivity properties
are weak forms of multiple transitivity.  In each half of the following
definition, version (a) is a standard notion and obviously implies
ordinary transitivity, whereas (b) and (c) do not.  These properties
are unaffected by reversal of the order on $L$ or the orientation on
$C$.

\begin{defn} \label{DEFMULTTRANS}
\begin{rm}
Let $N$ be a nearly ordered structure, $G \leq \fs{Aut}(N)$, and $n
\geq 1$.

(1) Let $N=\pair{L}{<}$ or $N=\fs{ED}(L)$.

$\,\,\,\,\,\,\,\,\,\,$(a) We call $\pair{N}{G}$ {\em (exactly)}
$n${\em -o-transitive} (``{\it o}'' for ``order'') if for all $a_1 < a_2 <
\cdots <
a_n$ and $b_1 < b_2 < \cdots < b_n$ in $L$, there exists $g \in
G$ such that $g(a_i) = b_i$ for $i = 1,\ldots,n$.

$\,\,\,\,\,\,\,\,\,\,$(b) We call $\pair{N}{G}$ {\em approximately}
$n${\em -o-transitive} if for all $a_1 < \ldots < a_n$ in $L$ and
$L$-intervals $I_1 < \ldots < I_n$, there exists $g \in G$ such
that $g(a_i) \in I_i$ for $i = 1,\ldots,n$.

$\,\,\,\,\,\,\,\,\,\,$(c) We call $\pair{N}{G}$
$n${\em -interval-transitive} if for all $L$-intervals $I_1 < \ldots <
I_n$ and $J_1 <
\ldots < J_n$, there exists $g \in G$ such that $g(I_i) \cap
J_i \neq \emptyset$ for $i = 1,\ldots,n$.

(2) Let $N = C = \fs{CR}(L)$ or $N = \fs{EO}(C)$.

$\,\,\,\,\,\,\,\,\,\,$(a) We call $\pair{N}{G}$ {\em (exactly)}
$n${\em -o-transitive} (``{\it o}'' for ``orientation'') if whenever
$\fs{Cr}(a_1,\ldots,a_n)$ and  $\fs{Cr}(b_1,\ldots,b_n)$ in $C$, there
exists $g \in G$ such that $g(a_i) = b_i$ for $i = 1,\ldots,n$.
(When $n = 2$ the notation $\fs{Cr}(a_1,a_2)$ is undefined, and we
interpret $\fs{Cr}(a_1,a_2)$ to mean simply that $a_1$ and $a_2$ are
distinct.  We interpret $\fs{Cr}(a_1)$ to be automatically true.)

$\,\,\,\,\,\,\,\,\,\,$(b) We call $\pair{N}{G}$ {\em approximately}
$n${\em -o-transitive} if for all $a_1, \ldots ,a_n$ in $C$ and
$C$-intervals $I_1,\ldots,I_n$ such that $\fs{Cr}(a_1,\ldots,a_n)$ and
$\fs{Crs}(I_1,\ldots,I_n)$, there exists $g \in G$ such that
$g(a_i) \in I_i$ for $i = 1,\ldots,n$.  (We interpret  $\fs{Crs}(I_1,I_2)$
to mean that $I_1$ and $I_2$ are disjoint, and interpret
$\fs{Crs}(I_1)$ to be
automatically true.)

$\,\,\,\,\,\,\,\,\,\,$(c) We call $\pair{N}{G}$
$n${\em -interval-transitive} if for all $C$-intervals $I_1,\ldots,I_n$
and
$J_1,\ldots,J_n$ such that $\fs{Crs}(I_1,\ldots,I_n)$ and
$\fs{Crs}(J_1,\ldots,J_n)$, there exists $g \in G$ such that
$g(I_i) \cap J_i \neq \emptyset$ for $i = 1,\ldots,n$.

\end{rm}
\end{defn}

Each of these six forms of multiple transitivity is
enjoyed by $\pair{N}{\fs{Opp}(G)}$ iff it is enjoyed by $\pair{N}{G}$,
except in (1) when $n=1$ (ordinary transitivity), and in (2)
when $n \leq 2$.

Clearly (a) $\Rightarrow$ (b) $\Rightarrow$ (c) within (1) and also
within (2).  At least in (1) these implications are not reversible for
$n=2$,
even when restricted to locally moving $\pair{N}{G}$; see
Examples \ref{WR2} and \ref{WR1}.

\begin{question}
\label{QREV}
\begin{rm}
Are these non-reversibilities also valid in (2) for $n=2$?  For $n=3$?
\end{rm}
\end{question}

For each form, $n$-transitive $\Rightarrow$ $(n-1)$-transitive, and
there is a notion of {\em highly transitive} (meaning $n$-transitive
for all $n$), referred to as highly {\it o-}transitive, etc.  We find in
Section~\ref{TRANSPROPS} that in the locally moving case the
situation is as follows:  In (1) 3-transitivity implies high transitivity
for each of the three forms, but in view of Example~\ref{WR1}
2-transitivity does not (though it
does in (a) and (b) if $L$ is Dedekind complete).   In (2)
4-transitivity implies high transitivity for each of the forms (as does
3-transitivity in (a) and (b) when $C$ is Dedekind complete).

In the locally moving case the distinction between ``approximate''
and ``interval'' transitivity matters only at the lower levels:
3-interval-transitivity implies high approximate transitivity in (1),
as
does 4-interval transitivity in (2).  See Section~\ref{TRANSPROPS}.

Let $\pair{N}{G}$ be 2-interval-transitive.  If any one
$g \in G- \sngltn{\fs{Id}}$ is bounded then $\pair{N}{G}$ is locally
moving; see
Proposition~\ref{IT}.  (Among 2-{\it o-}transitive linear permutation
groups, the term {\em pathological} has often been used to refer to
those lacking nonidentity bounded elements.)

A generic reconstruction theorem for linear permutation groups will
read like this, with the blank to be filled in by some multiple
transitivity property:

\bigskip
\noindent {\bf Theorem A$\,$}\begin{it}
Let $\pair{L_i}{G_i}$ be a locally moving linear permutation group
and assume that $\pair{L_i}{G_i}$ is
\ \lower 0.2mm \hbox{\rule{2.5cm}{.2mm}},
 $i=1,2$.
Suppose that $\iso{\alpha}{G_1}{G_2}$.
Then there exists a unique monotonic bijection
$\fnn{\tau}{\barL_1}{\barL_2}$
which induces $\alpha$, that is,
$\alpha(g) = \tau \circ g \circ \tau\inverse$ for all $g \in G_1$.
\end{it}
\bigskip

It follows from Theorem A that every group isomorphism
$\iso{\alpha}{G_1}{G_2}$ either preserves or reverses pointwise
order.

There is no hope, even under very strong multiple transitivity
hypotheses, of sharpening the conclusion to guarantee that $\tau$
preserves (rather than possibly reversing) the order on the
$\barL_i$'s, because for $L_1 = \pair{{\Q}}{<}$ and
$G_1 = \fs{Aut}(L_1)$, $G_1$ has automorphisms induced by
order-reversing automorphisms of $L_1$.  Nor is there hope of
guaranteeing that  $\iso{\tau}{L_1}{L_2}$ (without Dedekind
completions) because for
$L_2 = \pair{{\bar{\Q}} - {\Q}}{<}$
and $G_2 = \fs{Aut}(L_2)$,
$G_1 \cong G_2$.

Let $N$ be a nearly ordered structure. We define $N^{\smallr}$,
the {\em relaxed} version of $N$, as follows.
If $N = \pair{L}{<}$ or $N = \fs{ED}(L)$, then
$N^{\smallr} = \fs{ED}(L)$.
If $N = \fs{CR}(L)$ or $N = \fs{EO}(L)$, then
$N^{\smallr} = \fs{EO}(L)$.
In all four cases, $\overline{N^{\smallr}} = \bar{N}^{\smallr}$.

For an arbitrary type
\ \lower 0.2mm \hbox{\rule{0.5cm}{0.2mm}} \
of nearly ordered structure the generic
theorem would read:

\bigskip
\noindent {\bf Theorem B\,}\begin{it}
Let $\pair{N_1}{G_1}$ and $\pair{N_2}{G_2}$ be locally moving
nearly ordered permutation groups of  type
\ \lower 0.2mm \hbox{\rule{0.5cm}{0.2mm}},
and
assume that $\pair{N_i}{G_i}$ is
\ \lower 0.2mm \hbox{\rule{2.5cm}{.2mm}}, $i=1,2$.
Suppose that $\iso{\alpha}{G_1}{G_2}$.
Then there exists a unique
$\iso{\tau}{\barN_1^{\smallr}}{\barN_2^{\smallr}}$
such that $\tau$ induces $\alpha$, that is,
$\alpha(g) = \tau \circ g \circ \tau\inverse$ for all $g \in G_1$;
or equivalently, such that
$\iso{\tau \cup \alpha}
{\pair{\barN_1^{\smallr}}{G_1}}
{\pair{\barN_2^{\smallr}}{G_2}}$.
\end{it}
\bigskip

For a single nearly ordered structure $\pair{N}{G}$, reconstruction
permits us to identify $\fs{Aut}(G)$ with the normalizer of $G$ in
$\fs{Aut}(\barN^{\smallr})$, that is, in the group
$\fs{Homeo}(\barN)$ of homeomorphisms of $\barN$:

\bigskip
\noindent {\bf Theorem C\,}\begin{it}
Let $\pair{N}{G}$ be a locally moving nearly ordered permutation
group, and assume that $\pair{N}{G}$ is
\ \lower 0.2mm \hbox{\rule{2.5cm}{.2mm}}.
Then every automorphism $\alpha$ of $G$ is conjugation by a
unique automorphism  $\tau$ of $\barN^{\smallr}$, so that
\[
\fs{Aut}(G) =
 \fs{Norm}_{\mbox{\scriptsize {\fs{Aut}}}(\barN^{\mbox
{\smallr}})}(G) =
 \fs{Norm}_{\mbox{\scriptsize {\fs{Homeo}}}(\barN)}(G).
\]
\end{it}

We recapitulate some of the history of such reconstruction theorems,
focusing here on how to fill in the blank in Theorem C.  In each result
mentioned here the hypotheses force $\pair{N}{G}$ to be locally
moving, and to be highly transitive (for whatever form of multiple
transitivity the result assumes), though this is not always obvious.
All these results
will follow from the more general reconstruction theorems of the
present paper, some of which will have hypotheses that do not force
high transitivity.

We begin with the linear and monotonic cases.  Here
\begin{eqnarray*}
& \fs{Bdd}(G) \eqdf \setm{g \in G}{g \mbox{ has bounded support}}
\mbox{, and}\\*
& \fs{Bdd}(L) \eqdf \fs{Bdd}(\fs{Aut}(L)) =
\fs{Bdd}(\fs{Aut}(\fs{ED}(L))).
\end{eqnarray*}

N.G. Fine and G.E. Schweigert \cite{FS} proved Theorem C for the full
automorphism group $G =\fs{Aut}(\fs{ED}({\R}))$, {\R} the chain of
real numbers (that is, for $G$ the group of all homeomorphisms of
{\R}).
J.V. Whittaker \cite{Wh} proved it for any $G$ such that
$\fs{Bdd}({\R}) \leq G \leq \fs{Aut}(\fs{ED}({\R}))$,
and his arguments apply when {\R} is replaced by any chain $L$ for
which $\fs{Aut}(L)$ is 2-{\it o-}transitive.  This more general version
was also established explicitly in different ways by D. Eisenbud
\cite{Ei} (under additional hypotheses) and by E.B. Rabinovich
\cite{Ra}.

W.C. Holland \cite{Ho1} proved an analogue of Theorem C for
lattice-ordered permutation groups $\pair{L}{G}$, whose definition
we now
give.  When any $G \leq \fs{Aut}(L)$ is ordered pointwise ($f \leq g$
iff $f(x) \leq g(x)$ for all $x \in L$), $G$ becomes a {\em partially
ordered group} (meaning that the order is preserved by
multiplication on either side).  When $G = \fs{Aut}(L)$ this gives a
{\em lattice} order, with the supremum and infimum of a pair of
elements formed pointwise.  More generally, $\pair{L}{G}$ is called a
{\em lattice-ordered permutation group} if $G$ is a sublattice of
$\fs{Aut}(L)$, that is, if $G$ is closed under pointwise suprema and
infima.  Holland proved that when $\pair{L}{G}$ is a lattice-ordered
permutation group which is 2-{\it o-}transitive and has a nonidentity
bounded element, then Theorem C holds for those automorphisms
$\alpha$ of $G$ which preserve order.  He did {\em not} have to
assume that $\fs{Bdd}(L) \leq G$.

S.H. McCleary \cite{McHOMEOS} proved Theorem C (for {\em all}
group
automorphisms $\alpha$) for linear permutation groups which are
3-{\it o-}transitive and have a positive bounded element.
(An element $g$ is {\em positive} if $g(x) \geq x$ for all $x \in L$
and $g \neq \fs{Id}$.)  M. Giraudet \cite{Gi}
sharpened McCleary's argument to obtain a first order version of
Theorem B dealing with elementarily equivalent groups, a matter to
which we return shortly.

R. Bieri and R. Strebel \cite{BS} achieved the next weakening of
hypotheses.  They proved Theorem C for linear permutation groups
which are 3-interval-transitive, and which have a nonidentity
bounded
element and a positive element (but these properties need not be
conjoined in a
single element).  (The 3-interval-transitivity was stated in a
different form which is equivalent in the presence of a nonidentity
bounded element.)  Actually, they were dealing with intervals of
{\R},
but that was not essential for their development of Theorem C.  Also,
they had an additional assumption about suprema and infima in
$\barL$ of supports of elements of $G$, namely that
\[
\fs{sup}(\fs{supp}(g_1)) = \fs{inf}(\fs{supp}(g_2)) \mbox{ for some }
g_1,g_2 \in G,
\]
though actually their methods can be used to obtain Theorem C even
without that assumption.  We shall revisit their results in
Section~\ref{APPLICATIONS}.

M.Rubin \cite{Ru} developed results on local movability (some of
which are reproduced in Section~\ref{LMG}) and used them to prove
Theorem C for the {\em full} automorphism groups $G=\fs{Aut}(L)$
and $G=\fs{Aut}(\fs{ED}(L))$ when $\pair{L}{G}$ is
2-{\it o-}transitive.  Like Giraudet, he got a first order version of
Theorem
B\@.  However, those methods do not work in the generality we are
dealing with here.

Our present results are more general yet. We get Theorems B and C
for linear permutation groups which are 2-interval-transitive
(Theorem~\ref{LRT}ff) and
for monotonic permutation groups which are 3-interval-transitive
(Theorem~\ref{MNRT}).

We infer these results from the existence of first order
interpretations of the nearly ordered permutation groups
(now acting on certain dense subchains of $\barN$)
in their respective groups.
The use of a first order interpretation rather than a second order
interpretation
allows us to deduce the extra fact that if these groups $G_1$ and
$G_2$
are elementarily equivalent,
then so are the corresponding nearly ordered permutation groups.

We also get Theorems B and C under the hypothesis of
{\em span-transitivity} (Theorem~\ref{SPANRT}), meaning that
for all distinct $a,b \in \barN$ and
all bounded $\barN$-intervals $I \subseteq J$,
there exists  $g \in G$ such that
either $g(a) \in I$ and $g(b) \not\in J$,
or $g(b) \in I$ and $g(a) \not\in J$.
Here the interpretations are in second order logic, so
we cannot deduce the result about elementary equivalence.

The history of the circular and monocircular cases is briefer.  In
\cite{Wh} Whittaker also proved Theorem C for $G =
\fs{Aut}(\fs{EO}([0,1)))$, the group of all homeomorphisms of the
real circle {\C}.  His proof also establishes the theorem for $G =
\fs{Aut}(\fs{CR}([0,1)))$, the group of orientation-preserving
homeomorphisms of {\C}.  In \cite{Ru} Rubin dealt also with the
full automorphism groups $G = \fs{Aut}(\fs{CR}(L))$ and $G =
\fs{Aut}(\fs{EO}(L))$, proving Theorem C (and a first order version of
Theorem B) for both cases when $\pair{C}{G}$ is 3-{\it o-}transitive.

Recently Matt Brin asked us whether we knew of any more general
result for the circle.  We thank him for making us write this paper,
and for his enthusiasm for the enterprise.  (His question also helped
bring about \cite{Ru}.)  Brin uses our Corollary~\ref{1GR}, together
with its analogue for circular permutation groups, in
\cite{Br}, and again in \cite{BG}.

Our present results in the circular and monocircular cases parallel
those in the linear and monotonic cases.  We get Theorems B and C
for circular permutation groups
which are 3-interval-transitive (Theorem~\ref{CRT}) and for
monocircular permutation groups which are 4-interval-transitive
(Theorem~\ref{MCRT}), using first order interpretations.

We also get these theorems under the hypothesis of
{\em nest-transitivity} (Theorem~\ref{NESTRT}), meaning that
for all bounded $\barN$-intervals $K$ and
$I \subseteq J$ with $I$ and $J$ having no common endpoints, there
exists $g \in G$
such that
$I \subseteq g(K) \subseteq J$,
that is, $g(K)$ is nested between $I$ and $J$.
The interpretations are in second order logic.

The paper has been arranged so as to accommodate a variety of
readers.  We offer now a guide to the rest of the paper, with
suggestions for readers who want to read only selected parts.

Section \ref{OBJECTS}  provides a bit of basic background and
offers some context that is helpful but not explicitly used in the
sequel.

Section \ref{PRE}, after settling some preliminaries, surveys
some simplicity results about the derived group $\fs{Bdd}'(G)$, again
as helpful context rather than for explicit use.  Finally, for a nearly
ordered permutation group $\pair{N}{G}$, it introduces the Boolean
Algebra of regular open sets of $\barN$ on which $G$ will act.

Section \ref{OVERVIEW} deals with the linear case, the simplest of
the four cases.  It provides an informal and not too technical proof of
the
Linear Reconstruction Theorem~\ref{LRT} (cf. Theorem A of the
introduction), but with two huge gaps.  For a reader wanting to
understand just the linear case, guidance is offered on what else to
read to fill in the gaps.

Section \ref{LMG} is a technical section establishing (in
Theorem~\ref{ET}) that certain relationships between
automorphisms of
a Boolean algebra can be re-expressed in the first order language of
groups.  There is guidance on how to read just a selected subset
of the section.

Section \ref{INTERPS}, also technical, provides model theoretic
underpinnings and definitions of several kinds of ``interpretations''.
It uses the results of Section~\ref{LMG} to establish that all ``local
movement systems'' for a given group $G$ are ``the same''
(Theorem~\ref{LMSRT}).  Results involving second order
interpretation, which are needed only for the second order
reconstruction in Section~\ref{FURTHER}, are starred so that they can
be skipped on first reading.

Section \ref{RECON} gives the statement of the Linear Reconstruction
Theorem~\ref{LRT} (cf. Theorem B in the introduction), deduces
numerous corollaries, and deduces \ref{LRT} from the more technical
Linear Interpretation Theorem~\ref{LIT}. Then it does the same for
the other three types of nearly
ordered permutation groups.  Finally it proves \ref{LIT} and its
analogues.

Section \ref{FURTHER} deals with second order reconstruction, as
discussed above.  It can be skipped without detriment to
understanding the rest of the paper, especially if the reader is willing
to look up occasional definitions there.

Section \ref{APPLICATIONS} further develops certain ideas from
Section~\ref{RECON} and deduces various theorems from
papers cited in our introduction, but this time with weaker
hypotheses.  At the end is a discussion of the uses made of our
results in \cite{Br} and \cite{BG}.  Section~\ref{APPLICATIONS} can
be skipped in the same fashion as Section~\ref{FURTHER}.

Section \ref{PARAMETERS} deduces results of Giraudet \cite{Gi} on
interpretations with a parameter, now with weaker hypotheses.  It
can also be skipped in the above fashion.

Section \ref{LATTICE} discusses lattice-ordered permutation groups.

Section \ref{TRANSPROPS} discusses relationships among the various
transitivity properties, and is often referred to from earlier in the
paper.  It is nearly self-contained, and much of it could be read
immediately after the introduction.

Section \ref{EXAMPLES} offers numerous examples.

\medskip

Open questions are stated as they arise.  Chief among these are
whether in the circular and monocircular cases, assuming local
movability, the degree of multiple transitivity sufficient to give a
reconstruction theorem, or to guarantee high interval-transitivity,
can be reduced (Questions \ref{QCIR3TO2}, \ref{QMONOCIR4TO3}, and
\ref{QHIGH}).

A matter of notation:
The symbols ``$<$'' and ``$\leq$'' are not being used just as strict and
non-strict versions of the same order.
For the order on the Boolean algebras we consistently use the symbol
``$\leq$'', whereas for the order on the
chains $L$ we almost always use the symbol ``$<$'', with the few
exceptions being clear from context.
We also use ``$\leq$'' for the pointwise order on the group, on those
few occasions when we need to refer to it.

To understand this paper it is essential that the reader draw many,
many figures.
\newpage

%

\section{The objects of study}
\label{OBJECTS}

We now elaborate on our assertion that the equal direction relation
on a chain $L$ is essentially the same as the usual {\em betweenness
relation}
\[
\fs{Bet} = \fs{Bet}^L \eqdf
\setm{\trpl{x_1}{x_2}{x_3} \in L^3}{(x_1 < x_2 < x_3)  \vee  (x_1 >
x_2 > x_3)}.
\]
As with \fs{Ed}, $\fs{Bet}^{L^*}$ coincides with $\fs{Bet}^L$, and $L$
and $L^*$ are the only two linear orders inducing $\fs{Bet}^L$.
This ensures that \fs{Ed} and \fs{Bet} can be recovered from each
other, but we need a sharper statement.  The definitions of \fs{Ed}
and \fs{Bet} are valid for arbitrary chains $L$ (not necessarily dense
and without endpoints), and so is the following proposition.

The relations \fs{Ed} and \fs{Bet} are bi-interpretable in the
following sense:

\begin{prop}\
\label{BETW}

{\rm (1)} There is a first order formula
$\varphi_{\fss{Bet}}(x_1,x_2,x_3)$ in the language $\{ \fs{Ed} \}$
such that for any chain $L$ and any
$a_1, a_2, a_3 \in L$:
\[
\fs{ED}(L) \models \varphi_{\fss{Bet}}[a_1,a_2,a_3]
\mbox{ \ iff \ }
\fs{Bet}^L (a_1,a_2,a_3).
\]

{\rm (2)}  There is a first order formula
$\varphi_{\fss{Ed}}(x_1,x_2;x_3,x_4)$ in the language $\{ \fs{Bet} \}$
such that for any chain $L$ and any
$a_1,a_2,a_3,a_4 \in L$:
\[
\pair{L}{\fs{Bet}} \models \varphi_{\fss{Ed}}[a_1,a_2;a_3,a_4]
\mbox{ \ iff \ }
\fs{Ed}^L (a_1,a_2;a_3,a_4).
\]
\end{prop}

\noindent
{\bf Proof }
(1) $\varphi_{\fss{Bet}}(x_1,x_2,x_3) \eqqdf
\fs{Ed}(x_1,x_2;x_2,x_3)$
is as required.

(2) $\varphi_{\fss{Ed}}(x_1,x_2;x_3,x_4) \eqqdf
(x_1 \neq x_2) \ \wedge \ (x_3 \neq x_4) \ \wedge$
\newline
$(\fs{Bet}(x_1,x_2,x_3) \rightarrow \fs{Bet}(x_2,x_3,x_4)) \ \wedge \
(\fs{Bet}(x_1,x_3,x_2) \rightarrow \fs{Bet}(x_1,x_3,x_4)) \ \wedge
\newline
(\fs{Bet}(x_3,x_1,x_2) \rightarrow \neg\fs{Bet}(x_4,x_3,x_1))  \
\wedge
\newline
((x_3 = x_1) \rightarrow \neg\fs{Bet}(x_4,x_3,x_2)) \ \wedge\
((x_3 = x_2) \rightarrow \fs{Bet}(x_1,x_3,x_4))$
\smallskip
\newline
is as required.
\hfill\qed
\medskip

For circularly ordered sets $C$, the analogue of the betweenness
relation is the quaternary {\em separation relation}:
\[
\fs{Sep} = \fs{Sep}^C \eqdf
\setm{\vecx \in C^4}{\fs{Cr}(x_1,x_2,x_3,x_4) \vee
\fs{Cr}(x_4,x_3,x_2,x_1)}.
\]
As with \fs{Eo}, $\fs{Sep}^{C^*}$ coincides with $\fs{Sep}^C$, and $C$
and $C^*$ are the only two circular orders inducing $\fs{Sep}^C$.
This can be seen directly, or deduced from the following proposition
(whose proof can be adapted to arbitrary $C$).

The relations \fs{Eo} and \fs{Sep} are bi-interpretable in the
following sense:

\begin{prop}\
\label{SEP}

{\rm (1)} There is a first order formula
$\varphi_{\fss{Sep}}(x_1,x_2,x_3,x_4)$ in the language $\{ \fs{Eo}
\}$
such that for any circle $C$ and any
$a_1, a_2, a_3,a_4 \in C$:
\[
\fs{EO}(C) \models \varphi_{\fss{Sep}}[a_1,a_2,a_3,a_4]
\mbox{ \ iff \ }
\fs{Sep}^C (a_1,a_2,a_3,a_4).
\]

{\rm (2)}  There is a first order formula
$\varphi_{\fss{Eo}}(x_1,x_2,x_3;y_1,y_2,y_3)$ in the language $\{
\fs{Sep} \}$
such that for any circle $C$ and any
$\veca,\vecb \in C^3$:
\[
\pair{C}{\fs{Sep}} \models \varphi_{\fss{Eo}}[\veca,\vecb]
\mbox{ \ iff \ }
\fs{Eo}^C (\veca,\vecb).
\]
\end{prop}

\noindent
{\bf Proof }(for dense circular orders).

(1) This follows from the definition of $\fs{Sep}^C$.

(2) We define the following formulas:
\smallskip
\newline
$
\varphi_{\fss{Eo}}'(\vecx,\vecz) \, \eqqdf
\begin{displaystyle}
\mbox{\large(}
\bigwedge_{i=1}^3 \fs{Sep}(x_1,x_2,x_3,z_i) \mbox{\large)}
\, \wedge \, \fs{Sep}(x_2,z_1,z_2,z_3)
\, \wedge \, \fs{Sep}(x_1,x_3,z_1,z_3),
\end{displaystyle}
\newline
\varphi_{\fss{Eo}}''(\vecx,\vecz)  \, \eqqdf
\varphi_{\fss{Eo}}'(x_1,x_2,x_3;\vecz)  \ \vee \
   \varphi_{\fss{Eo}}'(x_2,x_3,x_1;\vecz)  \ \vee \
   \varphi_{\fss{Eo}}'(x_3,x_1,x_2;\vecz),
\smallskip
\newline
\varphi_{\fss{Eo}}(\vecx,\vecy)  \,\eqqdf
\exists \vecz \,
 (\varphi_{\fss{Eo}}''(\vecx,\vecz)  \ \wedge \
  \varphi_{\fss{Eo}}''(\vecy,\vecz)).
$
\smallskip
\newline
Then $\varphi_{\fss{Eo}}$ is as required.
\hfill\qed
\medskip

We remark that there does exist a (very long) quantifier free
formula which works in (2) for all (not necessarily dense) circular
orders.

Though we shall make no use of the rest of this section (except for
the discussion of pointwise order), it provides helpful context.

Various finite sets of first order axioms have been given to
characterize those ternary relations which are betweenness relations,
and obviously axioms can be added to specify denseness and the
absence of endpoints.  In view of Proposition~\ref{BETW}, the theory
of equal direction structures is also finitely first order axiomatizable.

P. Cameron \cite[Theorem 4.2ff]{Ca} gives a finite first order
axiomatization to characterize those quarternary relations which are
separation relations, and again an axiom can be added to specify
denseness.  Hence the theory of equal orientation structures is also
finitely first order axiomatizable.

Now we address axiomatizations of some of the classes of {\em
groups} with which we are dealing.  Consider the class
\[
\setm{G}{\exists  \ L \mbox{ such that } \pair{L}{G}
\mbox{ is a linear permutation group}}.
\]
Here the class is the same regardless of whether density and absence
of endpoints are required.  As observed by M. Giraudet and F. Lucas
\cite[p.82]{GL}, P. Conrad's classic characterization of this class of
groups provides for the class an infinite first order axiomatization in
the language of group theory.

Consider also
\[
\setm{G}{\exists \ L \mbox{ such that } \pair{\fs{ED}(L)}{G}
\mbox{ is a monotonic permutation group}}.
\]
Again the class is the same regardless of whether density is required.
Giraudet and Lucas show that this class is also first order
axiomatizable in group theoretic language.

Finally, the following class is finitely axiomatizable (Giraudet
\cite[Corollary 1-7]{Gi}):  The class of groups $G$ for which there
exists $L$ such that $\pair{L}{G}$ is a 2-{\it o-}transitive
lattice-ordered permutation group having a nonidentity element of
bounded support.

For linear permutation groups $\pair{L}{G}$, the {\em pointwise
order} on $G$ is defined by setting
$f \leq g$ iff $f(x) \leq g(x)$ for all $x \in L$.
As mentioned in the introduction, with this order $G$ becomes a
partially ordered group.

For monotonic permutation groups $\pair{\fs{ED}(L)}{G}$,
the pointwise order is defined in the same way.  This order depends
on whether $\fs{ED}(L)$ is based on $L$ or $L^*$; it is reversed by
reversal of the order on $L$.
Here the pointwise order yields what is
called a {\em half-ordered group}, a notion introduced in \cite{GL}
by Giraudet and Lucas:

\begin{defn} \label{HALFORDERED}
\begin{rm}\
(1) A {\em half-ordered group} consists of a group $G$ together with
a  nontrivial partial order ``$\leq$'' on $G$ such that for any $h \in
G$:

$\,\,\,\,\,\,\,\,\,\,$(a) For all $g_1 \leq g_2 \in G$, $g_1h \leq g_2h$;
and

$\,\,\,\,\,\,\,\,\,\,$(b) Either for all $g_1 \leq g_2 \in G$, $hg_1 \leq
hg_2$, or else for
all $g_1 \leq g_2 \in G$, $hg_1 \geq hg_2$.

(2) The {\em increasing part} of $G$ consists of those $h \in G$ such
that $hg_1 \leq hg_2$ for all $g_1 \leq g_2$; it is a subgroup of index
$\leq 2$.

(3) $G$ is a {\em half-lattice-ordered group} if its increasing part
forms a lattice.

(4) $\pair{\fs{ED}(L)}{G}$ is a {\em half-lattice-ordered permutation
group} if $\fs{Opp}(G)$ is a sublattice of $\fs{Aut}(L)$.
\end{rm}
\end{defn}

Cameron \cite[Theorem 6.1]{Ca} provides an interesting perspective
on the four classes of nearly ordered permutation groups.  Let
$\pair{X}{G}$ be a permutation group on an (unordered) infinite set
$X$.  $\pair{X}{G}$ is {\em n-homogeneous} if for any two
$n$-element sets $S_1,S_2 \subseteq X$, there exists $g \in G$ such
that
$g(S_1)=S_2$.  $n$-homogeneity $\Rightarrow$ ($n-1)$-homogeneity
\cite{Ca}.

Cameron shows that if $\pair{X}{G}$ is highly homogeneous, and if
for some $r$ $\pair{X}{G}$ is $r$-transitive (in the usual unordered
sense) but not $(r+1)$-transitive, then $X$ can be equipped with a
relation which makes it a nearly ordered structure $N$ with $G \leq
\fs{Aut}(N)$.

More fully, he shows that one of the following four cases must obtain:

(1) $r = 1$ and there is a linear order on $X$ (unique up to reversal)
giving a linear permutation group $\pair{N}{G}$, and $\pair{N}{G}$ is
2-{\it o-}transitive.  (Here 4-homogeneity suffices.)

(2) $r = 2$ and the action of $G_T$ (setwise stabilizer) on a 3-point
set $T$ contains a 2-cycle, and there is a unique equal direction
structure on $X$ giving a monotonic permutation group
$\pair{N}{G}$, and $\pair{N}{G}$ is 3-{\it o-}transitive.
(Here 4-homogeneity suffices; and by 3-homogeneity, the description
of $G_T$ is independent of $T$.)

(3) $r = 2$ and the action of $G_T$ on a 3-point set $T$ contains a
3-cycle, and there is a circular order (unique up to reversal) on $X$
giving a circular permutation group $\pair{N}{G}$, and $\pair{N}{G}$
is 3-{\it o-}transitive.
(Again 4-homogeneity suffices.)

(4) $r = 3$ and there is a unique equal orientation structure on $X$
giving a monocircular permutation group $\pair{N}{G}$, and
$\pair{N}{G}$ is 4-{\it o-}transitive.
(Here 5-homogeneity suffices.  The action of $G_T$ on a 4-point set
$T$, ordered by  $Cr$ for ``the'' circular order on which the equal
orientation structure is based, contains all 8 cyclic and reverse-cyclic
permutations.)

The statements about {\it o-}transitivity are not mentioned explicitly
in \cite{Ca}, but they follow from the descriptions of $G_T$.  Nor are
the uniqueness claims, but they follow easily from the rest of the
information.  Consider (3), for example.  Given $a,b,c \in X$ and the
choice of $Cr(a,b,c)$ rather than $Cr(a,c,b)$, this choice can be
transferred to any other three points by 3-homogeneity, and
uniquely so because of the description of $G_T$ and the lack of
3-transitivity.
\newpage

%

\section{Preliminaries}
\label{PRE}

Let $N$ be a nearly ordered structure.  We equip the universe
$|\barN|$ of $\barN$ with the topology generated by the set of
$\barN$-intervals (the {\em order-topology} on $\barN$).  The first
two propositions are obvious.

\begin{prop}\label{DENSESUB}
Let $\pair{N}{G}$ be a nearly ordered permutation group, and
suppose that  $M$ is a dense $G$-invariant subset of $\barN$.  Then

{\rm (1)} $G$ acts faithfully on $M$, and the substructure
$\pair{M}{G}$ of $\pair{\barN}{G}$ is a nearly ordered permutation
group of the same type as $\pair{N}{G}$.

{\rm (2)} If $\pair{N}{G}$ is locally moving, so is $\pair{M}{G}$.

{\rm (3)} If $\pair{N}{G}$ is $n$-interval-transitive for a given $n$,
so is $\pair{M}{G}$.
\end{prop}

However, (3) fails for the other versions of multiple transitivity in
Definition~\ref{DEFMULTTRANS}, even when $\pair{N}{G}$ is locally
moving  and $M$ consists of a single orbit of $\pair{\barN}{G}$; see
Examples \ref{WR1} and \ref{BACKANDFORTH}.
This is one way in which interval
transitivity is more ``robust'' than the other versions.  Another way
is part (2) of the following proposition, which also fails for
2-{\it o-}transitivity (see Example~\ref{WR2VARIANT}), and (we
conjecture) for
approximate 2-{\it o-}transitivity.  We shall refer to these attributes
as ``robustness under change of $G$-invariant subset'' and
``robustness under enlargement of the group''.

\begin{prop}\label{ENLARGEGROUP}
Let $\pair{N}{G}$ be a nearly ordered permutation group, let $G \leq
G_2 \leq \fs{Aut}(\barN)$, and let $G_2(N) \eqdf \setm{g_2(a)}{g_2
\in G_2, a \in N}$.  Then:

{\rm (1)} If $\pair{N}{G}$ is locally moving, so is
$\pair{G_2(N)}{G_2}$.

{\rm (2)} If $\pair{N}{G}$ is $n$-interval-transitive for a given $n$,
so is $\pair{G_2(N)}{G_2}$.

\end{prop}

We mention here a rather weak transitivity property (robust in both
the above senses) which is obviously implied by
2-interval-transitivity.  This property was mentioned in \cite{DS} for
the linear
case, and was known there as {\em feeble 2-transitivity}.
\begin{defn}
\label{DEFINCLTRANS}
\begin{rm}
We call a nearly ordered permutation group $\pair{N}{G}$ {\em
inclusion-transitive} if for any bounded $N$-intervals $I$ and $J$
there exists $g \in G$ such that $g(I) \subseteq J$.  Requiring that $g
\in \fs{Opp}(G)$ yields an equivalent definition.  (Given $I$ and $J$,
pick a bounded interval $K$ of $N$ such that $I \cup J \subseteq K$,
having first replaced $J$ by a subinterval of itself if necessary to
make this possible.  Picking $g \in G$ such that $g(K) \subseteq J$,
we have $g^2 \in \fs{Opp}(G)$ and $g^2(I) \subseteq g^2(K) \subseteq
g(J) \subseteq g(K) \subseteq J$.)
\end{rm}
\end{defn}

Note that although the chain $L$ in a linear permutation group is by
definition dense and without endpoints, the definition of
inclusion-transitivity makes sense without these restrictions on $L$,
and in fact implies them, so that in the presence of
inclusion-transitvity they become redundant (when $|L| > 2$).
This applies also to the other three types of nearly ordered
permutation groups.

\begin{defn} \label{}
\begin{rm}\

(1) For $G$ a group and $g,h \in G$, let
$g^h \eqdf h g h\inverse$.

(2) Let $K^{\fss{NO}}$ denote the class of all locally moving nearly
ordered permutation groups.
\end{rm}
\end{defn}

\begin{prop}\label{IT}
Let $\pair{N}{G}$ be a nearly ordered permutation group which is
inclusion-transitive.  Then

{\rm (1)} If $G$ has a nonidentity bounded element, then
$\pair{N}{G}$ is locally moving, that is,
$\pair{N}{G} \in K^{\fss{NO}}$.

{\rm (2)} The centralizer $Z_G(B)$ of any conjugacy class $B$ of
bounded elements is trivial.

{\rm (3)} Each orbit of $\pair{\barN}{G}$ is dense in
$\barN$, that is, $\pair{\barN}{G}$ is approximately 1-transitive.

{\rm (4)} Proposition~\ref{DENSESUB} applies to each nonempty
$G$-invariant subset $M$ of $\barN$.
\end{prop}

\noindent
{\bf Proof }
(1) and (2) hold because $\fs{supp}(g^h) = h(\fs{supp}(g))$. (3) and
(4) are clear.
\hfill\qed

\begin{prop}\label{BDD}
Suppose $\pair{N}{G}$ is locally moving, and let $\sngltn{\fs{Id}}
\neq H \unlhd G$.  Then

{\rm (1)} $H$ has a nonidentity bounded element.

{\rm (2)} If $\pair{N}{G}$ is inclusion-transitive, then $\pair{N}{G}$
is locally moving.
 \end{prop}

\noindent
{\bf Proof } (1) Pick $h \in H - \sngltn{\fs{Id}}$, and pick a bounded
$N$-interval
$I$ such that $h(I) \cap I = \emptyset$ and such that $h(I) \cup I$ is
bounded.  Then pick $p \in G - \sngltn{\fs{Id}}$ such that
$\fs{supp}(p) \subseteq I$.  Then ${\fs{Id}}\neq hph^{-1}\cdot p^{-1}
\in H$ and $\fs{supp}(p^hp^{-1})$ is bounded.

(2) follows from (1) since $\fs{supp}(k^g) = g(\fs{supp}(k))$.
\hfill\qed
\medskip

In results such as Theorem B of the
introduction, once the existence
of $\tau$ is assured, its uniqueness will follow from:

\begin{prop}\label{UNIQ}
Let $\pair{N_1}{G_1}$, $\pair{N_2}{G_2} \in K^{\fss{NO}}$.  Suppose
that $\iso{\alpha}{G_1}{G_2}$.
Then there is at most one
$\iso{\tau}{\barN_1^{\smallr}}{\barN_2^{\smallr}}$
which induces $\alpha$.
\end{prop}

\noindent
{\bf Proof }
It suffices to show that $\fs{Z}_{\mbox{\scriptsize
{\fs{Aut}}}(\barN^{\mbox
{\smallr}})}(G)$
is trivial.  Let $h \in \fs{Aut}(\barN^{\smallr})$, and suppose $h
\neq\fs{Id}$.  Pick an $\barN$-interval $I$ such that $h(I) \cap I =
\emptyset$, and pick $g \in G$ such that $\fs{supp}(g) \subseteq I$.
Then $\fs{supp}(g^h) \subseteq h(I)$, so $g^h \neq g$.
\hfill\qed

\begin{defn}
\label{}
\begin{rm}
Let $\pair{N}{G}$ be a nearly ordered permutation group, with $N =
\pair{L}{<}$ or $N = \fs{ED}(L)$.

(1) Recall that
\begin{eqnarray*}
\fs{Bdd}(G)  \eqdf  \setm{g \in G}{\fs{supp}(g)
\mbox{ is bounded in }\pair{L}{<}}, \end{eqnarray*}
a normal subgroup of $G$.

(2) Let
\begin{eqnarray*}
\fs{Lft}(G)  \eqdf  \setm{g \in G}{\fs{supp}(g) \mbox{ is bounded
above in } \pair{L}{<}}, \\*
\ \fs{Rt}(G)  \eqdf  \setm{g \in G}{\fs{supp}(g) \mbox{ is bounded
below in } \pair{L}{<}}.
\end{eqnarray*}
Both are subnormal subgroups of $G$ (normal in $G$ if $N =
\pair{L}{<}$), and
$\fs{Bdd}(G) = \fs{Lft}(G) \cap \fs{Rt}(G)$.  For $N = \fs{ED}(L)$,
$\fs{Lft}(G)$ and $\fs{Rt}(G)$ depend on whether $\fs{ED}(L)$ is
based on $L$ or $L^*$, and
are interchanged by reversal of the order on $L$.  They are normal
subgroups of $\fs{Opp}(G)$, and conjugation by any order-reversing
$g \in G$ interchanges them.
\end{rm}
\end{defn}

\begin{defn}
\label{DEFBDD}
\begin{rm}
Let $\pair{N}{G}$  be a nearly ordered permutation group, with $N =
C = \fs{CR}(L)$ or $N =  \fs{EO}(C)$.  We define
$\fs{Bdd}(G)$ to be
the (normal) subgroup of $G$ generated by the set of bounded
elements.
\end{rm}
\end{defn}

We prepare to survey some simplicity results for the derived group
$\fs{Bdd}'(G)$, and in certain cases for $\fs{Bdd}(G)$ as well.  As in
Section~\ref{OBJECTS}, though we shall not make much use of these
results,
they provide a helpful context for our reconstruction theorems.

Let $\fs{Bdd}^{(n)}(G)$ denote the $n^{\mbox{\scriptsize th}}$
derived group of $\fs{Bdd}(G)$.

\begin{prop}\label{BDDDER}
Suppose $\pair{N}{G}$ is locally moving.  Then for each $n$,
\newline
$\pair{N}{\fs{Bdd}^{(n)}(G)}$ is locally moving.
\end{prop}

\noindent
{\bf Proof }
It suffices to show this for $n=1$.  Let $I$ be any $N$-interval.  Pick
$b \in \fs{Bdd}(G) - \sngltn{\fs{Id}}$ such that $\fs{supp}(b)
\subseteq I$, then pick an $N$-interval $J$ such that $b(J) \cap J =
\emptyset$, and then pick $c \in \fs{Bdd}(G) - \sngltn{\fs{Id}}$ such
that $\fs{supp}(c) \subseteq J$.  We have $bcb^{-1}c^{-1} \in
\fs{Bdd}'(G) - \sngltn{\fs{Id}}$ and $\fs{supp}(bcb^{-1}c^{-1})
\subseteq I$.  Therefore $\pair{N}{\fs{Bdd}'(G)}$ is locally moving.
\hfill\qed
\medskip

For linear/monotonic permutation groups, we focus on some results
of G. Higman \cite{Hi}.
Higman deals with permutation groups $\pair{X}{B}$, with $X$ an
infinite set,
having the property that for all $f,g,h \in B$ with $h \neq {\fs Id}$,
there exists $k \in B$ for which $k(Y)$ is unstable with respect to
$h$ in the sense that $k(Y) \cap hk(Y) = \emptyset$, where $Y =
\fs{supp}(f) \cup \fs{supp}(g)$.

$\pair{N}{G}$ satisfies Higman's condition if it is inclusion-transitive
and all its elements are bounded (and these conditions force
$\pair{N}{G}$ to be locally moving).  Conversely, for locally moving
$\pair{N}{G}$, Higman's condition forces $\pair{N}{G}$ to be
inclusion-transitive with all its elements bounded.  (Clearly all the
elements must be bounded since some of them are.  Also, given
bounded $N$-intervals $I$ and $J$, pick $f,g \in G - \sngltn{\fs{Id}}$
such that
$$
\fs{supp}(f) < I < \fs{supp}(h) \mbox{ \,and\, } \fs{supp}(h)
\subseteq J,
$$
and pick $k$ as in Higman's condition; then $k(I) \subseteq J$.)

Higman's condition, then, will be applied to $\fs{Bdd}(G)$, and we
will phrase his condition in terms of inclusion-transitivity.  We will
see in Corollary~\ref{USED} that for any 3-interval-transitive linear
or monotonic $\pair{N}{G}$ having a nonidentity bounded element,
$\pair{N}{\fs{Bdd}(G)}$ is inclusion-transitive.  (We do not know
whether 3-interval-transitivity could be reduced to
2-interval-transitivity.)

\begin{theorem}\label{HIGMAN}
{\rm (Higman)} \
Let $\pair{N}{G}$ be a linear or monotonic permutation group having
a nonidentity bounded element.  Then

{\rm (1)} If $\pair{N}{\fs{Bdd}(G)}$ is inclusion-transitive (for
example, if $\pair{N}{G}$ is 3-interval-transitive), then
$\fs{Bdd}'(G)$ is simple, $\fs{Bdd}'(G)$ is the smallest nontrivial
normal subgroup of $\fs{Bdd}(G)$, and $\fs{Bdd}''(G) = \fs{Bdd}'(G)$.

{\rm (2)} In (1), $\fs{Bdd}'(G)$ is also the smallest nontrivial {\em
subnormal subgroup of} $G$.

{\rm (3)} If we assume only that $\pair{N}{G}$ is
inclusion-transitive, then $\fs{Bdd}'(G)$ is still the smallest nontrivial
normal
subgroup of $G$, and $\fs{Bdd}''(G) = \fs{Bdd}'(G)$.
\end{theorem}

\noindent
{\bf Proof } What Higman actually proves is (1).  We address the rest.

(2) follows from (1).  For by (1), $\fs{Bdd}'(G)$ is simple, and by
Proposition~\ref{IT} $Z_G(\fs{Bdd}'(G))$ is trivial.  Hence
$\fs{Bdd}'(G)$ is contained in every nontrivial normal subgroup $H$
of $G$, then in every nontrivial normal subgroup of $H$, etc.

Now we deal with (3), assuming inclusion-transitivity only for
$\pair{N}{G}$.  For any $f,g \in \fs{Bdd}(G)$ and $h \in G -
\sngltn{\fs{Id}}$, there exists $k \in G$ for which $k(Y) \cap hk(Y) =
\emptyset$, where $Y = \fs{supp}(f) \cup \fs{supp(g)}$.
Higman's proof of (1) then shows that $fgf^{-1}g^{-1}$ is a product of
conjugates of $h$, so that $\fs{Bdd}'(G)$ is contained in every
nontrivial normal subgroup of $G$; and shows also that $\fs{Bdd}''(G)
= \fs{Bdd}'(G)$.  By Proposition~\ref{IT}, $\fs{Bdd}'(G)$ is nontrivial.

We mention that parts of (2) and (3) are also pointed out in
\cite{DGM}.
\hfill\qed
\medskip

There are many results establishing the simplicity of $\fs{Bdd}'(G)$
in particular examples, but here we are surveying more general
results.

\medskip

When $G = \fs{Aut}(N)$ in Theorem~\ref{HIGMAN}, more can be said.
Part (1) of the following theorem is due to Higman
\cite{Hi}, and (2) to R. Ball and Droste \cite{BD}.  (3) is well known in
the linear case (see \cite{GL}, for example); the monotonic case is due
to McCleary \cite {McSIMPLE}.

\begin{theorem}\label{HIGMANBDD}
For $N=L$ or $N= \fs{ED}(L)$, let $\pair{N}{\fs{Aut}(N)}$ be
2-{\it o-}transitive (hence highly {\it o-}transitive).  Then

{\rm (1)} $\fs{Bdd}'(N) = \fs{Bdd}(N)$.  Therefore
$\fs{Bdd}(N)$ is simple, and is the smallest nontrivial subnormal
subgroup of $\fs{Aut}(N)$.

{\rm (2)} Every subnormal subgroup of $\fs{Aut}(L)$ is normal.

{\rm (3)} Suppose $L$ has a countable subset which is coinitial and
cofinal in $L$.  Then $\fs{Aut}(\fs{ED}(L))$, $\fs{Aut}(L)$,
$\fs{Lft}(L)$, $\fs{Rt}(L)$, and $\fs{Bdd}(L)$ are the only nontrivial
subnormal subgroups of $\fs{Aut}(\fs{ED}(L))$; and $\fs{Aut}(L)$,
$\fs{Lft}(L)$, $\fs{Rt}(L)$, and $\fs{Bdd}(L)$ are the only nontrivial
subnormal subgroups of $\fs{Aut}(L)$.
\end{theorem}

\noindent
{\bf Note } Droste and Shortt \cite{DS} relax the requirement in (1)
that the group be $\fs{Aut}(N)$, proving (1) for a still quite
restrictive class of $\pair{N}{G}$'s.

\medskip

We turn briefly to circular/monocircular permutation groups
$\pair{{\C}}{G}$, where $\C$ is the real circle.  D. B. A. Epstein
\cite{Ep}
gives conditions which force $G'$ to be simple (and also force
$\fs{Bdd}(G) = G$).  J. Schreier and S. Ulam \cite{SU} prove that
$\fs{Aut}({\C})$ is simple.  (Here $\fs{Bdd}'({\C}) =
\fs{Bdd}({\C}) = \fs{Aut}({\C})$.)  It follows that
$\fs{Aut}(\fs{EO}({\C}))$ and $\fs{Aut}({\C})$ are the only
nontrivial subnormal subgroups of $\fs{Aut}(\fs{EO}({\C}))$.

\bigskip

Finally, we introduce the Boolean algebra of regular open subsets of
$|\barN|$, where $N$ is a nearly ordered structure.
Because $\barN$ is Dedekind complete, every open $U \subseteq
\barN$ can be written uniquely as the union $U = \bigcup \CI_U$ of
a pairwise disjoint set $\CI_U$ of (nonempty open)
$\barN$-intervals (recall Definition~\ref{DEFINT}).
The members of $\CI_U$ are the connected components of $U$.

Let $U \subseteq \barN$.  Let $\fs{int}(U)$ denote the interior of $U$
and $\fs{cl}(U)$ the closure.  $U$ is said to be {\em regular open} in
$|\barN|$  if $U = \fs{int}(\fs{cl}(U))$.  For any $U \subseteq \barN$,
$\fs{int}(\fs{cl}(U))$ is itself regular open.

For any open $U$, $\fs{int}(\fs{cl}(U))$ consists of those points $x \in
\barN$ such that there is an $\barN$-interval $J_x$ which contains
$x$  and contains no $\barN$-interval $K$ for which $K \cap U =
\emptyset$ (so that $J_x \subseteq \fs{cl}(U)$).

For $\barN = \barL$ or $\barN = \fs{ED}(\barL)$:  All
$\barL$-intervals are regular open.  An open set $U$ is regular open
iff for all
distinct $I,J \in \CI_U$ there exists an $\barL$-interval $K$ such that
$K \cap U = \emptyset$ and $K$ is between $I$ and $J$ (that is,
either $I< K< J$ or $J< K< I$).

For $\barN = \barC$ or $\barN = \fs{EO}(\barC)$:  All
$\barC$-intervals are regular open except those of the form $\barC_x
=  \setm{z \in \barC}{z \neq x}$.  An open set $U$ is regular open iff
for
all distinct $I,J \in \CI_U$ there exist $\barC$-intervals $K_1$ and
$K_2$ such that $K_1 \cap U = \emptyset = K_2 \cap U$ and
$\fs{Crs}(I,K_1,J,K_2)$.

We shall deal extensively with the Boolean algebra of regular open
subsets of $|\barN|$.

We regard a Boolean algebra (BA) as a partially ordered set
$\pair{B}{\leq^B}$. The operations and constants of $B$ are denoted
by
$+^B$, $\cdot^B$, $-^B$, $0^B$ and $1^B$.
For a subset $C \subseteq B$, the supremum and infimum of $C$ in
$B$,
if they exist, are denoted by $\sum^B C$ and $\prod^B C$.
The superscript $B$ is usually omitted.
The language of Boolean algebras, which in our setting is just the
language of posets \{$\leq$\}, is denoted by $\CL^{\fss{PST}}$.

A Boolean algebra $B$ is {\em complete} if every subset of $B$ has
a supremum and an infimum, and is {\em atomless} if
$B - \sngltn{0}$
has no minimal elements.

It is well known and easy to prove that for any topological space $X$,
the poset $\pair{\fs{Ro}(X)}{\subseteq}$ of regular open subsets is a
complete BA.  For $N$ a nearly ordered structure, let $\barB(N)$
denote $\fs{Ro}(|\barN|)$.

For $g \in \fs{Aut}(N)$ let $g^B$ be the automorphism that $g$
induces on $\barB(N)$, that is, $g^B(U) = \setm{g(u)}{u \in U}$.
The map $g \mapsto g^B$ is an embedding of $\fs{Aut}(N)$ in
$\fs{Aut}(\barB(N))$.  We thus regard $\fs{Aut}(N)$ as a subgroup of
$\fs{Aut}(\barB(N))$.

We record the above information as

\begin{prop}\label{BN}
Let $N$ be a nearly ordered structure.
Then $\pair{\barB(N)}{\subseteq}$
is a complete atomless Boolean algebra.
For $U,V \in \barB(N)$:
\begin{eqnarray*}
& U + V = \fs{int}(\fs{cl}(U \cup V)),  \;\; U \cdot V = U \cap V,
\;\; -U = \fs{int}(|\barN| - U), \\
& \sum \setm{U_i}{i \in I} = \fs{int}(\fs{cl}(\bigcup_{i \in I} U_i)),
\mbox{ and }
\prod \setm{U_i}{i \in I} = \fs{int}(\bigcap_{i \in I} U_i).
\end{eqnarray*}
Also, $\fs{Aut}(N) \leq \fs{Aut}(\barB(N))$.
\end{prop}
\newpage

%

\section{An overview of the linear case}
\label{OVERVIEW}

Here we give an overview of the proof of the Reconstruction
Theorem in the linear case.  This is part (2) of Theorem~\ref{LRT},
also stated here for convenience.  Its hypotheses guarantee local
movability; see Proposition~\ref{IT}.  Actually, what we present here
is an informal proof that is fairly complete except for two huge gaps.

\begin{theorem}\label{OVERVIEWTHM}
Let $\pair{L_i}{G_i}$ be a linear permutation group which is
2-interval-transitive and has a nonidentity bounded element,
$i=1,2$.
Suppose that $\iso{\alpha}{G_1}{G_2}$.
Then there exists a unique monotonic bijection
$\fnn{\tau}{\barL_1}{\barL_2}$
which induces $\alpha$, that is,
\[
\alpha(g) = \tau \circ g \circ \tau\inverse \mbox{ for all } g \in G_1.
\]
\end{theorem}

Until the very end of the proof, we deal with just one $\pair{L}{G}$.
Let $g \in G$.  Even more important to us than the open set
$\fs{supp}(g)$ will be the regular open set
$\fs{int}(\fs{cl}(\fs{supp}(g))) \in \barB(L)$,
which we denote by $\fs{var}(g)$.
Let
\[
\fs{Var}(L,G) \eqdf
\setm{\fs{var}(g)}{g \in G} \subseteq \barB(L).
\]
We view $G$ as acting on $\barB(L)$ and also on $\fs{Var}(L,G)$,
and because $\pair{L}{G}$ is locally moving even this latter action is
faithful.

For $g \in G$, let $\fs{Ep}(\fs{var}(g))$ consist of the endpoints in
$\barL$ of the convex hull of $\fs{var}(g)$, or equivalently of
$\fs{supp}(g)$.  Let
\begin{eqnarray*}
& \fs{Ep}(N,G) \eqdf \bigcup _{g \in G} \fs{Ep}(\fs{var}(g)),  \\*
& \fs{EP}(N,G) \eqdf \pair{\fs{Ep}(N,G)}{\fs{Ed}^{\barL}},
\end{eqnarray*}
the restriction to $\fs{Ep}(L,G)$ of the equal direction relation on
$\barL$.  We emphasize that $\fs{EP}(L,G)$ is an equal direction
structure, not a chain.  $\fs{Ep}(L,G)$ is dense in $\barL$ since
$\pair{L}{G}$ is locally moving; also, $\fs{Ep}(L,G)$ is $G$-invariant.
So $G$ acts
faithfully on $\fs{EP}(L,G)$.

We want to ``reconstruct''  $\fs{EP}(L,G)$ , and the action
$\pair{\fs{EP}(L,G)}{G}$, from the group $G$ alone; or put another
way, to ``interpret'' $\pair{\fs{EP}(L,G)}{G}$ in $G$ using the language
$\CL^{\fss{GR}}$ of groups.

We begin with the {\em set} $\fs{Ep}(L,G)$.  The plan is to view an $x
\in \fs{Ep}(L,G)$ as, in effect, the set of all $g \in G$ such that $x =
\fs{inf}(\fs{var}(g))$ or $x = \fs{sup}(\fs{var}(g))$.
We start building a library of notions that can be expressed by
formulas in $\CL^{\fss{GR}}$.

The crux of the proof is that disjointness of $\fs{var}$'s (or
equivalently, disjointness of $\fs{supp}$'s), that is,
\[
\fs{var}(g_1) \cdot \fs{var}(g_2) = 0,
\]
can be expressed in $\CL^{\fss{GR}}$!  This will be established in
Corollary~\ref{DSJNT}, and is one of the two gaps referred to above.
(A shortcut to \ref{DSJNT} will be mentioned after
Theorem~\ref{ET}.)

The advantage of $\fs{var}$'s over $\fs{supp}$'s is that containment
of $\fs{var}$'s can be expressed in terms of disjointness:
\[
\fs{var}(g_1) \subseteq \fs{var}(g_2)
\]
can be expressed as
\[
\mbox{for every } h \in G, \
\fs{var}(h) \cdot \fs{var}(g_2) = 0 \ \Rightarrow \
\fs{var}(h) \cdot \fs{var}(g_1) = 0.
\]
This is valid because if $\fs{var}(g_1) \not\subseteq \fs{var}(g_2)$
then $\fs{var}(g_1)$ contains a nonempty open set $U$ disjoint from
$\fs{var}(g_2)$ since $\fs{var}(g_2)$ is regular open, whereupon
local movability provides an $h \in G - \sngltn{\fs{Id}} $ with
$\fs{var}(h)
\subseteq U$. Then of course
\[
\fs{var}(g_1) = \fs{var}(g_2)
\]
can also be expressed in $\CL^{\fss{GR}}$.

Let $U_1,U_2 \in \fs{Var}(L,G)$.  We say that $U_1$ and $U_2$ are
{\em segregated} if their convex hulls in $\barL$ are disjoint, that is,
if \[
U_1 < U_2 \mbox{ \ or \ } U_2 < U_1 \mbox{ \ or \ }  U_1 =
\emptyset \mbox{ \ or \ }  U_2 = \emptyset.
\]
This can be expressed as
\[
\neg(\exists f)((f(U_1) \cdot U_2 \neq 0) \wedge
(f(U_2) \cdot U_1 \neq 0)).
\]

This works because of 2-interval-transitivity, and will be the sole
use of that hypothesis (aside from permitting us to assume only the
existence of a single bounded element).  Although the above
expression isn't in $\CL^{\fss{GR}}$, we see from the fact that
$h(\fs{var}(g)) = \fs{var}(g^h)$ that the statement
\[
\fs{var}(g_1) \mbox{ and } \fs{var}(g_2) \mbox{ are segregated}
\]
can indeed be expressed in $\CL^{\fss{GR}}$.

To interpret $\fs{Ep}(L,G)$, we consider those ordered pairs
$\pair{U_1}{U_2}$ such that $U_1$ and $U_2$ are segregated and
nonempty.  Such pairs will be used to ``represent'' the points in
$\fs{Ep}(L,G)$, and will be called {\em representatives}.
We consider also the ``point'' function $\fs{Pt}$  from the set of
representatives onto $\fs{Ep}(L,G)$, given by
\[
\fs{Pt}(U_1,U_2) = \left\{
\begin{array}{ll}
\fs{sup}(U_1) & \mbox{if $U_1 < U_2$,} \\*
\fs{inf}(U_1)  & \mbox{if $U_1 > U_2$.}
\end{array}
\right.
\]
$U_1$ is the important argument here; $U_2$ merely discriminates
between $\fs{inf}( U_1)$ and $\fs{sup}(U_1)$.

We need to express in $\CL^{\fss{GR}}$ when two representatives
yield the same point, that is, when $\fs{Pt}(U_1,U_2) =
\fs{Pt}(V_1,V_2)$; and to interpret the equal direction structure
$\fs{EP}(L,G)$ we need to capture the betweenness relation on
$\fs{Ep}(L,G)$.  We deal with these matters in Lemma~\ref{DDLN};
this is the second gap in the overview.  (An impatient reader could
fill in this gap by reading the {\em proof} (not the statement) of
Lemma~\ref{DDLN}, taking ``T'' to be $\fs{Var}(L,G)$, and armed with
the understanding that $\CL^{\fss{SEG}}$ means
$\dbltn{\leq}{\fs{Seg}}$, the language of posets augmented by a
relation symbol \fs{Seg} to express the segregation relation.)

The action of $G$ on $\fs{EP}(L,G)$ can also be expressed in
$\CL^{\fss{GR}}$ since for any representative
$\pair{\fs{var}(g_1)}{\fs{var}(g_2)}$ we have
\begin{eqnarray*}
&    & h(\fs{Pt}(\fs{var}(g_1),\fs{var}(g_2)) \\*
&=  & \fs{Pt}(h(\fs{var}(g_1),h(\fs{var}(g_2)) \\*
&=  & \fs{Pt}(\fs{var}(g_1^h),\fs{var}(g_2^h).
\end{eqnarray*}

Finally, we return to a consideration of two linear permutation
groups $\pair{L_1}{G_1}$ and $\pair{L_2}{G_2}$.We are given
$\iso{\alpha}{G_1}{G_2}$.  We define
\[
\iso{\tau_{EP}}{\fs{EP}(L_1,G_1)}{\fs{EP}(L_2,G_2)}
\]
in the natural way, by setting
\[
\tau_{EP}(\fs{Pt}(\fs{var}(g),\fs{var}(h)) =
\fs{Pt}(\fs{var}(\alpha(g)),\fs{var}(\alpha(h)).
\]
Then $\tau_{EP}$ is an isomorphism of the equal direction structures,
and
\[
\iso{\tau_{EP} \cup \alpha}
{\pair{\fs{EP}(L_1,G_1)}{G_1}}
{\pair{\fs{EP}(L_2,G_2)}{G_2}}.
\]

Since $\fs{Ep}(L_i,G_i)$ is dense in $\barL_i$, extending $\tau_{EP}$
to $\fnn{\tau}{\barL_1}{\barL_2}$
gives the desired monotonic bijection $\tau$.  Proposition~\ref{UNIQ}
guarantees the uniqueness of $\tau$.
\newpage

%

\section{Locally moving groups}
\label{LMG}

In this section we prove the Expressibility Theorem for Locally
Moving Groups. This theorem appears also in \cite{Ru}.

\begin{defn} \label{}
\begin{rm}
Let $B$ be a Boolean algebra.

(1) A subset $D \subseteq B$ is {\em dense} in $B$ if for every
$a \in B - \sngltn{0}$ there exists $d \in D - \sngltn{0}$
such that $d \leq a$.

(2) Let $a \in B$.
$B \rest a \eqdf \setm{b \in B}{b \leq a}$.

(3) Let $g \in \fs{Aut}(B)$ and $a \in B$.
Then $g \rest a \eqdf g \rest (B \rest a)$.

(4) Suppose $B$ is complete, and let $g \in \fs{Aut}(B)$.

$\,\,\,\,\,\,\,\,\,\,$(a) Let $\fs{fix}(g) \eqdf \sum \setm{a \in B}{g
\rest a = \fs{Id}}$.

$\,\,\,\,\,\,\,\,\,\,$(b) Let $\fs{var}(g) \eqdf - \fs{fix}(g)$.

$\,\,\,\,\,\,\,\,\,\,$(c) For $G \leq \fs{Aut}(B)$, let
$\fs{Var}(B,G) \eqdf \setm{\fs{var}(g)}{g \in G}$.

(5) Let $G$ a group and let $f,g \in G$.

$\,\,\,\,\,\,\,\,\,\,$(a) $[f,g] \eqdf f g f\inverse g\inverse = f^g
f\inverse$.

$\,\,\,\,\,\,\,\,\,\,$(b) $C_G (f) \eqdf \{ h \in G \mid [h,f] = \fs{Id}\}$.
We abbreviate $C_G(f)$ by $C(f)$.

\end{rm}
\end{defn}

We consider in detail the special case in which
 $B = \barB(N)$ for some nearly ordered structure $N$ and
$g \in \fs{Aut}(N) \leq \fs{Aut}(\barB(N))$ (see
Proposition~\ref{BN}).
For $a \in \barB(N)$, the expression ``$g\rest a = \fs{Id}$'' in (4a) has
an equivalent meaning in which $a \subseteq \barN$ and \fs{Id} is
the identity map on $a$.  From this second point of view, $\fs{fix}(g)$
is the largest regular open subset $U$ of $\barN$ such that $g\rest U
= \fs{Id}$ because
$\sum \setm{a_i}{i \in I} = \fs{int}(\fs{cl}(\bigcup_{i \in I} a_i))$.
Let $\fs{fps}_{\barN}(g) \eqdf \setm{x \in \barN}{g(x) = x}$,
the fixed point set of $g$ in $\barN$, which is a closed subset of
$\barN$.

\begin{prop}
\label{VARBN}
Let $N$ be a nearly ordered structure, and let $g \in \fs{Aut}(N)
\leq \fs{Aut}(\barB(N))$.  Then

{\rm (1)} $\fs{fix}(g) = \fs{int}(\fs{fps}_{\barN}(g))$,  the union of
the interiors of the nonempty connected components of
$\fs{fps}_{\barN}(g)$.

{\rm (2)}  $\fs{var}(g) = \fs{int}(\fs{cl}(\fs{supp}(g)))$,  the union of
the interiors of the connected components of $|\barN| - \fs{fix}(g)$.

{\rm (3)} For $h \in \fs{Aut}(N)$, $\fs{fix}(g^h) = h(\fs{fix}(g))$ and
$\fs{var}(g^h) = h(\fs{var}(g))$.
\end{prop}

\noindent
{\bf Proof }
We have already established (1).  For (2), note that
\[
\fs{var}(g) = -\fs{fix}(g) = \fs{int}(|\barN| - \fs{fix}(g))  =
\fs{int}(\fs{cl}(\fs{supp}(g))).
\]
(3) is clear.
\hfill\qed
\medskip

In particular, the present definition of $\fs{var}(g)$ has as a special
case the definition given in Section \ref{OVERVIEW}.

\begin{defn} \label{}
\begin{rm}
Let $\pair{N}{G}$ be a nearly ordered permutation group.
Let
\[
\fs{Var}(N,G) \eqdf \fs{Var}(\barB(N),G) =
\setm{\fs{var}(g)}{g \in G} \subseteq \barB(N),
\]
a $G$-invariant subset of $\barB(N)$.
\end{rm}
\end{defn}

For $\pair{\fs{ED}(L)}{G}$, and for any $g \in G - \fs{Opp}(G)$,
$\fs{fps}_{\barN}(g)$ is a singleton and $\fs{var}(g) = \barL$.
$\fs{Var}(\fs{ED}(L,G)$ thus coincides with $\fs{Var}(L,\fs{Opp}(G))$
except that $\barL$ is adjoined as a (possibly new) member if
$\fs{Opp}(G) \neq G$.
Similar remarks apply to $\pair{\fs{EO}(L)}{G}$, with
$\fs{fps}_{\barN}(g)$ a doubleton and $\fs{var}(g) = \barC$ when $g
\in G - \fs{Opp}(G)$.

\medskip

We return to general Boolean algebras.
Note that  the following definition of a locally moving subgroup $G$
of $\fs{Aut}(B)$ specializes for $G \leq \fs{Aut}(N) \leq
\fs{Aut}(\barB(N))$ to the definition in \ref{FIRSTDEFLM}.

\begin{defn}
\label{}
\begin{rm}\

(1) Let $B$ be an atomless complete BA and
$G \leq \fs{Aut}(B)$.
$G$ is a {\em locally moving subgroup} of $\fs{Aut}(B)$
if for every $a \in B - \sngltn{0}$ there
exists $g \in G - \sngltn{\fs{Id}}$ such that $\fs{var}(g) \leq a$;
that is, if $\fs{Var}(B,G)$ is dense in $B$.
The pair $\pair{B}{G}$ is then called a {\em local movement system}.

(2) A group $G$ is a {\em locally moving group}
if there is an atomless complete BA $B$ such that
$G$ is isomorphic to a locally moving subgroup of $\fs{Aut}(B)$.
\end{rm}
\end{defn}

Rubin showed in \cite{Ru} that there exists a sentence $\varphi_{\fss{LM}}$
in the first order language of group theory such that for every group $G$:
\[
G \models \varphi_{\fss{LM}} \, \mbox{ iff \,$G$ is a locally moving group}.
\]
However, we shall not make use of that in this paper.

\begin{theorem}
\label{ET}
{\rm (The Expressibility Theorem)} \
There are first order formulas $\varphi_{Eq}(x,y)$,
$\varphi_{\leq}(x,y)$, and
$\varphi_{Ap}(x,y,z)$ in the language $\CL^{\fss{GR}}$
of groups such that for every atomless complete Boolean algebra $B$
and locally
moving subgroup $G \leq \fs{Aut}(B)$ the following holds:
For all $f, g, h \in G$,
\begin{eqnarray*}
&  G \models \varphi_{Eq} [f,g] \mbox{ \ iff \ }
    \fs{var}(f) = \fs{var}(g), \\*
&  G \models \varphi_{\leq}[f,g] \mbox{ \ iff \ }
    \fs{var}(f) \leq \fs{var}(g), \\*
&  G \models \varphi_{Ap} [f,g,h] \mbox{ \ iff \ }
     f(\fs{var}(g)) = \fs{var}(h).
\end{eqnarray*}

\end{theorem}

The Expressibility Theorem, restated in Theorem~\ref{INTERPVA} in
terms of interpretations, will provide the basis for our reconstruction
theorems.
The proof of the Expressibility Theorem will be completed after
Proposition~\ref{VLE}. For applications to nearly ordered
permutation groups $\pair{N}{G}$ we need this theorem only for the
special case $B = \barB(N)$.  There the proof is far more transparent.
For a treatment of just the special case, the suggested approach is to
read Definition~\ref{DEFVLE} and then read from Lemma~\ref{NGEQ}
to the end of the proof of Theorem~\ref{ET}.  (The results through
Lemma~\ref{COMMUTATOR} are clear in the special case.)   For {\em
linear} permutation groups, it is possible to skip even
Lemma~\ref{NGEQ}; see Note~\ref{NOTEGUIDE}.

We begin the development of the general case.
It may be helpful to keep in mind that every BA is isomorphic to a
field of sets.
First we establish a few results not involving local movability.

\begin{lemma}
\label{DOT}
Let $B$ be a BA.  Let $g \in \fs{Aut}(B)$ and $a \in B$.

{\rm (1)} If $g(a) \neq a$,  then
there exists $0 \neq b \leq a$ such that $g(b) \cdot b = 0$.

{\rm (2)} Let $B$ be complete.

$\,\,\,\,\,\,\,\,\,\,${\rm (a)} If $g(b) \cdot b = 0$, then
$b \leq \fs{var}(g)$.

$\,\,\,\,\,\,\,\,\,\,${\rm (b)} If $g(a) \neq a$, then $a \cdot
\fs{var}(g) \neq 0$.

\end{lemma}

\noindent
{\bf Proof }
(1) If $a-g(a) \neq 0$, then $a-g(a)$ serves as the required
$b$.  Otherwise $g(a) \geq a$, and by hypothesis $g(a) \neq a$.
Let $c = g(a) -a$ and $b = g^{-1} (c) \leq a$.
Then $b$ is as required.

(2a) If $g(b) \cdot b = 0$, then $b \cdot \fs{fix}(g) = 0$,
so $b \leq \fs{var}(g)$.

(2b) This follows from (1) and (2a).
\hfill\qed

\begin{lemma}
\label{}
Let $B$ be a BA, $A \subseteq B$, and $b \in B$.  Then the following
are equivalent:

{\rm (1)}  For every upper bound $u$ of $A$, $b \leq u$.

{\rm (2)} For every $0 \neq b_1 \leq b$ there exists $a \in
A$ such that $b_1 \cdot a \neq 0$.

\end{lemma}

\noindent
{\bf Proof }
(1) $\Rightarrow$ (2) Let $0 \neq b_1 \leq b$,
and suppose by way of contradiction that for every $a \in A$,
$b_1 \cdot a = 0$.
Then $-b_1$ is an upper bound of $A$.
Hence $b \leq -b_1$, a contradiction.

(2) $\Rightarrow$ (1)  Let $u$ be an upper bound of $A$.
Let $c = b-u$.  Then $c \cdot u = 0$, so $c \cdot a = 0$ for all $a \in
A$.  Since $c \leq b$, $c = 0$.  Hence $b \leq u$.
\hfill\qed

\begin{prop}
\label{VSUP}
Let $B$ be a complete BA. Let $b \in B$ and $g \in \fs{Aut}(B)$.

{\rm (1)}  $b \leq \fs{fix}(g)$ iff $g \rest b = \fs{Id}$.

{\rm (2)}  $b \leq \fs{var}(g)$ iff for all $0 \neq c \leq b,\ g \rest c
\neq \fs{Id}$.

{\rm (3)} $\fs{var}(g) = \sum \setm{a \in B}{a \cdot g(a) = 0}$.

\end{prop}

\noindent
{\bf Proof }
(1) Let $b \leq \fs{fix}(g)$. Suppose $g(c) \neq c$ for some $c \leq
b$.
By Lemma~\ref{DOT}, $c \cdot \fs{var}(g) \neq 0$, and thus $c
\not\leq \fs{fix}(g)$.
This contradiction establishes one direction of (1a), and the other
direction is obvious.

(2) We have $b \leq \fs{var}(g) =
-\fs{fix}(g)$ iff $b \cdot \fs{fix}(g) = 0$ iff
for all $0 \neq c \leq b$, $g \rest c \neq \fs{Id}$ by (1).

(3) Let $A = \setm{a \in B}{g(a) \cdot a = 0}$.
For every $a \in A$, $a \leq \fs{var}(g)$,
and thus $\sum A \leq \fs{var}(g)$.
Conversely, we show that $\fs{var}(g) \leq \sum A$.
Let $0 \neq b_1 \leq \fs{var}(g)$.
Then $g \rest b_1 \neq \fs{Id}$.
Pick $c \leq b_1$ such that $g(c) \neq c$, and then use (1)
to pick $0 \neq a \leq c$ such that $g(a) \cdot a = 0$.
Thus $a \in A$ and $b_1 \cdot a \neq 0$.  By the lemma,
$\fs{var}(g) \leq \sum A$.
\hfill\qed

\begin{prop}
\label{VARCONJ}
Let $B$ be a complete BA, and let $f,g,h \in \fs{Aut}(B)$.

{\rm (1)} $\fs{var}(f\inverse) = \fs{var}(f)$ and
$\fs{var}(f^h) = h(\fs{var}(f))$.

{\rm (2)} $\fs{var}(fg) \leq \fs{var}(f) + \fs{var}(g)$.

{\rm (3)} If $\fs{var}(f) \cdot \fs{var}(g) = 0$, then

$\,\,\,\,\,\,\,\,\,\,${\rm (a)} $\fs{var}(fg) = \fs{var}(f) + \fs{var}(g)$,

$\,\,\,\,\,\,\,\,\,\,${\rm (b)} $[f,g] = \fs{Id}$.

{\rm (4)} Let $G \leq \fs{Aut}(B)$.  Then $\fs{Var}(B,G)$ is a
$G$-invariant subset of $\barB(N)$.
\end{prop}

\noindent
{\bf Proof }

(1) This holds for {\em var} because it holds for {\em fix}.

(2) By Proposition~\ref{VSUP} $\fs{fix}(f) \cdot \fs{fix}(g) \leq
\fs{fix}(fg)$,
which is equivalent to (2).

(3) Assume $\fs{var}(f) \cdot \fs{var}(g) = 0$.

$\,\,\,\,\,\,\,\,\,\,$(a) It suffices to show that
$\fs{fix}(fg) \leq \fs{fix}(f) \cdot \fs{fix}(g)$.
Let $b \leq \fs{fix}(fg)$, so that $(fg) \rest b = \fs{Id}$.
Then $g \rest b = f\inverse \rest b$.
If $g \rest b = \fs{Id}$, then
$b \leq \fs{fix}(f\inverse) \cdot \fs{fix}(g) =
\fs{fix}(f) \cdot \fs{fix}(g)$, as desired.
If not, then $g(c) \cdot c = 0$ for some $0 \neq c \leq b$,
and then by Lemma~\ref{DOT} we have $c \leq \fs{var}(g) \cdot
\fs{var}(f\inverse)$, a contradiction.

$\,\,\,\,\,\,\,\,\,\,$(b) Suppose there exists $a \in B$ such that
$f\inverse g\inverse fg(a) \neq a$.
Pick $0 \neq b \leq a$ such that
$f\inverse g\inverse fg(b) \cdot b = 0 =
g\inverse f\inverse gf(b)$.
We may assume with no loss of generality that $g(b) \neq b$.
Pick $0 \neq c \leq b$ such that $g(c) \cdot c = 0$, making $c \leq
var(g)$, and thus also making $g(c) \leq g(\fs{var}(g)) = \fs{var}(g)$.
Since $\fs{var}(f) \cdot \fs{var}(g) = 0$,
$c \leq \fs{fix}(f) = \fs{fix}(f\inverse)$ and $g(c) \leq \fs{fix}(f)$.
Therefore $f\inverse g\inverse fg(c) = c$, which is a contradiction
since $0 \neq c \leq b$ and $f\inverse g\inverse fg(b) \cdot b = 0$.

(4)  This follows from (1).
\hfill\qed

\begin{lemma}
\label{ALLOFF}
Let $B$ be a complete BA and let $g_1,\dots,g_n \in \fs{Aut}(B)$.
Let $0 \neq a \leq \displaystyle {\prod_{i=1}^n} \fs{var}(g_i )$.
Then there exists $0 \neq b \leq a$ such that
$b \cdot \displaystyle {\sum_{i=1}^n} g_i (b) = 0$.
\end{lemma}

\noindent
{\bf Proof }
The proof is by induction on $n$, using Proposition~\ref{VSUP}(2)
and Lemma \ref{DOT}(1).
\hfill\qed

\begin{defn} \label{DEFVLE}
\begin{rm}
Let $\pair{B}{G}$ be a local movement system.
For $a \in B$ let $\fs{vle}_G(a)
\eqdf \setm{g \in G}{\fs{var}(g) \leq a}$.
We abbreviate $\fs{vle}_G(a)$ by $\fs{vle}(a)$.
\end{rm}
\end{defn}

For $f \in G$, we want to capture the set $\fs{vle}(\fs{var}(f))$ in
$\CL^{\fss{GR}}$.  This will be accomplished in Proposition~\ref{VLE}.

\begin{lemma}
\label{COMMUTATOR}
Let $\pair{B}{G}$ be a local movement system.

{\rm (1)} Let $f \in G$ and $0 \neq a \leq \fs{var}(f)$.
Then there exists $g \in \fs{vle}(a)$ such that $[g,f] \neq \fs{Id}$.

{\rm (2)} Let $0 \neq a \in B$ and let $n \geq 1$.  Then there exists
$h \in \fs{vle}(a)$ such that $h^n \neq \fs{Id}$.

{\rm (3)} Let $f,g \in G$ and $0 \neq a \leq \fs{var}(f) \cdot
\fs{var}(g)$.
Then there exists $h \in \fs{vle}(a)$ such that $[f^h ,g] \neq \fs{Id}$.

\end{lemma}

\noindent
{\bf Proof }
(1) Let $0 \neq a \leq \fs{var}(f)$.
Pick $0 \neq b \leq a$ such that $b \cdot f(b) = 0$, and pick $g \in
\fs{vle}(b) - \sngltn{\fs{Id}}$.
Since $f(\fs{var}(g)) \neq \fs{var}(g), [g,f] \neq \fs{Id}$.

(2) We prove by induction on $n$ that there are
$h \in \fs{vle}(a)$ and $0 \neq b \leq a$
such that $b, h(b),\dots ,h^n (b)$ are pairwise disjoint.
For $n=1$, let $h \in \fs{vle}(a) - \sngltn{\fs{Id}}$.
Then by Lemma~\ref{DOT} there exists $0 \neq b \leq \fs{var}(h)$
such that
$b \cdot h(b) = 0$.  Suppose the claim is true for $n$.  Let $h
\in \fs{vle}(a)$ and $0 \neq b \leq a$ be as assured by the
induction hypothesis.  If $h^{n+1} \rest b \neq \fs{Id}$,
pick $0 \neq c \leq b$ such that $c \cdot h^{n+1} (c) = 0$,
and $h$ and $c$ are as required.
Suppose now that $h^{n+1} \rest b = \fs{Id}$.
Pick $k \in \fs{vle}(b) - \sngltn{\fs{Id}}$
and $0 \neq c \leq b$ such that $c \cdot k(c) = 0$.
Let $g = kh$.  Then for $i \leq n$, $g^i (c) = h^i (c) \leq h^i (b)$
since $b,\dots,h^n(b)$ are pairwise disjoint.
Hence $c,\dots ,g^n (c)$ are pairwise disjoint.
Moreover $g^{n+1} (c) \cdot c = 0$ because
$g^{n+1} (c) = (kh)^{n+1} (c) = kh (kh)^n (c) = kh(h^n (c)) =
k(c)$.   Since $g^{n+1} (c) \leq
b$, we have $g^{n+1} (c) \cdot g^i (c) = 0$ for $i=1, \dots ,n$.

(3) Let $0 \neq a \leq \fs{var}(f) \cdot \fs{var}(g)$.
If $[f,g] \neq \fs{Id}$, then $h = \fs{Id}$ is as required.
Thus we may assume that $[f,g] = \fs{Id}$.
Use Lemma~\ref{ALLOFF} to pick $0 \neq b \leq a$ such that
$b \cdot (f(b ) + g(b )) = 0$.

{\bf Case 1:}
$b \cdot \fs{var}(fg) = 0$.  Then pick
$h \in \fs{vle}(b)$ such that $h^2 \neq \fs{Id}$, and pick
$c \leq b$  such that $h^2(c) \neq c$.  Then since
$f(c) \cdot \fs{var}(h) =0 = c\cdot \fs{var}(fg)$ we have
\[
gf^h (h(c)) = ghfh^{-1} h(c) = ghf(c) = gf (c) = fg(c) = c;
\]
whereas since $\fs{var}(h) \cdot g(b) = 0 = h(c) \cdot \fs{var}(fg)$
we have
\[
f^h g(h(c)) = hfh^{-1} gh(c) = hfgh(c) = h^2(c).
\]
Hence $g$ and $f^h$ do not commute.

{\bf Case 2:} $b \cdot \fs{var}(fg) \neq 0$.  Then pick
$0 \neq b_1 \leq b$ such that $b_1 \cdot fg(b_1) = 0$.
Pick $h \in \fs{vle}(b_1) - \sngltn{\fs{Id}}$, and pick $c \leq b_1$
such that
$h(c) \neq c$.  Then as above
\[
gf^h (h(c)) = gf(c);
\]
whereas
\[
f^h g(h(c)) =  hfgh(c) = fgh(c) = gf(h(c)).
\]
Again $g$ and $f^h$ do not commute.
\hfill\qed
\medskip

The next lemma deals with automorphisms of an arbitrary BA,
and its proof does not use any of the previous results.

\begin{lemma}
\label{NGEQ}
Let $B$ be any BA and let $G = \fs{Aut}(B)$.
Let $k_0 ,\dots ,k_n \in {\Z}$, $f \in G$, and
$a\in B$.  Suppose that $f^{k_0} (a),\dots ,f^{k_n} (a)$ are
pairwise disjoint.  Then for all $h_1 ,\dots ,h_n \in C(f)$
and $0 \neq b \leq a$,
$$ \sum_{i=1}^n h_i (a) \not\geq \sum_{i=0}^n f^{k_i} (b) .$$
\end{lemma}

\noindent
{\bf Proof } The proof is by induction on $n$.  First we treat $n=1$.
If $f^{k_0} (b) +
f^{k_1} (b) \leq h_1 (a)$, then $0 \neq b \leq f^{-k_0} (h_1
(a)) \cdot f^{-k_1} (h_1 (a))$.  Then since $h_1 \in C(f)$,
$0 \neq f^{-k_0} (a) \cdot f^{-k_1} (a)$.
Apply $f^{k_0 +k_1}$ to both sides.
Then $f^{k_1} (a) \cdot f^{k_0} (a) \neq 0$, a contradiction.

Suppose the claim is true for $n$.  Let $k_0 ,\dots ,k_{n+1}$,
$h_1 ,\dots ,h_{n+1}$, and $a$ and $b$ be as in the lemma.  If
$h_{n+1} (a) \cdot \displaystyle {\sum_{i=0}^{n+1}} f^{k_i} (b)
= 0$, the induction hypothesis would give
$$
\sum_{i=1}^{n+1} h_i (a) \not\geq \sum_{i=0}^{n+1} f^{k_i}(b).
$$
Suppose otherwise, and pick $j$ such that
$h_{n+1}(a) \cdot f^{k_j}(b) \neq 0$.
We may assume that $j = n+1$, so that $c \eqdf h_{n+1}(a) \cdot
f^{k_{n+1}}(b) \neq 0$.
Let $b' = f^{-k_{n+1}} (c)$.
Since $0 \neq b' \leq b \leq a$, the $n=1$ case tells us that
$h_{n+1} (a) \not\geq f^{k_0} (b' ) + f^{k_{n+1}} (b' ) =
f^{k_0} (b' ) + c$.  Since $h_{n+1} (a) \geq c$,
$h_{n+1}(a) \not\geq f^{k_0} (b')$.

For $0 \leq i \leq n$, we  define
$b_i$ inductively so that
\[
b \geq b' \geq b_0 \geq b_1 \geq \dots \geq b_n \neq 0 \mbox{ \
and \ } h_{n+1} (a) \cdot f^{k_i} (b_i) = 0.
\]
Let $b_0 = f^{-k_0} (f^{k_0} (b' ) - h_{n+1} (a)) \neq 0$.
Then $b_0$ is as specified above.  Suppose $b_i \neq 0$  has been
defined as specified, with $i<n$.  We have $h_{n+1}(a) \not\geq
f^{k_{i+1}}(b_i ) + f^{k_{n+1}}
(b_i)$, and $h_{n+1}(a) \geq c = f^{k_{n+1}}(b') \geq f^{k_{n+1}}(b_i)$.
Let $b_{i+1} = f^{-k_{i+1}} (f^{k_{i+1}} (b_i )-h_{n+1} (a)) \neq 0$.
Then $b_{i+1}$ is as specified.

We obtain that $h_{n+1}(a) \cdot
\displaystyle {\sum_{i=0}^n} f^{k_i}(b_n ) = 0$.
By the external induction hypothesis,
$$\sum_{i=1}^n h_i(a) \not\geq \sum_{i=0}^n f^{k_i}(b_n)$$
and so
$$
\sum_{i=1}^{n+1} h_i(a) \not\geq \sum_{i=0}^n f^{k_i} (b_n).
$$
A fortiori
$$
\hspace{4.8cm}
\sum_{i=1}^{n+1} h_i (a) \not\geq \sum_{i=0}^{n+1} f^{k_i} (b).
\hspace{4.5cm}
\qed
$$

\begin{lemma}
\label{ALMOST}
Let
\begin{eqnarray*}
& \varphi (f,\hat{f} ) \ \eqqdf \
\forall g \mbox{\Large(} \mbox{\large(}[g,f] \neq \fs{Id}
\mbox{\large)} \rightarrow \\*
& \mbox{\large(}\exists h_1,h_2 \in C(\hat{f} )
\mbox{\large)} \mbox{\large(} ([g,h_1 ,h_2 ] \neq \fs{Id})
\wedge
([[g,h_1 ,h_2 ], \hat{f} ] = \fs{Id}) \mbox{\large)} \mbox{\Large)}.
\end{eqnarray*}
For every local movement system $\langle B,G\rangle$ and all
$f,\hat{f} \in G$:
\[
\fs{var}(f) \cdot \fs{var}(\hat{f}) = 0  \ \Rightarrow \  G \models
\varphi [f,\hat{f} ] \ \Rightarrow \
\fs{var}(f) \cdot \fs{var}(\hat{f}^{12} ) = 0.
\]
\end{lemma}

\noindent
{\bf Proof }
Let $\pair{B}{G}$ be a local movement system and let $f,\hat{f} \in
G$.

For the first implication, suppose that
$\fs{var}(f) \cdot \fs{var}(\hat{f}) = 0$.
Let $g \in G$ be such that $[g,f] \neq \fs{Id}$.
Then $\fs{var}(f) \cdot \fs{var}(g) \neq 0$ by
Proposition~\ref{VARCONJ}.
We use Lemma~\ref{DOT} to pick
$0 \neq a \leq \fs{var}(f)$
such that $g(a) \cdot a = 0$, and we pick $h_1 \in
\fs{vle}(a) - \sngltn{\fs{Id}}$.
Then we pick
$0 \neq b \leq \fs{var}(h_1 )$ such that
$h_1(b) \cdot b = 0$, and pick $h_2 \in
\fs{vle}(b) - \sngltn{\fs{Id}}$.
We have
\[
0 \neq \fs{var}(h_2) \leq b \leq \fs{var}(h_1) \leq a \leq \fs{var}(f).
\]
Hence $h_1,h_2 \in C(\hat{f})$.
Since $g(a) \cdot a = 0, h_1^g$ commutes with $h_1$ and $h_2$.  We
have
\begin{eqnarray*}
& \fs{var}([[g,h_1],h_2]) = \fs{var}(h_2^{[g,h_1]} \cdot h_2\inverse) =
\fs{var}(h_2^{h_1^gh_1\inverse} \cdot h_2\inverse) = \\*
& \fs{var}(h_2^{h_1\inverse} \cdot h_2\inverse) \leq
\fs{var}(f).
\end{eqnarray*}
Therefore $[g,h_1,h_2] \in C(\hat{f})$.  Moreover, since
$\fs{var}(h_2^{h_1\inverse}) \cdot \fs{var}(h_2\inverse) = 0$,
Proposition~\ref{VARCONJ} guarantees that $[g,h_1,h_2] \neq
\sngltn{\fs{Id}}$.
We have shown that $G \models \varphi [f,\hat{f} ]$.

For the second implication, suppose that $a \eqdf \fs{var}(f) \cdot
\fs{var}(\hat{f}^{12}) \neq 0$.
Then for $i = 1,\ldots,4$,
$a \leq \fs{var}(\hat{f} ^i )$.
We use Lemma~\ref{ALLOFF} to pick
$0 \neq b \leq a$ such that
$b \cdot \displaystyle {\sum_{i=1}^4} \hat{f}^i (b) = 0$.
Hence $\{ \hat{f}^i (b) \mid 0 \leq i \leq 4\}$ is pairwise disjoint.
Using Lemma~\ref{COMMUTATOR}(1), we pick $g \in \fs{vle}(b)$
such that $[g,f] \neq \fs{Id}$.

Let $h_1 ,h_2 \in C(\hat{f} )$ be such that
$g_2 \eqdf [[g,h_1 ],h_2 ] \neq \fs{Id}$.
We show that $[g_2 ,\hat{f}] \neq \fs{Id}$.
Now
\begin{eqnarray*}
& g_2 = [g \cdot (g^{-1} )^{h_1} ,h_2 ] = g(g^{-1} )^{h_1} \cdot
((g(g^{-1} )^{h_1} )^{-1} )^{h_2} = \\*
& g(g^{-1} )^{h_1} \cdot (g^{h_1} g^{-1} )^{h_2} = g(g^{-1}
)^{h_1} g^{h_2 h_1} (g^{-1} )^{h_2} .
\end{eqnarray*}
So $0 \neq \fs{var}(g_2 ) \leq b + h_1(b) + h_2 h_1(b) +
h_2(b)$.
Pick $h \in \{ \fs{Id}, h_1 ,h_2 h_1 ,h_2 \}$  such that
$c \eqdf \fs{var}(g_2) \cdot h(b) \neq 0$.
Since $\{ \hat{f} ^i (b) \mid 0 \leq i \leq 4\}$ is pairwise disjoint
and $h \in C(\hat{f} )$,
$\{ \hat{f}^i (h(b)) \mid 0 \leq i \leq 4\}$ is also pairwise disjoint.

Suppose by way of contradiction that $[g_2 ,\hat{f} ] = \fs{Id}$.
Then for $i = 0,\ldots,4$,
$(g_2 )^{\hat{f}^i} = g_2$.  Hence
\begin{eqnarray*}
& \displaystyle{\sum_{i=0}^4 \hat{f} ^i (c) \leq
   \sum_{i=0}^4 \hat{f}^i (\fs{var}(g_2 )) =
   \sum_{i=0}^4 \fs{var}(g_2^{\hat{f}^i} ) =
   \fs{var}(g_2) \leq} \\
& b + h_1(b) + h_2 h_1(b) + h_2(b) \leq \\
& h^{-1} (h(b)) + h_1 h^{-1}(h(b)) + h_2 h_1 h^{-1} (h(b)) +
     h_2 h^{-1} (h(b)).
\end{eqnarray*}
Since $h^{-1}$, $h_1 h^{-1}$,
$h_2 h_1 h^{-1}$, and $h_2 h^{-1}$ commute with $\hat{f}$ and $0
\neq c \leq h(b)$,
this contradicts Lemma~\ref{NGEQ}.  (Alternately, see the following
note.)  Therefore  $[g_2,\hat{f}] \neq \fs{Id}$.
We have shown that $G \not\models \varphi [f,\hat{f} ]$.
\hfill\qed

\begin{note}
\label{NOTEGUIDE}
\begin{rm}
For the special case of a local movement system $\pair{\barB(L)}{G}$,
with $\pair{L}{G}$ a linear permutation group, it is illuminating to
see how to avoid the above reference to Lemma~\ref{NGEQ} by using
additional care when choosing $b$ in the proof of the second
implication. For $\hat{f} \in G$ and $x \in \fs{supp}(\hat{f})$, the
convex hull of
$\setm{\hat{f}^n(x)}{n \in {\Z}}$ is called a {\em supporting interval}
of $\hat{f}$.  Already we had, among other things,
$0 \neq b \leq \fs{var}(\hat{f})$ and $b\cdot \hat{f}(b) \neq 0$; now
we arrange also that $b$ be an $\barL$-interval contained in a single
supporting interval of $\hat{f}$ (which guarantees that $b$ cannot
contain both $x$ and $\hat{f}(x)$ for any $x \in \barL$).  Then since
$h \in C(\hat{f})$, the same is true of $h(b)$, though the supporting
interval of $\hat{f}$ may be different; and likewise for
$h^{-1} (h(b))$, $h_1 h^{-1}(h(b))$, $h_2 h_1 h^{-1} (h(b))$, and
$h_2 h^{-1} (h(b)$.  For $0 \neq c \leq h(b)$ and $x \in c$,
$\displaystyle{\sum_{i=0}^4 \hat{f} ^i (c)}$ contains
$\setm{\hat{f}^i(x)}{0\leq i \leq 5}$, whence the contradiction.

\end{rm}
\end{note}

\begin{defn}
\label{}
\begin{rm}
Let $\pair{B}{G}$ be a local movement system and $f \in G$.

(1) Let $D (f) \eqdf \{ \hat{f} \in G \mid G \models
\varphi [ f,\hat{f} ]\}$.

(2) Let $V(f) \eqdf C(\{ \hat{f} ^{12} \mid \hat{f} \in D (f)\} )$.
\end{rm}
\end{defn}

\begin{prop}
\label{VLE}
Let $\pair{B}{G}$ be a local movement system.
For every $f \in G$, \,$V(f) = \fs{vle}(\fs{var}(f))$.
\end{prop}

\noindent
{\bf Proof }
First we show that
$\fs{vle}(\fs{var}(f)) \subseteq V(f)$.
Let $g \in \fs{vle}(\fs{var}(f))$. Let $\hat{f} \in D(f)$.
Then $\fs{var}(f) \cdot \fs{var}(\hat{f} ^{12} )= 0$,
and hence $\fs{var}(g) \cdot \fs{var}(\hat{f} ^{12} ) = 0$.
Hence $[g, \hat{f} ^{12} ] = \fs{Id}$.

Next we show that
$V(f) \subseteq \fs{vle}(\fs{var}(f))$.  Let $g \in G -
\fs{vle}(\fs{var}(f))$.
Then $a \eqdf \fs{var}(g) - \fs{var}(f) \neq 0$.
We use Lemma~\ref{COMMUTATOR}(2) to pick $\hat{f} \in
\fs{vle}(a)$ such
that $\hat{f} ^{12} \neq \fs{Id}$.
Let $0 \neq a_1 \eqdf \fs{var}(\hat{f}^{12}) \leq
\fs{var}(\hat{f}) \leq \fs{var}(g)$.  Using
Lemma~\ref{COMMUTATOR}(3)
we pick $h \in \fs{vle}(a_1)$
such that $[(\hat{f} ^{12} )^h ,g] \neq \fs{Id}$.
Now $\fs{var}(f) \cdot \fs{var}(\hat{f}) = 0$
and thus $\fs{var}(f) \cdot \fs{var}(\hat{f}^h) = 0$.
It follows that $\hat{f} ^h \in D (f)$.  Hence $g \notin V(f)$.
\hfill\qed
\bigskip

\noindent
{\bf Proof of Theorem~\ref{ET} }
Let $\varphi_{\leq}(f,g) \ \eqqdf \ V(f) \subseteq V(g)$.
Clearly, $\varphi_{\leq}(f,g)$ can be written as a first order formula
in $\CL^{\fss{GR}}$.
If $\fs{var}(f) \leq \fs{var}(g)$, then by the above proposition,
$G \models \varphi_{\leq}[f,g]$.
If on the other hand, $\fs{var}(f) \not\leq \fs{var}(g)$, then since
$G$ is a locally moving subgroup of $\fs{Aut}(B)$, there exists $h \in
\fs{vle}(\fs{var}(f)) - \fs{vle}(\fs{var}(g))$.
By the above proposition, $V(f) \not\subseteq V(g)$.
So $\varphi_{\leq}(f,g)$ is as required.

Let $\varphi_{Eq} (f,g) \ \eqqdf \
\varphi_{\leq}(f,g) \wedge \varphi_{\leq}(g,f)$
and $\varphi_{Ap} (f,g,h) \ \eqqdf \ \varphi_{Eq} (g^f ,h)$.
Clearly $\varphi_{Eq}$ and $\varphi_{Ap}$ are as required.
\hfill\qed
\medskip

The Expressibility Theorem yields

\begin{cor}\label{DSJNT}
There is a first order formula $\varphi_{Dsjnt}^{Gr}(x,y)$ in
$\CL^{\fss{GR}}$
 such that for every local movement system $\pair{B}{G}$ the
following holds:
For all $f, g \in G$,
\begin{eqnarray*}
&  G \models \varphi_{Dsjnt}^{Gr} [f,g] \mbox{ \ iff \ }
    \fs{var}(f) \cdot \fs{var}(g) =0. \\*
\end{eqnarray*}

\end{cor}

\vspace{-4 mm}
\noindent
{\bf Proof }
The proof, though very easy, does involve expressing the disjointness
of $\fs{var}(f)$ and $\fs{var}(g)$ in terms of the {\em poset}
$\pair{\fs{Var}(G)}{\leq}$, whereas disjointness is a BA notion.  This
brings up our first use of the following proposition, from which the
corollary follows immediately.
\hfill\qed

\begin{prop}
\label{PST}
For every first order formula $\varphi(\vecx)$ in the language
$\CL^{\fss{PST}} = \{ \leq \}$ of Boolean algebras,
there is a first order formula $\psi_{\varphi}(\vecx)$ in
$\CL^{\fss{PST}}$ such that for any atomless BA $B$, any dense
subset $D$ of $B$, and any $\veca \in D^{\fss{length}(\vecx)}$:
\[
B \models \varphi[\veca] \mbox{ \ iff \ }
\pair{D}{\leq} \models \psi_{\varphi}[\veca].
\]
In particular, there are formulas
$\varphi_{\fss{dsjnt}}(u,v)$, $\varphi_{\fss{lesum}}(u,v,w)$,
$\varphi_{\fss{sumle}}(u,v,w)$, and $\varphi_{\fss{cmp}}(u,v)$
that represent in the above sense the formulas
$u \cdot v = 0$, $u \leq v + w$, and $u + v \leq w$, and $v= -u$.
\end{prop}

\noindent
{\bf Proof for the special cases }
The point here is that not only must
$\psi_{\varphi}(\vecx)$ be in $\CL^{\fss{PST}}$,
but also the bound variables of $\psi_{\varphi}(\vecx)$
are to be quantified only over $D$.  As we need only the special cases
mentioned above, we forego proving the proposition and deal just
with the special cases.  Let
\newline
\smallskip
$\varphi_{\fss{dsjnt}}(u,v) \eqqdf
(\forall w_1,w_2)
\mbox{\Large(}\mbox{\large(}
(\bigwedge_{i=1}^2 (w_i \leq u \, \wedge \, w_i \leq
v)\mbox{\large)} \rightarrow
(w_1 = w_2)\mbox{\Large)}$,
\newline
\smallskip
$\varphi_{\fss{lesum}}(u,v,w) \eqqdf
(\forall x)\mbox{\Large(}\mbox{\large(}
\varphi_{\fss{dsjnt}}(x,v) \wedge
\varphi_{\fss{dsjnt}}(x,w)\mbox{\large)} \rightarrow
\varphi_{\fss{dsjnt}}(x,u)\mbox{\Large)}$,
\newline
\smallskip
$\varphi_{\fss{sumle}}(u,v,w) \eqqdf
u \leq w \ \wedge \ v \leq w$,
\newline
\smallskip
$\varphi_{\fss{cmp}}(u,v) \eqqdf
\varphi_{\fss{dsjnt}}(u,v) \,\wedge \,
(\forall w_1,w_2)
\mbox{\large(}
(\bigwedge_{i=1}^2 \varphi_{\fss{sumle}}(u,v,w_i) \rightarrow
(w_1 = w_2)\mbox{\large)}$.
\newline
\smallskip
It is left to the reader to verify that
these formulas are as required.
\hfill\qed
\medskip

For later use we mention two other properties of $\fs{Var}(B,G)$:

\begin{prop}
\label{BSUP}
Let $\pair{B}{G}$ be a local movement system, and let $b\in B$.
Then
$b = \sum \setm{\fs{var}(g)}{g\in G \mbox{ and } \fs{var}(g) \leq
b}$.
\end{prop}

\noindent
{\bf Proof }
$\fs{Var}(B,G)$ is a dense subset of an atomless complete BA.
\hfill\qed

\begin{defn}\label{DEFDDNS}
\begin{rm} $\:$
Let $B$ be a BA and let $T \subseteq B$.  We say that $T$ is {\em
doubly dense} in $B$ if
for all $a_1,a_2 \in B - \sngltn{0}$
there exists $t \in T$ such that $t \leq a_1 + a_2$ and
$t \cdot a_1, t \cdot  a_2 \neq 0$.
\end{rm}
\end{defn}

Note that if $T$ is doubly dense in $B$,
then $T$ is dense in $B$.
Also, $B$ is doubly dense in itself.

\begin{prop}\label{DDNS}
Let $\pair{B}{G}$ be a local movement system.  Then  $\fs{Var}(B,G)$
is doubly dense in $B$.
\end{prop}

\noindent
{\bf Proof }
Let $a_1,a_2 \in B - \sngltn{0}$.
Since $G$ is a locally moving subgroup, we may pick $g_i \in
\fs{vle}_G(a_i) - \sngltn{\fs{Id}}$, $i = 1,2$, and may arrange also
that  $\fs{var}(g_1) \cdot \fs{var}(g_2) = 0$.
Let $g = g_1 \circ g_2$.
By Proposition~\ref{VARCONJ} $\fs{var}(g) =\fs{var}(g_1) +
\fs{var}(g_2) $.  Hence $\fs{var}(g) \leq a_1 + a_2$ and
$\fs{var}(g) \cap a_i \neq 0$ for $i = 1,2$.
\hfill\qed
\newpage

%

\section{Model theory and interpretations}
\label{INTERPS}

Let $M$ be a structure in the sense of model theory.
$|M|$ denotes the universe of $M$, and $\CL(M)$ the language of
$M$.
$\CL(M)$ may include $n$-place relation symbols $R$ and $n$-place
function symbols $F$ (where $n \in \omega - \sngltn{0}$), and
individual constant symbols $c$.
If $P$ is any symbol of $\CL(M)$,
$P^M$ denotes the interpretation of $P$ in $M$.

The convention in model theory is that if $F$ is an $n$-place
function symbol in $\CL(M)$, then $\fnn{F^M}{|M|^n}{|M|}$.
We will often define ``functions'' $f$ on $|M|^n$ (where $n \in \omega
- \sngltn{0}$) which are  actually only {\em partial} functions, that
is, for which
$\fs{Dom}(f) \subseteq |M|^n$ and $\fs{Rng}(f) \subseteq |M|$.
When we do so, we regard the partial function $f$ as an abbreviation
for a pair consisting of

(1) The $n$-place function $f^{\rm tot}:|M|^n \raro |M|$
which extends $f$ and for which $f(a_1,\ldots,a_n) = a_1$ when
$(a_1,\ldots,a_n) \not\in \fs{Dom}(f)$, and

(2) The $n$-place relation $\fs{Dom}(f)$.
\newline
Thus, without adding or losing information, we have a structure in
the classical sense.

Sometimes we build a structure from others already defined, with
1-place relations to distinguish the original universes, and perhaps
adding new relations, functions, and constants:

\begin{defn}
\label{}
\begin{rm}
Let $M_1,\ldots,M_t$ be structures, and let
$\CL_i = \CL(M_i)$.
Suppose that these structures are compatible in the sense that for
any symbol $P \in \CL_i \cap \CL_j$, $P^{M_i}\rest (|M_i| \cap |M_j|)
= P^{M_j}\rest (|M_i| \cap |M_j|)$.
Let $U = \bigcup_{i = 1}^t |M_i|$.
Let $\tldR_1, \ldots, \tldR_k$ be relations on $U$,
let $\tldF_1, \ldots, \tldF_l$ be functions with $\tldF_j:
U^{n_j} \raro U$, and
let $\tldc_1, \ldots, \tldc_m \in U$.
We define the structure
\[
M \eqdf (M_1, \ldots, M_t;
\tldR_1, \ldots, \tldR_k, \tldF_1, \ldots, \tldF_l, \tldc_1, \ldots,
\tldc_m) {\mbox :}
\]

(1) $|M| = U$.

(2) $\CL(M) =
(\bigcup_{i = 1}^t \CL_i) \cup
\fsetn{R_1}{R_k} \cup
\fsetn{F_1}{F_l} \cup
\fsetn{c_1}{c_m} \cup
\fsetn{P_1}{P_t}$.

(3) For every $P \in \bigcup_{i = 1}^t \CL_i$,
$P^M = \bigcup \setm{P^{M_i}}{P \in \CL(M_i)}$, where the
convention about partial functions is employed when $P$ is a
function symbol.

(4) $R_i^M = \tldR_i$ for $i=1, \ldots, k$,
and similarly for the $F_i$'s and the $c_i$'s.

(5) $P_i^M = |M_i|$ for $i=1, \ldots, t$.  (Usually the context will make
explicit mention of the $P_i$'s unnecessary.)

\end{rm}
\end{defn}

Our first application of this construction is to actions, in general and
in the special case of local movement systems.

\begin{defn}
\label{}
\begin{rm}\

(1) Let $M$ be a structure and $\pair{G}{\circ}$ a subgroup of
$\fs{Aut}(M)$.
Then $\fs{ACT}(M,G) \eqdf (M, G; \fs{Ap})$,
where $\fnn{\fs{Ap}}{G \times |M|}{|M|}$ is the application function,
that is, $\fs{Ap}(g,a) = g(a)$.

(2) If $K$ is a class of structures all in the same language $\CL$, then
$\CL(K)$ denotes $\CL$.

(3) Recall that for a local movement system $\pair{B}{G}$,
$\fs{Var}(B,G)$
is a dense $G$-invariant subset of $B$.
The action of $G$ on $\fs{Var}(B,G)$ is faithful, and for $g \in G$
we identify $g\rest \fs{Var}(B,G)$ with $g$.
Let $\fs{VACT}(B,G)$ denote the substructure of $\fs{ACT}(B,G)$
whose universe is $\fs{Var}(B,G) \cup G$.
For contrast we write $\fs{BACT}(B,G) \eqdf \fs{ACT}(B,G)$.

(4) Let $K^{\fss{LM}}$ be the class of local movement systems
$\pair{B}{G}$.  Let
\smallskip
\newline
$K^{\fss{LM}}_{\fss{GR}} \eqdf
\setm{G}{\mbox{there exists } B \mbox{ such that }
\pair{B}{G} \in K^{\fss{LM}}}$,
\smallskip
\newline
$K^{\fss{LM}}_{\fss{VACT}} \eqdf
\setm{\fs{VACT}(B,G)}{\pair{B}{G} \in K^{\fss{LM}}}$,
\smallskip
\newline
$K^{\fss{LM}}_{\fss{BACT}} \eqdf
\setm{\fs{BACT}(B,G)}{\pair{B}{G} \in K^{\fss{LM}}}$.

\smallskip

(5) Let $\CL^{\fss{VACT}} \eqdf \CL(K^{\fss{LM}}_{\fss{VACT}} ) =
\{\leq,\circ,\fs{Ap},P_1,P_2\}$.
Let $\CL^{\fss{BACT}} \eqdf \CL(K^{\fss{LM}}_{\fss{BACT}})$.  Note
that $\CL^{\fss{BACT}}$ = $\CL^{\fss{VACT}}$.
\end{rm}
\end{defn}

For $g \in G$ we can capture $\fs{var}(g)$ first order, either in
$\fs{VACT}(B,G)$ or in $\fs{BACT}(B,G)$, and we can capture
$\fs{Var}(B,G)$ in $\fs{BACT}(B,G)$:

\begin{prop}
\label{PHIVAR}
Let $\pair{B}{G}$ be a local movement system.  Let
\[
\varphi_{var}(g,u)  \eqqdf
\forall w \mbox{\large(} u \leq w
\leftrightarrow (\forall a
(a \cdot g(a) = 0 \raro w \geq a)\mbox{\large)},
\]
re-expressed as a first order formula in $\CL^{\fss{VACT}}$ =
$\CL^{\fss{BACT}}$.  Then

{\rm (1)} For $g \in G$ and $v \in \fs{Var}(B,G)$:
\[
\fs{VACT}(B,G) \models \varphi_{var}[g,v] \mbox{ \ iff \ }
\fs{var}(g) = v.
\]

{\rm (2)} For $g \in G$ and $b \in B$:
\[
\fs{BACT}(B,G) \models \varphi_{var}[g,b] \mbox{ \ iff \ }
\fs{var}(g) = b.
\]

{\rm (3)} For $b \in B$:
\[
\fs{BACT}(B,G) \models \exists g (\varphi_{var}(g,b)) \mbox{ \ iff \ }
b \in \fs{Var}(B,G).
\]
\end{prop}

\noindent
{\bf Proof }
That  $\varphi_{var}(g,u)$ can be re-expressed in $\CL^{\fss{VACT}}$
is a consequence of Proposition~\ref{PST}.

(1) For any $g \in G$,
\[
\fs{var}(g) = \sum \setm{a \in \fs{Var}(B,G)}{a \cdot g(a) = 0}
\]
 by Propositions \ref{VSUP} and \ref{BSUP}.  (1) follows.

The proof of (2) is a subset of the proof of (1), and (3) is obvious.
\hfill\qed
\medskip

\smallskip

We conclude reconstruction results from the notion of an
interpretation.
The simplest such notion is a ``first order interpretation".
This is defined in \ref{INTERP}, and its key properties are stated in
Proposition~\ref{ELEM}.
A detailed example involving $K^{\fss{LM}}_{\fss{VACT}}$ is given in
Theorem~\ref{INTERPVA}; another is given in Lemma~\ref{DDLN}.

\begin{defn}
\label{}
\begin{rm}\

(1) If $\gf(x_1,\ldots,x_n)$ is a first order formula
in $\CL(M)$,
$\gf(M^n) \eqdf
\setm{\vec{a}\in |M|^n\ }{\ M\models \gf[\vec{a}]}$.

(2) If $\gf(\vecx,\vecy)$ is a first order formula in $\CL(M)$,
with $\ff{length}(\vecx) = k$ and $\ff{length}(\vecy) = l$,
then $\gf(M^k,M^l) \eqdf
\setm{\pair{\veca}{\vecb}}{M\models\gf[\veca,\vecb]}$.
\end{rm}
\end{defn}

When dealing with interpretations, we regard $m$-place functions
as
\newline
$(m+1)$-place relations and constants as 1-place relations, and
omit explicit reference to the function and constant symbols.

\begin{defn}
\label{INTERP}
\begin{rm}
Let $K$ and $K^*$ be classes of structures in the languages
$\CL$ and $\CL^*$ respectively,  and let $\CR\cntd K\times K^*$.

A {\em first order interpretation} ({\em FO-interpretation}) of
$K^*$ in $K$
relative to $\CR$ consists of first order formulas $\gf_U(\vecx)$ and
$\gf_{Eq}(\vecx,\vecy)$ in
$\CL$ (with $n \eqdf \ff{length}(\vecx) = \ff{length}(\vecy)$),
and for each $m$-place relation symbol $R$ in $\CL^*$
a first order formula $\gf_R(\vecx^1,\ldots,\vecx^m)$ in $\CL$
(with $\ff{length}(\vecx^j) = n$),
such that:

For each $\pair{M}{M^*}\in \CR$ there is a mapping
$\mu:\gf_U(M^n)\stackrel{onto}{\raro}|M^*|$ (called an {\em
interpreting mapping}
of $M^*$ in $M$)
satisfying

(1) For all $\veca,\vecb\in\gf_U(M^n)$,
$\mu(\veca)=\mu(\vecb)$ iff $M\models\gf_{Eq}[\veca,\vecb]$.

(2) For all $R$ in $\CL^*$ and
$\veca^1,\ldots,\veca^m\in\gf_U(M^n)$,\newline
$\trpl{\mu(\veca^1)}{\ldots}{\mu(\veca^m)}\in R^{M^*}$ iff
$M\models\gf_R[\veca^1,\ldots,\veca^m]$.

When such an interpretation
$\langle\gf_U,\,\gf_{Eq},\,\setm{\gf_R}{R \mbox{ is a relation
symbol in } \CL^*}
\rangle$
exists, $K^*$ is said to be {\em first order interpretable}
({\em FO-interpretable}) in $K$ relative to $\CR$.

\end{rm}
\end{defn}

We write $M_1 \equiv M_2$ if $M_1$ and $M_2$ are elementarily
equivalent, that is,
if $M_1$ and $M_2$ are in the same language
and satisfy the same first order sentences.
We write $M_1 \cong M_2$ if $M_1$ and $M_2$ are isomorphic, and
$\iso{\rho}{M_1}{M_2}$ means that $\rho$ is an isomorphism
between $M_1$ and $M_2$.
$\fs{Aut}(M)$ denotes the automorphism group of $M$.

\begin{prop}
\label{ELEM}
Suppose that $K^*$ is FO-interpretable in $K$ relative to $\CR$.
Let $\pair{M_i}{M_i^*} \in \CR$, $i = 1,2$.

{\rm (1)}  If $M_1\equiv M_2$, then $M_1^*\equiv M_2^*$.

{\rm (2)} If $M_1\cong M_2$, then $M_1^*\cong M_2^*$.

\end{prop}

{\bf Proof } The proof is easy and is left to the reader.  (Actually, the
proof is contained in the proofs of Propositions \ref{REEXPRESS} and
\ref{ISOEXT}.)
\hfill\qed
\medskip

In this paper we shall be especially interested in the following kind of
interpretation:

\begin{defn}
\label{DEFFOSTR}
\begin{rm}
Let $K$ and $K^*$ denote classes of structures in the
languages $\CL$ and
$\CL^*$ respectively, and let $\CR \subseteq K \times K^*$.

(1) We call $\trpl{K}{K^*}{\CR}$ a {\em subuniverse system}
if for every $\pair{M}{M^*} \in \CR$, $|M| \subseteq |M^*|$.

(2)  Let $\trpl{K}{K^*}{\CR}$ be a subuniverse system.
A {\em first order strong interpretation}
({\em FO-STR-interpretation}) of $K^*$ in $K$ relative to $\CR$
is an FO-interpretation
together with a first order formula $\gf_{\fss{Imap}}(\vecx,x)$ in
$\CL$, in which it is asked that the required interpreting mappings
$\mu:\gf_U(M^n)\stackrel{onto}{\raro}|M^*|$ not only satisfy conditions (1)
and (2) of
Definition~\ref{INTERP}, but also that they be {\em strong}, that is, that
they satisfy the additional condition:
\[
\setm{\pair{\veca}{\mu(\veca)}}{\mu(\veca) \in |M|} =
\gf_{\fss{Imap}}(M^n,M).
\]
\noindent (The restriction $\mu \rest \setm{\veca}{\mu(\veca) \in
|M|}$, whose range is $|M|$, is to be captured by
$\gf_{\fss{Imap}}$; or again, $\mu$ is to extend the function
$\gf_{\fss{Imap}}(M^n,M)$.)

When such an interpretation $\trpl{\gf_{\fss{Imap}}}{\gf_U,
\gf_{Eq}}{\ldots}$ exists, we say that $K^*$ is  {\em first order
strongly interpretable}
({\em FO-STR-interpretable}) in $K$ relative to $\CR$.

(3) A subuniverse system $\trpl{K}{K^*}{\CR}$ is a {\em first order
spanning system} if
for some $k$
there exists a first order formula $\gf^*_{\fss{Spn}}(x_1,\ldots,x_k,x)$
in the language $\CL^*$ (called a {\em spanning formula}) such that
for every $\pair{M}{M^*}\in \CR$,
$\gf^*_{\fss{Spn}}((M^*)^k,M^*)$
is a function whose domain is contained in $|M|^k$ (not merely in
$|M^*|^k$)
and whose range contains $|M^*| - |M|$.
\end{rm}
\end{defn}

\begin{prop}
\label{ISOEXT}
Let $\trpl{K}{K^*}{\CR}$ be a subuniverse system, and
suppose that $K^*$ is FO-STR-interpretable in $K$ relative to $\CR$.
Let $\pair{M_i}{M_i^*} \in \CR$, $i = 1,2$,
and let $\iso{\rho}{M_1}{M_2}$.  Then

{\rm (1)} There exists $\iso{\rho^*}{M_1^*}{M_2^*}$ such that
$\rho^* \supseteq \rho$.

{\rm (2)} If $\trpl{K}{K^*}{\CR}$ is a first order spanning system,
there exists a unique $\iso{\rho^*}{M_1^*}{M_2^*}$ such that
$\rho^* \supseteq \rho$.
\end{prop}

\noindent
{\bf Proof } (1) Let
$\trpl{\gf_{\fss{Imap}}(\vecx,x)}{\gf_U(\vecx),\gf_{Eq}(\vecx,\vecy)
}{\ldots}$
be an FO-STR-interpretation of $K^*$ in $K$,
and let $n = \ff{length}(\vecx)$.
For $i = 1,2$,
let $U_i = \gf_U((M_i)^n)$, $E_i = \gf_{Eq}((M_i)^n,(M_i)^n)$,
and $V_i = U_i/E_i$.
Let $\mu_i$ be a strong interpreting mapping of
$M^*_i$ in $M_i$.
$\mu_i$ induces a bijection
$\tldgt_i$ from $V_i$ to $|M_i^*|$,
and $\rho$ induces a bijection
$\tilde{\rho}$ from $V_1$ to $V_2$. Let
$\rho^* = \tldgt_2 \circ \tilde{\rho} \circ \tldgt_1\inverse$.
Clearly $\iso{\rho^*}{M^*_1}{M_2^*}$,
and we show that $\rho^* \supseteq \rho$.

Let $a \in |M_1|$,
and pick $\veca \in U_1$ such that $\mu_1(\veca) = a$.
$M_1 \models \gf_{\fss{Imap}}[\veca,a]$,
so $M_2 \models
\gf_{\fss{Imap}}[\rho(\veca),\rho(a)]$.
Also
$M_2 \models \gf_{\fss{Imap}}[\rho(\veca),\mu_2(\rho(\veca))]$.
Since $\gf_{\fss{Imap}}$
defines a function, $\rho(a) = \mu_2(\rho(\veca))$.
Therefore
\[
\rho^*(a) =
\tldgt_2 \circ \tilde{\rho} \circ \tldgt_1\inverse(a) =
\tldgt_2(\tilde{\rho}(\veca/E_1)) =
\tldgt_2(\rho(\veca)/E_2) =
\mu_2(\rho(\veca)) =
\rho(a).
\]

(2) Let $\gf^*_{\fss{Spn}}(\vecx,x)$
be a spanning formula.
Let $\iso{\rho'}{M_1^*}{M_2^*}$ be any isomorphism extending
$\rho$.
Let $a \in |M_1^*| - |M_1|$, and pick
$\veca \in |M_1|^k$ such that
$M_1^* \models \gf^*_{\fss{Spn}}[\veca,a]$.
Then $M_2^* \models \gf^*_{\fss{Spn}}[\rho'(\veca),\rho'(a)]$, and
this uniquely determines $\rho'(a)$.
\hfill\qed
\medskip

For some of our results, we need the more general notion of {\em
second order interpretability}.  Formulas in second order logic are
permitted to have second order variables (for $n$-place relations or
functions) as well as first order variables.
Most of this paper can be read without any dependence on the
second order results, and we mark these results with asterisks.

\begin{defn}
\label{}
$\!\!${\rm *}
\begin{rm}
If in the definition of FO-interpretability we allow
$\gf_U, \gf_{Eq}, \ldots$ to be second order formulas,
then $K^*$ is said to be {\em second order interpretable}
({\em SO-interpretable}) in $K$ relative to $\CR$.  (The range of the
interpreting map $\mu$ is still to be $|M^*|$ even though its {\em
arguments} are permitted to be second order.)

The notions of {\em second order strong interpretability}
({\em SO-STR-interpretability}) and {\em second order spanning
system} are defined similarly, with the last free variable of
$\gf_{\fss{Imap}}$ and of $\gf^*_{\fss{Spn}}$ required to be first
order.

In this paper the order of $\gf^*_{\fss{Spn}}$ never matters, and we
often speak just of a {\em spanning system} (without reference to
first or second order).
\end{rm}
\end{defn}

\begin{prop}
\label{2ORD}
$\!\!${\rm *}
Proposition~\ref{ELEM} (part (2)) and Proposition~\ref{ISOEXT}
remain valid when ``first order'' is replaced by ``second order''.
\end{prop}

The next theorem is the promised restatement of the Expressibility
Theorem in terms of interpretations.  From it the reconstruction
theorems for nearly ordered permutation groups will eventually
follow.
Note that the interpretation in $G$ of $\fs{VACT}(B,G)$ is first order,
whereas that of $\fs{BACT}(B,G)$ is second order.
In the statement of such results, $\CR$ will often be clear from
context and won't be mentioned.  (Here
$\CR = \setm{\pair{G}{\fs{VACT}(B,G)}}{\pair{B}{G} \in K^{LM}}$ in
(1)).

\begin{theorem}
\label{INTERPVA}
The following interpretational results hold:

{\rm (1)} $K^{\fss{LM}}_{\fss{VACT}}$ is FO-STR-interpretable in
$K^{\fss{LM}}_{\fss{GR}}$, and
$\pair{K^{\fss{LM}}_{\fss{GR}}}{K^{\fss{LM}}_{\fss{VACT}}}$ is a first
order spanning system.

$\!\!\!${\rm *(2)} $K^{\fss{LM}}_{\fss{BACT}}$ is
SO-STR-interpretable in
$K^{\fss{LM}}_{\fss{GR}}$, and
$\pair{K^{\fss{LM}}_{\fss{GR}}}{K^{\fss{LM}}_{\fss{BACT}}}$ is a
second order spanning system.
\end{theorem}

\noindent
{\bf Proof }
(1) The Expressibilty Theorem makes the FO-STR-interpretability
intuitively clear, and the formulas in $\CL^{\fss{GR}}$ it provides,
which we denote here
by $\varphi_{Eq}^{Ex}$, $\varphi_{\leq}^{Ex}$, and
$\varphi_{Ap}^{Ex}$, are the building blocks for the technical details.
We spell these details out here, though
never again will we be quite so meticulous.

The plan is to represent elements of $\fs{VACT}(B,G)$ by pairs
$\pair{g_1}{g_2}$ of elements of $G$, using {\em all} ordered pairs.
The interpreting mapping $\mu:|G|^2 \stackrel{onto}{\raro}
|\fs{VACT}(B,G)|$
will be given by
\[
\mu(g_1,g_2) = \left\{
\begin{array}{ll}
g_1 & \mbox{if $g_1 = g_2$,} \\*
\fs{var}(g_1)  & \mbox{if $g_1 \neq g_2$.}
\end{array}
\right.
\]
Accordingly, the interpretation consists of:

\smallskip
\noindent
$\varphi_{U}(x_1,x_2) \eqqdf x_1=x_1$.

\smallskip
\noindent
$\varphi_{Eq}(\vecx,\vecy) \eqqdf (x_1=x_2=y_1=y_2) \vee
\mbox{\large(} (x_1\neq x_2) \wedge (y_1\neq y_2) \wedge
\varphi_{Eq}^{Ex}(x_1,y_1) \mbox{\large)}$.

\smallskip
\noindent
$\varphi_{Imap}(\vecx,x) \eqqdf x_1=x_2=x$.

\smallskip
\noindent
$\varphi_{\leq}(\vecx,\vecy) \eqqdf (x_1 \neq x_2) \wedge (y_1
\neq y_2) \wedge
\varphi_{\leq}^{Ex}(x_1,y_1)$.

\smallskip
\noindent
$\varphi_{\circ}(\vecx,\vecy,\vecz) \eqqdf (x_1=x_2) \wedge
(y_1=y_2)
\wedge (z_1=z_2) \wedge (x_1 \circ y_1 = z_1)$.

\smallskip
\noindent
$\varphi_{Ap}(\vecx,\vecy,\vecz) \eqqdf (x_1=x_2) \wedge (y_1
\neq y_2)
\wedge (z_1 \neq z_2) \wedge \varphi_{Ap}^{Ex}(x_1,y_1,z_1)$.

\smallskip
\noindent
$\varphi_{P_1}(x_1,x_2) \eqqdf x_1 \neq x_2$.

\smallskip
\noindent
$\varphi_{P_2}(x_1,x_2) \eqqdf x_1=x_2$.

\smallskip

Finally, $\pair{K^{\fss{LM}}_{\fss{GR}}}{K^{\fss{LM}}_{\fss{VACT}}}$
is a spanning system because the formula $\varphi_{var}$ of
Proposition~\ref{PHIVAR} serves as the required $\varphi_{Spn}^*$.

\medskip

(2) Let $\fs{ACT}(B,G) \in K^{\fss{LM}}_{\fss{ACT}}$.  This time the
plan is to represent
members of $B$ by subsets of $G$.  If
$S \subseteq G$ then $S$ will represent
$\sum \setm{\fs{var}(g)}{g \in S}$; every member of $B$ has this
form by Proposition~\ref{BSUP}.  (Technically, the interpretation mapping
$\mu$ will be a function of two first order variables and one second
order variable.)  The building blocks of the interpretation are the
following three formulas (which can be written in $\CL^{\fss{GR}}$
with the aid of the formula $\varphi_{\fss{Dsjnt}}^{Gr}$ of
Proposition~\ref{DSJNT}):

(a) $\psi_{\leq}(S,T)$
is the formula which says:
For all $h \in G$,
if $\fs{var}(h) \cdot \fs{var}(g) = 0$
for all $g \in T$,
then $\fs{var}(h) \cdot \fs{var}(g) = 0$ for all $g \in S$.

\smallskip

(b) $\psi_{\fss{Eq}}(S,T) \eqqdf
\psi_{\leq}(S,T) \wedge \psi_{\leq}(T,S)$.

\smallskip

(c) $\psi_{\fss{Ap}}(f,S,T) \eqqdf
\psi_{\fss{Eq}}(S^f,T)$,
where $S^f = \setm{g^f}{g \in S}$.

\smallskip
\noindent
The rest of the proof of interpretability parallels that of (1).

Using the first order spanning formula developed above for
$\pair{K^{\fss{LM}}_{\fss{GR}}}{K^{\fss{LM}}_{\fss{VACT}}}$,
it is easy to construct a second order spanning formula for
$\pair{K^{\fss{LM}}_{\fss{GR}}}{K^{\fss{LM}}_{\fss{BACT}}}$.
\break
\rule{1pt}{0pt}\hfill\qed
\medskip

Usually in formulas such as $\varphi_{Spn}^*(g,v)$  in the proof
above, we omit explicit mention of ``$g \in G$'' and ``$v \in
\fs{Var}(B,G)$'' and the like, as the meaning will be clear from
context.

Note the following relationship between automorphisms of
$\fs{ACT}(M,G)$ and automorphisms of $G$. We use this relationship
in deducing reconstruction results from interpretability results.

\begin{lemma}
\label{ISO}
Let $\fnn{\sigma}{|\fs{ACT}(M_1,G_1)|}{|\fs{ACT}(M_2,G_2)}|$.
Then\newline
$\iso{\sigma}{\fs{ACT}(M_1,G_1)}{\fs{ACT}(M_2,G_2)}$ iff
$\iso{\sigma \rest |M_1|}{M_1}{M_2}$,
$\fs{Rng}(\sigma \rest G_1) = G_2$,
and for every $g \in G_1$,
$\sigma(g) = (\sigma \rest |M_1|) \circ g \circ (\sigma \rest
|M_1|)\inverse$.
\end{lemma}

We exhibit a reconstruction theorem here.  Most of part (2) is
contained in \cite{Ru}.  However, the reconstruction theorems for
nearly ordered permutation groups will be deduced from
Theorem~\ref{INTERPVA} rather than from Theorem~\ref{LMSRT}.

\begin{theorem}
\label{LMSRT}
{\rm (Reconstruction Theorem for Local Movement Systems)} \
Let $\pair{B_i}{G_i}$ be a local movement system, $i=1,2$.

{\rm (1)} If $G_1 \equiv G_2$ then
$\fs{VACT}(B_1,G_1) \equiv
    \fs{VACT}(B_2,G_2)$.

{\rm (2)} Suppose that $\iso{\alpha}{G_1}{G_2}$.  Then

$\,\,\,\,\,\,\,\,\,\,${\rm (a)} There exists a unique isomorphism
$\iso{\tau}{B_1}{B_2}$ such that $\tau$ induces $\alpha$, that is,
$\alpha(g) = \tau \circ g \circ \tau\inverse
\mbox{ for all } g \in G_1$;
or equivalently, such that
$\iso{\tau \cup \alpha}
{\fs{BACT}(B_1,G_1)}
{\fs{BACT}(B_2,G_2)}$.

$\,\,\,\,\,\,\,\,\,\,${\rm (b)} There exists a unique isomorphism
$\iso{\tau_V}{\fs{Var}(B_1,G_1)}{\fs{Var}(B_2,G_2)}$ (as partially ordered
sets) such that
$\iso{\tau_V \cup \alpha}
{\fs{VACT}(B_1,G_1)}
{\fs{VACT}(B_2,G_2)}$.
\end{theorem}

\noindent
{\bf Proof }
The theorem follows from Theorem~\ref{INTERPVA} with the aid of
Propositions \ref{ELEM} and \ref{2ORD}, Lemma~\ref{ISO}, and
finally Proposition ~\ref{PHIVAR}.
\hfill\qed

\begin{defn}
\label{}
\begin{rm}
For a class $K$ of structures, let
\[
K^{\fss{ACTS}} \eqdf
\setm{\fs{ACT}(M,G)}{M \in K \mbox{ and } G \leq \fs{Aut}(M)}.
\]
Suppose that
$K^*$ is FO-STR-interpretable (or even SO-STR-interpretable) in $K$
relative to $\CR$, and that $\trpl{K}{K^*}{\CR}$ is a spanning system.
For $\pair{M}{M^*} \in \CR$ and $g \in \fs{Aut}(M)$, let
$g^{M,M^*}$ denote the unique member of $\fs{Aut}(M^*)$
extending $g$ (as provided by Propositions \ref{ISOEXT} and
\ref{2ORD}).
The function $g \mapsto g^{M,M^*}$ is an embedding of
$\fs{Aut}(M)$
in $\fs{Aut}(M^*)$.
Thus we may regard every $G \leq \fs{Aut}(M)$ as a subgroup of
$\fs{Aut}(M^*)$.  Let
\[\CR^{\fss{ACTS}} \eqdf
\setm{\pair{\fs{ACT}(M,G)}{\fs{ACT}(M^*,G)}}
{\pair{M}{M^*} \in \CR \mbox{ and } G \leq \fs{Aut}(M)}.
\]\end{rm}
\end{defn}

\begin{prop}\
\label{ACTS}

{\rm (1)} Let $K^*$ be  FO-STR-interpretable in $K$ relative to $\CR$,
and suppose that $\trpl{K}{K^*}{\CR}$ is a spanning system.
Then $(K^*)^{\fss{ACTS}}$ is FO-STR-interpretable in
$K^{\fss{ACTS}}$ relative to $\CR^{\fss{ACTS}}$,
and $\trpl{K^{\fss{ACTS}}}{(K^*)^{\fss{ACTS}}}{\CR^{\fss{ACTS}}}$
is a spanning system.

$\!\!\!${\rm *(2)} Part (1) holds also for SO-STR-interpretability.
\end{prop}

\noindent
{\bf Proof } The proof is easy.
\hfill\qed
\medskip

For applications of the next proposition, see Corollary~\ref{FOP} and
Lemma~\ref{LEMPARAM}.

\begin{prop}
\label{REEXPRESS}\

{\rm (1)} Suppose that $K^*$ is FO-STR-interpretable in $K$ relative
to $\CR$.
Then for every first order formula $\varphi(z_1,\ldots,z_m)$ in
$\CL(K^*)$ whose free variables are all in $\CL(K)$,
there is a first order formula $\psi_{\varphi}(z_1,\ldots,z_m)$ in
$\CL(K)$
such that for every $\pair{M}{M^*} \in \CR$ and $\vecb \in |M|^m$:
$M^* \models \varphi[\vecb]$ iff $M \models \psi_{\varphi}[\vecb]$.

$\!\!\!${\rm *(2)} Part (1) remains valid when $K^*$ is only
SO-STR-interpretable in $K$ (with $\varphi$ still first order), but
$\psi_{\varphi}$ is then only second order (though its {\em free}
variables are first order).
\end{prop}

\noindent
{\bf Proof }
(1) Let $\trpl{\gf_{\fss{Imap}}(\vecx,x)}
{\gf_U(\vecx),\gf_{Eq}(\vecx,\vecy)}{\ldots}$
be an FO-STR-interpretation of $K^*$ in $K$ relative to $\CR$, and let
$n = \ff{length}(\vecx)$.  It is easy to see that there is a first order
formula $\varphi^{\vee}(\vecx^1,\ldots,\vecx^m)$ in $\CL(K)$
such that for every $\pair{M}{M^*} \in \CR$, every interpreting map
$\mu:\gf_U(M^n) \raro |M^*|$, and all
$\veca^1,\ldots,\veca^m\in\gf_U(M^n)$:
\[
M^* \models \varphi[\mu(\veca^1),\ldots,\mu(\veca^m)] \mbox{ \ iff
\ }
M\models \varphi^{\vee}[\veca^1,\ldots,\veca^m].
\]
Then
\[
\psi_{\varphi}(z_1,\ldots,z_m) \eqqdf
(\exists \vecx^1,\ldots,\vecx^m)
(\bigwedge_{i=1}^m \gf_{\fss{Imap}}(\vecx^i,z_i) \wedge
\varphi^{\vee}(\vecx^1,\ldots,\vecx^m))
\]
is as required.

(2) The proof of (2) is similar.
\hfill\qed
\medskip

In the next proposition we describe a general
argument which is used in
proving the reconstruction theorems of this paper (see Theorem~\ref{BIG},
for example).

\begin{prop}
\label{FORMALIZE}\

{\rm (1)} Suppose that $K^*$ is FO-STR-interpretable in $K$ relative
to $\CR$.
Let $\varphi(u)$ be a first order formula in $\CL(K^*)$, and let
$\CL' \subseteq \CL(K^*)$.
For $M^* \in K^*$, let $\hat{M}^*$ be the substructure of $M^*$
whose universe is $\varphi(M^*)$, and let
$M^{*}{'}$ be the reduct of $\hat{M}^*$ to the language $\CL'$.
Suppose that for every $\pair{M}{M^*{'}} \in \CR$, $|M| \subseteq |M^*{'}|$.
Then $\setm{M^*{'}}{M^* \in K^*}$ is FO-STR-interpretable in $K$
relative to $\setm{\pair{M}{M^*{'}}}{\pair{M}{M^*} \in \CR}$.

$\!\!\!${\rm *(2)} Part (1) remains valid when
FO-STR-interpretability is replaced throughout by
SO-STR-interpretability
(with $\varphi$ still first order).
\end{prop}

\noindent
{\bf Proof }
This proposition follows from the previous one.
\hfill\qed

\begin{prop}
\label{TRANS}
The following transitivity properties hold:

{\rm (1) (a)} Suppose that $K^*$ is FO-interpretable in $K$ relative to
$\CR$,
and $K^{**}$ is FO-interpretable in $K^*$ relative to $\CR^*$.
Then
$K^{**}$ is FO-interpretable in $K$ relative to $\CR \circ \CR^*$.

$\,\,\,\,\,\,\,\,${\rm (b)} The same holds for FO-STR-interpretablity.

$\!\!\!${\rm *(2) (a)} In {\rm (1a)}, if $K^*$ is SO-interpretable in $K$
and $K^{**}$ is FO-interpretable in $K^*$, then $K^{**}$
is SO-interpretable in $K$.

$\,\,\,\,\,\,\,\,${\rm (b)} The same holds for SO-STR-interpretability.
\end{prop}

\noindent
{\bf Proof }
We give a proof of (1b); this will include a proof of (1a). Let
$$\qdrpl{\gf_{\fss{Imap}}(\vecx,x)}
{\gf_U(\vecx)}{\gf_{Eq}(\vecx,\vecy)}{\ldots}$$
be an FO-STR-interpretation of $K^*$ in $K$ relative to $\CR$,
and let $n = \ff{length}(\vecx)$.
Let $$\qdrpl{\gf_{\fss{Imap}}^*(x_1,\ldots,x_{n^*},x)}
{\varphi_U^*(x_1,\ldots,x_{n^*})}
{\gf_{Eq}(x_1,\ldots,x_{n^*};y_1,\ldots,y_{n^*})}{\ldots}$$
be an FO-STR-interpretation of $K^{**}$ in $K^*$ relative to $\CR^*$.

We define an interpretation
$\qdrpl{\varphi'_{\fss{Imap}}}{\varphi'_U}{\varphi'_{\fss{Eq}}}
{\ldots}$ of $K^{**}$ in $K$.
Let
$$
\varphi'_U(\vecx^1,\ldots, \vecx^{n^*})
\eqqdf
\bigwedge_{i = 1}^{n^*} \varphi_U(\vecx^i) \wedge
\varphi_U^{*\vee}(\vecx^1,\ldots, \vecx^{n^*}),
$$
where $\varphi_U^{*\vee}$ is defined from $\varphi_U^*$ as in the
proof of Proposition \ref{REEXPRESS}.
Let
$$
\varphi'_{\fss{Eq}}
(\vecx^1,\ldots, \vecx^{n^*};\vecy^1,\ldots, \vecy^{n^*}) \eqqdf
\varphi_{\fss{Eq}}^{*\vee}
(\vecx^1,\ldots, \vecx^{n^*};\vecy^1,\ldots, \vecy^{n^*}).
$$
We skip the definition of the formulas which interpret the relation
symbols.

We show that the required formula
$\gf_{\fss{Imap}}'$ is given by
\begin{eqnarray*}
& \gf_{\fss{Imap}}'(\vecx^1,\ldots,\vecx^{n^*},x)
\eqqdf \\
& \begin{displaystyle}\bigwedge_{i=1}^{n^*}
\end{displaystyle}\gf_U(\vecx^i) \wedge
\exists \vecx
\mbox{\large(}\gf_{\fss{Imap}}(\vecx,x) \wedge
\gf_{\fss{Imap}}^{*\vee}(\vecx^1,\ldots,\vecx^{n^*},\vecx)
\mbox{\large)}.
\end{eqnarray*}
Let $\pair{M}{M^*} \in \CR$ and
$\pair{M^*}{M^{**}} \in \CR^*$.
Let $\mu$ and $\mu^*$ be strong interpreting mappings of $M^*$ in
$M$
and of $M^{**}$ in $M^*$, respectively.
We define an interpreting map $\mu'$ of $M^{**}$ in $M$.
Let
$\veca^1\cat\dots\cat\veca^{n^*} \in \varphi'_U(M^{nn^*})$.
Then $\veca^i \in \varphi_U(M^{n})$ for $i = 1, \dots n^*$, and
$M \models \varphi_U^{*\vee}[\veca^1,\ldots,\veca^{n^*}]$.
Let $b_i = \mu(\veca^i)$.  Then
$M^* \models \varphi^*_U[b_1,\ldots b_{n^*}]$. We define
$\mu'(\veca^1,\ldots,\veca^{n^*}) = \mu^*(b_1,\ldots,b_{n^*})$.
We leave it to the reader to check that $\mu'$ is an
interpreting mapping.

Now we show that $\mu'$ is strong.
Let
$\veca^1\cat\dots\cat\veca^{n^*} \in \varphi'_U(M^{nn^*})$
and suppose that $a \eqdf \mu'(\veca^{1}, \dots, \veca^{n^*}) \in
|M|$.
Pick $\veca$ such that $a = \mu(\veca)$; then
$M \models \varphi_{\fss{Imap}}[\veca,a]$.
Let $b_i = \mu(\veca^i)$.  Then $\mu^*(b_1,\ldots,b_{n^*}) =
\mu'(\veca^1,\ldots,\veca^{n^*}) = a$,
and hence $M^* \models \varphi^*_{\fss{Imap}}[b_1,\ldots
b_{n^*},a]$, so that
$M \models
\varphi_{\fss{Imap}}^{*\vee}[\veca^1,\ldots \veca^{n^*},\veca]$.
Obviously, $M \models \varphi_U[\veca^i]$ for $i = 1, \ldots n^*$.
Therefore $M \models
\varphi'_{\fss{Imap}}[\veca^1,\ldots,\veca^{n^*},a]$.

Conversely, if
$M \models \varphi'_{\fss{Imap}}[\veca^1,\ldots,\veca^{n^*},a]$
then $a = \mu'(\veca, \dots, \veca^{n^*})$ and $a \in |M|$;
this we leave this to the reader to check.
\hfill\qed
\medskip

 Caution:  The notion of spanning system is not transitive.
\newpage

%

\section{The core reconstruction results}
\label{RECON}

Here we present our core reconstruction results.  We state the results
for linear permutation groups, and observe that most of them follow
easily from the Linear Interpretation Theorem~\ref{LIT}.  Then we
give brief parallel statements for the other three types of nearly
ordered permutation groups.  Next we prove \ref{LIT}, and then we
modify the proof to deal with the other cases.  Finally, we consider
interrelations among the four cases.

First we establish some notation common to all four cases.

Let $\pair{N}{G}$ be a nearly ordered
permutation group.  As in Proposition~\ref{BN},
$G \leq \fs{Aut}(N) \leq \fs{Aut}(\barB(N))$.
Always our hypotheses will guarantee that $\pair{N}{G}$ is locally
moving, that is, that $\pair{N}{G} \in K^{\fss{NO}}$, the class of all
locally moving nearly ordered permutation groups.

\begin{defn}
\label{}
\begin{rm}
For each $n \geq 2$:

(1) $K^{\fss{LNn}}$ denotes the class of $n$-interval-transitive linear
permutation groups $\pair{L}{G}$ which have a nonidentity bounded
element, or equivalently, which are locally moving.

(2) $K^{\fss{MNn}}$ denotes the analogue for  monotonic permutation
groups.

(3) $K^{\fss{CRn}}$ denotes the analogue for  circular permutation
groups.

(4) $K^{\fss{MCn}}$ denotes the analogue for  monocircular
permutation groups.
\end{rm}
\end{defn}

Our core reconstruction results will be for the classes
$K^{\fss{LN2}}$, $K^{\fss{MN3}}$,  $K^{\fss{CR3}}$, and
$K^{\fss{MC4}}$.

We associate with each locally moving $\pair{N}{G}$ a crucial
structure $\fs{EP}(N,G)$ on which $G$ will act.  It is by ``interpreting
the action of $G$ on $\fs{EP}(N,G)$ in the language $\CL^{\fss{GR}}$ of
groups'' that we prove our reconstruction theorems.  The matters
discussed here are unaffected by reversal of the order or orientation
on $N$.

\begin{defn} \label{DEFEPLN}
\begin{rm}
Let $\pair{N}{G} \in K^{\fss{NO}}$, with $N = L$ or $N = \fs{ED}(L)$.

(1) For a regular open set $U \in \barB(N)$, let $\fs{Ep}(U)$ consist of
the endpoints in
$\barL$ of the convex hull of $U$, of which there are at
most two.

(2) Let
\begin{eqnarray*}
& \fs{Ep}(N,G) \eqdf \bigcup _{g \in G} \fs{Ep}(\fs{var}(g)),  \\*
& \fs{EP}(N,G) \eqdf \pair{\fs{Ep}(N,G)}{\fs{Ed}^{\barL}}.
\end{eqnarray*}
$\fs{EP}(N,G)$ is only an equal direction structure even if $N$ is a
chain.  The structures $\fs{EP}(\fs{ED}(L),G)$ and
$\fs{EP}(L,\fs{Opp}(G))$
coincide. \end{rm}
\end{defn}

\begin{defn} \label{DEFEPCR}
\begin{rm}
Let $\pair{N}{G}\in K^{\fss{NO}}$, with $N = C = \fs{CR}(L)$ or $N =
\fs{EO}(C)$.

(1) For $U \in \barB(N)$, let $\fs{Ep}(U)$ consist of the endpoints
(boundary points) of the connected components of $-U$.  Thus for $g
\in G$, $\fs{Ep}(\fs{var}(g))$ consists of the endpoints of the
connected components of $\fs{fix}(g)$, or
equivalently, of the {\em nonsingleton} connected components of
$\fs{fps$_{\barN}$}(g)$.   (Recall Proposition \ref{VARBN}.)
$\fs{Ep}(\fs{var}(g))$ contains two points for each of these
components (except when $g = \fs{Id}$).

(2) Let
\begin{eqnarray*}
& \fs{Ep}(N,G) \eqdf \bigcup _{g \in G} \fs{Ep}(\fs{var}(g)),  \\*
& \fs{EP}(N,G) \eqdf \pair{\fs{Ep}(N,G)}{\fs{Eo}^{\barC}}.
\end{eqnarray*}
$\fs{EP}(N,G)$ is only an equal orientation structure even if $N$ is a
circular order.  $\fs{EP}(\fs{EO}(C),G)$  and $\fs{EP}(C,\fs{Opp}(G))$
coincide.
\end{rm}
\vspace{-2.5mm}
\end{defn}

In all four cases $\fs{Ep}(N,G)$ is $G$-invariant, and the assumption
that $\pair{N}{G}$ is locally moving ensures that $\fs{Ep}(N,G)$ is
dense in $\barN$.  The faithful action on $\fs{EP}(N,G)$ is a locally
moving monotonic (or monocircular) permutation group which for
any $n$ is $n$-interval-transitive if $\pair{N}{G}$ is.

\begin{defn}
\label{}
\begin{rm}
Let $\pair{N}{G} \in K^{\fss{NO}}$.

(1) Let $\fs{EPACT}(N,G) \eqdf \fs{ACT}(\fs{EP}(N,G),G)$.

(2) Let $\fs{DCACT}(N,G) \eqdf \fs{ACT}(\barN^{\smallr},G)$.

\end{rm}
\end{defn}

\begin{defn}
\label{}
\begin{rm}
Let $K \subseteq K^{\fss{NO}}$.

(1) Let $K_{EPACT} \eqdf \setm{\fs{EPACT}(N,G)}{\pair{N}{G} \in K}$.

(2) Let $K_{DCACT} \eqdf \setm{\fs{DCACT}(N,G)}{\pair{N}{G} \in K}$.

(3) Let $K_{GR} \eqdf \setm{G}{\mbox{there exists } N \mbox{such
that } \pair{N}{G}
\in K}$.

\end{rm}
\end{defn}

When $\pair{L}{G}$ is a linear permutation group and we speak of
$G$ as a partially ordered group, the order will always be understood
to be the pointwise order, and similarly in the monotonic case for
``half-ordered group'' (Definition~\ref{HALFORDERED}).

\medskip

Now we state the {\bf linear reconstruction results}.
We begin with the Linear Reconstruction Theorem, which has as
immediate consequences all the results through
Corollary~\ref{CUTDOWN}.
Its hypotheses do not imply 3-interval-transitivity or approximate
2-{\it o-}transitivity (Example~\ref{WR1VARIANT}).
In Section~\ref{FURTHER} we establish some of these same
conclusions under different hypotheses.

\begin{theorem}\label{LRT}
{\rm (Linear Reconstruction Theorem)} \
Let $\pair{L_i}{G_i}$ be a linear permutation group which is
2-interval-transitive and has a nonidentity bounded element,
$i=1,2$.

{\rm (1)} If $G_1 \equiv G_2$ then
\[\fs{ACT}(\fs{EP}(L_1,G_1),G_1) \equiv
    \fs{ACT}(\fs{EP}(L_2,G_2),G_2).
\]

{\rm (2)} Suppose that $\iso{\alpha}{G_1}{G_2}$.
Then there exists a unique monotonic bijection
$\fnn{\tau}{\barL_1}{\barL_2}$
which induces $\alpha$, that is,
\[
\alpha(g) = \tau \circ g \circ \tau\inverse \mbox{ for all } g \in G_1;
\]
or put another way, there exists a unique
$\iso{\tau}{\barL_1^{\smallr}}{\barL_2^{\smallr}}$ such that
\[
\iso{\tau \cup \alpha}
{\pair{\barL_1^{\smallr}}{G_1}}
{\pair{\barL_2^{\smallr}}{G_2}}.
\]
\end{theorem}

\noindent
{\bf Note }
The map $\fnn{\tau}{\barL_1}{\barL_2}$ is a map between chains,
whereas the map $\iso{\tau}{\barL_1^{\smallr}}{\barL_2^{\smallr}}$
is an isomorphism between equal direction structures
(or equiv\-a\-lently---because of Dedekind completeness---a
homeomorphism).
Also, (2) implies that
\[
\iso{\tau\rest \fs{Ep}(L_1,G_1) \cup \alpha}
{\fs{ACT}(\fs{EP}(L_1,G_1),G_1)}
{\fs{ACT}(\fs{EP}(L_2,G_2),G_2)}.
\]

\begin{cor}
\label{2GRORBITS}
In Theorem~\ref{LRT}, $\tau = \tau_\alpha$ induces a bijection
$M \mapsto \tau(M)$ between the set of orbits of
$\pair{\barL_1}{G_1}$ and that of
$\pair{\barL_2}{G_2}$.
Either $\tau \rest M$ preserves order for each orbit $M$ (if
$\fnn{\tau}{\barL_1}{\barL_2}$
preserves order) or $\tau \rest M$ reverses order for each $M$ (if
$\tau$ reverses order).
In the former case
$$
\iso{(\tau \rest M) \cup \alpha}
{\pair{M}{G_1}}
{\pair{\tau(M)}{G_2}},
$$
and in the latter
$$
\iso{(\tau \rest M) \cup \alpha}
{\pair{M}{G_1}}
{\pair{\tau(M)^*}{G_2}}.
$$
\end{cor}

\smallskip

\noindent
{\bf Note }
$\alpha$ either preserves pointwise order on $G_1$ and $G_2$
(if $\tau$ preserves order) or reverses
it (if $\tau$ reverses order).  Hence if $\iso{\alpha}{G_1}{G_2}$ as
partially ordered groups (under the pointwise order), then $\iso{\tau \cup
\alpha}{\pair{\barL_1}{G_1}}{\pair{\barL_2}{G_2}}$.
\medskip

As mentioned in the introduction, there is also a one-group version
of Theorem~\ref{LRT}
which permits us to  identify $\fs{Aut}(G)$ with the normalizer of
$G$ in $\fs{Aut}(\barL^{\smallr})$:

\begin{cor}
\label{1GR}
Let $\pair{L}{G}$ be a linear permutation group which is
2-interval-transitive and has a nonidentity bounded element.  Then
every automorphism $\alpha$ of $G$ is conjugation by a unique
automorphism  $\tau$ of $\barL^{\smallr}$, so that
\[
\fs{Aut}(G) =
\fs{Norm}_{\mbox{\scriptsize {\fs{Aut}}}(\barL^{\mbox
{\smallr}})}(G) =
\fs{Norm}_{\mbox{\scriptsize {\fs{Homeo}}}(\barL)}(G).
\]
\end{cor}

\noindent
{\bf Note }
$\tau\rest\fs{Ep}(L,G) \cup \alpha$ is an automorphism of
$\fs{ACT}(\fs{EP}(L,G),G)$.

\begin{cor}
\label{1GRORBITS}
In Corollary~\ref{1GR}, $\tau = \tau_\alpha$ induces a permutation
$M \mapsto \tau(M)$ of the set of orbits of
$\pair{\barL}{G}$, and either $\tau \rest M$ preserves order for
each orbit $M$ or else $\tau \rest M$ reverses order for each $M$.
In the former case
$$
\iso{(\tau \rest M) \cup \alpha}
{\pair{M}{G}}
{\pair{\tau(M)}{G}},
$$
and in the latter
$$\iso{(\tau \rest M) \cup \alpha}
{\pair{M}{G}}
{\pair{\tau(M)^*}{G}}.
$$
\end{cor}

\smallskip

\noindent
{\bf Note }
Here $\alpha$ either preserves or reverses the pointwise order on
$G$.  We denote by {\it o-}$\fs{Aut}(G)$ the group of automorphisms
of $G$ which preserve the pointwise order, that is, which are
automorphisms of the partially ordered group $G$.
{\it o-}$\fs{Aut}(G) \leq
\fs{Aut}(\barL)$, and the index of {\it o-}$\fs{Aut}(G)$ in
$\fs{Aut}(G)$ is $\leq 2$.

\begin{cor}
\label{CHAR}
Let $\pair{L}{G} \in K^{\fss{LN2}}$.  Then

{\rm (1)} $\fs{Bdd}(G)$ and $\fs{Bdd}'(G)$ are characteristic
subgroups of $G$.

{\rm (2)} Every automorphism $\alpha$ of $G$ either preserves or
interchanges $\fs{Lft}(G)$ and $\fs{Rt}(G)$.

{\rm (3)} $\fs{Aut}(G) \leq
\fs{Aut}\mbox{\large(}\fs{Bdd}(G)\mbox{\large)} \leq
\fs{Aut}\mbox{\large(}\fs{Bdd}'(G)\mbox{\large)} \leq
\fs{Aut}(\barL^{\mbox{\smallr}})$, and each automorphism of $G$
is determined by its
restriction to $\fs{Bdd}'(G)$.

{\rm (4)} {\it o-}$\fs{Aut}(G) \leq
\fs{Aut}\mbox{\large(}\fs{Lft}(G)\mbox{\large)}$; and
$\fs{Aut}\mbox{\large(}\fs{Lft}(G)\mbox{\large)} \leq
\fs{Aut}(\barL)$ provided
$\fs{Lft}(G) \neq \fs{Bdd}(G)$.
The same holds for
$\fs{Aut}\mbox{\large(}\fs{Rt}(G)\mbox{\large)}$.
\end{cor}

When can we sharpen Corollary~\ref{1GR} by replacing $\barL$ by
$L$ in its conclusion, so that we get $\fs{Aut}(G) \leq
\fs{Aut}(L^{\smallr})$?  Obviously this is valid when $L$ is already
Dedekind complete, but when else?

\begin{cor}
\label{CUTDOWN}
Let $\pair{L}{G} \in K^{\fss{LN2}}$, and suppose that $\pair{L}{G}$ is
transitive.  Then
\[
\fs{Aut}(G) =
\fs{Norm}_{\mbox{\scriptsize {\fs{Aut}}}(L^{\mbox
{\smallr}})}(G)
\]
iff there is no orbit $M \neq L$ of $\pair{\barL}{G}$
with $\pair{L}{G} \cong \pair{M}{G}$ or $\pair{L}{G} \cong
\pair{M^*}{G}$.
\end{cor}

\noindent
{\bf Proof }
For every $\alpha \in \fs{Aut}(G)$,
$\iso{(\tau_\alpha \rest L) \cup \alpha}
{\pair{L}{G}}
{\pair{\tau_\alpha(L)}{G}}$
(or $\pair{\tau_\alpha(L)^*}{G}$) by Corollary~\ref{1GRORBITS}; this
establishes sufficiency of the condition.

Conversely, if $\iso{\sigma \cup \beta}
{\pair{L}{G}}{\pair{M}{G}}$ (or $\pair{M^*}{G}$) for some orbit $M
\neq L$ of
$\pair{\barL}{G}$, then $\beta \in
\fs{Aut}(G)$ and $\tau_\beta(L) = M \neq L$ since $\tau_\beta =
\sigma$ by the
uniqueness of $\tau_\beta$ in Theorem~\ref{LRT}.
\hfill\qed
\medskip

This obviously applies to $\pair{\R}{\fs{Aut}({\R})}$, and almost as
obviously to
\newline
$\pair{\Q}{\fs{Aut}({\Q})}$.
For these and other applications, see Theorems \ref{RQ} and
\ref{ALPHA}.
This line of thought will be continued in Theorem~\ref{ORBITCORR}.

\medskip

Theorem~\ref{LRT} will be deduced from Theorem~\ref{LIT}(1),
which says that for $\pair{L}{G} \in  K^{\fss{LN2}}$ the action of $G$
on $\fs{EP}(L,G)$
can be recovered from $G$ using first order formulas.  Part (2) of
Theorem~\ref{LIT} is included here for purposes of comparison.
However, we shall not really use such second order results until
Section~\ref{FURTHER}, and we continue marking them with
asterisks.

We have not included any statement in Theorem~\ref{LIT} about
$\pair{K^{\fss{LN2}}_{\fss{GR}}}
{K^{\fss{LN2}}_{\fss{EPACT}}}$ or
$\pair{K^{\fss{LN2}}_{\fss{GR}}}
{K^{\fss{LN2}}_{\fss{DCACT}}}$
being a spanning system.  In fact both {\em are} spanning systems,
but rather than showing this we achieve the same end by the use of
Proposition~\ref{UNIQ}.

\begin{theorem} {\rm (Linear Interpretation Theorem)}
\label{LIT}

{\rm (1)} $K^{\fss{LN2}}_{\fss{EPACT}}$ is FO-STR-interpretable in
$K^{\fss{LN2}}_{\fss{GR}}$.

$\!\!\!${\rm *(2)} $K^{\fss{LN2}}_{\fss{DCACT}}$ is
SO-STR-interpretable in
$K^{\fss{LN2}}_{\fss{GR}}$.
\end{theorem}

\noindent
{\bf Deduction of Theorem~\ref{LRT} from Theorem~\ref{LIT}(1) }
Let $\pair{L_1}{G_1}$  and $\pair{L_2}{G_2} \in K^{LN2}$.

(1) This follows from Proposition~\ref{ELEM}.

(2) Let $\iso{\alpha}{G_1}{G_2}$.  By Proposition~\ref{ISOEXT},
$\alpha$ can be extended to an isomorphism
\[
\iso{\alpha^*}{\fs{ACT}(\fs{EP}(L_1,G_1),G_1)}
{\fs{ACT}(\fs{EP}(L_2,G_2),G_2)}.
\]
Then $\alpha^*\rest\fs{Ep}(L,G)$ can be extended to an isomorphism
$\iso{\tau}{\barL_1^{\smallr}}{\barL_2^{\smallr}}$.
Lemma \ref{ISO} certifies that the two ways of stating part (2) are
equivalent, and that $\tau$ induces $\alpha$.  Finally, by
Proposition~\ref{UNIQ}, $\tau$ is the {\em only} isomorphism
inducing $\alpha$.

Observe that part (2) of \ref{LRT} could have been deduced from
part (2) of \ref{LIT}.
\hfill\qed
\bigskip

\noindent
{\bf Deduction of Theorem \ref{LIT}(2) from Theorem~\ref{LIT}(1) }
Let $\pair{L}{G} \in K^{\fss{LN2}}$.  Then $\barL^{\smallr} =
\overline{\fs{EP}(L,G)}$ by the density of $\fs{Ep}(L,G)$ in
$\barL^{\smallr}$.  Using the classical Dedekind cut construction, the
points of $\overline{\fs{EP}(L,G)}$ can be represented by pairs of
subsets of $\fs{Ep}(L,G)$.  The equal direction relation
$\fs{Ed}^{\barL}$ and the action of $G$ are easily captured, and
Theorem~\ref{TRANS} completes the proof.
\newline
\rule{0pt}{1pt}\hfill\qed
\medskip

The next result says that for $\pair{L}{G} \in K^{\fss{LN2}}$, certain
important notions expressible in the language of monotonic
permutation groups can in fact be expressed in the language of
groups.

Let $\pair{\fs{ED}(L)}{G}$ be a monotonic permutation group.  A
nonidentity element $g \in G$ is called {\em signed} if, in
$\pair{L}{<}$,
either
$g(x) \geq x$ for all $x \in L$
or else $g(x) \leq x$ for all $x \in L$; that is, if
\[
\fs{Ed}(g(x_1),x_1; g(x_2),x_2) \mbox{ for all } x_1,x_2 \in
\fs{supp}(g) \cap L.
\]
Signed elements are necessarily in $\fs{Opp}(G)$.
An element $g \in \fs{Opp}(G)$ has to be signed if it has just one
supporting interval (see Note \ref{NOTEGUIDE} for the definition).
Two signed elements $g_1$ and $g_2$ have the {\em same sign} if
\[
\fs{Ed}(g_1(x_1),x_1; g_2(x_2),x_2) \mbox{ for all } x_i \in
\fs{supp}(g_i) \cap L.
\]
Equivalent definitions of these notions can be obtained by replacing
$L$ by any dense subset of $\barL$.

\begin{cor}
\label{FOP}
For every first order formula $\varphi(x_1,\ldots,x_n)$ in the
language
$\CL^{\fss{MNPG}}$ of monotonic permutation groups whose free
variables are all group variables, there is a first
order formula $\psi _\varphi(x_1,\ldots,x_n)$ in the language
$\CL^{\fss{GR}}$ of groups such that for every $\pair{L}{G} \in
K^{\fss{LN2}}$ and $\vecg \in |G|^n$:
\[
\fs{ACT}(\fs{EP}(L,G),G) \models \varphi [\vecg] \: \mbox{ iff } \:\: G
\models \psi_\varphi [\vecg].
\]
In particular each of the following can be expressed by a first order
formula in $\CL^{\fss{GR}}$:

{\rm (1)} $g$ is signed.

{\rm (2)} $g_1$ and $g_2$ have the same sign.

{\rm (3)} $\fs{Supp}(g)$ is bounded.

{\rm (4)} $\fs{Supp}(g)$ is bounded either above or below, that is, $g
\in
\fs{Lft}(G) \cup \fs{Rt}(G)$.

{\rm (5)} $g_1$ and $g_2$ are both in $\fs{Lft}(G)$ or both in
$\fs{Rt}(G)$.

{\rm (6)} $\fs{ACT}(\fs{EP}(L,G),G)$ is transitive; or more generally,
$\fs{ACT}(\fs{EP}(L,G),G)$ has precisely $n$ orbits (where
$n$ is a given natural number).

\end{cor}

\noindent
{\bf Note}*
There is also a second order version of this corollary in which
$\fs{EP}(L,G)$ is replaced by $\barL^{\smallr}$ and in which for first
order $\varphi$ the formula $\psi_{\varphi}$ is second order.  By
(6), there is a second order sentence in $\CL^{\fss{GR}}$
expressing, for {\em transitive} $\pair{L}{G} \in K^{LN2}$, that $L$ is
Dedekind complete.

\bigskip

\noindent
{\bf Deduction of Corollary \ref{FOP} from Theorem~\ref{LIT} }
In view of
Proposition~\ref{REEXPRESS}, Theorem~\ref{LIT} tells us that the
desired formula $\psi _\varphi(x_1,\ldots,x_n)$ exists.
\hfill\qed

\begin{cor}
\label{BDDNOTEE}
Let $\pair{L}{G} \in K^{\fss{LN2}}$.  Then $G \not\equiv
\fs{Bdd}(G)$ unless $G=\fs{Bdd}(G)$, and similarly for $\fs{Lft}(G)$
and for $\fs{Rt}(G)$.
\end{cor}

\noindent
{\bf Deduction of Corollary \ref{BDDNOTEE} from
Corollary~\ref{FOP} }
``$G = \fs{Bdd}(G)$'' can be expressed first order in $\CL^{GR}$,
establishing the result for \fs{Bdd}.  So can ``$G = \fs{Lft}(G)$ or $G =
\fs{Rt}(G)$'' by (4) and (5) of Corollary~\ref{FOP}.  Now if $G \equiv
\fs{Lft}(G)$ but $G \neq \fs{Lft}(G)$,
then $G = \fs{Rt}(G)$, making $\fs{Lft}(G) = \fs{Bdd}(G)$.  Now the
result for \fs{Bdd} gives the result for \fs{Lft}. \hfill\qed
\medskip

Corollary~\ref{CHAR} also has a two-group version, of which the next
corollary is a first order analogue.  The proof of Corollary~\ref{LR}
will be
given after
Theorem~\ref{BIG}.

\begin{cor}
\label{LR}
Let $\pair{L_i}{G_i} \in K^{\fss{LN2}}$, $i=1,2$.  Suppose $G_1 \equiv
G_2$.

{\rm (1)} Then
\[
\fs{ACT}(\fs{EP}(L_1,G_1),\fs{Bdd}(G_1))  \equiv
\fs{ACT}(\fs{EP}(L_2,G_2),\fs{Bdd}(G_2)).
\]

{\rm (2)} Let $M_i^{\fss{Lft}}$ be the structure
$\fs{ACT}(\fs{EP}(L_i,G_i),\fs{Lft}(G_i))$, and similarly for
$M_i^{\fss{Rt}}$. Then
\[
M_1^{\fss{Lft}} \equiv M_2^{\fss{Lft}} \mbox{ and } M_1^{\fss{Rt}}
\equiv M_2^{\fss{Rt}},
\mbox{ or }
M_1^{\fss{Lft}} \equiv M_2^{\fss{Rt}} \mbox{ and } M_1^{\fss{Rt}}
\equiv M_2^{\fss{Lft}}.
\]
\end{cor}

Caution:  $\pair{L_i}{G_i} \in K^{\fss{LN2}}$ does not imply that
$\pair{L_i}{\fs{Lft}(G_i)} \in
K^{\fss{LN2}}$, or even that it has dense orbits in $\barL$; see
Example~\ref{WR1}.

\begin{defn}
\label{LFOR}
\begin{rm}
Let $K$ be a class of locally moving linear permutation groups.

{\rm (1)}  We say that {\em First Order Reconstruction holds for $K$}
if Theorem~\ref{LRT} through Corollary~\ref{LR} hold
when $K^{LN2}$ is replaced by $K$.

$\!\!\!${\rm *(2)} We say that {\em Second Order Reconstruction
holds for $K$} if First Order Reconstruction holds except for:
Theorem~\ref{LIT}(1), Corollary~\ref{FOP}
(only its Note remains valid), and
statements dealing with elementary equivalence.
\end{rm}
\end{defn}

In this terminology we have

\begin{theorem}
\label{LNRT}
{\rm (Linear Reconstruction Package)} \
First Order Reconstruction holds for linear permutation groups
$\pair{L}{G}$ which are
2-interval-transitive and have a nonidentity bounded element.
\end{theorem}

Here 2-interval-transitivity cannot be reduced to
1-interval-transitivity
(or to 1-{\it o-}transitivity); this does not guarantee even
Second Order Reconstruction.  See
Example~\ref{WRS2}.

We shall find other instances in which First Order Reconstruction
holds---see Theorems \ref{NESTRT}(1) and \ref{SPANRT}(1).  For
instances of Second Order Reconstruction, see
Theorems \ref{NESTRT}(2) and \ref{SPANRT}(2).

\bigskip

We turn now to reconstruction of {\bf monotonic permutation
groups}
 $\pair{\fs{ED}(L)}{G}$.  Under appropriate hypotheses, the above
results (including the accompanying notes) all carry over to the
monotonic case.  In Theorem~\ref{LRT},
$\fnn{\tau}{\barL_1}{\barL_2}$
is still a map between chains, and
$\iso{\tau}{\barL_1^{\smallr}}{\barL_2^{\smallr}}$ a map between
equal direction structures.  (Here $L$ is the chain on which
$\fs{ED}(L)$ is based.   Reversal of the order on $L$ interchanges
$\fs{Lft}(G)$ and $\fs{Rt}(G)$, but the matters
below, except for the half-order, are otherwise unaffected.)
Of course $\barL^{\smallr} = \overline{L^{\smallr}} = \fs{ED}(\barL) =
\overline{\fs{ED}(L)}$.

An orbit $M$ of $\pair{\fs{ED}(\barL)}{G}$ is either a single orbit of
$\pair{\barL}{\fs{Opp}(G)}$ or the union $M = M_1 \cup M_2$ of two
such orbits (with $M_2 \cong M_1^*$), and in the latter case $M$
must be symmetric (that is, $M \cong M^*$).  If $G$ contains an
order-reversing permutation, all orbits of $\pair{\fs{ED}(\barL)}{G}$
are symmetric.

For $\pair{\fs{ED}(L)}{G}$, the pointwise order on $G$ produces a
half-ordered group (see Definition~\ref{HALFORDERED}).  This order
is reversed by reversal of the order on $L$.

We deal with monotonic permutation groups $\pair{\fs{ED}(L)}{G}$
which are 3-interval-transitive and have a nonidentity bounded
element.  Unlike the hypotheses of the Linear Reconstruction
Theorem, these hypotheses imply {\em high} interval-transitivity
(Theorem~\ref{3H}), which in turn implies high approximate
order-transitivity for $\pair{\fs{ED}(\barL)}{G}$
(Proposition~\ref{HAT}).

\begin{theorem}
\label{MNRT}
{\rm (Monotonic Reconstruction Package)} \
For monotonic permutation groups $\pair{\fs{ED}(L)}{G}$ which are
3-interval-transitive and have a nonidentity bounded element, the
linear results Theorem~\ref{LRT} through Corollary~\ref{LR}
(with $K^{\fss{LN2}}$ replaced by $K^{\fss{MN3}}$ and ``partially
ordered group'' by ``half-ordered group'') all hold verbatim.

In particular, every automorphism $\alpha$ of $G$ is conjugation by
a unique automorphism  $\tau$ of $\barL^{\smallr}$, so that
\[
\fs{Aut}(G) =
\fs{Norm}_{\mbox{\scriptsize {\fs{Aut}}}(\barL^{\mbox
{\smallr}})}(G) =
\fs{Norm}_{\mbox{\scriptsize {\fs{Homeo}}}(\barL)}(G).
\]

Moreover: $\fs{Opp}(G)$ is a characteristic subgroup of $G$ and
$\fs{Aut}(G) \leq \fs{Aut}\mbox{\large(}\fs{Opp}(G)\mbox{\large)}
\leq \fs{Aut}\mbox{\large(}\fs{Bdd}(G)\mbox{\large)}$; ``$g$
preserves order'' can be expressed in $\CL^{\fss{GR}}$
in Corollary~\ref{FOP}ff; Corollary~\ref{BDDNOTEE} holds for
$\fs{Opp}(G)$; and
part (1) of Corollary~\ref{LR} holds for
$\fs{ACT}(\fs{EP}(\fs{ED}(L),G),\fs{Opp}(G))$.
\end{theorem}

\noindent
{\bf Note }
When $L$ is Dedekind complete, $\pair{\fs{ED}(L)}{G}$ satisfies the
hypotheses of Theorem~\ref{MNRT} provided it is {\em
approximately} 2-{\it o-}transitive and has a nonidentity bounded
element (Theorem~\ref{3H}).

\begin{defn}
\label{MFOR}
\begin{rm}
In the monotonic case,
``First Order Reconstruction'' refers to the results listed in
Theorem~\ref{MNRT}, and ``Second Order Reconstruction'' involves
the same exceptions as in Definition~\ref{LFOR}.
\end{rm}
\end{defn}

Again the (analogue of the) Linear Interpretation Theorem~\ref{LIT}
is the fundamental result, and the others all follow from it just as in
the linear case.

Here 3-interval-transitivity cannot be reduced to
2-interval-transitivity
(or to 2-{\it o-}transitivity); this does not guarantee even
Second Order Reconstruction. See
Example~\ref{WRSN}.

\bigskip

Next we consider reconstruction of {\bf circular permutation groups}
$\pair{C}{G}$, where $C = \fs{CR}(L)$.  Again under appropriate
hypotheses and with a few shifts in terminology, the results from the
linear case almost all carry over.  Of course $\barC^{\smallr}=
\overline{C^{\smallr}} = \fs{EO}(\barC) = \overline{\fs{EO}(C)}$.

One difference is that in the circular case $\fs{Lft}(G)$ and
$\fs{Rt}(G)$ are not defined.  Also we no longer have a ``pointwise
order'' on $G$.
Nevertheless, we define ``signed'' elements in the language
$\CL^{\fss{MCPG}}$ of monocircular permutation groups
$\pair{\fs{EO}(C)}{G}$.  We call an element $g \in G$ {\em signed} if
$g^2 \neq \fs{Id}$ and
\[
\fs{Eo}(g^2(x_1), g(x_1), x_1; g^2(x_2), g(x_2), x_2) \mbox{ for all }
x_1, x_2 \in \fs{supp}(g^2) \cap C.
\]
Signed elements are necessarily in $\fs{Opp}(G)$.
An element $g \in \fs{Opp}(G)$ such that $g^2 \neq \fs{Id}$
is necessarily signed if it has just one supporting interval.  (By a
{\em supporting interval} is meant a $\barC$-interval maximal with
respect to being contained in $\fs{supp}(g)$.)
We say that signed elements $g_1$ and $g_2$ have the {\em same
sign} if
\[
\fs{Eo}(g_1^2(x_1), g_1(x_1), x_1; g_2^2(x_2), g_2(x_2), x_2) \mbox{
for all } x_i \in C \cap \fs{supp}(g_i^2).
\]

The circular reconstruction results deal with circular permutation
groups which are 3-interval-transitive and have a nonidentity
bounded element.  Presumably these hypotheses do not imply
high interval-transitivity.  This would parallel the linear (but not the
monotonic) case.

\begin{theorem}\label{CRT}
{\rm (Circular Reconstruction Package)} \
For circular permutation groups $\pair{C}{G}$ which are
3-interval-transitive and have a nonidentity bounded element, the
linear results Theorem~\ref{LRT} through Corollary~\ref{FOP} all
hold with the obvious changes---that is, references to
$\fs{Lft}(G)$ and $\fs{Rt}(G)$
and to pointwise order are to be omitted, and
$L$ is to be replaced by $C$, ``order'' by ``orientation'', ``chain'' by
``circle'',
``monotonic'' by ``monocircular'', ``equal direction structure'' by
``equal orientation structure'', $\CL^{\fss{MNPG}}$ by
$\CL^{\fss{MCPG}}$, and $K^{\fss{LN2}}$ by $K^{\fss{CR3}}$.

In particular,
every automorphism $\alpha$ of $G$ is conjugation by a unique
automorphism  $\tau$ of $\barC^{\smallr}$, so that
\[
\fs{Aut}(G) =
\fs{Norm}_{\mbox{\scriptsize {\fs{Aut}}}(\barC^{\mbox
{\smallr}})}(G) =
\fs{Norm}_{\mbox{\scriptsize {\fs{Homeo}}}(\barC)}(G).
\]
\end{theorem}

Just as in the previous cases,  the analogue of Theorem~\ref{LIT}
yields all the other results.
Note that  Corollaries \ref{BDDNOTEE} and \ref{LR} do not have
circular analogues, as
the set of bounded elements is not a subgroup of $G$.

\begin{defn}
\label{CFOR}
\begin{rm}
In the circular case,
``First Order Reconstruction'' refers to the results listed in
Theorem~\ref{CRT}, and ``Second Order Reconstruction'' involves the
same exceptions as in Definition~\ref{LFOR}.
\end{rm}
\end{defn}

\begin{question}
\label{QCIR3TO2}
\begin{rm}
Can 3-interval-transitivity be reduced to 2-interval-transitivity in
Theorem~\ref{CRT}, perhaps at the expense of obtaining Second
Order Reconstruction rather than First?
\end{rm}
\end{question}

\smallskip

Finally we turn to reconstruction of {\bf monocircular permutation
groups}
$\pair{\fs{EO}(C)}{G}$, where $C$ is the circle on which $\fs{EO}(C)$ is
based; the matters below are unaffected by reversal of orientation on
$C$.
Under appropriate hypotheses, the circular results carry over to the
monocircular case.
The orbits of $\pair{\fs{EO}(\barC)}{G}$ are related to those of
$\pair{\barC}{\fs{Opp}(G)}$ as those of $\pair{\fs{ED}(\barL)}{G}$
are to those of $\pair{\barL}{\fs{Opp}(G)}$.

We deal with monocircular permutation groups $\pair{\fs{EO}(C)}{G}$
which are 4-interval-transitive and have a nonidentity bounded
element.  Paralleling the monotonic case, these hypotheses imply
{\em high} interval-transitivity, which in turn implies high
approximate orientation-transitivity for $\pair{\fs{EO}(\barC)}{G}$.

\begin{theorem}
\label{MCRT}
{\rm (Monocircular Reconstruction Package)} \
For monocircular permutation groups $\pair{\fs{EO}(C)}{G}$ which
are
4-interval-transitive and have a nonidentity bounded element,
Theorem~\ref{CRT}
(with $K^{\fss{CR3}}$ replaced by $K^{\fss{MC4}}$) holds verbatim.

In particular,
every automorphism $\alpha$ of $G$ is conjugation by a unique
automorphism  $\tau$ of $\barC^{\smallr}$, so that
\[
\fs{Aut}(G) =
\fs{Norm}_{\mbox{\scriptsize {\fs{Aut}}}(\barC^{\mbox
{\smallr}})}(G) =
\fs{Norm}_{\mbox{\scriptsize {\fs{Homeo}}}(\barC)}(G).
\]

Moreover: $\fs{Opp}(G)$ is a characteristic subgroup of $G$;
``$g$ preserves orientation'' can be expressed in $\CL^{\fss{GR}}$
in Corollary~\ref{FOP}ff; Corollary~\ref{BDDNOTEE} holds for
$\fs{Opp}(G)$; and part (1) of Corollary~\ref{LR} holds for
$\fs{ACT}(\fs{EP}(\fs{EO}(C),G),\fs{Opp}(G))$.
\end{theorem}

\noindent
{\bf Note }
When $C$ is Dedekind complete, $\pair{\fs{EO}(C)}{G}$ satisfies the
hypotheses of Theorem~\ref{MCRT} provided it is {\em
approximately} 3-{\it o-}transitive and has a nonidentity bounded
element (Theorem~\ref{3H}).
\medskip

Just as in the previous cases,  the analogue of Theorem~\ref{LIT}
yields all the other results.

\begin{defn}
\label{MCFOR}
\begin{rm}
In the monocircular case,
``First Order Reconstruction'' refers to the results listed in
Theorem~\ref{MCRT}, and ``Second Order Reconstruction'' involves
the same exceptions as in Definition~\ref{LFOR}.
\end{rm}
\end{defn}

\begin{question}
\label{QMONOCIR4TO3}
\begin{rm}
Can 4-interval-transitivity be reduced to 3-interval-transitivity in
Theorem~\ref{MCRT}?
\end{rm}
\end{question}

\bigskip

\noindent
{\bf The proofs}

\smallskip

We begin the proofs with some matters common to all four cases.
Recall that $\fs{Var}(N,G)$ is an abbreviation for
$\fs{Var}(\barB(N),G)$.

\begin{defn}
\label{}
\begin{rm}
Let $\pair{N}{G} \in K^{NO}$, which ensures that $G$ acts faithfully
on $\fs{Var}(N,G)$.

(1) Let $\fs{VACT}(N,G) \eqdf
\fs{ACT}(\fs{Var}(N,G),G) = \fs{VACT}(\barB(N),G)$.

$\!\!\!${\rm *(2)} Let $\fs{BACT}(N,G) \eqdf
\fs{ACT}(\barB(N),G) = \fs{BACT}(\barB(N),G)$.

\end{rm}
\end{defn}

\begin{defn}
\label{}
\begin{rm}
Let $K \subseteq K^{\fss{NO}}$.

(1) Let $K_{\fss{VACT}} \eqdf
\setm{\fs{VACT}(N,G)}{\pair{N}{G} \in K}$.

$\!\!\!${\rm *(2)} Let $K_{\fss{BACT}} \eqdf
\setm{\fs{BACT}(N,G)}{\pair{N}{G} \in K}$.

\end{rm}
\end{defn}

As a special case of Theorem~\ref{LMSRT} we have

\begin{cor}
\label{NOIT}\

{\rm (1)}  $K^{\fss{NO}}_{\fss{VACT}}$ is FO-STR-interpretable in
$K^{\fss{NO}}_{\fss{GR}}$.

$\!\!\!${\rm *(2)} $K^{\fss{NO}}_{\fss{BACT}}$ is
SO-STR-interpretable in
$K^{\fss{NO}}_{\fss{GR}}$.

\end{cor}

\smallskip
\noindent
{\bf The linear case}

\smallskip

We embark on the proof of Theorem~\ref{LIT}(1), which will be
concluded after Theorem~\ref{BIG}.  It may help to keep in mind the
overview presented in Section~\ref{OVERVIEW}.

\begin{defn}\label{}
\begin{rm}
Let $\pair{L}{<}$ be a dense chain without endpoints.

(1) For $U \subseteq \barL$, let $\fs{ch}^{\barL}(U)$ denote the
convex hull of $U$ in $\barL$. We abbreviate $\fs{ch}^{\barL}(U)$ by
$\fs{ch}(U)$.

(2) Let $U_1,U_2 \subseteq \barL$.  We say that $U_1$ and $U_2$
are {\em segregated} if $\fs{ch}(U_1) \cap \fs{ch}(U_2) =
\emptyset$, that is, if $U_1 < U_2$ or $U_2 < U_1$ or $U_1 =
\emptyset$ or $U_2 = \emptyset$.  Let
\[
\fs{Seg}^{L}  \eqdf \setm{\pair{U_1}{U_2} \in
\barB(L)^2}{\fs{ch}(U_1) \cap \fs{ch}(U_2) = \emptyset}.
\]

(3) Let  $U_1,U_2,U_3 \subseteq \barL$.  We write
$\fs{Bets}(U_1,U_2,U_3)$ if each $U_i \neq \emptyset$, and if
$\fs{Bet}(x_1,x_2,x_3)$ whenever each $x_i \in U_i$ (or
equivalently, if
$x_1< x_2< x_3$ whenever each $x_i \in U_i$, or else $x_3< x_2< x_1$
whenever each $x_i \in U_i$).

(4) These notions are unaffected by reversal of the order on $L$, and
thus are also applicable to $\fs{ED}(L)$.
\end{rm}
\end{defn}

We now show that for $\pair{L}{G} \in K^{\fss{LN2}}$ the relation
$\fs{Seg}^{L}$ on $\fs{Var}(L,G)$, or even on $\barB(L)$, is definable
by a first order formula in $\CL^{\fss{VACT}}$.  This will be our sole
use of
2-interval-transitivity (aside from permitting us to assume only the
existence of a single bounded element).
Recall that $\fs{var}(g) = \fs{int}(\fs{cl}(\fs{supp}(g)))$.  For a
restatement of Proposition~\ref{LSEG} in terms of interpretability,
see Proposition~\ref{LSEGTRANSLATION}.

\begin{prop}
\label{LSEG}
Let
\[
\varphi_{\fss{Seg}}^{\fss{LN2}}(u_1,u_2) \eqqdf
\neg(\exists g)\mbox{\large(}(g(u_1) \cdot u_2 \neq 0) \wedge
(g(u_2) \cdot u_1 \neq 0)\mbox{\large)},
\]
re-expressed as a first order formula in $\CL^{\fss{VACT}} =
\CL^{\fss{BACT}}$.  Let $\pair{L}{G} \in K^{\fss{LN2}}$.  Then

{\rm (1)} For any $U_1,U_2 \in \fs{Var}(L,G)$:
\[
\fs{VACT}(L,G) \models \varphi_{\fss{Seg}}^{\fss{LN2}}[U_1,U_2] \:
\mbox{ iff } \: \fs{Seg}^{L}(U_1,U_2).
\]

$\!\!\!${\rm *(2)} For any $U_1,U_2 \in \barB(L)$:
\[
\fs{BACT}(L,G) \models \varphi_{\fss{Seg}}^{\fss{LN2}}[U_1,U_2] \:
\mbox{ iff } \: \fs{Seg}^{L}(U_1,U_2).
\]

\end{prop}

\noindent
{\bf Proof } $\varphi_{\fss{Seg}}^{\fss{LN2}}$ can be re-expressed as
specified because of Proposition~\ref{PST}.

(1) Let $\pair{L}{G} \in K^{\fss{LN2}}$
and $U_1,U_2 \in \fs{Var}(L,G)$.
Obviously if
$\fs{Seg}^{L}(U_1,U_2)$
then
$\fs{VACT}(L,G) \models \varphi_{\fss{Seg}}^{\fss{LN2}}[U_1,U_2]$.

Suppose that
$\fs{Seg}^{L}(U_1,U_2)$ does not hold.
We may assume that there are open intervals
$I < J < K$ of $L$ such that $I,K \subseteq U_1$ and $J \subseteq
U_2$.
2-interval-transitivity provides $g \in G$ such that $g(I) \cap J \neq
\emptyset$
and $g(J) \cap K \neq \emptyset$, and thus
$g(U_1) \cap U_2 \neq \emptyset$ and $g(U_2) \cap U_1 \neq
\emptyset$.
Hence $\fs{VACT}(L,G) \not\models
\varphi_{\fss{Seg}}^{\fss{LN2}}[U_1,U_2]$.

(2) The proof is virtually identical to that of (1).
\hfill\qed
\medskip

Theorem~\ref{BIG} will say that for any class $K$ of locally moving
linear permutation groups $\pair{L}{G}$, in order to capture
$\fs{EPACT}(L,G)$ in the group $G$ it is enough to capture the
segregation relation in $\fs{VACT}(L,G)$.   Applied to $K^{\fss{LN2}}$
in combination with Proposition~\ref{LSEG}, Theorem~\ref{BIG} will
give Theorem~\ref{LIT}; later it will be applied to other $K$'s.  Most
of its proof is contained in Lemma~\ref{DDLN}, in which the group
$G$ does not play any role.  This lemma will deal with doubly dense
subsets $T$ of $\barB(L)$
(see Definition~\ref{DEFDDNS}).
It will say that one can capture in the structure
$\fs{SEG}(L,T) \eqdf \trpl{T}{\leq^{\barB(L)}}{\fs{Seg}^L}$ a certain
dense subset of $\barL$, namely the set $\fs{Ep}(L,T)$ of  endpoints
of convex hulls of members of $T$, along with the equal direction
structure $\fs{EP}(L,T)$ induced on it by $\barL$.
When specialized to $T = \fs{Var}(L,G)$, $\fs{EP}(L,T)$ will give
$\fs{EP}(L,G)$.

\begin{defn}\label{DEFLNNOGR}
\begin{rm} Let $\pair{L}{<}$ be a dense chain without endpoints and
let $T$ be a dense subset of $\barB(L)$.

(1) Let
\begin{eqnarray*}
& \fs{SEG}(L,T) \eqdf \trpl{T}{\leq^{\barB(L)}}{\fs{Seg}^L},  \\*
& \CL^{\fss{SEG}} \eqdf \CL(\fs{SEG}(L,T)) = \{ \leq,\fs{Seg} \}.
\end{eqnarray*}

(2) Let
\begin{eqnarray*}
& \fs{Ep}(L,T) \eqdf \bigcup_{U \in T} \fs{Ep}(U),  \\*
& \fs{EP}(L,T) \eqdf \pair{\fs{Ep}(L,T)}{\fs{Ed}^{\barL}}, \\*
& \CL^{\fss{EP}} \eqdf \CL(\fs{EP}(L,T)) = \{\fs{Ed} \}.
\end{eqnarray*}
$\fs{Ep}(L,T)$ is dense in $\barL$ because $T$ is dense in
$\barB(L)$.

\smallskip
Now we prepare to represent the points of $\fs{Ep}(L,T)$
by pairs of members of $T$:

\smallskip

(3) Let $\fs{Rep}(L,T) \eqdf
\setm{\pair{U_1}{U_2} \in T^2}{
\fs{Seg}^{L}(U_1,U_2) \mbox{ and }
U_1,U_2 \neq \emptyset}$.

(4) Let $\fnn{\fs{Pt}^{L,T}}{\fs{Rep}(L,T)}{\fs{Ep}(L,T)}$ be the
surjection defined by
\[
\fs{Pt}^{L,T}(U_1,U_2) = \left\{
\begin{array}{ll}
\fs{sup}(U_1) & \mbox{if $U_1 < U_2$,} \\*
\fs{inf}(U_1)  & \mbox{if $U_1 > U_2$.}
\end{array}
\right.
\]

(5) Let
\begin{eqnarray*}
& \fs{SEGEP}(L,T) \eqdf (\fs{SEG}(L,T),\fs{EP}(L,T);\fs{Pt}^{L,T}), \\*
& \CL^{\fss{SEGEP}} \eqdf \CL(\fs{SEGEP}(L,T)) =
 \{ \leq,\fs{Seg}, \fs{Ed} , \fs{Pt} , P_1 , P_2, Q \}.
\end{eqnarray*}
The universe of $\fs{SEGEP}(L,T)$ is $T \cup \fs{Ep}(L,T)$, and $P_1$
and $P_2$ are
1-place relation symbols distinguishing $T$ and $\fs{Ep}(L,T)$.
$Q$ is a 2-place relation symbol, interpreted as
$\fs{Dom}(\fs{Pt}^{L,T})$, which according to our convention is added
when the partial function $\fs{Pt}^{L,T}$ is transformed into a total
function.

\end{rm}
\end{defn}

Because $\fs{EP}(L,T)$ is an equal direction structure rather than a
chain, nothing in Definition~\ref{DEFLNNOGR} is affected by reversal
of the order on $L$.  Hence these notions (and Lemma~\ref{DDLN}
and its proof) are also applicable to $\fs{ED}(L)$.

\begin{defn}
\label{DEFDDLN}
\begin{rm}\

(1) Let $K^{\fss{DDLN}}$ be the class of all pairs $\pair{L}{T}$ such
that $L$ is a dense chain without endpoints and $T$ is a doubly
dense subset of $\barB(L)$.

(2) Let $K^{\fss{DDLN}}_{\fss{SEG}} \eqdf
\setm{\fs{SEG}(L,T)}{\pair{L}{T} \in K^{\fss{DDLN}}}$.

(3) Let $K^{\fss{DDLN}}_{\fss{SEGEP}} \eqdf
\setm{\fs{SEGEP}(L,T)}{\pair{L}{T} \in K^{\fss{DDLN}}}$.

\end{rm}
\end{defn}

\begin{lemma}
\label{DDLN}
$K^{\fss{DDLN}}_{\fss{SEGEP}}$
is FO-STR-interpretable in
$K^{\fss{DDLN}}_{\fss{SEG}}$, and
\newline
$\pair{K^{\fss{DDLN}}_{\fss{SEG}}}{K^{\fss{DDLN}}_{\fss{SEGEP}}}$
is a spanning system.
\end{lemma}

\noindent
{\bf Proof }
We need to capture several notions by first order formulas in
$\CL^{\fss{SEG}}$.

First we capture ``betweenness'' of sets $U \in T$.  Let
\begin{eqnarray*}
& \varphi_{\fss{Bets}}(u_1,u_2,u_3) \eqqdf  \\*
&
\begin{displaystyle}
\mbox{\Large(} \! \bigwedge_{1 \leq i \leq3} \!u_i \neq 0
\mbox{\Large)}  \wedge
\mbox{\Large(} \!\bigwedge_{1 \leq i < j \leq 3} \!\! \fs{Seg}(u_i,u_j)
\,\,\wedge \,\, \exists w \mbox{\large (}(w \leq u_1 + u_3)
\wedge \neg\fs{Seg}(w,u_2)\mbox{\large )} \mbox{\Large)}.
\end{displaystyle}
\end{eqnarray*}
Then for any $\pair{L}{T} \in K^{\fss{DDLN}}$, double density
guarantees that for any $U_1,U_2,U_3 \in T$:
\[
\fs{SEG}(L,T) \models \varphi_{\fss{Bets}}[U_1,U_2,U_3]
\mbox{ \ iff \ } \fs{Bets}^{\barL}(U_1,U_2,U_3).
\]
An application of Proposition~\ref{PST} allows us to re-express
$\varphi_{\fss{Bets}}$ in $\CL^{\fss{SEG}}$.

Obviously we can capture $\fs{Rep}(L,T)$.  Just let
\[
\varphi_{\fss{Rep}}(u_1,u_2) \eqqdf
\fs{Seg}(u_1,u_2) \wedge  (u_1 \neq 0) \wedge (u_2 \neq 0).
\]
Then for any $\pair{L}{T} \in K^{\fss{DDLN}}$ and $U_1,U_2 \in T$:
\[
\fs{SEG}(L,T) \models \varphi_{\fss{Rep}}[U_1,U_2]   \mbox{ \ iff \ }
\pair{U_1}{U_2} \in \fs{Rep}(L,T).
\]

We want a first order formula $\varphi_{\fss{Eqp}}(\vecu,\vecv)$ in
$\CL^{\fss{SEG}}$
such that for every $\pair{L}{T} \in K^{\fss{DDLN}}$
and $\vecU, \vecV \in T^2$:
\begin{eqnarray*}
& \fs{SEG}(L,T) \models  \varphi_{\fss{Eqp}}[\vecU,\vecV] \mbox{ \
iff \ } \\*
& \vecU \in \fs{Rep}(L,T) \mbox{ and } \vecV \in \fs{Rep}(L,T)
\mbox{ and }
\fs{Pt}^{L,T}(\vecU) = \fs{Pt}^{L,T}(\vecV).
\end{eqnarray*}

For $U_1,U_2,U'_2 \subseteq \barL$ with $U_1 \neq \emptyset$, we
say that $U_2$ and $U'_2$ are on the {\em same side} of $U_1$ if
\begin{eqnarray*}
& U_2 < U_1 \mbox{ (or $U_2 = \emptyset$)} \mbox{ and } U'_2 <
U_1
\mbox{ (or $U'_2 = \emptyset$)}
\mbox{, or} \\*
& U_2 > U_1 \mbox{ (or $U_2 = \emptyset$)} \mbox{ and } U'_2 >
U_1   \mbox{ (or $U'_2 = \emptyset$)}
  \mbox{;}
\end{eqnarray*}
or equivalently, if
\[
\fs{Seg}^{L}(U_1,U_2) \mbox{, } \fs{Seg}^{L}(U_1,U'_2) \mbox{, and }
U_1 \mbox{ is not between } U_2 \mbox{ and } U'_2.
\]
Clearly this can be captured in $\CL^{\fss{SEG}}$.

For $\pair{U_1}{U_2} \in \fs{Rep}(L,T)$, let
\[
\fs{Sd}(U_1,U_2) \eqdf
\setm{U'_2 \in T}{U_2 \mbox{ and } U'_2 \mbox{ are on the same
side of }U_1},
\]
which we call the $U_2${\em -side} of $U_1$.
Now for $\vecU,\vecV \in \fs{Rep}(L,T)$:
\[
\fs{Pt}^{L,T}(\vecU) = \fs{Pt}^{L,T}(\vecV)   \mbox{ \ iff \ }
\mbox{\large(}\fs{Sd}(\vecU) = \fs{Sd}(\vecV) \mbox{\large)}
\,\vee\, \mbox{\large(}
\fs{Sd}(\vecU) = \fs{Bcmp}(\fs{Sd}(\vecV))\mbox{\large)},
\]
where the {\em Boolean complement} $\fs{Bcmp}(\fs{Sd}(\vecV))$
means
\[
\setm{W \in T}{W \cdot V'_2 = \emptyset \mbox{ for all } V'_2 \in
\fs{Sd}(V_1,V_2)}.
\]
 From this it is easy to see that the desired first order
$\varphi_{\fss{Eqp}}$ can be constructed.

We can also capture in $\CL^{\fss{SEG}}$ the restriction to
\fs{Ep}(L,T) of the betweenness relation $\fs{Bet}^{\barL}$ on
$\barL$.  We seek a first order formula
$\varphi_{\fss{Bet}}(\vecu,\vecv,\vecw)$ in $\CL^{\fss{SEG}}$ such
that for any $\pair{L}{T} \in K^{\fss{DDLN}}$ and any
$\vecU,\vecV,\vecW \in \fs{Rep}(L,T)$:
\[
\fs{SEG}(L,T) \models \varphi_{\fss{Bet}}[\vecU,\vecV,\vecW]
\mbox{ \ iff \ }
\fs{Bet}^{\barL}(\fs{Pt}^{L,T}(\vecU),\fs{Pt}^{L,T}(\vecV),\fs{Pt}^{L,T}(
\vecW)).
\]
Let $\fs{Sd}_1(\vecU) = \fs{Sd}(\vecU)$ and
      $\fs{Sd}_2(\vecU) = \fs{Bcmp}(\fs{Sd}(\vecU))$.  Then
\begin{eqnarray*}
& \fs{Bet}^{\barL}
(\fs{Pt}^{L,T}(\vecU),\fs{Pt}^{L,T}(\vecV),\fs{Pt}^{L,T}(\vecW))
\mbox{ iff}  \\*
& {\displaystyle \bigvee} \setm{\fs{Sd}_i(\vecU) \subsetneqq
 \fs{Sd}_j(\vecV) \subsetneqq  \fs{Sd}_k(\vecW)}
{i,j,k \in \dbltn{1}{2}}.
\end{eqnarray*}
 From this it is clear that the desired $\varphi_{\fss{Bet}}$  can be
constructed.  Then, because of the bi-interpretability of
$\fs{Bet}^{\barL}$ and $\fs{Ed}^{\barL}$ (Proposition~\ref{BETW}),
there is a first order formula
$\varphi_{\fss{Ed}}(\vecx_1,\vecx_2,\vecx_3,\vecx_4)$ in
$\CL^{\fss{SEG}}$ which similarly captures the restriction to
$\fs{Ep}(L,T)$ of the equal direction structure $\fs{ED}(\barL)$.

The formulas $\varphi_{Rep}$,
$\varphi_{Eqp}$, and
$\varphi_{Ed}$ provide the building blocks for the interpretation.
The interpreting mapping $\mu:T^2 \stackrel{onto}{\raro} T \cup
\fs{Ep}(L,T)$
will be given by
\[
\mu(U,V) = \left\{
\begin{array}{ll}
U & \mbox{if $U=V$,} \\*
\fs{Pt}(U,V)  & \mbox{if $\pair{U}{V} \in \fs{Rep}(L,T)$.}
\end{array}
\right.
\]
Some of the formulas of the interpretation are:

\smallskip
\noindent
$\varphi_{U}(u_1,u_2) \eqqdf (u_1=u_2) \vee
\varphi_{Rep}(u_1,u_2)$.

\smallskip
\noindent
$\varphi_{Eq}(\vecu,\vecv) \eqqdf (u_1=u_2=v_1=v_2) \vee
\varphi_{Eqp}(\vecu,\vecv)$.

\smallskip
\noindent
$\varphi_{Imap}(\vecu,v) \eqqdf u_1=u_2=v$.

\smallskip
\noindent
$\varphi_{Pt}(\vecu,\vecv,\vecw) \eqqdf (u_1=u_2) \wedge (v_1
= v_2)
\wedge  \varphi_{Eqp}(u_1,v_1;\vecw)$.

\smallskip

Finally, we have a spanning system because $\fs{Pt}^{L,T}$ is a
surjection onto $|\fs{SEGEP}(L,T)| - |\fs{SEG}(L,T)|$.
\hfill\qed
\smallskip

\begin{prop}
\label{LEP}
Let $\pair{L}{G}$ be a locally moving linear permutation group.
Then

{\rm (1)} $\fs{EP}(L,\fs{Var}(L,G)) = \fs{EP}(L,G)$.

$\!\!\!${\rm *(2)} $\fs{EP}(L,\barB(L)) = \barL^{\smallr}$.

\end{prop}

\noindent
{\bf Proof }
These are the special cases of $\fs{EP}(L,T)$ arising from
Definition~\ref{DEFLNNOGR} when
$T = \fs{Var}(L,G)$ and when $T = \barB(L)$.
\hfill\qed
\medskip

The following notational pattern will be used repeatedly:

\begin{defn}
\label{PATTERN}
\begin{rm}
Let $K \subseteq K^{\fss{NO}}$.

(1) For any structure $\fs{STRUC}(N,G)$ defined
for each $\pair{N}{G} \in K$, let $K_{\fss{STRUC}} \eqdf
\setm{\fs{STRUC}(N,G)}{\pair{N}{G} \in K}$.

(2) If also $G \leq \fs{Aut}(\fs{STRUC}(N,G))$ for each $\pair{N}{G}
\in K$,  let
$K_{\fss{STRUCACT}} \eqdf
\setm{\fs{ACT}(\fs{STRUC}(N,G),G)}{\pair{N}{G} \in K}$.

\end{rm}
\end{defn}

\begin{defn} \label{}
\begin{rm}
Let $K$ be a class of locally moving linear permutation
groups, and let $\pair{L}{G} \in K$.

(1) (a) Let $\fs{VSEG}(L,G) \eqdf \fs{SEG}(L,\fs{Var}(L,G)) =
\newline
\trpl{\fs{Var}(L,G)}{\leq ^{\barB(L)}}{\fs{Seg}^L}$.
Note that $\fs{ACT}(\fs{VSEG}(L,G),G) =
\newline
(\fs{VACT}(L,G);\fs{Seg}^L)$.

$\,\,\,\,\,\,\,\,\,\,$(b) Let $\fs{VSEGEP}(L,G) \eqdf
\fs{SEGEP}(L,\fs{Var}(L,G)) =
\newline
(\fs{VSEG}(L,G),\fs{EP}(L,G);\fs{Pt}^{L,\fss{Var}(L,G)})$.

$\,\,\,\,\,\,\,\,\,\,$(c) $K_{\fss{VSEG}}$, $K_{\fss{VSEGACT}}$,
$K_{\fss{VSEGEP}}$, and $K_{\fss{VSEGEPACT}}$ are defined in
accord with Definition~\ref{PATTERN}.

$\!\!\!$*(2) (a) Let $\fs{BSEG}(L,G) \eqdf
\fs{SEG}(L,\barB(L))$.

$\,\,\,\,\,\,\,\,\,\,$(b) Let $\fs{BSEGDC}(L,G) \eqdf
\fs{SEGEP}(L,\barB(L))) =
\newline
(\fs{BSEG}(L,G),\barL^{\smallr};\fs{Pt}^{L,\barB(L)})$.

$\,\,\,\,\,\,\,\,\,\,$(c) $K_{\fss{BSEG}}$, $K_{\fss{BSEGACT}}$,
$K_{\fss{BSEGDC}}$, and $K_{\fss{BSEGDCACT}}$ are defined in accord
with Definition~\ref{PATTERN}.
 \end{rm}
\vspace{-2.5mm}
\end{defn}

\begin{prop}
\label{LSEGTRANSLATION}\

{\rm (1)} $K^{\fss{LN2}}_{\fss{VSEGACT}}$
is FO-STR-interpretable in $K^{\fss{LN2}}_{\fss{VACT}}$.

$\!\!\!${\rm *(2)} $K^{\fss{LN2}}_{\fss{BSEGACT}}$
is FO-STR-interpretable in $K^{\fss{LN2}}_{\fss{BACT}}$.

\end{prop}

\noindent
{\bf Proof }
This is just a translation of Proposition~\ref{LSEG}.
\hfill\qed
\medskip

The next theorem (which includes Theorem~\ref{LIT} as a special
case) presents results which hold for classes of locally moving linear
permutation groups $\pair{L}{G}$ in which the segregation relation
$\fs{Seg}^{L}$ can be captured first order in $\fs{VACT}(L,G)$, or
even in $\fs{BACT}(L,G)$:

\begin{theorem}
\label{BIG}
Let $K$ be a class of locally moving linear permutation groups.

{\rm (1)} Suppose there is a first order formula
$\varphi_{\fss{Seg}}(u_1,u_2)$ in
$\CL^{\fss{VACT}}$ such that for any $\pair{L}{G} \in K$  and
$U_1,U_2 \in
\fs{Var}(L,G)$:
\[
\fs{VACT}(L,G) \models \varphi_{\fss{Seg}}[U_1,U_2]
\mbox{ \ iff \ }
\fs{Seg}^{L}(U_1,U_2).
\]
Then First Order Reconstruction holds for $K$ (see
Definition~\ref{LFOR}).

$\!\!\!${\rm *(2)} Suppose there is a {\em first} order formula
$\varphi_{\fss{Seg}}(u_1,u_2)$ in
$\CL^{\fss{BACT}}$ such that for any $\pair{L}{G} \in K$ and
$U_1,U_2 \in
\barB(L)$:
\[
\fs{BACT}(L,G) \models \varphi_{\fss{Seg}}[U_1,U_2]
\mbox{ \ iff \ }
\fs{Seg}^{L}(U_1,U_2).
\]
Then Second Order Reconstruction holds for $K$.

$\!\!\!${\rm *(3)}  The hypothesis in (1) implies that in (2).
\end{theorem}

\noindent
{\bf Proof }
(1) We know that
\[      K_{\fss{VACT}} \mbox{ is FO-STR-interpretable in }
K_{\fss{GR}}  \]
since this holds for any $K \subseteq K^{\fss{NO}}$ by
Corollary~\ref{NOIT}.
\[
K_{\fss{VSEGACT}} \mbox{ is FO-STR-interpretable in }
K_{\fss{VACT}}
\]
by the hypothesis about $\varphi_{\fss{Seg}}$.
By Lemma~\ref{DDLN} (with $T = \fs{Var}(L,G))$ and Propositions
\ref{DDNS} and \ref{LEP},
$K_{\fss{VSEGEP}}$ is FO-STR-interpretable in  $K_{\fss{VSEG}}$,
and $\pair{K_{\fss{VSEG}}}{K_{\fss{VSEGEP}}}$ is a spanning system.
Then by Proposition~\ref{ACTS},
\[
K_{\fss{VSEGEPACT}} \mbox{ is FO-STR-interpretable in }
K_{\fss{VSEGACT}}.
\]
By the transitivity of FO-STR-interpretations
(Proposition~\ref{TRANS}),
\[
K_{\fss{VSEGEPACT}} \mbox{ is FO-STR-interpretable in }
K_{\fss{GR}}. \]
Retaining only the information about the endpoint action (see
Proposition~\ref{FORMALIZE}),
\[
K_{\fss{EPACT}} \mbox{ is FO-STR-interpretable in } K_{\fss{GR}}.
\]

Replacement of $K^{LN2}$ by $K$ in Theorem~\ref{LIT}(1) gives a
generalization of \ref{LIT}(1), and the results
mentioned in (1) all follow from the generalization just as they did
for $K^{\fss{LN2}}$ (except for Theorem~\ref{LIT}(2), which follows
from (2) and (3) of the present theorem, and except for
Corollary~\ref{LR}, which is proved below).

(2) The proof parallels that of (1), using the second order parts of the
cited results.
The application of Proposition~\ref{TRANS}, establishing that
$K_{\fss{BSEGDCACT}}$ is SO-STR-interpretable in $K_{\fss{GR}}$,
is valid because only the interpretation of $K_{\fss{VACT}}$ in
$K_{GR}$ is second order, and that is why the formula
$\varphi_{\fss{Seg}}$ in (2) was assumed to be first order.

(3) For $U_1,U_2 \in \barB(L)$, $\fs{Seg}^{L}(U_1,U_2)$ iff
$\fs{Seg}^{L}(U'_1,U'_2)$ for all  $U'_1,U'_2 \in \fs{Var}(L,G)$
with $U'_i \leq U_i$, because of Proposition~\ref{BSUP}.
In view of Proposition~\ref{PHIVAR}(3), from a formula
$\varphi'_{\fss{Seg}}$ as in (1) it is easy to construct a formula
$\varphi_{\fss{Seg}}$ as in (2).
\hfill\qed
\bigskip

\noindent
{\bf Proof of Theorem~\ref{LIT} }
Because of Proposition~\ref{LSEG}, Theorem \ref{LIT} is a special
case of Theorem~\ref{BIG}.
\hfill\qed
\bigskip

\noindent
{\bf Proof of Corollary~\ref{LR}}

(1) In view of Corollary~\ref{FOP}, there is a first order formula
$\psi_{\fss{Bdd}}(g)$ in $\CL^{\fss{GR}}$ expressing that $g \in
\fs{Bdd}(G)$.

Now let $\varphi$ be any first order sentence in $\CL^{\fss{MNPG}}$.
Then
\[
\fs{ACT}(\fs{EP}(L_i,G_i),\fs{Bdd}(G_i)) \models \varphi \mbox{ \ iff
\ }
\fs{ACT}(\fs{EP}(L_i,G_i),G_i) \models \varphi^{\fss{Bdd}},
\]
where $\varphi^{\fss{Bdd}}$ is the sentence in $\CL^{\fss{MNPG}}$
obtained by using $\psi_{\fss{Bdd}}(g)$ to modify $\varphi$ so that
for each group variable $g$ of $\varphi$ the quantification is over
$\fs{Bdd}(G)$
rather than over $G$.  Then (1) follows from Theorem~\ref{LRT}(1).

(2) Let $\psi_{\fss{1way}}(g)$ and $\psi_{\fss{Smway}}(g,h)$ be first
order formulas in $\CL^{\fss{GR}}$, with $\psi_{\fss{1way}}(g)$
expressing that $g \in \mbox{\Large(}\mbox{\large(}\fs{Lft}(G) \cup
\fs{Rt}(G)\mbox{\large)} - \fs{Bdd}(G)\mbox{\Large)}$,
and $\psi_{\fss{Smway}}(g,h)$ expressing that $g$ and $h$ are both
in $\fs{Lft}(G)$ or both in $\fs{Rt}(G)$ (membership in $\fs{Bdd}(G)$
is permitted).

Again let $\varphi$ be any first order sentence in
$\CL^{\fss{MNPG}}$.  For $g \in G$ such that $\psi_{\fss{1way}}(g)$,
let $\varphi^{\fss{1way}}(g)$ be the formula in $\CL^{\fss{MNPG}}$
obtained by modifying $\varphi$ so that its group variables $h$ are
quantified over $\setm{h \in G}{\psi_{\fss{Smway}}(g,h)}$.

Then for $i=1,2$:
\begin{eqnarray*}
& \mbox{\large(}M_i^{\fss{Lft}} \models \varphi \mbox{ \,or }
M_i^{\fss{Rt}} \models \varphi \mbox{\large)} \mbox{ \ iff \ }
\\*
& \mbox{\Large(}
\exists g\mbox{\large(} \psi_{\fss{1way}}(g) \wedge
\varphi^{\fss{1way}}(g)\mbox{\large)}
\ \vee \  \mbox{\large(} \forall g(\psi_{\fss{1way}}(g) \raro
\psi_{\fss{Bdd}}(g))
\, \wedge \,\varphi^{\fss{Bdd}} \mbox{\large)}\mbox{\Large)}.
\end{eqnarray*}
Moreover,
\begin{eqnarray*}
&  \mbox{\large(}M_i^{\fss{Lft}} \models \varphi \mbox{ \,and \,}
M_i^{\fss{Rt}} \models \varphi \mbox{\large)} \mbox{ \ iff \ }
\\*
& \mbox{\LARGE(} \mbox{\Large(}
\exists g \exists h \mbox{\large(} \psi_{\fss{1way}}(g) \wedge
\psi_{\fss{1way}}(h) \wedge \neg \psi_{\fss{Smway}}(g,h) \wedge
\varphi^{\fss{1way}}(g) \wedge \varphi^{\fss{1way}}(h)
\mbox{\large)} \mbox{\Large)}
\ \vee \\*
& \mbox{\Large(}
\exists g\mbox{\large(} \psi_{\fss{1way}}(g) \wedge
\forall h
( \psi_{\fss{1way}}(h) \raro \psi_{\fss{Smway}}(g,h))
\,\wedge\, \varphi^{\fss{1way}}(g) \mbox{\large)}  \wedge
\varphi^{\fss{Bdd}} \mbox{\Large)}
\ \vee  \\*
& \mbox{\large(} \forall g(\psi_{\fss{1way}}(g) \raro
\psi_{\fss{Bdd}}(g))
\, \wedge \, \varphi^{\fss{Bdd}} \mbox{\large)}\mbox{\LARGE)}.
\end{eqnarray*}

{\bf Case 1:} $M_1^{\fss{Lft}} \not\equiv M_2^{\fss{Lft}}$.  We claim
that $M_1^{\fss{Lft}} \equiv M_2^{\fss{Rt}}$.  First, pick a sentence
$\varphi'$ holding in $M_1^{\fss{Lft}}$ but not in $M_2^{\fss{Lft}}$.
Now let $\varphi$ be any sentence holding in $M_1^{\fss{Lft}}$.
Then $\varphi \wedge \varphi'$ holds in $M_2^{\fss{Lft}}$ or in
$M_2^{\fss{Rt}}$
(by Theorem~\ref{LRT}(1) and what was shown above),
whence $\varphi \wedge \varphi'$ holds in $M_2^{\fss{Rt}}$,
and thus $\varphi$ holds in $M_2^{\fss{Rt}}$.
The claim follows.  A parallel argument then shows that
$M_2^{\fss{Lft}} \equiv M_1^{\fss{Rt}}$.

{\bf Case 2:} $M_1^{\fss{Lft}} \equiv M_2^{\fss{Lft}}$.
This case is left to the reader.
\hfill\qed
\medskip

\medskip

\noindent
{\bf The monotonic case}

\smallskip

Now we treat the monotonic case.  The proof of Theorem~\ref{MNRT}
is almost identical to the proof in the linear case.  There are only two
real differences.
The first is that we need Proposition~\ref{USED} from
Section~\ref{TRANSPROPS}, which says that if a locally moving
monotonic permutation group $\pair{\fs{ED}(L)}{G}$ is
3-interval-transitive then $\pair{\fs{ED}(L)}{\fs{Bdd}(G)}$ is
2-interval-transitive.  The other difference
occurs in capturing the segregation relation in
Proposition~\ref{LSEG}.

\begin{defn}
\label{}
\begin{rm}
Let $\pair{N}{G} \in K^{\fss{NO}}$, and let $g \in G$.  We call $g$ a
{\em non-unit} if
$\fs{var}(g) \neq \barN$, that is, if there exists $f \in G -
\sngltn{\fs{Id}}$ such that
$\fs{var}(f) \cdot \fs{var}(g) = \emptyset$.
In view of Proposition~\ref{DSJNT}, this can be expressed in a first
order formula $\varphi_{\fss{NU}}(u)$ in $\CL^{\fss{GR}}$.  Clearly
\[
\fs{Supp}(g) \mbox{ is bounded } \Rightarrow \varphi_{\fss{NU}}(g)
\Rightarrow g \in \fs{Opp}(G).
\]
\end{rm}
\end{defn}

The monotonic version of Proposition~\ref{LSEG}(1) is:

\begin{prop}
\label{MSEG}
Let
\[
\varphi_{\fss{Seg}}^{\fss{MN3}}(u_1,u_2) \eqqdf
\neg(\exists g)\mbox{\large(}(g(u_1) \cdot u_2 \neq 0) \wedge
(g(u_2) \cdot u_1 \neq 0) \wedge \varphi_{\fss{NU}}(g)
\mbox{\large)},
\]
re-expressed as a first order formula in $\CL^{\fss{VACT}}$.  Let
$\pair{\fs{ED}(L)}{G} \in K^{\fss{MN3}}$.  Then for any $U_1,U_2 \in
\fs{Var}(L,G)$:
\[
\fs{VACT}(\fs{ED}(L),G) \models
\varphi_{\fss{Seg}}^{\fss{MN3}}[U_1,U_2] \: \mbox{ iff } \:
\fs{Seg}^{L}(U_1,U_2).
\]
\end{prop}

\noindent
{\bf Proof }
Let $\pair{\fs{ED}(L)}{G} \in K^{\fss{MN3}}$.
Then $\pair{L}{\fs{Bdd}(G)} \in K^{\fss{LN2}}$  by
Proposition~\ref{USED}.
Now let $U_1,U_2 \in \fs{Var}(\fs{ED}(L),G)$.

If $\fs{Seg}^{L}(U_1,U_2)$
then
$\fs{VACT}(\fs{ED}(L),G) \models
\varphi_{\fss{Seg}}^{\fss{MN3}}[U_1,U_2]$
since every non-unit $g$ lies in $\fs{Opp}(G)$.

Conversely, if $\fs{Seg}^{L}(U_1,U_2)$ fails then applying
Proposition~\ref{LSEG} to
\newline
$\pair{L}{\fs{Bdd}(G)} \in K^{\fss{LN2}}$
guarantees that
\[
g(U_1) \cap U_2 \neq \emptyset \mbox{ and } g(U_2) \cap U_1 \neq
\emptyset
\]
for some $g \in \fs{Bdd}(G)$, so that
$\varphi_{\fss{Seg}}^{\fss{MN3}}[U_1,U_2]$ fails.
\hfill\qed
\medskip

The rest of the proof parallels the linear proof; see
Definition~\ref{DEFLNNOGR}ff, which justifies the use of
Lemma~\ref{DDLN}.

\begin{theorem}
\label{BIGMN}
The analogue of Theorem~\ref{BIG}
holds for monotonic permutation groups, where now
``Reconstruction'' refers to Definition~\ref{MFOR}.
\end{theorem}

Theorem~\ref{MNRT} is a special case.

\medskip
\medskip

\noindent
{\bf The circular case}

\smallskip

Next we treat the circular case, which requires more extensive
modification of the linear proof.

\begin{defn}\label{}
\begin{rm}
Let $C = \fs{CR}(L)$ be a circularly ordered set, where $L$ is a dense
chain having at most one endpoint.  Let $U_1,U_2,U_3,U_4 \in
\barB(C)$.

(1) We say that $U_1$ and $U_2$ are {\em segregated} if
there exist disjoint $\barC$-intervals $I_1$ and $I_2$ such that $U_j
\subseteq I_j$ for $j=1,2$ (or if some $I_j = \emptyset$).  Let
\[
\fs{Seg}^{C}  \eqdf \setm{\pair{U_1}{U_2} \in \barB(C)^2}{U_1
\mbox{ and } U_2 \mbox{ are segregated}}.
\]

(2) Recall that $Crs(U_1,U_2,U_3)$ means that
$Cr(x_1,x_2,x_3)$ whenever each $x_j \in U_j$, and each $U_j \neq
\emptyset$.
We write
$Crs^{-}(U_1,U_2,U_3)$ to mean that
\newline
$Crs(U_3,U_2,U_1)$.

(3) We write $Crs^{\pm}(U_1,U_2,U_3)$ to mean that either
$Crs(U_1,U_2,U_3)$ or $Crs^{-}(U_1,U_2,U_3)$, or equivalently, that
each $U_j \neq \emptyset$ and there exist pairwise disjoint
$\barC$-intervals $I_1,I_2,I_3$ such that $U_j \subseteq I_j$ for
$j=1,2,3$.

(4) We write
$\fs{Seps}(U_1,U_2,U_3,U_4)$ if $\fs{Sep}(x_1,x_2,x_3,x_4)$
whenever $x_j \in U_j$ for all $j$, and if
each $U_j \neq \emptyset$.  (The relation \fs{Sep} was defined in
Section~\ref{OBJECTS}.)  This implies that there exist pairwise
disjoint $\barC$-intervals $I_1,I_2,I_3,I_4$ such that $U_j \subseteq
I_j$ for $j=1,2,3,4$.

(5) These notions are unaffected by reversal of the orientation on $C$
(except that $Crs(U_1,U_2,U_3)$ and $Crs^{-}(U_3,U_2,U_1)$
 are interchanged), and thus are also applicable to $\fs{EO}(L)$.
\end{rm}
\end{defn}

We now show that for $\pair{C}{G} \in K^{\fss{CR3}}$ the relation
$\fs{Seg}^{C}$ on $\fs{Var}(C,G)$ is definable by a first order formula
in $\CL^{\fss{VACT}}$.  This will be our sole use of
3-interval-transitivity.  Here we omit the statement of the obvious
analogue
involving $\barB(C)$.

\begin{prop}
\label{CSEG}
There is a first order formula
$\varphi_{\fss{Seg}}^{\fss{CR3}}(u_1,u_2)$
in $\CL^{\fss{VACT}}$ having the following property.
For any $\pair{C}{G} \in K^{\fss{CR3}}$ and any $U_1,U_2 \in
\fs{Var}(C,G)$:
\[
\fs{VACT}(C,G) \models \varphi_{\fss{Seg}}^{\fss{CR3}}[U_1,U_2] \:
\mbox{ iff } \: \fs{Seg}^{C}(U_1,U_2).
\]
\end{prop}

\noindent
{\bf Proof }
Let $\pair{C}{G} \in K^{\fss{CR3}}$ and let $U_1,\dots,U_4 \in
\fs{Var}(C,G)$.

Now let
\begin{eqnarray*}
& \varphi_{\fss{Crs}^{\pm}}^{\fss{VA}}(u_1,u_2,u_3) \eqqdf
   \varphi(u_1,u_2,u_3) \wedge \varphi(u_2,u_3,u_1) \wedge
   \varphi(u_3,u_1,u_2),  \\
& \mbox{where } \varphi(u_1,u_2,u_3)  \eqqdf (\bigwedge_{i=1}^3
u_i\neq 0)
   \wedge \\*
& \neg (\exists g)
   \mbox{\large(} g(u_1) \cdot u_1 \neq 0 \, \wedge \,
   g(u_2) \cdot u_3 \neq 0 \, \wedge \,
   g(u_3) \cdot u_2 \neq 0 \mbox{\large)}.
\end{eqnarray*}
We claim that
\begin{eqnarray*}
& \fs{VACT}(C,G) \models
\varphi_{\fss{Crs}^{\pm}}^{\fss{VA}}[U_1,U_2,U_3] \: \mbox{ iff } \:
\fs{Crs}^{\pm}(U_1,U_2,U_3).
\end{eqnarray*}

Obviously the implication from right to left holds.

Now suppose that $\fs{VACT}(C,G) \models
\varphi_{\fss{Crs}^{\pm}}^{\fss{VA}}[U_1,U_2,U_3]$, but that
\newline
$\fs{Crs}^{\pm}(U_1,U_2,U_3)$ does not hold.
Then $U_1,U_2,U_3$ are nonempty; also they are pairwise disjoint
(else $g=\fs{Id}$ would give a contradiction).
Since
\newline
$\fs{Crs}^{\pm}(U_1,U_2,U_3)$ does not hold, there cannot exist
pairwise disjoint $\barC$-intervals containing $U_1$, $U_2$,
$U_3$.  Hence without loss of generality
there exists a $\barC$-interval $K_1 \subseteq U_1$ such that,
letting $I_1$ be the largest $\barC$-interval which contains $K_1$
and is disjoint from $U_2 \cup U_3$, we have $I_1 \not\supseteq
U_1$.

Pick a $\barC$-interval $K'_1 \subseteq U_1 - I_1$.
There must exist $\barC$-intervals $K_1 \subseteq U_2$ and
$K_3 \subseteq U_3$ such that $\fs{Seps}(K_1,K_2,K'_1,K_3)$.  We
may assume that $\fs{Crs}(K_1,K_2,K'_1,K_3)$, the other case being
similar.  Then
$\fs{Crs}(K_1,K_2,K_3)$ and $\fs{Crs}(K'_1,K_3,K_2)$.
3-interval-transitivity provides $g\in G$ such that
$g(K_1) \cap K'_1 \neq \emptyset$,
$g(K_2) \cap K_3 \neq \emptyset$,
and $g(K_3) \cap K_2 \neq \emptyset$, which is a contradiction.

Now $\fs{Seps}(U_1,U_2,U_3,U_4)$ is captured by
\begin{eqnarray*}
& \varphi_{\fss{Seps}}^{\fss{VA}}(u_1,u_2,u_3,u_4) \eqqdf
   \varphi_{\fss{Crs}^{\pm}}^{\fss{VA}}(u_1,u_2,u_3) \wedge
   \varphi_{\fss{Crs}^{\pm}}^{\fss{VA}}(u_1,u_4,u_3)  \wedge \\*
& (\exists v)\mbox{\large(}(0 \neq v \leq  u_2+u_4) \wedge
   \neg
\varphi_{\fss{Crs}^{\pm}}^{\fss{VA}}(u_1,v,u_3)\mbox{\large)}.
\end{eqnarray*}

Then $\fs{Seg}^C$ is captured by
\begin{eqnarray*}
& \varphi_{\fss{Seg}}^{\fss{Cr}}(u_1,u_2) \eqqdf
   \neg  (\exists v_1,v_3 \leq  u_1)
            (\exists v_2,v_4 \leq  u_2)
\varphi_{\fss{Seps}}^{\fss{VA}}(v_1,v_2,v_3,v_4).
\end{eqnarray*}

Finally, all these formulas can be re-expressed in
$\CL^{\fss{VACT}}$ because of Proposition~\ref{PST}.
\hfill\qed
\medskip

Again there will be an analogue of
Theorem~\ref{BIG} saying that in order to capture
$\fs{EPACT}(C,G)$ in the group $G$ it is enough to capture the
segregation relation in $\fs{VACT}(C,G)$.  The language
$\CL^{\fss{SEG}} =  \{ \leq,\fs{Seg} \}$ will be the same as before.
Once again we deal with doubly dense subsets $T$ of $\barB(L)$, and
capture a
dense set $\fs{Ep}(C,T)$ of ``endpoints'', representing these points by
triples of members of $T$.

\begin{defn}
\label{}
\begin{rm}
Let $U_1,U_2 \in \barB(C)$.

(1) Suppose $U_1$ and $U_2$ are nonempty and segregated.
Let $I$ be the intersection of all $\barC$-intervals which contain
$U_1$ and are disjoint from $U_2$.  Then $I$ is a $\barC$-interval.
Let $I = (x,y)$.  We denote $x$ by
$\fs{cw}(U_1;U_2)$, the clockwise end of $U_1$
relative to $U_2$, and $y$ by $\fs{ccw}(U_1;U_2)$, the
counterclockwise end of $U_1$ relative to $U_2$. When $U$ is a
bounded $\barC$-interval, we abbreviate $\fs{cw}(U;-U)$ by
$\fs{cw}(U)$, and dually for  $\fs{ccw}(U;-U)$.

(2) For $U_1 \neq \emptyset$ let
\begin{eqnarray*}
\fs{Ep}(U_1) \eqdf
& \setm{\fs{cw}(U_1;U_2)}{U_2 \neq \emptyset
\mbox{ and } \fs{Seg}^C(U_1,U_2)}
\,\, \cup
\\*
& \setm{\fs{ccw}(U_1;U_2)}{U_2 \neq \emptyset
\mbox{ and } \fs{Seg}^C(U_1,U_2)};
\end{eqnarray*}
and let $\fs{Ep}(\emptyset) = \emptyset$.
Note that $\fs{Ep}(U_1)$ consists of the endpoints
in $\barC$ of the connected components of $-U_1$.
\end{rm}
\end{defn}

\begin{defn}\label{DEFCRNOGR}
\begin{rm} Let $T$ be a dense subset of $\barB(C)$. We
define:
\newline
\smallskip
$\fs{SEG}(C,T) \eqdf \trpl{T}{\leq^{\barB(C)}}{\fs{Seg}^C}.
\newline
\smallskip
\fs{Ep}(C,T) \eqdf \bigcup_{U\in T}\fs{Ep}(U).
\newline
\smallskip
\fs{EP}(C,T) \eqdf \pair{\fs{Ep}(C,T)}{\fs{Eo}^{\barC}},
   \mbox{ with } \CL(\fs{EP}(C,T)) = \{\fs{Eo} \}.
\newline
\smallskip
\fs{Rep}(C,T) \eqdf
   \setm{\trpl{U_1}{U_2}{U_3} \in (T-\sngltn{\emptyset})^3}
   {\fs{Crs}^{\pm}(U_1,U_2,U_3)}.
\newline
\smallskip
\fnn{\fs{Pt}^{C,T}}{\fs{Rep}(C,T)}{\fs{Ep}(C,T)} \mbox{ is the
   surjection defined by }
\newline
\smallskip
\ \ \ \ \ \ \ \ \ \ \ \ \ \  \fs{Pt}^{C,T}(U_1,U_2,U_3) = \left\{
\begin{array}{ll}
\fs{ccw}(U_1;U_2) & \mbox{if $\fs{Crs}(U_1,U_2,U_3)$,} \\*
\fs{cw}(U_1;U_2)  & \mbox{if $\fs{Crs}^{-}(U_1,U_2,U_3)$.}
\end{array}
\right. \\
\fs{SEGEP}(C,T) \eqdf (\fs{SEG}(C,T),\fs{EP}(C,T);\fs{Pt}^{C,T}).$
\end{rm}
\end{defn}

Because $\fs{EP}(C,T)$ is an equal orientation structure, nothing
in Definition~\ref{DEFCRNOGR} is affected by reversal of the
orientation on $C$.

\begin{defn}\label{DEFIU}
\begin{rm}
For a {\em representative}  $\vecU =  \trpl{U_1}{U_2}{U_3} \in
\fs{Rep}(C,T)$:

(1) Let $I_{\vecU}$ be the largest
$\barC$-interval which contains $U_2$ and is disjoint from
$U_1 \cup U_3$.

(2) Let $x_{\vecU} \eqdf  \fs{Pt}^{C,T}(\vecU) \in \barC$.
\end{rm}
\end{defn}

We clarify the landscape with an easy lemma, whose verification is
left to the reader.

\begin{lemma}
\label{LEMIU}
Let $\vecU =  \trpl{U_1}{U_2}{U_3} \in
\fs{Rep}(C,T)$.  Then there is a connected component $J$ of $-U_1$
such that:

{\rm (1)} $U_2,U_3 \subseteq J$, $I_{\vecU} \subseteq J$,
$x_{\vecU} \in \fs{Ep}(J)$, and $x_{\vecU} \in \fs{Ep}(I_{\vecU})$.

{\rm (2) (a)} If $\fs{Crs}(\vecU)$ then $I_{\vecU} = (x_{\vecU},y)
\subsetneqq J$, where $y = \fs{cw}(U_3;U_2)$.

$\,\,\,\,\,\,\,\,\,${\rm (b)} If $\fs{Crs}^-(\vecU)$ then $I_{\vecU} =
(y,x_{\vecU})
\subsetneqq J$, where $y = \fs{ccw}(U_3;U_2)$.
\end{lemma}

The definitions of $K^{\fss{DDCR}}$, $K^{\fss{DDCR}}_{\fss{SEG}}$, and
$K^{\fss{DDCR}}_{\fss{SEGEP}}$ parallel the linear definitions in
Definition~\ref{DEFDDLN}.

\begin{lemma}
\label{DDCR}
$K^{\fss{DDCR}}_{\fss{SEGEP}}$
is FO-STR-interpretable in
$K^{\fss{DDCR}}_{\fss{SEG}}$, and
\newline
$\pair{K^{\fss{DDCR}}_{\fss{SEG}}}{K^{\fss{DDCR}}_{\fss{SEGEP}}}$
is a spanning system.
\end{lemma}

\noindent
{\bf Proof }
We capture several notions in $\CL^{\fss{SEG}}$,
some of which were captured above in $\CL^{\fss{VACT}}$.
Some of the verifications will be left to the reader.
We begin with the notion of ``representative''.  Let
\begin{eqnarray*}
& \varphi_{\fss{Rep}}(u_1,u_2,u_3) \eqqdf
    \varphi(u_1,u_2,u_3) \wedge
    \varphi(u_2,u_3,u_1) \wedge
    \varphi(u_3,u_1,u_2),  \\*
& \mbox{where } \varphi(u_1,u_2,u_3) \eqqdf
\mbox{\large(} \! \bigwedge_{i=1}^{3} \!u_i \neq \emptyset
\mbox{\large)}  \wedge
\forall v \mbox{\large (}(v \leq u_1 + u_3)
\raro \fs{Seg}(u_2,v)\mbox{\large )}.
\end{eqnarray*}
We show that for any $\pair{C}{T} \in K^{\fss{DDCR}}$
and $U_1,U_2,U_3 \in T$:
\[
\fs{SEG}(C,T) \models \varphi_{\fss{Rep}}[U_1,U_2,U_3]
\mbox{ \ iff \ } \trpl{U_1}{U_2}{U_3} \in \fs{Rep}(C,T).
\]
Obviously the implication from right to left holds.
Suppose that
$\fs{SEG}(C,T) \models \varphi_{\fss{Rep}}[U_1,U_2,U_3]$ holds, but
that $\trpl{U_1}{U_2}{U_3} \not\in \fs{Rep}(C,T)$.
Then each $U_i \neq \emptyset$, and we may assume that
there exist
$x_1,x_2 \in U_2$, $y \in U_1$, and $z \in U_3$ such that
$\fs{Sep}(x_1,y,x_2,z)$.
Let $I,J$ be $C$-intervals such that
$y \in I \subseteq U_1$ and $z \in J \subseteq U_3$.
The double density of $T$ guarantees the
existence of $V \in T$ such that $V \leq I + J$ and
$V \cap I, V \cap J \neq \emptyset$.
Hence $\fs{Seg}^{C}(U_2,V)$ does not hold, and thus
$\fs{SEG}(C,T) \not\models \varphi_{\fss{Rep}}[U_1,U_2,U_3]$, a
contradiction.

The formula $\varphi_{\fss{Seps}}(u_1,u_2,u_3,u_4)$,
built from  $\varphi_{\fss{Rep}}$ just as
$\varphi_{\fss{Seps}}^{\fss{VA}}$
was built from
$\varphi_{\fss{Crs}^{\pm}}^{\fss{VA}}$ in the proof of
Lemma~\ref{CSEG},
captures $\fs{Seps}(U_1,U_2,U_3,U_4)$.

Our next goal is to find a formula
$\varphi_{\fss{Eqp}}(\vecu,\vecv)$
which says that $x_{\vecU} = x_{\vecV}$.

For $U \in \barB(C)$, let  $T\rest U \eqdf \setm{V \in T}{V \leq U}$.
For $\vecU \in \fs{Rep}(C,T)$,
we verify that the 4-place relation
``$V \in T \rest I_{\vecU} - \sngltn{\emptyset}$'' can
be expressed by a first order formula
$\varphi_{\fss{Mid}}(v,u_1,u_2,u_3)$ in $\CL^{\fss{SEG}}$.  Let
$$
\fs{Ext}(U_2;U_1,U_3) \eqdf
\setm{V \in T}
{\trpl{U_1}{V}{U_3} \in \fs{Rep}(C,T) \,\wedge
\neg \fs{Seps}(U_1,V,U_3,U_2)}.
$$
Then $\fs{Ext}(U_2;U_1,U_3) =
T \rest I_{\vecU} - \sngltn{\emptyset}$,
and obviously there is a first order formula
in $\CL^{\fss{SEG}}$
expressing the relation ``$V \in \fs{Ext}(U_2;U_1,U_3)$''.

For $\vecU = \trpl{U_1}{U_2}{U_3} \in \fs{Rep}(C,T)$,
$$
\fs{Sm1pt}(U_1,U_2,U_3) \eqdf
\setm{\trpl{U_1}{W_2}{W_3} \in \fs{Rep}(C,T)}
{x_{\trpl{U_1}{U_2}{U_3}} = x_{\trpl{U_1}{W_2}{W_3}}},
$$
the set of representatives having the same $U_1$ and the same $x$
as $\vecU$.  Now for
$\vecU' = \trpl{U_1}{U'_2}{U'_3} \in \fs{Rep}(C,T)$,
it is easily seen with the aid of Lemma~\ref{LEMIU} that the
following are equivalent:

(a) $x_{\vecU} = x_{\vecU'}$.

(b) $I_{\vecU} \subseteq I_{\vecU'}$ or $I_{\vecU} \supseteq
I_{\vecU'}$.

(c) $T \rest I_{\vecU} - \sngltn{\emptyset} \ \subseteq\
T \rest I_{\vecU'} - \sngltn{\emptyset}$
\ \ or\ \
$T \rest I_{\vecU} - \sngltn{\emptyset}\ \supseteq\
T \rest I_{\vecU'} - \sngltn{\emptyset}$.

Since ``$V \in T \rest I_{\vecU} - \sngltn{\emptyset}$''
is expressible by the formula $\varphi_{\fss{Mid}}(v,u_1,u_2,u_3)$,
we conclude that
``$\vecU' \in \fs{Sm1pt}(\vecU)$''
is expressible by a first order formula in $\CL^{\fss{SEG}}$.

We claim that for $\vecU, \vecV \in \fs{Rep}(C,T)$:
\begin{eqnarray*}
&  x_{\vecU} = x_{\vecV} \,\,\, \mbox{ iff} \\*
\mbox{(*) \,\,\,\,\,\,\,\,\,\,\,\,}
& \left(\rule{0pt}{4mm}
   \forall \vecU' \in \fs{Sm1pt}(\vecU)\right)
   \left(\rule{0pt}{4mm}
   \forall \vecV' \in \fs{Sm1pt}(\vecV)\right)
   \left(\rule{0pt}{4mm}
   \forall W_1,W_2 \in T - \sngltn{\emptyset}\right) \\*
&  \left(\rule{0pt}{4mm}
\exists X \in T \rest I_{\vecU'} - \sngltn{\emptyset}\right)
   \left(\rule{0pt}{4mm}
   \exists Y \in T \rest I_{\vecV'} - \sngltn{\emptyset}\right)
   \neg\fs{Seps}(X,W_1,Y,W_2).
\end{eqnarray*}

Suppose first that
$x_{\vecU} = x_{\vecV}$.
Let $\vecU' \in \fs{Sm1pt}(\vecU)$,
$\vecV' \in \fs{Sm1pt}(\vecV)$,
and $W_1,W_2 \in T - \sngltn{\emptyset}$.
Then $x \eqdf x_{\vecU} = x_{\vecV}$ is a common endpoint of
$I_{\vecU'}$ and $I_{\vecV'}$.
We assume that $I_{\vecU'} = (x,y)$ for some $y \in \barC$, the
other case being the dual.  If also $I_{\vecV'} = (x,y')$ for some $y'
\in \barC$, then there exists $Z \in
(T \rest I_{\vecU'} - \sngltn{\emptyset)} \cap
(T \rest I_{\vecV'} - \sngltn{\emptyset)}$,
and obviously $\fs{Seps}(Z,W_1,Z,W_2)$ does not hold.
If instead $I_{\vecV'} = (y',x)$ for some $y' \in \barC$,
then the condition
$$
 \left(\rule{0pt}{4mm}
\forall X \in T \rest I_{\vecU'} - \sngltn{\emptyset}\right)
   \left(\rule{0pt}{4mm}
   \forall Y \in T \rest I_{\vecV'} - \sngltn{\emptyset}\right)
   \fs{Seps}(X,W_1,Y,W_2)
$$
would require $\fs{Crs}(I_{\vecV'},W_i,I_{\vecU'})$ either for $i=1$
or for $i=2$, which is impossible.

Conversely, suppose that
$x_{\vecU} \neq x_{\vecV}$.
Choose $\barC$-intervals $I'_1,I'_2,J_1,J_2$ such that
$x_{\vecU} \in I'_1$,
$x_{\vecV} \in I'_2$, and
$\fs{Seps}(I'_1,J_1,I'_2,J_2)$.
Let $I_1 = I'_1 \cap I_{\vecU}$ and
$I_2 = I'_2 \cap I_{\vecV}$.
Use the density of $T$ to choose
$W_i \in T \rest J_i - \sngltn{\emptyset}$ for $i = 1,2$.
Now choose $U'_2$ and $U'_3$ in $T \rest I_1 - \sngltn{\emptyset}$
such that
\[
\left\{
\begin{array}{ll}
\fs{Crs}(U_1,U'_2,U'_3) &
\mbox{if $\fs{Crs}(U_1,U_2,U_3)$,} \\*
\fs{Crs}^{-}(U_1,U'_2,U'_3) & \mbox{if
$\fs{Crs}^{-}(U_1,U_2,U_3)$.} \end{array}
\right.
\]
Then $\vecU' \eqdf \trpl{U_1}{U'_2}{U'_3} \in \fs{Sm1pt}(\vecU)$.
Since $I_{\vecU'} \subseteq I_1$ we have
$T \rest I_{\vecU'} \subseteq T \rest I_1$.
Choose $V'_2, V'_3 \in T \rest I_2  - \sngltn{\emptyset}$ similarly,
so that $T \rest I_{\vecV'} \subseteq T \rest I_2$.
Let $X \in T \rest I_{\vecU'} - \sngltn{\emptyset}$ and
$Y \in T \rest I_{\vecV'} - \sngltn{\emptyset}$.
Then $X \subseteq I'_1$, $W_1 \subseteq J_1$, $Y \subseteq I'_2$,
and
$W_2 \subseteq J_2$,
so $\fs{Seps}(X,W_1,Y,W_2)$ holds.

Since (*) can be expressed first order in $\CL^{\fss{SEG}}$, there is
indeed a first order formula
$\varphi_{\fss{Eqp}}(\vecu,\vecv)$ in $\CL^{\fss{SEG}}$
such that for every
$\pair{C}{T} \in K_{\fss{DDCR}}$ and
$\vecU, \vecV \in \fs{Rep}(C,T)$:
\[
\fs{SEG}(C,T) \models
\varphi_{\fss{Eqp}}[\vecU,\vecV]
\mbox{ \ iff \ }
x_{\vecU} = x_{\vecV}.
\]

Now we want to capture by a first order formula in $\CL^{\fss{SEG}}$
the restriction
$\fs{Eo}^{\barC} \rest \fs{Ep}(C,T)$ of the equal orientation relation.
By Proposition~\ref{SEP} it suffices to capture
$\fs{Sep}^{\barC} \rest \fs{Ep}(C,T)$ instead.
To do this, let
$\varphi_{\fss{Sep}}(\vecu^1,\vecu^2,\vecu^3,\vecu^4)$
be the formula which says the following:
There exist $\vecv^1,\ldots,\vecv^4$ such that

(1) $\bigwedge_{i = 1}^4(x_{\vecv^i} = x_{\vecu^i})$, and

(2) $\mbox{\large(} \forall w_1 \in T \rest I_{\vecv^1} -
\sngltn{\emptyset}\mbox{\large)},\ldots,
\mbox{\large(}\forall w_4 \in T \rest I_{\vecv^4} -
\sngltn{\emptyset}\mbox{\large)}
\fs{Seps}(w_1,\ldots,w_4)$.
\newline
Then for every $\pair{C}{T} \in K_{\fss{DDCR}}$ and
$\vecU = \trpl{\vecU^1}{\ldots}{\vecU^4} \in (\fs{Rep}(C,T))^4$:
\[
\fs{SEG}(C,T) \models \varphi_{\fss{Sep}}[\vecU]
\mbox{ \ iff \ }
\fs{Sep}(x_{\vecU^1},\ldots,x_{\vecU^4}) \mbox{ holds}.
\]

The rest of the proof parallels that of Lemma~\ref{DDLN}.
\hfill\qed
\medskip

As in the linear case:

\begin{theorem}
\label{BIGCR}
The analogue of Theorem~\ref{BIG} holds for circular permutation
groups, where now ``Reconstruction'' refers to Definition~\ref{CFOR}.
\end{theorem}

Theorem~\ref{CRT} is a special case.

\medskip
\medskip

\noindent
{\bf The monocircular case}

\smallskip

The monocircular case is to the circular case as the monotonic case is
to the linear case.
We use Proposition~\ref{USED}, which says that if a locally moving
monocircular permutation group $\pair{\fs{EO}(C)}{G}$ is
4-interval-transitive, then $\pair{\fs{EO}(C)}{\fs{Bdd}(G)}$ is
3-interval-transitive, and even better that
$\pair{\fs{EO}(C)}{\fs{3Bdd}(G)}$
is 3-interval-transitive, where
$$\fs{3Bdd}(G)  \eqdf  \setm{g \in G}{g \mbox{ is the product of
3 elements of bounded support}}.$$

In the monocircular version of Proposition~\ref{CSEG},
$\varphi_{\fss{Seg}}^{\fss{MC4}}(u_1,u_2)$
is constructed exactly as was $\varphi_{\fss{Seg}}^{\fss{CR3}}$,
except that in the definition of
$\varphi_{\fss{Crs}^{\pm}}^{\fss{VA}}$
the requirement is added that $g$ be the product of three non-units
(see Proposition~\ref{MSEG}), which guarantees that $g \in
\fs{Opp}(G)$.

\begin{theorem}
\label{BIGMC}
The analogue of Theorem~\ref{BIG}
holds for monocircular permutation groups, where now
``Reconstruction'' refers to Definition~\ref{MCFOR}.
\end{theorem}

Theorem~\ref{MCRT} is a special case.

\medskip
\medskip

\noindent
{\bf Distinguishing the four types}

\smallskip

Finally, we show that the four types of nearly ordered permutation
groups that we have been dealing with can be distinguished first
order in the language of groups (even when the degrees of multiple
transitivity are relaxed a bit), of course with the exception that when
a monotonic permutation group $\pair{\fs{ED}(L)}{G}$ has
$\fs{Opp}(G) = G$ then $\pair{L}{G}$ is a linear permutation group,
and analogously for monocircular permutation groups.

\begin{defn}
\label{}
\begin{rm}\

(1) Let $K^{\fss{LN}}$ denote the class of locally moving linear
permutation groups, and analogously for $K^{\fss{MN}}$,
$K^{\fss{CR}}$, and $K^{\fss{MC}}$.

(2) For $K \subseteq K^{\fss{MN}}$, let
$K_{\fss{Orp}}$ denote the class of $\pair{N}{G} \in K$ such that
$G$ contains an order-reversing permutation,
and analogously for $K \subseteq K^{\fss{MC}}$.

(3) Let $K^{\fss{FOUR}} \eqdf K^{\fss{LN}} \cup
K^{\fss{MN2}}_{\fss{Orp}} \cup K^{\fss{CR2}} \cup
K^{\fss{MC3}}_{\fss{Orp}}$.
\end{rm}
\end{defn}

\begin{prop}
\label{FOUR}
There are first order sentences $\psi^{\fss{LN}}$,
$\psi^{\fss{MN2}}_{\fss{Orp}}$, $\psi^{\fss{CR2}}$,
and $\psi^{\fss{MC3}}_{\fss{Orp}}$ in the language $\CL^{\fss{GR}}$
of groups such that for any $\pair{N}{G} \in K^{\fss{FOUR}}$:

{\rm (1)} $G \models \psi^{\fss{LN}}$ iff $\pair{N}{G} \in
K^{\fss{LN}}$.

{\rm (2)} $G \models \psi^{\fss{MN2}}_{\fss{Orp}}$ iff $\pair{N}{G}
\in K^{\fss{MN2}}_{\fss{Orp}}$.

{\rm (3)} $G \models \psi^{\fss{CR2}}$ iff $\pair{N}{G} \in
K^{\fss{CR2}}$.

{\rm (4)} $G \models \psi^{\fss{MC3}}_{\fss{Orp}}$ iff $\pair{N}{G}
\in K^{\fss{MC3}}_{\fss{Orp}}$.

\noindent
Consequently for $\pair{N_1}{G_1}$, $\pair{N_2}{G_2} \in
K^{\fss{FOUR}}$ of different types, $G_1$ and $G_2$ cannot be
isomorphic or even elementarily equivalent.
\end{prop}

\noindent
{\bf Proof }
It suffices to distinguish these four classes of groups by first order
sentences in $\CL^{\fss{VACT}}$.  For then, since
$K^{\fss{NO}}_{\fss{VACT}}$ is FO-STR-interpretable in
$K^{\fss{NO}}_{\fss{GR}}$ by Corollary~\ref{NOIT},
Proposition~\ref{REEXPRESS} can be used to translate each of the four
sentences $\varphi$ in $\CL^{\fss{VACT}}$ to a corresponding
sentence $\psi_{\varphi}$ in $\CL^{\fss{GR}}$
such that for all $\pair{N}{G} \in K^{\fss{NO}}$:
$\fs{VACT}(N,G) \models \varphi \mbox{ \ iff \ } G \models
\psi_{\varphi}$.

The following serve as distinguishing sentences (and can be re-
expressed in $\CL^{\fss{VACT}}$ via Proposition~\ref{PST}):

\smallskip

(4)  $\varphi^{\fss{MC3}}_{\fss{Orp}} \eqqdf \newline
(\forall u_1,u_2,u_3 \neq 0)(\exists g)(g(u_1) \cdot u_1 \neq 0
\, \wedge \, g(u_2) \cdot u_3 \neq 0
\, \wedge \, g(u_3) \cdot u_2 \neq 0)$.

\smallskip

(3)  $\varphi^{\fss{CR2}} \eqqdf \newline
(\forall u \neq 0)(\forall v \neq 0)(\exists g)(\forall w)
(w \cdot u = 0 \raro g(w) \leq v) \, \wedge \, \neg
\varphi^{\fss{MC3}}_{\fss{Orp}}$.

\smallskip

(2)  $\varphi^{\fss{MN2}}_{\fss{Orp}} \eqqdf \newline
(\forall u_1,u_2 \neq 0)(\exists g)(g(u_1) \cdot u_2 \neq 0
\, \wedge \, g(u_2) \cdot u_1 \neq 0) \, \wedge \,
\neg \varphi^{\fss{MC3}}_{\fss{Orp}} \, \wedge \,
\neg \varphi^{\fss{CR2}}$.

\smallskip

(1)  $\varphi^{\fss{LN}} \eqqdf
\neg \varphi^{\fss{MC3}}_{\fss{Orp}} \, \wedge \,
\neg \varphi^{\fss{CR2}}  \, \wedge \,
\neg \varphi^{\fss{MN2}}_{\fss{Orp}}$.

\smallskip

\noindent
The verifications are left to the reader.
\hfill\qed
\newpage

%

\section{Second order reconstruction}
\label{FURTHER}

In Section \ref{RECON} we dealt with situations in which we could
obtain a first order interpretation of a dense substructure of
$\barN^{\smallr}$ in $G \leq \fs{Aut}(N)$.
We now consider different transitivity assumptions that yield only
second order interpretations.   These assumptions will involve
$\barN$-intervals in various ways.   Although the present
results will have first order analogues, those analogues will require
an extra hypothesis, and our primary interest here is in the second
order results.

\begin{defn}
\label{DEFNEST}
\begin{rm}\

(1) Let $N$ be a nearly ordered structure based on $L$, and let $(a,b)
\subseteq (c,d)$ be bounded
$\barN$-intervals. We write $(a,b) \subset^S (c,d)$
if $a \neq c$ and $b \neq d$.

(2) We call a nearly ordered permutation group $\pair{N}{G}$
{\em nest-transitive} if for all bounded $\barN$-intervals
$I,J,K$: If $I \subset^S J$, then there exists $g \in G$
such that
$I \subseteq g(K) \subseteq J$,
that is, $g(K)$ is nested between $I$ and $J$.
(It is not required that $g \in \fs{Opp}(g)$.)
Like interval-transitivity, nest-transitivity is robust in the senses of
Propositions \ref{DENSESUB} and \ref{ENLARGEGROUP}.
Clearly approximate 2-{\it o-}transitivity of $\pair{\barN}{G}$
implies nest-transitivity of $\pair{N}{G}$; for \mbox{$N = \pair{L}{<}$} and
for $N = \fs{CR}(L)$  the two notions are equivalent.
For each of the four types of nearly ordered permutation groups,
nest-transitivity implies inclusion-transitivity.

(3) Let $K^{\fss{NEST}}$ denote the class of nest-transitive nearly
ordered permutation groups (of all four types) which have a
nonidentity bounded element, or
equivalently, which are locally moving.

(4) Let  $K^{\fss{A2$\barN$}}$ denote the class of nearly
ordered permutation groups for which $\pair{\barN}{G}$ is
approximately 2-{\it o-}transitive, and which have a nonidentity
bounded element, or
equivalently, which are locally moving.  Note that
$K^{\fss{A2$\barN$}} \subseteq K^{\fss{NEST}} \subseteq
K^{\fss{2IT}}$.

(5) Let $K^{\fss{NEST-INT}}$ denote the class of nearly ordered
permutation groups $\pair{N}{G} \in K^{\fss{NEST}}$ for which there
exists $g \in G$ such that $\fs{var}(g)$ is a bounded
$\barN$-interval;
and similarly for $K^{\fss{A2$\barN$-INT}}$.
\end{rm}
\end{defn}

Even for lattice-ordered permutation groups $\pair{L}{G}$ which are
highly {\it o-}transitive and locally moving, the existence of $g \in G$
such that $\fs{var}(g)$ is a bounded $\barL$-interval is not
automatic, as
seen via an awesome example due to K.R.P. Pierce \cite{P}, which was
modified by A.M.W. Glass \cite{GLONPIERCE} to make it highly
{\it o-}transitive.

One of the two main theorems of this section says that
nest-transitivity (plus the
existence of a nonidentity bounded element) gives Second Order
Reconstruction (see Definitions \ref{LFOR},  \ref{MFOR}, \ref{CFOR},
and \ref{MCFOR}).  The additional assumption that there exists $g \in
G$ such that $\fs{var}(g)$ is a bounded $\barN$-interval gives First
Order Reconstruction.

In the linear case this gives nothing new
since First Order Reconstruction has already been established under
the weaker hypothesis of 2-interval-transitivity
(Theorem~\ref{LRT}ff).  In the monotonic case
nest-transitivity is weaker than the 3-interval-transitivity
hypothesis of Theorem~\ref{MNRT} (given a nonidentity bounded
element) because of Theorems \ref{3H}
and \ref{HAT}.
In the monocircular case
nest-transitivity is obviously weaker than the 4-interval-transitivity
hypothesis of Theorem~\ref{MCRT}.

\begin{question}
\label{QNEST}
\begin{rm}
In the presence of local movability:

(1) Is nest-transitivity strictly weaker than
3-interval-transitivity in the monotonic case?

(2) In the circular case, does 3-interval-transitivity imply
nest-transitivity, or conversely?

(3) Is nest-transitivity strictly weaker than 4-interval-transitivity in
the monocircular case?
\end{rm}
\end{question}

An assertion that {\em First Order Reconstruction} holds for a class
$K$ of
nearly ordered permutation groups will mean that for each of the
four types, First Order Reconstruction  holds for the class of all
members of $K$ of that type; and that there are first order
sentences in
$\CL^{\fss{GR}}$ distinguishing the types within $K$ in the sense of
Proposition \ref{FOUR}.  Similar remarks apply to {\em Second Order
Reconstruction}, the types now being distinguished by second order
sentences.

\begin{theorem}\
\label{NESTRT}

{\rm (1)} First Order Reconstruction holds for
$K^{\fss{NEST-INT}}$, and thus also for $K^{\fss{A2$\barN$-INT}}$.

{\rm (2)} Second Order Reconstruction holds for
$K^{\fss{NEST}}$, and thus also for $K^{\fss{A2$\barN$}}$.
\end{theorem}

The other main theorem is \ref{SPANRT}, which has even weaker
hypotheses, though it applies only in the linear/monotonic case.
\medskip

The plan of the proof of Theorem \ref{NESTRT} is as follows.
We show that the set of bounded
$\barN$-intervals is definable in $\fs{BACT}(N,G)$.  From
this it will follow that the segregation relation $\fs{Seg}^N$ is
also definable in $\fs{BACT}(N,G)$;
it is this part of the plan that we implement first.
Then we apply the machinery of Section \ref{RECON} to each of the
four types.  Finally, we distinguish the types.

\begin{defn}
\label{}
\begin{rm}
Let $N$ be a nearly ordered structure, and let $E \subseteq
\barB(N)$.
We say that $\pair{N}{E}$ is {\em interval-rich} if for all bounded
$\barN$-intervals $I \subset^S J$
there exists an $\barN$-interval $K \in E$
such that $I \subseteq K \subseteq J$.
Let $\fs{Int}^{\barN}$ be the 1-place relation
of being a bounded $\barN$-interval.
\end{rm}
\end{defn}

Note that $\fs{Var}(N,G)$ is interval-rich when $\pair{N}{G}$ is
nest-transitive and
there exists $g \in G$ such that $\fs{var}(g)$ is a bounded
$\barN$-interval, and that $\barB(N)$ is always interval-rich.

\begin{prop}
\label{BIGINT}
Let $K_1$ be a class of locally moving nearly ordered permutation
groups, all of the same type.

{\rm (1)} Suppose there is a first order formula
$\varphi_{\fss{Int}}(u)$ in
$\CL^{\fss{VACT}}$ such that for any $\pair{N}{G} \in K_1$  and $U
\in
\fs{Var}(N,G)$:
\[
\fs{VACT}(N,G) \models \varphi_{\fss{Int}}[U]
\mbox{ \ iff \ }
\fs{Int}^{\barN}(U).
\]
Suppose also that for any $\pair{N}{G} \in K_1$,
$\pair{N}{\fs{Var}(N,G)}$ is interval-rich.
Then First Order Reconstruction holds for $K_1$.

{\rm (2)} Suppose there is a first order formula
$\varphi_{\fss{Int}}(u)$ in
$\CL^{\fss{BACT}}$ such that for any $\pair{N}{G} \in K_1$
and $U \in \barB(N)$:
\[
\fs{BACT}(N,G) \models \varphi_{\fss{Int}}[U]
\mbox{ \ iff \ }
\fs{Int}^{\barN}(U).
\]
Then Second Order Reconstruction holds for $K_1$.
\end{prop}

\noindent
{\bf Proof } (1) Let
\begin{eqnarray*}
& \!\! \varphi_{\fss{Seg}}^1 (u,v) \eqqdf \\*
& \!\! \begin{displaystyle}
(\forall u_1,u_2 \leq u)
\mbox{\Large(} \bigwedge_{i=1}^2 u_i \neq 0  \rightarrow
\exists w \mbox{\large(}\varphi_{Int}(w) \,\wedge
\bigwedge_{i=1}^2 (w\cdot u_i \neq 0)
\,\wedge \, w \cdot v = 0
\mbox{\large)}\mbox{\Large)},
\end{displaystyle}
\end{eqnarray*}
and let
\[
\varphi_{\fss{Seg}}(u,v) \eqqdf
\varphi_{\fss{Seg}}^1 (u,v)
\, \wedge \, \varphi_{\fss{Seg}}^1 (v,u),
\]
re-expressed in $\CL^{\fss{VACT}}$ via Proposition \ref{PST}.
Let $\pair{N}{G} \in K_1$ and $U,V \in \fs{Var}(N,G)$.
By interval-richness,
$\varphi^1_{\fss{Seg}}[U,V]$ expresses that either $U = \emptyset$
or there exists a (not necessarily bounded) $\barN$-interval which
contains $U$ and is disjoint from $V$.  Therefore
\[
\fs{VACT}(N,G) \models \varphi_{\fss{Seg}}[U,V]
\mbox{ \ iff \ }
\fs{Seg}^{N}(U,V).
\]
Now (1) follows from Theorem \ref{BIG} and its analogues.

(2) The proof is similar, but uses part (2) of Theorem \ref{BIG}.
\hfill\qed
\medskip

Now we deal with the definability of bounded $N$-intervals.  For
$U,V \in \barB(N)$, we write
$U \sim V$ to mean that there exists $g \in G$ such that
$g(U) = V$.  The following lemma is obvious.

\begin{lemma}
\label{PHIBDD}
Let
\[
\varphi_{\fss{Bdd}}(u) \eqqdf
(\forall v \neq 0)(\exists u' \sim u)(u' \leq v).
\]
Suppose that $\pair{N}{G}$ is nest-transitive and locally moving.
Then

{\rm (1)}  For any $U \in \fs{Var}(N,G)$:
\[
\fs{VACT}(N,G) \models \varphi_{\fss{Bdd}}[U]
\mbox{ \ iff \ } U \mbox{ is bounded}.
\]

{\rm (2)} For any $U \in \barB(N)$:
\[
\fs{BACT}(N,G) \models \varphi_{\fss{Bdd}}[U]
\mbox{ \ iff \ } U \mbox{ is bounded}.
\]
\end{lemma}

\noindent
{\bf Proof of Theorem \ref{NESTRT} }
(1) Let
\begin{eqnarray*}
& \!\! \varphi_{\fss{Int}}(u) \eqqdf \varphi_{\fss{Bdd}}(u) \wedge
    u \neq 0  \wedge \\*
& \!\! \begin{displaystyle}
(\forall v_1,v_2,v_3 \sim u)
\mbox{\Large(}\mbox{\large(}
\bigwedge_{i=1}^3 (v_i \cdot u\neq 0 \wedge
v_i \cdot -u \neq 0) \mbox{\large)}
\rightarrow
\bigvee_{1 \leq i < j \leq 3}(v_i \cdot v_j \neq 0) \mbox{\Large)},
\end{displaystyle}
\end{eqnarray*}
re-expressed as a first order formula in $\CL^{\fss{VACT}}$ via
Proposition~\ref{PST}.
We claim that for all $\pair{N}{G} \in K^{\fss{NEST}}$ and $U \in
\fs{Var}(N,G)$:
\[
\fs{VACT}(N,G) \models \varphi_{\fss{Int}}[U]  \mbox{ \ iff \ }
\fs{Int}^{\barN}(U).
\]

The implication from right to left is clear since a bounded
$\barN$-interval has only two boundary points.

Conversely, suppose that $\fs{VACT}(N,G) \models
\varphi_{\fss{Int}}[U]$, so that $U$ is bounded, but that $U$ is not an
$\barN$-interval and thus has more than one connected component.

Assume temporarily that
$N$ has the form $\fs{CR}(L)$ or $\fs{EO}(L)$.
Pick a connected component $A$ of $-U$.
Let $s = \fs{cw}(U;A)$, the clockwise end of $U$ relative to $A$,
and let  $t = \fs{ccw}(U;A)$.

Pick connected components $B$ and $D$ of $U$ and $C$ of $-U$ such
that $\fs{Crs}(A,\sngltn{s},B,C,D,\sngltn{t})$.  Then pick
$a_1,a_2 \in A$,  $b_1,b_2,b_3,b_4 \in B$, $c_1,c_2,c_3,c_4 \in C$,
and $d_1,d_2 \in D$ such that
\[
\fs{Cr}(a_1,a_2,b_1,b_2,b_3,b_4,c_1,c_2,c_3,c_4,d_1,d_2).
\]
Let $I_1 \eqdf (a_2,b_1) \subset^S J_1 \eqdf (a_1,b_2)$,
$I_2 \eqdf (b_4,c_1) \subset^S J_2 \eqdf (b_3,c_2)$, and
$I_3 \eqdf (c_4,d_1) \subset^S J_3 \eqdf (c_3,d_2)$.
For $n=1,2,3$, pick $g_n \in G$ such that $I_n \subseteq g_n((s,t))
\subseteq J_n$ and let $V_n = g_n((s,t))$.  The $V_n$'s are pairwise
disjoint and each $V_n$ meets both $U$ and $-U$, so
$\fs{VACT}(N,G) \not\models \varphi_{\fss{Int}}[U]$,

If $N$ has the form $\pair{L}{<}$ or $\fs{ED}(L)$, the argument is the
same (picking $A$ to be the connected component of $-U$ which is
not bounded below), though the notation is slightly different.

Now apply Proposition~\ref{BIGINT} to complete the
argument for each type.

Finally, we distinguish the four types.  As in Proposition~\ref{FOUR},
it
suffices to distinguish them by first order sentences in
$\CL^{\fss{VACT}}$.  Let
$\pair{N}{G} \in K^{\fss{NEST-INT}}$.

With $\varphi_{\fss{Bdd}}(u)$ as in Lemma~\ref{PHIBDD},
$\fs{VACT}(N,G) \models
(\exists u \neq 0)(\neg\varphi_{\fss{Bdd}}(-U))$ in the linear and
monotonic cases but not in the other two cases.

Let $\varphi_{\fss{Seg}}(u,v)$ be the formula in the proof of
Proposition~\ref{BIGINT}.

For the linear and monotonic cases, let
\[
\varphi_{\fss{Orp}}(g) \eqqdf (\exists u,v)
\mbox{\large(} \varphi_{\fss{Seg}}(u,v) \,\wedge\,
g(u) \cdot v \neq 0  \,\wedge \,
g(v) \cdot u \neq 0 \mbox{\large)}.
\]
We claim that for any $\pair{N}{G} \in
K^{\fss{NEST-INT}}$ which is linear or monotonic,
and for any $g \in G$:
\[
\fs{VACT}(N,G) \models \varphi_{\fss{Orp}}[g]
\mbox{ \ iff \ } g \in G - \fs{Opp}(G).
\]
Clearly if $\fs{VACT}(N,G) \models \varphi_{\fss{Orp}}[g]$
then $g$ cannot preserve order.

Conversely, if $g$  reverses order
and thus fixes exactly one $x \in \barL$, pick $h \in G$ such that
$x < \fs{var}(h)$.  Let $V = g(\fs{var}(h)) < x$.

If $g^2(\fs{var}(h)) \cdot \fs{var}(h) \neq \emptyset$,
let $U = \fs{var}(h)$.
Then $V < x < U$, so $\fs{Seg}^N(U,V)$.
Also, $g(U) \cap V \neq \emptyset$ and
$g(V) \cap U \neq \emptyset$. So
$\fs{VACT}(N,G) \models \varphi_{\fss{Orp}}[g]$.

If $g^2(\fs{var}(h)) \cdot \fs{var}(h) = \emptyset$,
then by Proposition~\ref{VARCONJ},
$U \eqdf\fs{var}(h) + g^2(\fs{var}(h)) = \fs{var}(h \cdot h^{g^2})
\in \fs{Var}(N,G)$.
Again $V < x < U$, and $\fs{VACT}(N,G) \models
\varphi_{\fss{Orp}}[g]$.

For the circular and monocircular cases, let
$\varphi'_{\fss{Rep}}(u_1,u_2,u_3)$ be the first order formula in
$\CL^{\fss{VACT}}$ obtained from the formula
$\varphi_{\fss{Rep}}(u_1,u_2,u_3)$
in the proof of Lemma~\ref{DDCR} by replacing $\fs{Seg}(u_2,v)$ by
$\varphi_{\fss{Seg}}(u_2,v)$.  As in \ref{DDCR}, for any
$U_1,U_2,U_3 \in \fs{Var}(N,G)$:
\[
\pair{N}{\fs{Var}(N,G)} \models \varphi'_{\fss{Rep}}[U_1,U_2,U_3]
\mbox{ \ iff \ } \trpl{U_1}{U_2}{U_3} \in \fs{Rep}(N,\fs{Var}(N,G)).
\]
Let
\[
\varphi'_{\fss{Orp}}(g) \eqqdf (\exists u,v,w)
\mbox{\large(} \varphi'_{\fss{Rep}}(u,v,w) \,\wedge\,
g(u) \cdot v \neq 0  \,\wedge \,
g(v) \cdot u \neq 0  \,\wedge \,
g(w) \cdot w \neq 0 \mbox{\large)}.
\]
Then as in the linear and monotonic cases,
\[
\fs{VACT}(N,G) \models \varphi'_{\fss{Orp}}[g]
\mbox{ \ iff \ } g \in G - \fs{Opp}(G).
\]

Now it is easy to construct the four distinguishing first order
sentences in $\CL^{\fss{VACT}}$.

(2) The proof for any one type is similar to that of (1), but uses part
(2) of Proposition~\ref{BIGINT}.
The argument for distinguishing the types is also similar,
except that because
$K^{\fss{NO}}_{\fss{BACT}}$ is only SO-STR-interpretable
in $K^{\fss{NO}}_{\fss{GR}}$, the distinguishing sentences in
$\CL^{\fss{GR}}$ are only second order.
\hfill\qed
\medskip

\medskip

In the linear/monotonic case, there is yet another transitivity
assumption that leads to a reconstruction theorem.

\begin{defn}
\label{}
\begin{rm}
(1) Let $\pair{N}{G}$ be a linear or monotonic permutation group.
We call $\pair{N}{G}$ {\em span-transitive}
if for all distinct $a,b \in \barN$ and
all bounded $\barN$-intervals $I \subseteq J$,
there exists  $g \in G$ such that
either $g(a) \in I$ and $g(b) \not\in J$,
or $g(b) \in I$ and $g(a) \not\in J$.  Like interval-transitivity and
nest-transitivity, span-transitivity is robust in the senses of
Propositions \ref{DENSESUB} and \ref{ENLARGEGROUP}.

(2) Let $K^{\fss{SPAN}}$ denote the class of span-transitive locally
moving linear and monotonic permutation groups.

(3) Let $K^{\fss{SPAN-INT}}$ denote the class of all
$\pair{N}{G} \in K^{\fss{SPAN}}$ such that
$\fs{Var}(N,G)$ contains a bounded $\barL$-interval or an
$\barL$-ray.
\end{rm}
\end{defn}

Clearly nest-transitivity implies span-transitivity, as does
3-interval-tran\-si\-tivity.  2-interval-transitivity (or even
2-{\it o-}transi\-tivity) does not; see Example \ref{WR1}.

\begin{question}
\label{QSPAN}
\begin{rm}
In the presence of local movability:

(1) Is span-transitivity strictly weaker than
nest-transitivity?

(2) In the linear case, does
span-transitivity imply 2-interval-transitivity?

(3) In the monotonic case, is span-transitivity strictly weaker than
3-interval-transitivity?
\end{rm}
\end{question}

Some of the basic consequences of the other forms of multiple
transitivity with which we have dealt (for example,
inclusion-transitivity and the ability to deduce local movability from
the
existence of a single bounded element)
no longer hold, as far as we know.
At least we do still have:

\begin{prop}
\label{SPANDENSITY}
Let $\pair{N}{G}$ be span-transitive.  Then every orbit of
$\pair{\barN}{G}$ is dense in $\barN$.
\end{prop}

\noindent
{\bf Proof }
Let $a \in \barN$, and let $I$ be a bounded
$\barN$-interval.  Choosing
$\barN$-intervals $I_1 \subseteq J_1$ such that $a \in J_1 - I_1$,
we see that there exists $h \in G$ such that $b \eqdf h(a) \neq a$.
Now choose an $\barN$-interval $J \supseteq I$.  Then there exists
$g \in G$ such that either $g(a) \in I$ or $g(b) \in I$.  In both cases,
$G(a)$ meets $I$.
\hfill\qed
\medskip

Caution:  In the monotonic case, we don't know that the orbits of
$\pair{\barN}{\fs{Opp}(G)}$ have to be dense in $\barN$.
\medskip

Now for the analogue of Theorem \ref{NESTRT}:

\begin{theorem}\
\label{SPANRT}

{\rm (1)} First Order Reconstruction holds for
$K^{\fss{SPAN-INT}}$.

{\rm (2)} Second Order Reconstruction holds for
$K^{\fss{SPAN}}$.
\end{theorem}

\noindent
{\bf Proof }
(2) We start with the proof of (2), which is simpler and contains the
main new argument.  Let $N = L$ or $N = \fs{ED}(L)$.  The plan of the
proof is to show that the set of rays is
definable in $\fs{BACT}(N,G)$, whence the set of bounded
intervals is likewise definable, then to appeal to
Proposition~\ref{BIGINT}, and finally to distinguish the two types.

For $U,V \in \barB(L) = \barB(N)$, we write $U \comp V$ to mean
that $U
\subseteq V$ or $U \supseteq V$.
Let
\[
\varphi_{\fss{Ray}}(u) \eqqdf
(0 \neq u \neq 1) \wedge
(\forall g,h)
\mbox{\large(}u    \scomp g(u) \, \vee \,
u    \scomp h(u) \, \vee \,
g(u) \scomp h(u)\mbox{\large)},
\]
a first order formula in $\CL^{\fss{BACT}}$.
Let $\pair{N}{G} \in K^{\fss{SPAN}}$
and $U \in \barB(L)$.  We claim that:
\[
\fs{BACT}(N,G) \models \varphi_{\fss{Ray}}[U]
\mbox{ \ iff \  $U$ is an $\barL$-ray}.
\]
Clearly if $U$ is an $\barL$-ray, then
$\fs{BACT}(N,G) \models \varphi_{\fss{Ray}}[U]$.

Conversely, let $U \in \barB(L) - \dbltn{\emptyset}{\barL}$, and
suppose that $U$ is not an $\barL$-ray.
Replacing $U$ by $-U$ if necessary (permissible since
$\varphi_{\fss{Ray}}(u) \leftrightarrow \varphi_{\fss{Ray}}(-u)$),
we may assume that
some connected component of $U$ is a bounded $\barL$-interval
$I = (a,b)$.  Pick $c,d \in \barL$ such that
$c < a < b <  d$ and let $J = (c,d)$.
Pick $g \in G$
such that $g(a) \in I$ and $g(b) \not\in J$,
the other case being similar.
Since $g(a) \in I \subseteq \fs{int}(U)$ and $a$ is a boundary point
of $U$, $g(U)
\not\supseteq U$.  Then since also $g(b) \notin J$, $g(I)$ is an
$\barL$-interval meeting $-U$, so
$g(U) \not\subseteq U$.
Therefore $g(U) \ncomp U$.
Also, $I_1 \eqdf I \cap g(I)$ and
$J_1 \eqdf J \cup g(I)$
are $\barL$-intervals and $I_1 \subseteq I \subseteq J \subseteq
J_1$.
Pick $h \in G$
such that $h(a) \in I_1 \subseteq I$ and $h(b) \not\in J_1$, or
vice versa.  As with $g$, $h(U) \ncomp U$.
Finally, either
$h(a) \in I_1 \subseteq g(I)$,
making $h(U) \not\supseteq g(U)$ since $a$ is a boundary point of
$U$, and also $h(b) \notin J_1 \supseteq g(I)$, making
$h(U) \not\subseteq g(U)$ since $h(I)$ meets $-g(I)$;
or else vice versa.  Hence $g(U) \ncomp h(U)$.
Therefore $\fs{BACT}(N,G) \not\models \varphi_{\fss{Ray}}[U]$.

Clearly $U$ is a bounded $L$-interval iff $U \neq \emptyset$ and
$U$ is the intersection of some two $\barL$-rays $U$ and $V$ such
that $U \ncomp V$.  An application of Proposition~\ref{BIGINT}
completes the argument for each type.

The types are distinguished as in the linear/monotonic case of
\ref{NESTRT},
but using
\[
\varphi''_{\fss{Orp}}(g) \eqqdf (\exists u)
( \varphi_{\fss{Ray}}(u) \,\wedge\,
g(u) \nscomp u ).
\]
This concludes the proof of (2).
\smallskip

(1) Here is the plan of the proof of (1).  The set of rays
belonging to $\fs{Var}(N,G)$ may be
empty.  However,
we shall define a set $\fs{Vray}(N,G)$, with $\fs{Var}(N,G) \subseteq
\fs{Vray}(N,G) \subseteq \barB(N)$, which is $G$-invariant and {\em
does} contain a ray, and which is
FO-STR-interpretable in $\fs{VACT}(N,G)$. We shall then interpret in
$\fs{ACT}(\fs{Vray}(N,G),G)$
a dense subset of $\barL$, and finally interpret $\fs{Ep}(N,G)$.

Let $\pair{N}{G} \in K^{\fss{SPAN-INT}}$ and let $U,V \in
\fs{Var}(N,G)$.
We define the set
$\fs{Pray}(U,V) \subseteq \fs{Var}(N,G)$ by
$$
\fs{Pray}(U,V) \eqdf
\setm{W \in \fs{Var}(N,G)}{W \sim U \mbox{ and } W \not\supseteq
U
\mbox{ and } W \cap V \cap -U \neq \emptyset},
$$
and define the ``possible ray''
$$
\fs{pray}(U,V) \eqdf U + \sum  \fs{Pray}(U,V) \in \barB(L).
$$
By Lemma \ref{PST}, there is a first order formula
$\varphi_{\fss{Pray}}(u,v,w)$ which defines
$\fs{Pray}(U,V)$ in $\fs{VACT}(N,G)$.

The point of this is that when $I = (a,b) \in \fs{Var}(N,G)$ is a
bounded $\barL$-interval, then there exists $J
\in \fs{Var}(N,G)$ such that either $\fs{pray}(I,J) = (-\infty,b)$ or
dually, so that $\fs{pray}(I,J)$ is an $\barL$-ray.  To see this,
define $I$ to be
{\em left-expanding} if
$$
\setm{c < a}{\mbox{there exists } d \in (a,b) \mbox{ such that }
(c,d) \sim I}
$$
 is unbounded below.
{\em Right expanding} is defined dually.  By the span-transitivity of
$\pair{N}{G}$, $I$ is either left-expanding
or right-expanding, and without loss of generality we may assume
the former.  Pick
$J = (c,d) \sim I$, with
$c < a < d < b$.
Let $K = (e,f) \in \fs{Pray}(I,J)$.
Then $e < a$ since $K \cap J \cap -I \neq \emptyset$, whereupon $f <
b$
since $K \not\supseteq I$, making $K \subseteq (-\infty,b)$.
Hence $\fs{pray}(I,J) \subseteq (-\infty,b)$.
Since $I$ is left-expanding,
$\fs{pray}(I,J) \supseteq (-\infty,b)$.
Therefore  $\fs{pray}(I,J) = (-\infty,b)$.

\smallskip
Now we define:
\newline
$\fs{Vray}'(N,G) \eqdf \fs{Var}(N,G) \cup
\setm{\fs{pray}(U,V)}{U,V \in \fs{Var}(N,G)}$.
\newline
$\fs{Vray}(N,G) \eqdf \fs{Vray}'(N,G) \cup
\setm{-U}{U \in \fs{Vray}'(N,G)}$.
\newline
$\fs{VRAY}(N,G) \eqdf \pair{\fs{Vray}(N,G)}{\leq^B}$.
\newline
$\fs{VRACT}(N,G) \eqdf \fs{ACT}(\fs{VRAY}(N,G),G)$.
\newline
$K^{\fss{SPAN-INT}}_{\fss{VRACT}} \eqdf
\setm{\fs{VRACT}(N,G)}{\pair{N}{G} \in K^{\fss{SPAN-INT}}}$.
\smallskip

{\bf Claim 1}
$K^{\fss{SPAN-INT}}_{\fss{VRACT}}$ is FO-STR-interpretable in
$K^{\fss{SPAN-INT}}_{\fss{VACT}}$.

{\bf Proof } We let each $\vecU = \qdrpl{U}{V}{X}{X'}\in
(\fs{Var}(N,G))^4$
represent an element $R(\vecU) \in \fs{Vray}(N,G)$ in the following
way:
\newline
[1] If $X < X'$, then $\vecU$ represents $U$.
\newline
[2] If $X > X'$, then $\vecU$ represents $-U$.
\newline
[3] If $X = X'$, then $\vecU$ represents $\fs{pray}(U,V)$.
\newline
[4] If $X \ncomp X'$, then $\vecU$ represents $-\fs{pray}(U,V)$.
\newline
For $\vecU$ as above let $\fs{Comp}(\vecU)$ denote the appropriate
Boolean complement in $\fs{Var}(N,G)$:
\newline
$\fs{Bcmp}(U) = \setm{W \in \fs{Var}(N,G)}{W \cdot U =
\emptyset}$
if [1] holds.
\newline
$\fs{Bcmp}\Big(\fs{Bcmp}(U)\Big)$
if [2] holds.
\newline
$\fs{Bcmp}(\sngltn{U} \cup \fs{Pray}(U,V))$
if [3] holds.
\newline
$\fs{Bcmp}\Big(\fs{Bcmp}(\sngltn{U} \cup \fs{Pray}(U,V))\Big)$
if [4] holds.
\newline
Obviously there is a first order formula
$\varphi_{\fss{Comp}}(\vecu,w)$
in $\CL(K^{\fss{MN}}_{\fss{VACT}})$ which says that $w \in
\fs{Comp}(R(\vecu))$.
We define $\varphi_{\tiny{\leq}}(\vecu^1,\vecu^2)$
to be the formula which says that
$\fs{Comp}(\vecu^2) \subseteq \fs{Comp}(\vecu^1)$.
Now let
$\varphi_{\fss{Eq}}(\vecu^1,\vecu^2) \equiv
\varphi_{\tiny{\leq}}(\vecu^1,\vecu^2) \wedge
\varphi_{\tiny{\leq}}(\vecu^2,\vecu^1)$
and
$\varphi_{\fss{App}}(f,\vecu^1,\vecu^2) \equiv
\varphi_{\fss{Eq}}((\vecu^1)^f,\vecu^2)$.
\smallskip

Note that for $\pair{N}{G} \in K^{\fss{SPAN-INT}}$,
$$
\setm{b \in \barL}{(-\infty,b) \in \fs{Vray}(N,G)}
\mbox{ \,\Big(= }
\setm{b \in \barL}{(b,\infty) \in \fs{Vray}(N,G)}\Big)
$$
is dense in $\barL$.
For either $\fs{Var}(N,G)$ contains an  $\barL$-ray,
or $\fs{Var}(N,G)$ contains a bounded $\barL$-interval and then
$\setm{\fs{pray}(U,V)}{U,V \in \fs{Var}(N,G)}$ contains an
$\barL$- ray; and orbits of $\pair{\barN}{G}$ are dense by
Proposition \ref{SPANDENSITY}.

We define:
\newline
$\fs{Ray}(N,G) \eqdf \setm{U \in \fs{Vray}(N,G)}{U \mbox{ is an
$\barL$-ray}}$.
\newline
$\fs{Ptr}^N(U)$ denotes the endpoint of $U$ in $\barL$,
for $U \in \fs{Ray}(N,G)$.
\newline
$\fs{Epr}(N,G) \eqdf
\setm{\fs{Ptr}^N(U)}{U \in \fs{Ray}(N,G)}$.
\newline
$\fs{EPR}(N,G) \eqdf \pair{\fs{Epr}(N,G)}{\fs{Ed}^{\barN}}$.
\newline
$\fs{VREP}(N,G) \eqdf (\fs{VRAY}(N,G),\fs{EPR}(N,G);\fs{Ptr}^{N})$.
\newline
$\fs{VREPACT}(N,G) \eqdf \fs{ACT}(\fs{VREP}(N,G),G)$.
\newline
$K^{\fss{SPAN-INT}}_{\fss{VREPACT}} \eqdf
\setm{\fs{VREPACT}(N,G)}{\pair{N}{G} \in K^{\fss{SPAN-INT}}}$.
\smallskip

{\bf Claim 2 } $K^{\fss{SPAN-INT}}_{\fss{VREPACT}}$
is FO-STR-interpretable in $K^{\fss{SPAN-INT}}_{\fss{VRACT}}$.

{\bf Proof } We skip the formal details, and show that the required
formulas exist.  Every point $b \in \fs{Epr}(N,G)$ will be represented
by each of the rays $(-\infty,b)$ and $(b,\infty)$, both of which
belong to $\fs{Vray}(N,G)$.)

The first order formula $\varphi_{\fss{Ray}}(u)$ developed above in
the proof of part (2) defines the set
$\fs{Ray}(N,G)$ in $\fs{VRACT}(N,G)$, the proof being the same as
before.

The formula
$$
\varphi_{\fss{Eq}}(u,v) \eqqdf
(u = v) \vee (u = -v)
$$
holds for two rays
iff they represent the same point.

To define $\fs{Ed}$, it suffices to define $\fs{Bet}$
(Proposition \ref{BETW}).
A formula $\varphi_{\fss{Bet}}(u,v)$ that holds for the rays
$U$, $V$, $W$ iff
$\fs{Bet}(\fs{Ptr}(U),\fs{Ptr}(V),\fs{Ptr}(W))$ is given by
$$
\varphi_{\fss{Bet}}(u,v,w)
\eqqdf
\exists u' \exists v' \exists w'
\Big(\varphi_{\fss{Eq}}(u',u) \wedge \varphi_{\fss{Eq}}(v',v)
\wedge
\varphi_{\fss{Eq}}(w',w) \wedge
(u' <  v' < w') \Big).
$$
This proves Claim 2.
\smallskip

We define:
\newline
$\fs{EPRACT}(N,G) \eqdf \fs{ACT}(\fs{EPR}(N,G),G)$.
\newline
$K^{\fss{SPAN-INT}}_{\fss{EPRACT}} \eqdf
\setm{\fs{EPRACT}(N,G)}{\pair{N}{G} \in K^{\fss{SPAN-INT}}}$.
\smallskip

{\bf Claim 3 } $K^{\fss{SPAN-INT}}_{\fss{EPRACT}}$
is FO-STR-interpretable in $K^{\fss{SPAN-INT}}_{\fss{GR}}$.

{\bf Proof } Each of the classes
$K^{\fss{SPAN-INT}}_{\fss{GR}}$,
$K^{\fss{SPAN-INT}}_{\fss{VACT}}$,
$K^{\fss{SPAN-INT}}_{\fss{VRACT}}$,
$K^{\fss{SPAN-INT}}_{\fss{VREPACT}}$
is FO-STR-interpretable in its predecessor,
by Theorem \ref{INTERPVA}, Claim 1, and Claim 2.
By the transitivity of FO-STR-interpretations
(Proposition \ref{TRANS}),
$K^{\fss{SPAN-INT}}_{\fss{VREPACT}}$
is FO-STR-interpretable in $K^{\fss{SPAN-INT}}_{\fss{GR}}$.
Retaining only the information about the action on $\fs{Epr}$, we
conclude that
$K^{\fss{SPAN-INT}}_{\fss{REPACT}}$
is FO-STR-interpretable in $K^{\fss{SPAN-INT}}_{\fss{GR}}$.

There remains a minor issue.  To show that
$K^{\fss{SPAN-INT}}_{\fss{EPACT}}$
is FO-STR-interpretable in $K^{\fss{SPAN-INT}}_{\fss{GR}}$ (which is
part of First Order Reconstruction), we need to capture $\fs{Ep}(N,G)$
rather than $\fs{Epr}(N,G)$.  So now we capture $\fs{Ep}(N,G)$ in
$\fs{ACT}(\fs{VREP}(N,G),G)$.
Let $U \in \fs{Var}(N,G) - \sngltn{\emptyset}$.  The convex hull
of $U$ has at most two endpoints, corresponding to the at most two
nonempty maximal towers
$\CT$ of $\setm{R\in \fs{Ray}(N,G)}{R \supseteq U}$.
If $U$ is bounded above,
$$
\setm{x \in \barL}{x > U} = \fs{Bcmp}\Big(\fs{Bcmp}(\bigcup_{R \in
\CT}(-R))\Big)
$$
for one of these towers, and dually.  The rest is easy.

Finally, the types are distinguished just as in the proof of (2).
\hfill\qed

\begin{defn}
\label{WEAKSPANTRANS}
\begin{rm}
Let $\pair{N}{G}$ be linear or monotonic.
We say that
$\pair{N}{G}$ is {\em weakly span-transitive} if for all distinct
$a,b \in \barN$ and all bounded $\barN$-intervals $I$,
there exists $g \in G$ such that
either $g(a) \in I$ and $g(b) \not\in I$ or
$g(b) \in I$ and $g(a) \not\in I$.
Even weak span-transitivity has some bite---see the definition of {\em
convex semi-block} in the construction of Example~\ref{WRSN}.
Let $K^{\fss{WSPAN}}$ denote the class of all locally moving
weakly span-transitive $\pair{N}{G}$'s.
\end{rm}
\end{defn}

\begin{question}
\label{QWEAKSPAN}
\begin{rm}
For locally moving linear (resp., monotonic) permutation groups:

(1) Does Second Order Reconstruction hold for $K^{\fss{WSPAN}}$?

(2) Is weak span-transitivity strictly weaker than
span-transitivity? (Presumably yes.)
\end{rm}
\end{question}
\newpage

%



\section{Applications}
\label{APPLICATIONS}

In this section we further develop some of the ideas in Section
\ref{RECON},
and apply them to various nearly ordered permutation groups studied in
\cite{McHOMEOS}, \cite{BS}, \cite{Br}, and \cite{BG}.

A chain $L$ is called {\em homogeneous} if
$\pair{L}{\fs{Aut}(L)}$ is transitive, and we shall use analogous
terminology for our various flavors of multiple transitivity.  A
similar remark applies for circles.

We begin by mentioning a striking theorem due to Y. Gurevich
and W. C. Holland \cite{GH} about recognizing $\R$ and (almost) $\Q$
first order:

\begin{theorem}
\label{RR}
{\rm (Gurevich and Holland)} \
There are first order sentences $\varphi_{\R}$ and
$\varphi_{\Q}$ in
the language of groups such that for every chain $L$ for which
$\pair{L}{\fs{Aut}(L)}$ is transitive:

{\rm (1)} $\fs{Aut}(L) \models \varphi_{\R}$ iff $L \cong {\R}$.

{\rm (2)}  $\fs{Aut}(L) \models \varphi_{\Q}$ iff $L \cong {\Q}$ or $L
\cong \bar{\Q} - {\Q}$.
\end{theorem}

The heart of the proof deals with
2-{\it o-}homogeneous chains $L$.  For such $L$'s,
$\pair{L}{\fs{Aut}(L)} \in K^{\fss{LN2}}$ by Proposition~\ref{AUTLC},
and
we indicate here how our results could be
used in the
part of the proof dealing with capturing Dedekind completeness first
order.   Our Corollary~\ref{FOP}(6) tells us that
transitivity of $\fs{ACT}(\fs{EP}(L,G),G)$ can be expressed first order
in $\CL^{GR}$.  But for $G = \fs{Aut}(L)$ we have $\fs{EP}(L,G) =
\barL$ since
any pairwise disjoint set of elements of $\fs{Aut}(L)$ can be
``spliced'' together to make a single element of $\fs{Aut}(L)$.  Since
$\pair{L}{\fs{Aut}(L)}$ is transitive,
Dedekind completeness of $L$ is equivalent to transitivity of
$\fs{ACT}(\fs{EP}(L,\fs{Aut}(L)),\fs{Aut}(L))$, and so can be
expressed first order.

The rest of the proof, omitted here, involves showing that the
separability of $L$ and the existence of a countable orbit in
$\pair{\barL}{\fs{Aut}(L)}$ can be expressed by first order sentences in
the language of $\fs{ACT}(L,\fs{Aut}(L))$.

This proof does not add another item to Corollary~\ref{FOP}, as it
depends not just on having $\pair{L}{G} \in K^{\fss{LN2}}$, but also
on having $G = \fs{Aut}(L)$.

\medskip

The {\em outer automorphism group} $\fs{Outer}(G)$ of a group $G$
means the quotient of $\fs{Aut}(G)$ by $\fs{Inner}(G)$, the group of
inner automorphisms of $G$.
$G$ is called {\it complete} if all its
automorphisms are
inner and if it has trivial center, so that $\fs{Aut}(G) \cong
\fs{Inner}(G)$.  (Here the triviality of the center of $G$ will always
be obvious; our interest will be in the condition that all
automorphisms be inner.)

Recall that when First (or even Second) Order Reconstruction holds,
we have identified $\fs{Aut}(G)$ with
$\fs{Norm}_{\mbox{\scriptsize {\fs{Aut}}}(\barL^{\mbox
{\smallr}})}(G)$.  Under this identification, completeness of $G$ is
expresed as $\fs{Aut}(G) = G$.
In the next theorem we are dealing with
{\em full}
automorphism groups $G = \fs{Aut}(N)$.
In expressions such as
$\fs{Aut}\mbox{\large(}\fs{Aut}(\R)\mbox{\large)}$,  $\fs{Aut}(\R)$
is of course the group of order-automorphisms of the chain $\R$, and
$\fs{Aut}\mbox{\large(}\fs{Aut}(\R)\mbox{\large)}$ is that group's
automorphism group, which in the next theorem is identified with
$\fs{Aut}(\R^{\smallr})$.  No part
of this theorem is new; various parts are due to various people, as
was sorted out in the introduction.

\begin{theorem}
\label{RQ}\

{\rm (1)}  $\fs{Aut}\mbox{\large(}\fs{Aut}(\R)\mbox{\large)} =
\fs{Aut}\mbox{\large(}\fs{Bdd}(\R)\mbox{\large)} =
\fs{Aut}(\R^{\smallr})$, and
$\fs{Aut}(\R^{\smallr})$ is a complete group.

{\rm (2)}  $\fs{Aut}\mbox{\large(}\fs{Lft}(\R)\mbox{\large)} =
\fs{Aut}\mbox{\large(}\fs{Rt}({\R})\mbox{\large)} =
$ {\it o-}$\fs{Aut}\mbox{\large(}\fs{Aut}(\R)\mbox{\large)} =
\fs{Aut}(\R)$.

{\rm (3)}  $\fs{Aut}\mbox{\large(}\fs{Aut}(\C)\mbox{\large)} =
\fs{Aut}(\C^{\smallr})$, and
$\fs{Aut}(\C^{\smallr})$ is a complete group.

{\rm (4)} Parts (1), (2), and (3) also hold when $\R$ is replaced by
$\Q$ and $\C$ by $\C_{\Q}$.
\end{theorem}

\noindent
{\bf Proof }
This follows from Corollary~\ref{CUTDOWN} and its
analogues.  For (4) note that even in $\fs{Bdd}(\Q)$ the
irrationals form a single orbit, obviously not order-isomorphic to
$\Q$.
\hfill\qed
\medskip

A result such as
``$\fs{Aut}\mbox{\large(}\fs{Aut}(\R)\mbox{\large)} =
\fs{Aut}(\R^{\smallr})$''
does not add yet another result to the Linear Reconstruction  Package
(Theorem~\ref{LNRT}).  Rather, it is an application of that package,
valid because $\pair{\R}{\fs{Aut}(\R)}$ lies in some class $K$ for
which First Order Reconstruction holds, namely $K = K^{\fss{LN2}}$.

We continue here with a line of thought left dangling in
Corollary~\ref{CUTDOWN}.  The resulting theorem will have
Corollary~\ref{CUTDOWN} as a very special case.  Let $\pair{L}{G} \in
K^{\fss{LN2}}$, and
suppose that $\pair{L}{G}$ is transitive.
The ``tamest''
automorphisms $\alpha$ of $G$ are those for which the conjugator
$\tau = \tau_\alpha$ of Theorem~\ref{LRT} lies not just in
$\fs{Aut}(\barL^{\smallr})$ but
in $\fs{Aut}(L)$, that is, $\tau \in \fs{Norm}_{\fss{Aut}(L)}(G) =
\fs{Aut}(G) \cap \fs{Aut}(L)$ under the identification in
Corollary~\ref{1GR} of
$\fs{Aut}(G)$ with the normalizer of $G$ in
$\fs{Aut}(\barL^{\smallr})$.
We denote $\fs{Aut}(G) \cap \fs{Aut}(L)$ by
$\fs{Aut}_L(G)$, the group of automorphisms $\alpha \in
\fs{Aut}(G)$ such that $\tau_\alpha(L) = L$.
Almost as tame are the $\alpha$'s for which $\tau_\alpha \in
\fs{Aut}(L^{\smallr})$; that will be addresed in the monotonic
analogue of the present discussion.

We would like to describe  $\fs{Aut}(G)/\fs{Aut}_L(G)$.  However,
in our present degree of generality, $\fs{Aut}_L(G)$ need not be
normal in $\fs{Aut}(G)$, and we first describe instead a different
quotient.
Call an orbit
$M$ of $\pair{\barL}{G}$ {\em automorphic} if $\pair{L}{G} \cong
\pair{M}{G}$ or $\pair{L}{G} \cong \pair{M^*}{G}$.  (These
isomorphisms are not required to arise from $\tau_{\alpha}$'s,
though
it will turn out that in fact they do.  Also, we do not distinguish
between $\sigma$ and its unique extension to $\barL$.)  Let
$\fs{Autorb}(L,G)$ denote the set of automorphic orbits, one of which
is $L$ itself.  Let
\begin{eqnarray*}
\fs{Aut}_L^{\fss{Sgr}}(G) & \eqdf &  \setm{\alpha \in \fs{Aut}_L(G)}
{\tau_\alpha(M) = M \mbox{ for all } M \in \fs{Autorb}(L,G)} \\*
& \leq  & \fs{Aut}_L(G).
\end{eqnarray*}
Clearly $\fs{Aut}_L^{\fss{Sgr}}(G) \unlhd \fs{Aut}(G)$, and it is
$\fs{Aut}(G)  / \fs{Aut}_L^{\fss{Sgr}}(G)$ that we describe.

For each $M \in \fs{Autorb}(L,G)$, every isomorphism
$\iso{\sigma \cup \beta}
{\pair{L}{G}}{\pair{M}{G}}$ (or $\pair{M^*}{G}$) induces a
permutation $M_1 \mapsto \sigma(M_1)$ of $\fs{Autorb}(L,G)$.
(Certainly $\sigma$ permutes the set of {\em
all} orbits of $\pair{\barL}{G}$.  Moreover,
\newline
$\iso{\sigma \cup \sigma\beta \sigma^{-1}}
{\pair{M_1}{G}}{\pair{\sigma(M_1)}{G}}$ (or
$\pair{\sigma(M_1)^*}{G}$).
It follows easily that the induced map is a permutation of
$\fs{Autorb}(L,G)$.)  The group $\fs{Perm}(L,G)$ of permutations of
$\fs{Autorb}(L,G)$ induced in this way is transitive on
$\fs{Autorb}(L,G)$.  For each automorphism $\alpha$ of $G$,
$\iso{(\tau_{\alpha}\rest L) \cup \alpha}
{\pair{L}{G}}
{\pair{\tau_\alpha(L)}{G}}$ by Corollary~\ref{1GRORBITS},
and this induces an element $p \in \fs{Perm}(L,G)$.  Conversely,
each $p
\in \fs{Perm}(L,G)$ arises in this way, for $p$ arises from some
$\iso{\sigma \cup \beta}
{\pair{L}{G}}{\pair{M}{G}}$ (or $\pair{M^*}{G}$), with $\beta \in
\fs{Aut}(G)$, and $\beta$ is the
required automorphism of $G$ since $\tau_\beta = \sigma$ by the
uniqueness of $\tau_\beta$ in Theorem~\ref{LRT}.

\begin{theorem}
\label{ORBITCORR}
Let $\pair{L}{G} \in K^{\fss{LN2}}$, and suppose that $\pair{L}{G}$ is
transitive.

{\rm (1)} If $\pair{L}{G} \cong \pair{L^*}{G}$, then an isomorphism
$$
\fs{Aut}(G) / \fs{Aut}_L^{\fss{Sgr}}(G) \, \cong \,
\fs{Perm}(L,G) \times \{\pm 1\}
$$
is induced by the homomorphism $\alpha \mapsto
\pair{p}{\pm 1}$, where $p \in \fs{Perm}(L,G)$ is induced by
$\tau_\alpha
\cup \alpha$, and where the second coordinate is $+1$ if $\tau$
preserves order and $-1$ if $\tau$ reverses order.  $\fs{Aut}(G) /
\fs{Aut}_L^{\fss{Sgr}}(G)$
contains  two
cosets of $\fs{Aut}_L^{\fss{Sgr}}(G)$ for each $M \in
\fs{Autorb}(L,G)$.
Also, {\it o-}$\fs{Aut}(G) / \fs{Aut}_L^{\fss{Sgr}}(G) \cong
\fs{Perm}(L,G)$.

{\rm (2)} If $\pair{L}{G} \not\cong \pair{L^*}{G}$, then
$\tau_\alpha$ preserves
order for every $\alpha \in \fs{Aut}(G)$, and an isomorphism
$$
\fs{Aut}(G) / \fs{Aut}_L^{\fss{Sgr}}(G) \cong \fs{Perm}(L,G)
$$
is induced by the homomorphism $\alpha \mapsto p$ (with $p$ as
above).  $\fs{Aut}(G) / \fs{Aut}_L^{\fss{Sgr}}(G)$ contains one coset
of $\fs{Aut}_L^{\fss{Sgr}}(G)$ for each
$M \in \fs{Autorb}(L,G)$.  Also, {\it o-}$\fs{Aut}(G) /
\fs{Aut}_L^{\fss{Sgr}}(G)$ is the
group of cosets of $\fs{Aut}_L^{\fss{Sgr}}(G)$ corresponding to
order-preserving $\tau_\alpha$'s, a subgroup  of
$\fs{Aut}(G) / \fs{Aut}_L^{\fss{Sgr}}(G)$ of index $\leq 2$.
\end{theorem}

\noindent
{\bf Proof }
In each case we have $\tau_{\alpha_1\alpha_2}(L) =
\tau_{\alpha_1}(\tau_{\alpha_2}(L))$, so we do have a
homomorphism, and its kernel is $\fs{Aut}_L^{\fss{Sgr}}(G)$.

If (1) holds, then $\alpha \mapsto \pair{p}{\pm 1}$ is a
homomorphism from $\fs{Aut}(G)$ onto $\fs{Perm}(L,G) \times
\{\pm 1 \}$.  For since $\pair{L}{G}
\cong \pair{L^*}{G}$,
then for each $M \in \fs{Autorb}(L,G)$ there exist isomorphisms
$\pair{L}{G} \rightarrow \pair{M}{G}$  and also $\pair{L}{G}
\rightarrow
\pair{M^*}{G}$.

If (2) holds, then $\alpha \mapsto p$ is a
homomorphism from $\fs{Aut}(G)$ onto $\fs{Perm}(L,G)$.  For if
there were an $M$ with
isomorphisms $\pair{L}{G} \rightarrow \pair{M}{G}$ and also
$\pair{L}{G} \rightarrow \pair{M}{G}$,
then $\pair{L}{G} \cong
\pair{L^*}{G}$.
The rest is clear.
\hfill\qed
\medskip

Now we specialize to $G = \fs{Aut}(L)$, the full
automorphism group.
This makes $\fs{Aut}_L(G) = \fs{Inner}(G) \unlhd \fs{Aut}(G) $, and
$\fs{Aut}(G)/\fs{Aut}_L(G) =
\fs{Outer}(G)$.
We denote {\it o-}$\fs{Aut}(G) / \fs{Inner}(G)$ by
{\it o-}$\fs{Outer}(G)$; its index in ${Outer}(G)$ is $\leq 2$.

The condition that
``$\pair{L}{G} \cong \pair{L^*}{G}$'' becomes ``$L$ is symmetric''.
``$\pair{L}{G} \cong
\pair{M}{G}$ or $\pair{L}{G} \cong \pair{M^*}{G}$'' becomes
``$L \cong M$ or $L \cong M^*$ and the action of $\fs{Aut}(L)$ on
$M$ is all of $\fs{Aut}(M)$'', making the definition of ``automorphic
orbit'' more tangible..

Most important of all is that $\fs{Aut}_L^{\fss{Sgr}}(G)$ becomes
simply $\fs{Aut}_L(G)$.  For if $\tau_\alpha(L) = L$,
then $\tau_\alpha(M_1) = M_1$ for all $M_1 \in \fs{Autorb}(L,G)$.
This is because if $\tau_\alpha(L) = L$,
then since $G = \fs{Aut}(L)$ by hypothesis and since $\tau_\alpha$
preserves order, we have $\alpha \in
\fs{Aut}_L(G) = \fs{Inner}(G)$, so that $\tau_\alpha$ respects the
orbit
$M_1$ of $\pair{\barL}{G}$.

Also, the assumption that $G = \fs{Aut}(L)$ turns this discussion into
a discussion
about chains, specifically---in view of Proposition~\ref{AUTLC}---
about
2-{\it o-}homogeneous chains $L$, in which we happen to let $G =
\fs{Aut}(L)$.
Accordingly, we write $\fs{Autorb}(L)$ for $\fs{Autorb}(L,G)$ and
similarly for $\fs{Perm}(L)$.

Theorem~\ref{ORBITCORR} specializes to:

\begin{theorem}
\label{ORBITCORR2}
Let $L$ be a 2-{\it o-}homogeneous chain, and let $G = \fs{Aut}(L)$.

{\rm (1)} If $L$ is symmetric, then an isomorphism
$$
\fs{Outer}\mbox{\large(}\fs{Aut}(L)\mbox{\large)} = \fs{Outer}(G) =
\fs{Aut}(G)/\fs{Aut}_L(G) \, \cong \,
\fs{Perm}(L) \times \{\pm 1\}
$$
is induced by the homomorphism $\alpha \mapsto
\pair{p}{\pm 1}$, where $p \in \fs{Perm}(L)$ is induced by
$\tau_\alpha
\cup \alpha$, and where the second coordinate is $+1$ if $\tau$
preserves order and $-1$ if $\tau$ reverses order.  $\fs{Outer}(G)$
contains  two
cosets of $\fs{Aut}_L(G)$ for each $M \in \fs{Autorb}(L)$.
Also, {\it o-}$\fs{Outer}(G) \cong \fs{Perm}(L)$.

{\rm (2)} If $L$ is not symmetric, then $\tau_\alpha$ preserves
order for every $\alpha \in \fs{Aut}(G)$, and an isomorphism
$$
\fs{Outer}\mbox{\large(}\fs{Aut}(L)\mbox{\large)} = \fs{Outer}(G) =
\fs{Aut}(G)/\fs{Aut}_L(G) \cong \fs{Perm}(L)
$$
is induced by the homomorphism $\alpha \mapsto p$ (with $p$ as
above).  $\fs{Outer}(G)$ contains one coset of $\fs{Aut}_L(G)$ for
each
$M \in \fs{Autorb}(L)$.  Also, {\it o-}$\fs{Outer}(G)$ is the
group of cosets of $\fs{Aut}_L(G)$ corresponding to
order-preserving $\tau_\alpha$'s, a subgroup  of
$\fs{Outer}(G)$ of index $\leq 2$.
\end{theorem}

The last two theorems carry over to the other three types of nearly
ordered permutation groups.  The analogues of
Theorem~\ref{ORBITCORR2} all deal with 2-{\it o-}homogeneous
chains
and circles (see Proposition~\ref{AUTLC}).  In the monotonic case, if
$\pair{\fs{ED}(L)}{G}$ contains order-reversing permutations, then in
case (1) of both theorems
$\fs{Perm}(L,G) \times \{\pm 1 \}$ becomes just $\fs{Perm}(L,G)$
and there is only one coset per $M$, and
case (2) cannot occur.  {\it o-}$\fs{Outer}(G)$ is not defined in the
monotonic case.
As mentioned earlier, the monotonic case $\pair{\fs{ED}(L)}{G} =
\pair{L^{\smallr}}{G}$ can be used to deal with ``semi-tame''
automorphisms in the linear case since these coincide with ``tame''
automorphisms in the monotonic case.
The circular case requires no comment, and the monocircular case is
related to it as is the monotonic case to the linear case.

Theorem~\ref{ORBITCORR2} and its analogues have
Theorem~\ref{RQ}
as a very special case; these illustrate part (1).
$|\fs{Outer}\mbox{\large(}\fs{Aut}(L)\mbox{\large)}| = 2$ and
$\fs{Aut}(L^{\smallr})$ is a complete group, and similarly for circles
$C$.  They also have $\alpha$-sets and their circular analogues as
examples, again with no order-preserving outer automorphisms; see
Theorem~\ref{ALPHA} and Example~\ref{EXALPHA} and their linear
analogues.

Another illustration of (1), this time with
$\fs{Outer}\mbox{\large(}\fs{Aut}(L)\mbox{\large)}$ the
Klein 4-group and
$|\fs{Outer}\mbox{\large(}\fs{Aut}(L^{\smallr})\mbox{\large)}| = 2$,
comes
from an example due to Holland \cite{Ho2}; see (the linear analogue
of) our Example~\ref{EXOUTER}.

For a circular version of Holland's chain see Example~\ref{EXOUTER}.
Also see  Example~\ref{EXOUTERCTBL} and its linear analogue.

The semilong line $L$ illustrates case
(2) with {\it o-}$\fs{Outer}(\fs{Aut}(L)) =
\newline
\fs{Outer}(\fs{Aut}(L))$.

For other examples, including one illustrating case (2) with
\newline
{\it o-}$\fs{Outer}(\fs{Aut}(L)) \neq \fs{Outer}(\fs{Aut}(L))$, see
Example~\ref{EXALPHA} (with $C = C_{10}$).

\medskip

McCleary investigated automorphisms of linear and monotonic
permutation groups in \cite{McHOMEOS}.  (There, by definition, a
monotonic permutation group was required to include
order-reversing permutations.)  The hypotheses were that
$\pair{N}{G}$ is (exactly) 3-{\it o-}transitive and contains a positive
bounded element.  Here, in Section~\ref{RECON}, the hypotheses are
much weaker, with interval-transitivity instead of exact transitivity,
with the degree of transitivity reduced to 2 in the linear case, and
with the assumed bounded element not required to be positive.
Nevertheless we obtain from the weaker hypotheses the conclusions
reached in \cite{McHOMEOS}.  For
the full automorphism group
$\pair{L}{\fs{Aut}(L)}$, there is no real difference between the two
sets of hypotheses:  If $\pair{L}{\fs{Aut}(L)}$ is
2-interval-transitive
then it is 3-{\it o-}transitive on each of its
$\barL$-orbits and contains a positive bounded
element; see Proposition~\ref{AUTLC}.
This also holds in the monotonic case; just apply the linear result to
$\pair{\fs{ED}(L)}{\fs{Opp}(G)}$.

Many of these generalizations have already been dealt with, and for
all four types of nearly ordered permutation groups; these are among
the results in Section~\ref{RECON}.  Also, we mention here that
Theorems \ref{ORBITCORR} and
\ref{ORBITCORR2} could have been proved in
\cite{McHOMEOS} with the less powerful techniques available there,
but weren't.

We turn now to circular/monocircular analogues of results in
\cite{McHOMEOS} that did not appear in
Section~\ref{RECON}.  Ordinarily these will be about full
automorphism groups $G = \fs{Aut}(N)$, and then in the
linear/monotonic
case the present methods will give nothing not already contained in
\cite{McHOMEOS}.  We shall omit the statements of such
linear/monotonic results, and deal only with the analogues.  In these
analogues about $\pair{C}{\fs{Aut}(C)}$ and
$\pair{C^{\smallr}}{\fs{Aut}(C^{\smallr})}$, the hypothesis on $C$ will
be either 2-interval-homogeneity or 2-{\it o-}homogeneity.
By Proposition~\ref{AUTLC}, it follows that $\pair{C}{\fs{Aut}(C)} \in
K^{\fss{CR3}}$, so that First Order Reconstruction holds (see
Theorem~\ref{LNRT}).
No mention is made of
$\pair{C}{\fs{Bdd}(C)}$ since
$\fs{Bdd}(C) = \fs{Aut}(C)$; this applies also to other results later in
this section for non-full $\pair{N}{G}$'s.

The next theorem specializes several results in Section~\ref{RECON}
to the case of full automorphism groups.  Part (1) is an analogue of
Corollary 11 of
\cite{McHOMEOS}.  (2) is an
analogue of Corollary 14, and (3) of Corollary 15.  More
exactly, (3) is an analogue of an improved Corollary 15, improved by
omitting the monotonic case hypothesis that $L$ is symmetric.  This
improved Corollary 15 is
true for the same reason that Theorem~\ref{CIRCAT}(3) is true; see
the proof below.

In (1) $\barC^{\smallr}$ cannot be improved to $C^{\smallr}$;
$\fs{Aut}\mbox{\large(}\fs{Aut}(C)\mbox{\large)} \not\subseteq
\fs{Aut}(C^{\smallr})$ in Example~\ref{EXALPHA} (with $C = C_{10}$)
and in Example~\ref{EXOUTER}.
(2) has the circular parts of Theorem~\ref{RQ} as
special cases, as do Theorem~\ref{ALPHA} and
Example~\ref{EXALPHA}.  Example~\ref{EXALPHA} also provides a
good illustration of (3).

\begin{theorem}
\label{CIRCAT}\
Let $C$ be a 2-interval-homogeneous (hence highly
interval-homogeneous) circle.

{\rm (1)}
$\fs{Aut}\mbox{\large(}\fs{Aut}(C^{\smallr})\mbox{\large)}$ and
$\fs{Aut}\mbox{\large(}\fs{Aut}(C)\mbox{\large)}$ are trapped
between $\fs{Aut}(C^{\smallr})$ and
$\fs{Aut}(\barC^{\smallr})$:
$$
\fs{Aut}(C^{\smallr}) \leq
\fs{Aut}\mbox{\large(}\fs{Aut}(C^{\smallr})\mbox{\large)} \leq
\fs{Aut}\mbox{\large(}\fs{Aut}(C)\mbox{\large)} \leq
\fs{Aut}(\barC^{\smallr}).
$$

{\rm (2)$\,\,$(a)} Suppose there is no orbit $M \neq C$ of
$\pair{\barC}{\fs{Aut}(C)}$ such that $C \cong M$ and such that the
action of $\fs{Aut}(C)$ on $M$ is all of $\fs{Aut}(M)$.  Then
$|\fs{Outer}\mbox{\large(}\fs{Aut}(C)\mbox{\large)}| \leq 2$.

$\,\,\,\,\,\,\,\,\,\,${\rm (b)} Suppose there is no orbit $M \neq C$ of
$\pair{\barC}{\fs{Aut}(C)}$ such that $C \cong M$ or $C \cong M^*$,
and such that the
action of $\fs{Aut}(C)$ on $M$ is all of $\fs{Aut}(M)$.  Then
$\fs{Aut}(C)$ is a complete group.

$\,\,\,\,\,\,\,\,\,\,${\rm (c)} Suppose there is no orbit $M \neq C$ of
$\pair{\barC}{\fs{Aut}(C^{\smallr})}$ such that $C\cong M$ or
$C\cong M^*$, and
such that the action of $\fs{Aut}(C^{\smallr})$ on $M^{\smallr}$ is
all of
$\fs{Aut}(M^{\smallr})$.  Then $\fs{Aut}(C^{\smallr})$ is a
complete group.

{\rm (3)} $\fs{Aut}\mbox{\large(}\fs{Aut}(\barC)\mbox{\large)} =
\fs{Aut}(\barC^{\smallr})$ and
$\fs{Aut}(\barC^{\smallr})$ is a complete group.
\end{theorem}

\noindent
{\bf Proof } (1) follows immediately from Corollary~\ref{CHAR}, just
as in \cite{McHOMEOS}.

(2) follows immediately from Theorem~\ref{ORBITCORR2}, or
alternately from
Corollary~\ref{CUTDOWN} as in \cite{McHOMEOS}.

(3) The proof is like that of \cite[Corollary 15]{McHOMEOS}, except
that there is no need to establish exact 3-{\it o-}transitivity.  It is
sufficient to observe instead that 3-interval-transitivity carries over
from $\pair{C}{\fs{Aut}(C)}$ to $\pair{\barC}{\fs{Aut}(\barC)}$, from
which (3) follows.  (This is a good illustration of the advantages of
interval-transitivity over exact transitivity.)
\hfill\qed
\medskip

Let $L$ be a chain, and let $a \in \barL$.
The {\em left character} of $a$ means the cofinality of the chain
$\setm{x \in L}{x < a}$,  necessarily a regular cardinal number
$\aleph_{\alpha}$ (that is, all its cofinal subsets also have cardinality
$\aleph_{\alpha}$); and dually for {\em right character}.  $a$ is said
to have character $c_{\alpha\beta}$ if it has left character
$\aleph_{\alpha}$ and right character $\aleph_{\beta}$.
A {\em hole} in $L$ is an element of
$\barL - L$.

For a circle there are analogous definitions of clockwise
character, counterclockwise character, etc.

Part (1) of the next theorem is an analogue of a well known theorem
about chains, and (2) and (3) are analogues of \cite[Corollaries 16
and 17]{McHOMEOS}.

\begin{theorem}
\label{CIRCATCTBL}\
Let $C$ be a 2-{\it o-}homogeneous) circle.

{\rm (1)}  All holes in $C$ of countable character lie in the same
$\barC$-orbit of $\pair{C}{\fs{Aut}(C)}$.

{\rm (2)} Suppose the points of $C$ have countable character.  Then:

$\,\,\,\,\,\,\,\,\,\,${\rm (a)}
$|\fs{Outer}\mbox{\large(}\fs{Aut}(C)\mbox{\large)}| =
1, 2, \mbox{or } 4$ (and if $4$ then
$\fs{Outer}\mbox{\large(}\fs{Aut}(C)\mbox{\large)}$ is the Klein
4-group);
and $|\fs{Outer}\mbox{\large(}\fs{Aut}(C^{\smallr})\mbox{\large)}|
\leq 2$.

$\,\,\,\,\,\,\,\,\,\,${\rm (b)}
$\fs{Aut}\mbox{\large(}\fs{Aut}(C^{\smallr})\mbox{\large)} =
\fs{Aut}\mbox{\large(}\fs{Aut}(C)\mbox{\large)}$.

 {\rm (3)} Suppose all points and holes of $C$ have countable
character.
Then
$\fs{Aut}\mbox{\large(}\fs{Aut}(C)\mbox{\large)} =
\fs{Aut}(C^{\smallr})$ and $\fs{Aut}(C^{\smallr})$ is a
complete group.
\end{theorem}

We remark that in (3) the hypothesis that the holes also have
countable
character is in fact needed; see Example~\ref{EXOUTERCTBL}, in
which point characters are countable but $\fs{Aut}(C)$ has an
order-preserving outer automorphism.

\medskip

\noindent
{\bf Proof } (1) is proved just as in the linear case, constructing $g
\in \fs{Aut}(C)$ sending one hole of countable character to another
by splicing together maps between appropriate
$C$-intervals.

(2) Given (1), (2a) follows immediately from
Corollary~\ref{ORBITCORR2}; and the proof
of (2b) is like that of the analogue in \cite{McHOMEOS}.

(3) This proof is like that of \cite[Corollary 17]{McHOMEOS} but with
one change.  Now, rather than $\R$, it is $\C$ which we need to know
cannot be written as the disjoint union of two dense
2-{\it o-}homogeneous subcircles which are orientation-isomorphic
to each other (and similarly for any Dedekind complete
3-{\it o-}homogeneous circle $C$).  The proof of {\em this} is the
same as Holland's proof of the analogous fact for $\R$ \cite{Ho3}.
\hfill\qed
\medskip

We gather here some information on $\alpha$-sets (see
\cite{McALPHA}, for example).
An {\em $\alpha$-set} is a chain $L$ of cardinality
$\aleph_{\alpha}$ in which for any two (possibly empty) subsets
$U_1 < U_2$ of $L$ having cardinality less than
$\aleph_{\alpha}$, there exists $a \in L$ such that $U_1 < a < U_2$.
For an $\alpha$-set to exist, $\aleph_{\alpha}$ must be a regular
cardinal.  Conversely, let $\aleph_{\alpha}$ be a regular
cardinal. Then, assuming the Generalized Continuum Hypothesis,
there exists an $\alpha$-set $\LL_{\alpha}$, unique up to
order-isomorphism.  $\LL_{\alpha}$ has an order-reversing
permutation.  The points of $\LL_{\alpha}$ have character
$c_{\alpha\alpha}$, and the coinitiality and cofinality of
$\LL_{\alpha}$
are $\aleph_{\alpha}$.  All bounded and unbounded
$\LL_{\alpha}$-intervals $I$
are themselves $\alpha$-sets and thus order-isomorphic to each
other and to their reverses $I^*$. $\LL_{\alpha}$ is highly
{\it o-}homogeneous and
$\fs{Aut}(\LL_{\alpha})$ is
highly {\it o-}transitive on each of its
$\bar{\LL}_{\alpha}$-orbits.

The orbits of points and holes of
$\pair{\LL_{\alpha}}{\fs{Aut}(\LL_{\alpha})}$ are as follows:
One orbit is $\LL_{\alpha}$.  There is a single orbit of holes of
character $c_{\alpha\alpha}$; it is of cardinality
$2^{\aleph_{\alpha}}$, and thus is not order-isomorphic to
$\LL_{\alpha}$.  For each regular cardinal $\aleph_\beta <
\aleph_\alpha$, there is a single orbit of holes of character
$c_{\alpha\beta}$ and a single orbit of holes of character
$c_{\beta\alpha}$.
These orbits are preserved by $\fs{Aut}(\bar{\LL}_{\alpha})$, except
that
$\LL_{\alpha}$ and the orbit of holes of character
$c_{\alpha\alpha}$ coalesce into a single orbit
of $\fs{Aut}(\bar{\LL}_{\alpha})$.

Now we define the {\em circular $\alpha$-set} $\C_{\alpha}$ to be
the circularly ordered set based on $\LL_{\alpha}$.  ${\C}_0$ is the
familiar rational circle.  The bounded
$\C_\alpha$-intervals, totally ordered in the counterclockwise
direction, are $\alpha$-sets.
(For $\C_\alpha$-intervals which contain the hole at which the two
ends of $\LL_\alpha$ are joined, this is because the concatenation of
two $\alpha$-sets is an $\alpha$-set.)  This makes $\C_{\alpha}$
highly {\it o-}homogeneous and $\fs{Aut}(\C_{\alpha})$
highly {\it o-}transitive on each of its
$\bar{\C}_{\alpha}$-orbits.

The description of the orbits of points and holes of
$\pair{\C_{\alpha}}{\fs{Aut}(\C_{\alpha})}$ is analogous to that for
$\pair{\LL_{\alpha}}{\fs{Aut}(\LL_{\alpha})}$; in particular,
$\C_{\alpha}$ is not orientation-isomorphic to the orbit of holes of
character $c_{\alpha\alpha}$.

\begin{theorem} {\rm (GCH)} \
\label{ALPHA}\

{\rm (1)}
$\fs{Aut}\mbox{\large(}\fs{Aut}(\C_{\alpha})\mbox{\large)} =
\fs{Aut}(\C_{\alpha}^{\smallr})$, and
$\fs{Aut}(\C_{\alpha}^{\smallr})$ is a complete group.

{\rm (2)}
$\fs{Aut}\mbox{\large(}\fs{Aut}(\bar{\C}_{\alpha})\mbox{\large)} =
\fs{Aut}(\bar{\C}_{\alpha}^{\smallr})$, and
$\fs{Aut}(\bar{\C}_{\alpha}^{\smallr})$ is a complete group.
\end{theorem}

\noindent
{\bf Proof }
(1) follows from Corollary~\ref{CIRCAT}(2a), just as in
\cite{McHOMEOS}.
(2) is an instance of Corollary~\ref{CIRCAT}(3).
\hfill\qed

\begin{example} {\rm (GCH)} \
\label{EXALPHA}
In the circular 1-set ${\C}_1$, let $C_{10}$ be the set of holes of
character $c_{10}$ and let $C_{01}$ be its dual.
Let $C_{11}$ be the set of elements (including
points)  of
$\bar{\C}_1$ of character $c_{11}$.
Then
$$
\fs{Aut}(C_{10}) = \fs{Aut}(C_{01}) = \fs{Aut}(C_{10} \cup C_{01}) =
\fs{Aut}(C_{11})
=\fs{Aut}(\bar{\C}_1).
$$

{\rm (1)}
$C_{10}$ is highly {\it o-}homogeneous and not symmetric.
$C_{10}$ and $C_{01}$ coalesce into a single orbit of
$\fs{Aut}((C_{10} \cup C_{01})^{\smallr})$.

{\rm (2)}
$C_{11}$ is highly {\it o-}homogeneous
and symmetric.  $C_{10}$ and $C_{01}$,
the two orbits of $\fs{Aut}(C_{11})$ in $\barC_{11} - C_{11}$,
coalesce into a single orbit
$\barC_{11}$-orbit of
$\fs{Aut}(C_{11}^{\smallr})$.
\end{example}

As mentioned after Theorem~\ref{ORBITCORR},
Holland \cite{Ho2} constructed a symmetric
2-{\it o-}homogeneous chain $L$ such that $\fs{Aut}(L)$ has an
order-preserving outer automorphism.
In $\pair{\barL}{\fs{Aut(L)}}$ there is one orbit $M \neq L$ such
that $M \cong L$, and $M$ is automorphic.  By part (1) of
Theorem~\ref{ORBITCORR},
$\fs{Outer}(\fs{Aut}(L))$ is the Klein 4-group.

In preparation for the circular analogue, we mention a few details of
the construction.
Of course all $a \in L$ yield
(up to isomorphism) the same chain $\setm{x \in L}{x < a}$ and the
same chain $\setm{x \in L}{x > a}$, which we denote by
$\fs{Lftray}(L)$
and $\fs{Rtray}(L)$, respectively.  Holland describes a certain tower
of equivalence relations on $L$ ({\em not} respected by
$\fs{Aut}(L)$, and in it there is a largest proper equivalence relation
$R$.  The set of $R$-classes, under the order inherited from $L$, is
order-isomorphic to $\Z$.  Moreover, for each $R$-class $P$ and each
$a \in P$, $\setm{x \in P}{x < a} \cong \fs{Lray}(L)$ and dually.
Picking $a_1 < a_2$ such that $a_1$ and $a_2$ lie in consequtive
$R$-classes, we see that the $L$-interval $(a_1,a_2)$ is isomorphic to
the concatenation $\fs{Rtray}(L) \cat \fs{Lftray}(L)$.  By
2-{\it o-}homogeneity, this must then be true for every $L$-interval.

Now let $\dot{C}$ be the circularly ordered set based on Holland's
chain $L$.  $\dot{C}$ is
highly {\it o-}homogeneous by Proposition~\ref{WRAP}.  We apply
Theorem~\ref{ORBITCORR2} to get

\begin{example}
\label{EXOUTER}
$\dot{C}$ is highly {\it o-}homogeneous.
$\fs{Outer}\mbox{\large(}\fs{Aut}(\dot{C})\mbox{\large)}$ is the
Klein 4-group
and
$|\fs{Outer}\mbox{\large(}\fs{Aut}(\dot{C}^{\smallr})\mbox{\large)}
| = 2$.
\end{example}

Holland's chain has points of character $c_{11}$.  However, McCleary
\cite[Example 20]{McHOMEOS} modified the construction to
produce a symmetric 2-{\it o-}homogeneous chain $L_0$ for which
$\fs{Aut}(L_0)$ has an order-preserving outer automorphism, and in
which the points of $L_0$ have countable character.  (It cannot be
arranged that all $a \in \barL_0$ have countable character; see
Theorem~\ref{CIRCATCTBL}(3).)  This involved the $\alpha$-set
$\LL_\alpha$ ($\alpha > 0$), and an index chain which was taken to
be $\omega_0$ but should instead have been
$\omega_0^{\omega_0}$.  (For why $\omega_0$ doesn't work, see
the discussion of Holland's example in \cite[p.174]{G1}.)  As before,
for $a_1 < a_2$ with $a_1$ and $a_2$ from different $R$-classes, the
$L_0$-interval $(a_1,a_2)$ is isomorphic to $\fs{Rtray}(L_0) \cat
\fs{Lftray}(L_0)$; this is because deleting a
closed left ray $\setm{x}{x \leq a}$ from an $\alpha$-set leaves an
$\alpha$-set.

Now let $\ddot{C}$ be the circular order based on $L_0$.  As in the
previous example, we get

\begin{example} {\rm (GCH)} \
\label{EXOUTERCTBL}
$\ddot{C}$ is highly {\it o-}homogeneous, this time with countable
point character.
$\fs{Outer}\mbox{\large(}\fs{Aut}(\ddot{C})\mbox{\large)}$ is the
Klein 4-group,
and
$|\fs{Outer}\mbox{\large(}\fs{Aut}(\ddot{C}^
{\smallr})\mbox{\large)}|
= 2$.
\end{example}

Continuing with analogues of results in \cite{McHOMEOS},
we look at several relatively small groups acting on the real
circle $\C$ whose automorphism groups can be exactly ascertained.
We look first at groups involving differentiability, then at groups
involving piecewise linearity.

\begin{defn}\label{}
\begin{rm} \

(1) Let $D$ denote the group of diffeomorphisms $g$ of $\C$ ($C^1$
not required), that is, of those $g \in \fs{Aut}(\C^{\smallr})$ such
that $g$ and $g^{-1}$ are differentiable, or again, such that $g$ has
everywhere a nonzero derivative.

(2) Let $D^{\fss{1sd}}$ be the group of $g \in
\fs{Aut}(\C^{\smallr})$
such that $g$ and $g^{-1}$ are one-sided differentiable on each side.

(3) Let $D^{\fss{ccw}}$ be the group of order-preserving $g \in
\fs{Aut}(\C^{\smallr})$
such that $g$ and $g^{-1}$ are one-sided differentiable on the
counterclockwise side.

(4) Let \fs{PD} denote the group of $g \in \fs{Aut}({\C}^{\smallr})$
which are piecewise differentiable (and thus finitely piecewise
differentiable).
\end{rm}
\end{defn}

\begin{theorem}
\label{DIFF}\

{\rm (1)} $\fs{Aut}\mbox{\large(}\fs{Opp}(D)\mbox{\large)} = D$
and $D$ is a complete group.

{\rm (2)}
$\fs{Aut}\mbox{\large(}\fs{Opp}(D^{\fss{1sd}})\mbox{\large)} =
D^{\fss{1sd}}$ and $D^{\fss{1sd}}$ is a complete group.

{\rm (3)} $\fs{Aut}(D^{\fss{ccw}})$ is a complete group.

{\rm (4)} $\fs{Aut}\mbox{\large(}\fs{Opp}(\fs{PD})\mbox{\large)} =
\fs{PD}$ and \fs{PD} is a complete group.
\end{theorem}

\noindent
{\bf Proof }
The proofs are the same as those of the analogues in
\cite{McHOMEOS}.
\hfill\qed

\begin{defn}\label{}
\begin{rm} \
We view the real circle $\C$ as ${\C} = {\R}/{\Z}$, permitting us to
discuss piecewise linearity just as we would on $\R$.
Let \fs{PL} denote the group of $g \in \fs{Aut}({\C}^{\smallr})$
which are piecewise linear (and thus finitely piecewise
linear).
\end{rm}
\end{defn}

\begin{theorem}
\label{PL}\
 $\fs{Aut}\mbox{\large(}\fs{Opp}(\fs{PL})\mbox{\large)} = \fs{PL}$
and \fs{PL} is a complete group.
\end{theorem}

\noindent
{\bf Proof }
The proof is the same as that of the analogue in \cite{McHOMEOS}.
\hfill\qed

\bigskip

Bieri and Strebel \cite{BS} have conducted an extensive analysis of
groups of piecewise linear homeomorphisms on the real line $\R$.
We discuss here those of their results to which the present paper
makes a contribution.

\begin{defn}\label{GIAP}
\begin{rm}
Let $I$ be a nonsingleton closed interval of
$\R$ ($[c,d]$ with $c<d$, or a closed ray, or $\R$ itself), $P$ a
nontrivial subgroup of the multiplicative group of
positive reals, and $A$ a nontrivial subgroup of the additive group of
$\R$ which is closed under multiplication by elements of $P$.

(1) $G \eqdf G(I;A,P)$ denotes the group of all
PL-homeomorphisms $g:\R \rightarrow \R$ such that

$\,\,\,\,\,\,\,\,\,\,$(a) $\fs{Supp}(g) \subseteq I$,

$\,\,\,\,\,\,\,\,\,\,$(b) The slopes of $g$ are elements of $P$ (so $g$ is
increasing),

$\,\,\,\,\,\,\,\,\,\,$(c) $g$ has only finitely many break points
(discontinuities of slope),

$\,\,\,\,\,\,\,\,\,\,$(d) The break points of $g$ are elements of $A$,
and

$\,\,\,\,\,\,\,\,\,\,$(e) $g(A) = A$.

(2) $B \eqdf B(I;A,P)$ denotes the subgroup of $G$ consisting of
those $g$ for which $\fs{supp}(g)$ is bounded within $\fs{int}(I)$, so
that $B\rest I = \fs{Bdd}(G\rest I)$.

(3) $B'(I;A,P)$ denotes the derived group $B'$.
\end{rm}
\end{defn}

For example, $G([0,1]; {\Z}[\frac{1}{2}], \langle 2 \rangle)$ is one way
of defining the famous
Thompson Group $F$, an infinite 2-generator 2-relator group with
numerous important properties, discussed extensively in \cite{CFP}.
(${\Z}[\frac{1}{2}] \eqdf \setm{a/2^b}{a,b \in \Z}$.)  We say more
about $F$ below.

$G$ often fails to act transitively on $\fs{int}(I)$, but its orbits in
$\fs{int}(I)$ are dense.  By \cite[Corollary A4]{BS}, $G$ is highly
{\it o}-transitive on each of these orbits except in certain cases
involving $I = \R$.

For $G$ itself, the need for a notion less restrictive than {\em exact}
multiple transitivity arises only in the exceptional case.  However,
this need also arises for groups $H$ intermediate between $B'$ and
$G$.  The orbits in $\fs{int}(I)$ of such an $H$ are dense in
$\fs{int}(I)$.  From a more general result (to be discussed below)
involving 3-{\em interval}-transitivity, Bieri and Strebel
\cite[Theorem E4]{BS} deduce part (2) of the present
Theorem~\ref{LRT} for $\pair{\fs{int}(I_1)}{H_1}$ and
$\pair{\fs{int}(I_2)}{H_2}$, provided each $H_i$ contains a positive
element.  Here, besides deducing the entire First Order
Reconstruction package (see Theorem~\ref{LRT}), we find that the
positivity hypothesis can be discarded.

\begin{theorem}
\label{BSPL}\

{\rm (1)} Let $B'(I;A,P) \leq H \leq G(I;A,P)$.  Then
$\pair{\fs{int}(I)}{H}$ is highly interval-transitive and locally
moving, so that $\pair{\fs{int}(I)}{H} \in K^{\fss{LN2}}$.
Also, $\fs{EP}(\fs{int}(I),H) = A \cap \fs{int}(I)$.

{\rm (2)}  First Order Reconstruction holds for the class of all
$\pair{\fs{int}(I)}{H}$'s satisfying (1) (for various $I,A,P$).  In
particular, let $B'(I_i;A_i,P_i) \leq H_i \leq G(I_i;A_i,P_i)$ for
$i=1,2$.
Then:

$\,\,\,\,\,\,\,\,\,\,${\rm (a)} If $H_1 \equiv H_2$ then
\[
\fs{ACT}(A_1 \cap \fs{int}(I_1)),H_1) \equiv
\fs{ACT}(A_2 \cap \fs{int}(I_2)),H_2).
\]

$\,\,\,\,\,\,\,\,\,\,${\rm (b)} Suppose that $\iso{\alpha}{H_1}{H_2}$.
Then there exists a unique monotonic bijection
$\fnn{\tau}{\fs{int}(I_1)}{\fs{int}(I_2)}$ such that
\[
\alpha(h)\rest \fs{int}(I_2) = \tau \circ h\rest \fs{int}(I_1)  \circ
\tau\inverse \mbox{ for all } h \in H_1.
\]
Moreover,
\[
\iso{\tau\rest (A_1 \cap \fs{int}(I_1)) \cup \alpha}
{\fs{ACT}(A_1 \cap \fs{int}(I_1))}
{\fs{ACT}(A_2 \cap \fs{int}(I_2))}.
\]

\end{theorem}

\smallskip

The proof will be given below.

\medskip

Bieri and Strebel's more general result \cite[Theorem E3]{BS} says
approximately that part (2) of Theorem~\ref{LRT} holds for any
3-interval-transitive $\pair{L_1}{G_1}$ and $\pair{L_2}{G_2}$ which
contain nonidentity bounded elements and for which $L_i$ is an
open interval of $\R$.  This is a special case of Theorem~\ref{LRT},
which actually assumes only 2-interval-transitivity rather than
3-interval-transitivity.  We next discuss the ways in which their
result differs from the approximation just stated.

First, their form of 3-interval-transitivity is slightly different from
ours but is easily seen to be equivalent in the presence of a
nonidentity bounded element (since then our version implies high
interval-transitivity by Theorem~\ref{3H}).

The requirement that $L_1$ and $L_2$ be subintervals of $\R$ is not
essential for their development.  Nor, given the other hypotheses, is
an extra hypothesis (not assumed in our development) that
\[
\fs{sup}(\fs{supp}(g_1)) = \fs{inf}(\fs{supp}(g_2)) \mbox{ for some }
g_1,g_2 \in G.
\]
(Given $g_1 \in \fs{Lft}(G) - \sngltn{\fs{Id}}$, instead of a single
$g_2$ such that
$\fs{inf}(\fs{supp}(g_2)) =
\fs{sup}(\fs{supp}(g_1))$
one can use a {\em set} of $g_2$'s for which
$\fs{inf}\{\fs{inf}(\fs{supp}(g_2))\} =
\fs{sup}(\fs{supp}(g_1))$.)

Bieri and Strebel also assume the existence of a positive element,
which with present methods is unnecessary.

Finally, they actually assume about $\pair{L_1}{G_1}$ not the
existence of a nonidentity bounded element, but instead the ``extra''
hypotheses mentioned above.  However, the extra hypothesis,
together with 3-interval-transitivity, guarantees the existence of a
nonidentity bounded element; see Corollary~\ref{1WAYSUFF}.

The upshot of all this is that Bieri and Strebel's Theorem E3
(sharpened a bit) is a special case of Theorem~\ref{LNRT}.

\bigskip
\noindent
{\bf Proof of Theorem~\ref{GIAP} }
$\pair{\fs{Int}(I)}{H}$ satisfies the hypotheses of
Theorem~\ref{LNRT} just as in the proof of \cite[Theorem E4]{BS}.
Given Theorem~\ref{LRT} and its note, the rest is clear.
\hfill\qed
\bigskip

The present paper was inspired by a question asked us by Matt Brin
during his writing of \cite{Br}, at a time when he had already applied
results in
\cite{McHOMEOS} to several linear permutation groups acting on
$\R$.  These included:

(1) $\fs{PL}_2(\R)$, which is defined exactly as is
$G([0,1]; {\Z}[\frac{1}{2}], \langle 2 \rangle)$ but with the finitary
condition (c) omitted.

(2) $\widetilde{\fs{PL}}_2(\R)$, which is defined exactly as is
$\fs{PL}_2(\R)$ but with condition (b) relaxed to require only that
$g$ be monotonic and that the {\em absolute values} of its slopes be
elements of $P$.  $\fs{PL}_2(\R)$ has index 2 in
$\widetilde{\fs{PL}}_2(\R)$.

(3) $\fs{BPL}_2({\R}) \eqdf \fs{Bdd}(\fs{PL}_2(\R))$.

(4) The Thompson group $F$, now viewed as the group of elements
of $G({\R}; {\Z}[\frac{1}{2}], \langle 2 \rangle)$ such that there exist
$i,j \in \Z$ for which $g(x)$ is translation by $i$ for
$x$ near $+\infty$ and by $j$ for $x$ near $-\infty$.

\smallskip

Thus
$\fs{BPL}_2({\R}) \leq F \leq \fs{PL}_2({\R}) \leq
\widetilde{\fs{PL}}_2(\R)$, and all elements of $F$ are finitary.

Brin was also dealing with another of Thompson's groups $T$, a
closely related supergroup of $F$, whose importance stems partly
from its being an infinite simple group
having a finite presentation.  $T$ acts on the circle $\C = \R / \Z$,
and is the group of all PL-homeomorphisms $g:\C \rightarrow \C$
satisfying conditions (b) through(e) in the definition of
$G(\R; {\Z}[\frac{1}{2}], \langle 2 \rangle)$ but with ``increasing''
transmuted to ``orientation-preserving'' in (b).  Of course (c) is
automatic on the circle.

Finally, $\widetilde{T}$ is defined from $T$ as
$\widetilde{\fs{PL}}_2(\R)$ was from
$\fs{PL}_2(\R)$, and $T$ has index 2 in $\widetilde{T}$.

Brin asked us whether any of the general results in \cite{McHOMEOS}
for $\R$ had known analogues for $\C$, with the aim of applying
such results to $T$ and $\widetilde{T}$.  The answer ``not yet''
eventually
became ``yes''.  In \cite{Br} Brin quoted Theorem~\ref{LRT} (and its
circular analogue) from our present paper in the first step of his
analyses of $\fs{Aut}(T)$ and $\fs{Aut}(\widetilde{T})$.  (This was in
lieu of using the proofs he had already obtained for these two special
cases.)
What was essential about the improvement from \cite{McHOMEOS} to
the present paper was not the relaxation to approximate multiple
transitivity (each of the above groups is actually highly {\it o-
}transitive on a
dense orbit in {\R} or {\C}), but the ability to deal with groups acting
on circles.

In Brin's main theorem \cite[Theorem 1]{Br} he deals with groups $G$
for which
$\fs{BPL}_2({\R}) \leq G \leq \widetilde{\fs{PL}}_2(\R)$, and with
$G=T$ and $G = \widetilde{T}$.  He establishes that $\fs{Aut}(G)$ is
``small'' and that the conjugators $\tau_\alpha$ of
Corollary~\ref{1GR} are
unexotic.  (``Exotic'' means not piecewise linear.)

In the linear case $\fs{Aut}(G) \leq
\fs{Norm}_{\widetilde{\fss{PL}}_2(\R)}(G)$,
and for several special cases there is a sharper conclusion:

(1) $\fs{Aut}(\fs{BPL}_2({\R})) = \fs{Aut}(\fs{PL}_2({\R})) =
\widetilde{\fs{PL}}_2(\R)$, and $\widetilde{\fs{PL}}_2(\R)$ is a
complete group.

(2)  In the case of $F$, there is a short exact sequence
\[
1 \rightarrow F \rightarrow
\fs{Opp}(\fs{Norm}_{\widetilde{\fss{PL}}_2(\R)}(G)) \rightarrow T
\times T \rightarrow 1.
\]

In the circular case, $\fs{Aut}(T) = \widetilde{T}$ and
$\widetilde{T}$ is a complete group.

\smallskip

Later, in \cite{BG}, Brin and F. Guzman got sharply different results
for certain other generalized Thompson groups.  These groups $G$ act
on $\R$ or ${\R} / r\Z$.  Again the first step in analyzing
$\fs{Aut}(G)$ involves quoting our Theorem~\ref{LRT} (or its circular
analogue) to conclude that
$$
\fs{Aut}(G) = \fs{Norm}_{\fss{Homeo}({\R})}(G)  \mbox{ \, or \, }
\fs{Aut}(G) = \fs{Norm}_{\fss{Homeo}({\C})}(G).
$$
These other automorphism groups turn out to be ``large'' and {\em
do} have exotic elements.

\newpage

%

\section{Interpretations with parameters}
\label{PARAMETERS}

M. Giraudet sharpened \cite{McHOMEOS} by recasting various of its
results in terms of interpretations with parameters.  Here we follow
in her footsteps, but with the weaker hypotheses of the present
paper.  This procedure is built on

\begin{lemma}
\label{LEMPARAM}
Let $K$ be a class of locally moving linear permutation groups for
which First Order Reconstruction holds.  Then there exists a first
order formula $\psi_{Sgn}(p,x,y) \in \CL^{GR}$ such that for every
$\pair{L}{G} \in K$ and all $p,f,g \in G - \fs{Id}$:
\[
G \models \psi_{Sgn}[p,f,g] \mbox{ iff } \left\{
\begin{array}{ll}
p \in G^+ \mbox{ and } \fs{supp}(f) < \fs{supp}(g), \mbox{ or else} \\*
p \in G^- \mbox{ and } \fs{supp}(f) > \fs{supp}(g).\end{array}
\right.
\]
\end{lemma}

\medskip
\noindent
{\bf Proof }
The condition on nonsingleton $p,f,g$ is equivalent to:
$$
\fs{Ed}(x,p(x);y,z)  \mbox{ for all } x \in \fs{supp}(p), y \in
\fs{supp}(f), z \in \fs{supp}(g).
$$
Since the condition can be expressed first order in the language
$\CL^{\fss{MNPG}}$ of monotonic permutation groups, it can also be
expressed first order in $\CL^{\fss{GR}}$ by Corollary~\ref{FOP}.
\medskip
\hfill\qed

\begin{defn}\label{DEFEPL}
\begin{rm}
Let $\pair{L}{G} \in K^{\fss{LN}}$.

(1) Let
\begin{eqnarray*}
& \fs{Epl}(L,G) \eqdf
\setm{\fs{sup}(\fs{supp}(h))}{\fs{Id} \neq h \in \fs{Lft}(G)}
\mbox{, and} \\
& \fs{EPL}(L,G) \eqdf \pair{\fs{Epl}(L,G)}{\fs{Ed}^{\barL}}.
\end{eqnarray*}

(2) $\fs{Epr}(L,G)$ and $\fs{EPR}(L,G)$ are defined dually.

(3) The {\em pointwise order} on $G$, given by $f \leq g$ iff $f(x)
\leq g(x)$ for all $x \in L$ (or equivalently, for all $x \in \barL$), will
be denoted by ``$\leq$'', and its reverse by ``$\leq ^*$''.  For any dense
$G$-invariant subset $M$ of $\pair{\barL}{G}$,
the pointwise order for $\pair{M}{G}$ coincides with that for
$\pair{L}{G}$.  $\pair{G}{\leq}$ is a partially ordered group.

(4) $\fs{ACT}(\fs{EPL}(L,G),G,\leq)$ will denote
$\fs{ACT}(\fs{EPL}(L,G),G)$ with $G$ equipped with the (unambiguous)
pointwise order.
\end{rm}
\vspace{-2.5mm}
\end{defn}

\begin{defn}\label{}
\begin{rm}
Let $K$ be a class of locally moving linear permutation groups for
which First Order Reconstruction holds.  Then:

(1) Let $K_{\fss{EPLPWACT}} \eqdf
\setm{\fs{ACT}(\fs{EPL}(L,G),G,\leq)}
{\pair{L}{G} \in K}$.

(2) Let $K_{\fss{EPRPWACT}}$ be its dual.

(3) Let $K_{\fss{OGR}} \eqdf \setm{\pair{G}{\leq}}
{\mbox{there exists } L \mbox{ such
that } \pair{L}{G}
\in K}$.

(4) Let $\CL^{\fss{OGR}}$ denote the language of partially ordered
groups.
\end{rm}
\end{defn}

In the following formal definition, the class $K^*$ of
Definition~\ref{INTERP} is written instead as $K^{\times}$ so as to
avoid confusion with other meanings of the former.

\begin{defn}
\label{DEFPARAM}
\begin{rm}
Let $\trpl{K}{K^{\times}}{\CR}$ be a subuniverse system.
An {\em FO-STR-interpretation of $K^{\times}$ in
$K$ relative to $\CR$ with one parameter} consists of
a first order formula $\psi_{Par}(z)$ in $\CL(K)$ together with
an FO-STR-interpretation
$\trpl{\gf_{\fss{Imap}}}{\gf_U, \gf_{Eq}}{\ldots}$ of $K^{\times}$ in
$K^{\psi_{Par}}$ relative to $\CR^{\psi_{Par}}$, where
\begin{eqnarray*}
K^{\psi_{Par}} & \eqdf &
\setm{(M,a)}{M \in K \mbox{ and } M \models \psi_{Par}[a]},
\mbox{ and} \\*
\CR^{\psi_{Par}} & \eqdf &
\setm{\pair{(M,a)}{M^{\star}}}{\pair{M}{M^{\star}} \in \CR
\mbox{ and } (M,a) \in K^{\psi_{Par}}}.
\end{eqnarray*}
\end{rm}
\end{defn}

\begin{theorem} {\rm (Linear Interpretation Theorem with One
Parameter)}
\label{LITP}
Let $K$ be a class of locally moving linear permutation groups for
which First Order Reconstruction holds.  Then there is an
FO-STR-interpretation with one parameter $p$ of $K$ in $K_{GR}$,
with $\psi_{Par} = \psi_{Sgn}$,
which for each $\pair{L}{G} \in K$ interprets one of the following
structures:
\[
\left\{
\begin{array}{ll}
\fs{ACT}(\fs{EPL}(L,G),G,\leq) \mbox{ if $p$ is interpreted by some
} a \in G^+,\\*
\fs{ACT}(\fs{EPR}^*(L,G),G,\leq^*) \mbox{ if $p$ is interpreted by
some } a \in G^-.\end{array}
\right.
\]
\end{theorem}

\smallskip

\noindent
{\bf Proof }
The proof is Giraudet's proof with minor changes.  In order to
improve the theorem slightly by
interpreting more orbits, we change her formula $D \eqdf
\psi_{Bdd}$
that defines $\fs{Bdd}(G)$ to the formula $D \eqdf \psi_{LR}(g) \eqdf
\exists f(\psi_{Sgn}(p,g,f))$ that
defines either $\fs{Lft}(G) - \sngltn{\fs{Id}}$ or $\fs{Rt}(G) -
\sngltn{\fs{Id}}$, depending on the interpretation of the parameter.
Also, we replace the formula $\CO$ she used for defining the
order
by
$$
\CO(g,h) \eqqdf \forall k \mbox{\large (}\psi_{Sgn}(p,h,k)
\rightarrow
\psi_{Sgn}(p,g,k)\mbox{\large )} \wedge \,\exists f(\psi_{Sgn}(p,h,f)).
$$
$\fs{Lft}(G)$ is then definable with the parameter, and so since
$\fs{EP}(L,G)$ is definable in $\CL^{GR}$ without a parameter,
$\fs{EPL}(L,G)$ is definable {\em with} the parameter; and similarly
for
$\fs{Rt}(G)$ and $\fs{EPR}(L,G)$.  The rest of the proof is as in
\cite{Gi}.
\break
\rule{1pt}{0pt}\hfill\qed
\medskip

Giraudet's Corollaries 1--1 through 1--5 follow as in \cite{Gi}, but
now with weaker hypotheses.  Two of those corollaries have already
been mentioned here as Corollaries \ref{FOP}(6) and \ref{BDDNOTEE}.
We gather the others below in Corollary~\ref{CORPARAM}.
(The hypothesis ``$\fs{Lft}(g) \not\subseteq
\fs{Rt}(G)$'' in (3b) was inadvertently stated in \cite{Gi} as ``$G
\not\subseteq \fs{Rt}(G)$''.)

\begin{cor}
\label{CORPARAM}
Let $K$ be a class of locally moving linear permutation groups for
which First Order Reconstruction holds.  Then:

{\rm (1)}  $K_{\fss{EPLPWACT}}$ is FO-STR-interpretable in
$K_{OGR}$
without a parameter, as is $K_{\fss{EPRPWACT}}$.

{\rm (2)}  For any first order formula $\varphi(x_1,\dots,x_n)$ in
$\CL^{OGR}$
whose free variables are all group variables, there exist first order
formulas
$\psi_{\varphi}^1(x_1,\dots,x_n)$ and $\psi_{\varphi}^1(x_1,\dots,x_n)$
in $\CL^{GR}$ such that for every $\pair{L}{G} \in K$ containing a
positive element, and every $\vec{g} \in |G|^n$:
\begin{eqnarray*}
& G \models \psi_\varphi^1[\vec{g}] \mbox{ \,iff\, }
\mbox{\large (}M \models \varphi[\vec{g}] \mbox{ or }
M^* \models \varphi[\vec{g}]\mbox{\large )}, \mbox{ and}    \\*
& G \models \psi_\varphi^2[\vec{g}] \mbox{ \,iff\, }
\mbox{\large (}M \models \varphi[\vec{g}] \mbox{ and }
M^* \models \varphi[\vec{g}]\mbox{\large )ß}, \mbox{ \,\,\, }
\end{eqnarray*}
where $M = \fs{ACT}(\fs{EPL}(L,G),G,\leq)$
and $M^* = \fs{ACT}(\fs{EPR}^*(L,G),G,\leq^*)$.

{\rm (3)(a)} For $\pair{L}{G} \in K$ containing a positive element,
the
following are equivalent:

$\,\,\,\,\,\,\,\,\,\,\,\,\,\,$ {\rm (i)} There is a first order formula
$\varphi(x)$ in $\CL^{GR}$ (without a parameter) such that for all $g \in
G$:
$$
G \models \varphi[g] \mbox{ \, iff \,} g \in G^+.
$$

$\,\,\,\,\,\,\,\,\,\,\,$ {\rm (ii)} $(G,\leq) \not\equiv (G,\leq^*)$.

$\,\,\,\,\,\,\,\,\,$ {\rm (iii)}
$\fs{ACT}(\fs{EPL}(L,G),G,\leq) \not\equiv
\fs{ACT}(\fs{EPR}^*(L,G),G,\leq^*)$.

$\,\,\,\,\,\,${\rm (b)} There is a first order formula $\varphi_+(x)$
in
$\CL^{GR}$ (without a parameter) such that for every $\pair{L}{G}
\in
K$ for which $\fs{Lft}(G) \not\subseteq \fs{Rt}(G)$, and for every $g
\in \fs{Lft}(G)$:
$$
G \models \varphi_+[g] \mbox{ \, iff \,} g \in G^+.
$$
A similar statement holds with \fs{Lft} and \fs{Rt} interchanged.

{\rm (4)} For $\pair{L_1}{G_1},\pair{L_2}{G_2} \in K$, the following
are equivalent:

$\,\,\,\,\,\,\,\,\,\,${\rm (a)} $G_1 \equiv G_2$.

$\,\,\,\,\,\,\,\,\,\,${\rm (b)} $(G_1,\leq_1) \equiv (G_2,\leq_2)$ or
$(G_1,\leq_1) \equiv (G_2,\leq_2^*)$.

$\,\,\,\,\,\,\,\,\,\,${\rm (c)} $\fs{ACT}(\fs{EPL}(L_1,G_1),G_1,\leq_1)
\equiv \fs{ACT}(\fs{EPL}(L_2,G_2),G_2,\leq_2)$ or
\newline
$\,\,\,\,\,\,\,\fs{ACT}(\fs{EPL}(L_1,G_1),G_1,\leq_1) \equiv
\fs{ACT}(\fs{EPR}^*(L_2,G_2),G_2,\leq_2^*)$.
\end{cor}

There is a more general second order analogue of
Theorem~\ref{LITP} which applies to classes $K$ of linear
permutation groups for which Second Order Reconstruction holds
(see Definition~\ref{LFOR}).  In the lemma, $\psi_{Sgn}(p,x,y)$ is
second order, and the
theorem gives a SO-STR-interpretation with one parameter $p$.
$K_{\fss{EPLPWACT}}$ is SO-STR-interpretable in $K_{\fss{OGR}}$
without a parameter.
\newpage

\section{Lattice-ordered permutation groups}
\label{LATTICE}

Recall from the historical remarks in the introduction that when $L$
is a chain, $\pair{\fs{Aut}(L)}{\leq}$ is a lattice-ordered
group ({\em {\it l-}group}), where ``$\leq$'' is the pointwise order,
discussed in Definition~\ref{DEFEPL}.  The
lattice operations ``$\vee$'' and ``$\wedge$'' are pointwise,
that is, $(f \vee g)(x) = \fs{max}\{ f(x),g(x) \}$ and dually.

When $\pair{L}{G}$ is a linear permutation group $\trpl{L}{G}{\leq}$
is called an {\em
ordered permutation group}.  When $\pair{G}{\leq}$ is an
{\it l-}subgroup (simultaneously a subgroup and a sublattice) of
$\pair{\fs{Aut}(L)}{\leq}$, that is, when $G$ closed under pointwise
suprema and infima,
$\trpl{L}{G}{\leq}$  is called a {\em lattice-ordered permutation
group}.  Suprema and infima are then pointwise at all $a \in
\barL$,
so that $\trpl{\barL}{G}{\leq}$  is also a lattice-ordered permutation
group,
as is $\trpl{M}{G}{\leq}$  for any dense $G$-invariant subset $M$ of
$\pair{\barL}{G}$.

In this notation the note after Corollary~\ref{2GRORBITS} says for
$\pair{L_i}{G_i} \in K^{\fss{LN2}}$ that if
$\iso{\alpha}{\pair{G_1}{\leq_1}}
{\pair{G_2}{\leq_2}}$ then there exists a unique
$\iso{\tau}{\barL_1}{\barL_2}$ such that
$\iso{\tau \cup \alpha}{\trpl{L_1}{G_1}{\leq_1}}
{\trpl{L_2}{G_2}{\leq_2}}$.

There is an extensive literature on lattice-ordered permutation
groups (see \cite{G1} or \cite{McSURVEY}, for example).  However,
the present paper makes no contribution to this study because here
the weakening of hypotheses to 2-interval-transitivity offers no
extra generality---see the next proposition.
Local movability is not assumed here.  For a followup to this result
see Proposition~\ref{AUTLC}.

\begin{prop}
\label{NOHELP}\
Let $\trpl{L}{G}{\leq}$ be a lattice-ordered
permutation group.
If $\pair{L}{G}$ is 2-interval-transitive (or even
inclusion-transitive), then all orbits of $G$ in $\pair{\barL}{G}$ are
dense in $\barL$, and $G$ is highly {\it o-}transitive on each
of them.
\end{prop}

Strictly speaking, Proposition~\ref{NOHELP} is within the scope of
this paper only when $L$ is a dense chain, but see the remark after
Proposition~\ref{DEFINCLTRANS}.

\smallskip
\noindent
{\bf Proof }
Inclusion-transitivity guarantees density of orbits.  Now let $M$ be
an orbit of $\pair{\barL}{G}$.  The proof is by induction.
For the induction step, it suffices to show that for
$$
x_1 < \dots < x_{n-1} < x_n < y_n \mbox{ in $M$,}
$$
there exists $f \in G$ such that $f(x_i) = x_i$ for $i=1,\dots,n-1$, and
such that $f(x_n) = y_n$.  Pick $g \in G$ such that $g(x_n) = y_n$.
Replacing $g$ by $g \vee \fs{Id}$, we may assume that $g(x_i) \geq
x_i$ for $i=1,\dots,n-1$, as well as that $g(x_n) = y_n$.
Next pick bounded $L$-intervals $I$ and $J$ such that
$$
\fs{inf}(I) < x_1 < x_{n-1} < J < x_n < y_n < \fs{sup}(I).
$$
Now pick $h \in G$ such that $h(I) \subseteq J$.  Then $h^{-1}(J)
\supseteq I$, ensuring that $h^{-1}(x_i) < x_i$ for $i=1,\dots,n-1$ and
that $h^{-1}(x_n) > y_n$.  Then $f \eqdf g \wedge (h^{-1} \vee
\fs{Id})$ is as required.  Verification is left to the reader.
\hfill\qed
\medskip

We recall some standard terminology for linear permutation
groups, and for later use we apply it to all four types of nearly
ordered permutation groups.  First, by a {\em convex} subset of a
circle $C$ we
mean a subset $P \subseteq C$ which is (the intersection with $C$ of)
a $\barC$-interval (with endpoints of $P$ possibly adjoined if they
lie in $C$); singletons and $\emptyset$ are also considered convex.
A {\em convex} subset of an equal direction structure is simply a
convex subset of the underlying chain $L$ (or equally well of $L^*$),
and similarly for equal orientation structures.

For convenience we deal here just with the transitive case, which
guarantees that for any $\pair{N}{G}$-congruence, the congruence
classes are all images of each other under elements of $G$.

\begin{defn} \label{DEFOPRIM}
\begin{rm}
Let $\pair{N}{G}$ be a transitive nearly ordered permutation group.

(1) An {\em o-block} of $\pair{N}{G}$ is a convex block, that is, a
nonempty convex $P \subseteq N$ such that for all $g \in G$, $g(P) =
P$ or $g(P) \cap P = \emptyset$; and is {\em proper} if it is not a
singleton or $N$ itself.  Proper {\it o-}blocks cannot have endpoints
in $N$.

(2) An $\pair{N}{G}$-congruence is {\em convex} if all its congruence
classes are convex; and is {\em proper} if its classes are proper.

(3) $\pair{N}{G}$ is {\em o-primitive} if it has no proper
{\it o-}blocks, or equivalently, no proper convex congruences.
\end{rm}
\end{defn}

T. J. Scott showed in \cite{Sc1} that for a monotonic permutation
group $\pair{\fs{ED}(L)}{G}$, {\it o-}blocks of $\pair{L}{\fs{Opp}(G)}$
are {\it o-}blocks even of $\pair{\fs{ED}(L)}{G}$ and thus
{\it o-}primitivity of the former is equivalent to {\it o-}primitivity of
the latter, provided $\pair{L}{\fs{Opp}(G)}$ is a transitive
lattice-ordered permutation group.  The lattice-ordered proviso can
be weakened to require just that $\pair{L}{\fs{Opp}(G)}$ be {\em
coherent}, that is, for all $a<b \in L$ there exists $g \in \fs{Opp}(G)$
such that $g(a) = b$ and
$g(c) \geq c$ for all $c \in \barL$.  Coherence implies its dual.

For lattice-ordered permutation groups there is a detailed
description of the transitive {\it o-}primitive case, and it relates
nicely to the themes of this paper.  We state this description and
then define the key notion of periodic {\it o-}primitivity.

\begin{theorem} {\rm ({\it o-}Primitive Classification Theorem,
McCleary \cite{McPRIM})}
\label{CLASSIF}\
Every {\it o-}primitive transitive
lattice-ordered permutation group is of one of the following three
mutually exclusive types, and all transitive
lattice-ordered permutation groups of these types are
{\it o-}primitive:

{\rm (1)} 2-{\it o-}transitive (hence highly {\it o-}transitive).

{\rm (2)} The regular representation of a subgroup of the additive
reals.

{\rm (3)} Periodically {\it o-}primitive.

\end{theorem}

The periodically {\it o-}primitive case is best understood in terms of
the canonical example $\pair{{\R}}{G,\leq}$, where
$$
G = \setm{g \in \fs{Aut}({\R})}{g(a + 1) = g(a) + 1 \mbox{ for all } a
\in \R}.
$$
The {\em period}  here is the function $a \mapsto a+1$, and the
effect of any $g$ on any interval $[a,a+1)$ determines its effect ``one
period up'' on $[a+1,a+2)$, and similarly throughout all of $L$.

Now we consider a slightly more general kind of period $z$, lying in
$\fs{Aut}(\barL)$ but not necessarily in $\fs{Aut}(L)$.  A {\em
period} of a transitive lattice-ordered permutation group
$\pair{L}{\fs{Aut}(L)_{pw}}$ is an
element $z \in \fs{Aut}(\barL)^+$, one (hence each) of whose orbits
$\setm{z^n(a)}{n \in \Z}$ is coterminal in $\barL$, and which
generates (as a group) the centralizer $Z_{\fss{Aut}(\barL)}(G)$.
An {\it o-}primitive transitive lattice-ordered permutation group is
called {\em periodically {\it o-}primitive} if it has a period $z$
(except that isomorphic copies of the regular representation
$({\Z},{\Z})$ are not counted as
periodically {\it o-}primitive).

\smallskip

The Classification Theorem can be generalized to transitive coherent
ordered permutation groups.

Observe that in the regular and periodic cases of
Theorem~\ref{CLASSIF} there is no nonidentity bounded element.
Here is an alternate proof of Proposition~\ref{NOHELP} based on the
classification:

\smallskip
\noindent
{\bf Proof }
$\pair{L}{G,\leq}$ must be {\it o-}primitive, for if $P$ were a proper
{\it o-}block no $g \in G$ could map $P$ to a proper subinterval of
itself, violating inclusion-transitivity.

Now apply the Classification Theorem.  If $\pair{L}{G,\leq}$ were
periodic, no $g$ could map an interval $I= (a,z(a))$ to a proper
subinterval of itself since if $g(a) \in I$, then $g(z(a)) = z(g(a)) > z(a) >
I$.  A similar argument works for the regular representation of a
dense subgroup of the additive reals, and the discrete case was
disposed of just after the statement of the proposition.  Therefore
$\pair{L}{G,\leq}$ is highly {\it o-}transitive.  This argument works
equally well for other orbits $M$ of $\pair{\barL}{G}$.
\hfill\qed
\medskip

When $L$ is Dedekind complete and $\pair{L}{G}$ has a nonidentity
bounded element, ordinary transitivity is enough:

\begin{prop} {\rm (Holland \cite{Ho1})}
\label{}\
Let $\pair{L}{G,\leq}$ be a transitive lattice-ordered permutation
group.  If $L$ is Dedekind complete and $\pair{L}{G}$ contains a
nonidentity bounded element, then $\pair{L}{G}$ is highly
{\it o-}transitive.
\end{prop}

\noindent
{\bf Proof }
Proper {\it o-}blocks cannot have endpoints in $L$, so Dedekind
completeness guarantees {\it o-}primitivity.  In the regular and
periodic cases of Theorem~\ref{CLASSIF}, $\pair{L}{G}$ contains no
nonidentity bounded elements.
\hfill\qed
\newpage

%

\section{Transitivity properties}
\label{TRANSPROPS}

In this section we discuss relationships among various transitivity
properties of nearly ordered permutation groups $\pair{N}{G}$.  This
is basically a discussion of {\em small} degrees of multiple
transitivity (exact, approximate, and interval) because in the
circumstances of interest here, {\it n-}transitivity for relatively
modest $n$ turns out to imply high transitivity.  In all results in this
section, the transitivity hypotheses will be strong enough to
guarantee that each orbit of $\pair{\barN}{\fs{Opp}(G)}$ is dense in
$\barN$
(see Proposition~\ref {IT}).

For high transitivity, the approximate and interval forms coincide,
and carry over from $\pair{N}{G}$ to $\pair{\barN}{G}$:

\begin{prop}
\label{HAT}
If a nearly ordered permutation group $\pair{N}{G}$ is highly
interval-transitive, then $\pair{\barN}{G}$ is highly approximately
{\it o-}transitive.
\end{prop}

\noindent
{\bf Proof }
Let $\pair{L}{G}$ be a highly interval-transitive linear permutation
group.  Let $a_1 < \cdots < a_n$ in $\barL$, and let $J_1 < \cdots <
J_n$
be $\barL$-intervals.  For $i=1,\dots,n$, pick $L$-intervals $I_1' <
I_1''$ such that
\[
I_1' < a_1 < I_1'' < I_2' < a_2 < I_2'' < \cdots < I_n' < a_n < I_n'',
\]
and pick $L$-intervals $J_i' < J_i''$ within $J_i$.
Now pick $g \in G$ such that for $i=1,\dots,n$,
$g(I_i') \cap J_i' \neq \emptyset$ and $g(I_i'') \cap J_i'' \neq
\emptyset$, and thus $g(a_i) \in J_i$.

The proofs for the other three types are similar.
\hfill\qed
\medskip

The main results of this section will now be stated, and then we will
turn to their proofs.  Our primary interest is in approximate and
interval transitivity.  However, most of the results have analogues
for exact transitivity, and these will be collected in
Theorem~\ref{EXACT}.

What does it take to force high approximate {\it o-}transitivity
of $\pair{N}{G}$, or equivalently, of $\pair{\barN}{G}$?

\begin{theorem}
\label{3H}
Let $\pair{N}{G}$ be a nearly ordered permutation group having a
nonidentity bounded element.

{\rm (1)} In the linear/monotonic case:

$\,\,\,\,\,\,\,\,\,\,${\rm (a)} 3-interval-transitivity $\Rightarrow$
high approximate {\it o-}transitivity.

$\,\,\,\,\,\,\,\,\,\,${\rm (b)} When $N$ is Dedekind complete,
approximate 2-{\it o-}transitivity $\Rightarrow$
high approximate {\it o-}transitivity.

{\rm (2)} In the circular/monocircular case:

$\,\,\,\,\,\,\,\,\,\,${\rm (a)} 4-interval-transitivity $\Rightarrow$
high approximate {\it o-}transitivity.

$\,\,\,\,\,\,\,\,\,\,${\rm (b)} When $N$ is Dedekind complete,
approximate 3-{\it o-}transitivity $\Rightarrow$
high approximate {\it o-}transitivity.
\end{theorem}

Note that when the assumed degree $n$ of transitivity is reduced by
one for Dedekind complete $N$, the hypothesis involves {\em
approximate} (rather than {\em interval}) transitivity.
At least in (1b), interval-transitivity would not
suffice here; consider $\pair{\barL}{G}$, where  $\pair{L}{G}$ is as in
Example~\ref{WR1VARIANT}.
Also, at least in (1), the assumed degrees cannot be decreased
(Examples \ref{WR1VARIANT} and \ref{AUTRZ}, with the latter
acting on $\barL$). However, some relaxation of hypotheses is
possible; see Theorem~\ref{UPWARD}.

For lattice-ordered permutation groups, the hypothesis of
Theorem~\ref{3H} can be weakened; see Proposition~\ref{NOHELP}.

\begin{question}
\label{QHIGH}
\begin{rm}
In part (2) of Theorem~\ref{3H}:

(1) Can 4-interval-transitivity be reduced to 3-interval-transitivity?

(2) Can approximate 3-{\it o-}transitivity be reduced to approximate
2-{\it o-}transitivity?

(3) Can approximate 3-{\it o-}transitivity be relaxed to
3-interval-transitivity?
\end{rm}
\end{question}

It is a consequence of Theorem~\ref{3H} and Proposition~\ref{HAT}
that for a linear permutation group $\pair{L}{G}$ having a
nonidentity bounded element, the hypothesis of
approximate 3-{\it o-}transitivity is robust in the sense that
$\pair{M}{G}$ is also approximately 3-{\it o-}transitive for any
$G$-invariant subset $M$ of $\barL$.  Analogous remarks apply to
the other three types of nearly ordered permutation groups.

However, a linear permutation group $\pair{L}{G}$ having a
nonidentity bounded element can be (exactly) 2-{\it o-}transitive
and yet have an orbit $M$ in $\barL$ such that $\pair{M}{G}$ is not
even approximately 2-{\it o-}transitive
(Example~\ref{WR1}), and can be 3-{\it o-}transitive (hence
highly {\it o-}transitive, as we see in Theorem~\ref{EXACT}) and yet
fail to be even 2-{\it o-}transitive
on some orbit in $\barL$ (Example~\ref{BACKANDFORTH}).
The second of these phenomena occurs also for circular permutation
groups (Example~\ref{CIRCLE2}).

\begin{question}
\label{Q3OTIMP3ATORB}
\begin{rm}
Can a locally moving circular permutation group $\pair{C}{G}$ be
3-{\it o-}transitive and yet have an orbit $M$ in $\barC$ such that
$\pair{M}{G}$ is not even approximately 3-{\it o-}transitive?
\end{rm}
\end{question}

The next theorem explores consequences of high interval-transitivity
of $\pair{N}{G}$, which turns out to be inherited by nontrivial
subnormal subgroups $H \unlhd\!\unlhd\, G$.  Particularly
important here are $\fs{Bdd}(G)$, $\fs{Bdd}'(G)$, and in the
linear/monotonic case $\fs{Lft}(G)$ and $\fs{Rt}(G)$.  In the
linear/monotonic case,
$\fs{Bdd}'(G)$ is the smallest nontrivial subnormal subgroup of $G$
by Theorem~\ref{HIGMAN}.  In (2) and (3a) we consider data
confined to a bounded $N$-interval $K$.  For the
circular/monocircular case, we need to amplify Definition
\ref{DEFBDD}:

\begin{defn}
\label{NBDD}
\begin{rm}
Let $\pair{N}{G}$  be a nearly ordered permutation group, with $N =
C = \fs{CR}(L)$ or $N =  \fs{EO}(C)$.  Let
\begin{eqnarray*}
\fs{1Bdd}(G)  & \eqdf  & \setm{g \in G}{\fs{supp}(g) \subseteq I
\mbox{ for some bounded $C$-interval } I}, \mbox{and} \\*
\fs{nBdd}(G)  & \eqdf & \setm{g \in G}{g \mbox{ is the product of
$n$ elements of } \fs{1Bdd}(G)};
\end{eqnarray*}
these are normal subsets of $G$.
Recall that
\begin{eqnarray*}
\fs{Bdd}(G)  \eqdf  \sing{\fs{1Bdd}(G)},
\end{eqnarray*}
the (normal) subgroup of $G$ generated by $\fs{1Bdd}(G)$.
\end{rm}
\end{defn}

\begin{theorem}
\label{CONSEQ}
Let $\pair{N}{G}$ be a nearly ordered permutation group which is
highly interval-transitive.

{\rm (1)} Let $\{\fs{Id}\} \neq H \unlhd\!\unlhd\, G$.  Then

$\,\,\,\,\,\,\,\,\,\,${\rm (a)} $\pair{\barN}{H}$ is highly
approximately {\it o-}transitive.

$\,\,\,\,\,\,\,\,\,\,${\rm (b)} If $G$ has a bounded element, so does
$H$.

{\rm (2)} Suppose $\pair{N}{G}$ is linear or monotonic, and has a
nonidentity bounded element.  Let $a_1 < \ldots < a_n$ in $\barN$
and let $J_1 < \ldots < J_n$ be bounded $N$-intervals.  Let $K$ be any
bounded $N$-interval containing all the $a_i$'s and $J_i$'s.  Then
there exists $g \in G$ such that $g(a_i) \in \barJ_i$ for $i=1,\dots,n$
and such that $\fs{supp}(g) \subseteq K$.

{\rm (3)} Suppose $\pair{N}{G}$ is circular or monocircular, and has a
nonidentity bounded element.  Let $\fs{Cr}(a_1,\ldots,a_n)$ in
$\barN$, and let $J_1,\ldots,J_n$ be bounded $N$-intervals with
$\fs{Crs}(J_1,\ldots,J_n)$.

$\,\,\,\,\,\,\,\,\,\,${\rm (a)}  Let $K$ be any bounded $N$-interval
such that
\[
\fs{Cr}(\fs{cw}(K),a_1,\ldots,a_n,\fs{ccw}(K)) \mbox{ and }
\fs{Crs}(\fs{cw}(K),J_1,\ldots,J_n,\fs{ccw}(K)).
\]
Then there exists $g \in G$ such that $g(a_i) \in \barJ_i$ for
$i=1,\dots,n$ and such that $\fs{supp}(g) \subseteq K$.

$\,\,\,\,\,\,\,\,\,\,${\rm (b)}  There exists $g \in \fs{nBdd}(G)$
such that $g(a_i) \in \barJ_i$ for $i=1,\ldots,n$.
\end{theorem}

We state here two particular consequences of Theorems \ref{3H} and
\ref{CONSEQ} which are used in Section~\ref{RECON}:

\begin{cor}
\label{USED}
Suppose $\pair{N}{G}$ has a nonidentity bounded element.

{\rm (1)} In the linear/monotonic case,
if $\pair{N}{G}$ is 3-interval-transitive then
$\pair{\barN}{\fs{Bdd}(G)}$ is
highly approximately {\it o-}transitive.

{\rm (2)} In the circular/monocircular case,
if $\pair{N}{G}$ is 4-interval-transitive then
$\pair{\barN}{\fs{nBdd}(G)}$ is
approximately {\it n-o-}transitive.
\end{cor}

Each part of Theorems \ref{3H} and \ref{CONSEQ} and
Corollary~\ref{USED} has an analogue for {\em exact} multiple
transitivity.  For example, the analogue of \ref{USED}(1) states that in
the linear/monotonic case, if $\pair{N}{G}$ is 3-{\it o-}transitive and
has a nonidentity bounded element then $\pair{N}{\fs{Bdd}(G)}$ is
highly {\it o-}transitive.  (Of course we cannot say that
$\pair{\barN}{\fs{Bdd}(G)}$ is highly {\it o-}transitive, which is
impossible when $N \neq \barN$.)  This particular analogue was
already known \cite{McHOMEOS} (and it contains the analogue of
\ref{CONSEQ}(1a)), and the analogue of \ref{CONSEQ}(1b) follows
from Higman's Theorem \ref{HIGMAN}, but the other analogues may
be new.

\begin{theorem}
\label{EXACT}
For each part of Theorems \ref{3H} and \ref{CONSEQ} and
Corollary~\ref{USED}, the analogue for exact multiple transitivity is
valid.
\end{theorem}

However, a 2-{\it o-}transitive linear permutation group having a
nonidentity bounded element need not be 3-{\it o-}transitive
(Example~\ref{WR1}).

The last theorem blends interval-transitivity and exact transitivity.
In part (2) there is no exact transitivity {\em hypothesis} at all.

\begin{theorem}
\label{HOT}
Let $\pair{N}{G}$ be a nearly ordered permutation group which is
highly interval-transitive and has a nonidentity bounded element.

{\rm (1)} If $\pair{N}{G}$ is 2-{\it o-}transitive, then it is highly {\it
o-}transitive, and so is $\pair{N}{H}$ for any $\{\fs{Id}\} \neq H
\unlhd\!\unlhd\, G$.

{\rm (2)} In the linear/monotonic case, $\pair{\barN}{\fs{Bdd}(G)}$
is highly {\it o-}transitive on each of its own orbits.  Even better,
given any $a_1 < \cdots < a_n$ and $b_1 < \cdots < b_n$ in $\barN$ with
each $b_i$ in the same $\fs{Bdd}(G)$-orbit as $a_i$, there exists $g \in
\fs{Bdd}(G)$ such that $g(a_i) = b_i$ for $i = 1,\dots,n$.
\end{theorem}

\noindent
{\bf Proof of Theorem \ref{3H}}
(1a)
It suffices to treat the linear case since for monotonic
$\pair{\fs{ED}(L)}{G}$ the linear result can be applied to
$\pair{L}{\fs{Opp}(G)}$.  Let $\pair{L}{G}$ be a 3-interval-transitive
linear permutation group having a nonidentity bounded element.
Pick $p \in \fs{Bdd}(G)$ such that $p(d) < d$ for some $d \in
\fs{supp}(p)$.

By Proposition~\ref{HAT}, it suffices to show that $\pair{L}{G}$ is
highly interval-transitive.  Let $I_1 < \cdots < I_n$ and $J_1 < \cdots
< J_n$ be bounded $L$-intervals, with $n \geq 4$.  We need to show
there exists $g \in G$ such that $g(I_i) \cap J_i \neq \emptyset$ for
$i=1,\dots,n$.  Suppose by induction that there exists $h \in G$ such
that
$h(I_i) \cap J_i \neq \emptyset$ for $i=1,\dots,n-1$, and suppose
that $h(I_n) \cap J_n = \emptyset$ (else we are done).

{\bf Case 1:} $h(I_n) > J_n$.  We seek an $f \in G$ such that
$p^{f\inverse}(h(I_n)) \cap J_n \neq \emptyset$ and
$J_{n-1} < \fs{supp}(p^{f\inverse})$.
Then $g \eqdf p^{f\inverse}h$ will be as required.

Pick $L$-intervals $K_1 < K_2$ such that
\[
\fs{inf}(J_n) < K_1 < K_2 < \fs{sup}(J_n) < h(I_n).
\]
Next pick $L$-intervals $K_1' < K_2' < K_3'$ such that
\[
K_1' < \fs{inf}(\fs{supp}(p)) < p(d) < K_2' < K_3' < d.
\]
Now pick $f \in G$ such that
$f(K_1) \cap K_1' \neq \emptyset$,
$f(K_2) \cap K_2' \neq \emptyset$, and
$f(h(I_n)) \cap K_3' \neq \emptyset$.
Then $f\inverse pf(h(I_n)) \cap J_n \neq \emptyset$ since
\begin{eqnarray*}
&  \fs{inf}(f\inverse pf(h(I_n))) =
f\inverse pf(\fs{inf}(h(I_n))) <
f\inverse p(\fs{sup}(K_3')) \leq \\*
& f\inverse(p(d)) \leq
f\inverse(\fs{inf}(K_2')) <
\fs{sup}(K_2) \leq
\fs{sup}(J_n),
\end{eqnarray*}
\noindent
and since
\begin{eqnarray*}
& \fs{sup}(f\inverse pf(h(I_n))) =
f\inverse pf(\fs{sup}(h(I_n))) >
f\inverse p(\fs{inf}(K_3')) > \\*
& f\inverse(p(d)) >
f\inverse(\fs{sup}(K_1')) >
\fs{inf}(K_1)  \geq
\fs{inf}(J_n).
\end{eqnarray*}

\noindent
Moreover, $J_{n-1} < \fs{supp}(f\inverse pf)$ since
\[
J_{n-1} < \fs{inf}(K_1) < \fs{sup}(f\inverse (K_1')) =
f\inverse (\fs{sup}(K_1')) < f\inverse (\fs{supp}(p)).
\]

{\bf Case 2:} $h(I_n) < J_n$.  To reduce Case 2 to Case 1,
we seek an $f \in G$ such that
$(p\inverse)^{f\inverse}(h(I_n)) > J_n$ and
$h(I_{n-1}) < \fs{supp}((p\inverse)^{f\inverse})$.
Then $h_1 \eqdf (p\inverse)^{f\inverse} h$ will satisfy
Case 1.

Replacing $I_n$ by a subinterval of itself, we may assume that
$\fs{sup}(h(I_{n-1})) < \fs{inf}(h(I_n))$.
Pick $L$-intervals $K_1 < K_3$ such that
\[
h(I_{n-1}) < K_1 < h(I_n) < J_n < K_3.
\]
Pick $K_1' < K_2' < K_3'$ exactly as in Case 1.
Now pick $f \in G$ such that
$f(K_1) \cap K_1' \neq \emptyset$,
$f(h(I_n)) \cap K_2' \neq \emptyset$, and
$f(K_3) \cap K_3' \neq \emptyset$.  Then
\begin{eqnarray*}
&  \fs{sup}(f\inverse p\inverse f(h(I_n))) =
f\inverse p\inverse f(\fs{sup}(h(I_n))) >
f\inverse p\inverse(\fs{inf}(K_2')) \geq \\*
& f\inverse p\inverse p(d)) =
f\inverse(d) \geq
f\inverse(\fs{sup}(K_3')) >
\fs{inf}(K_3) \geq
\fs{sup}(J_n).
\end{eqnarray*}

\noindent
Again replacing $I_n$ by a subinterval of itself, we may assume that
\newline
$f\inverse p\inverse f(h(I_n)) > J_n$.  Moreover,
$h(I_{n-1}) < \fs{inf}(K_1) < \fs{supp}(f\inverse p\inverse f)$ since
$f(\fs{inf}(K_1)) < \fs{sup}(K_1') < \fs{supp}(p)$.

(2a) The proof is much like that of Case 1 of (1a).
Let $\pair{C}{G}$ be a 4-interval-transitive circular permutation
group having a nonidentity bounded element. Let $I_1,\dots,I_n$
and $J_1,\dots,J_n$ be $C$-intervals such that
$Crs(I_1,\dots,I_n)$ and $Crs(J_1,\dots,J_n)$,
with $n \geq 5$.  Suppose there exists $h \in G$ such that
$h(I_i) \cap J_i \neq \emptyset$ for $i=1,\dots,n-1$, but that $h(I_n)
\cap J_n = \emptyset$.
Replacing the $I_i$'s and $J_i$'s by subintervals of themselves, we
may assume that $h(I_i) = J_i$ for $i=1,\dots,n-1$, and that $h(I_n)$
has no endpoint in common with any $J_i$.  We also assume that
$Crs(J_n,h(I_n),J_1)$, the case with $Crs(J_{n-1},h(I_n),J_n)$ being
the dual.

Pick $C$-intervals $K_1,K_2,K_4$ such that
\[
Crs(\fs{cw}(J_n),K_1,K_2,\fs{ccw}(J_n),h(I_n),K_4,J_1).
\]
Pick $p \in G$ such that $\fs{supp}(p) \subseteq B$ for some
bounded $C$-interval $B$, and such that $Cr(p^2(d),p(d),d)$ for some
$d \in C$.
Pick $K_1', K_2', K_3', K_4'$ such that
\[
Crs(K_1', \fs{cw}(B), p(d), K_2', K_3', d,\fs{ccw}(B),K_4').
\]
Now pick $f \in G$ such that
$f(K_i) \cap K_i' \neq \emptyset$ for $i=1,2,4$,
and such that
$f(h(I_n)) \cap K_3' \neq \emptyset$.
As before,
$g \eqdf p^{f\inverse} h$ satisfies $g(I_i) \cap J_i \neq \emptyset$
for $i=1,\dots,n$.  Verification of this is left to the reader.

(1b) and (2b)  will be proved near the end of this section.
\hfill\qed
\medskip

In the foregoing proof of (1a), the assumption of a nonidentity
bounded element can be weakened a bit:

\begin{cor}
\label{1WAYSUFF}
Let $\pair{N}{G}$ be a linear or monotonic permutation group.

{\rm (1)} Suppose $\pair{N}{G}$ has a nonidentity element $p$
whose support is bounded above {\em or} below (that is, $p \in
\fs{Lft}(G) \cup
\fs{Rt}(G)$).  Then 3-interval-transitivity  $\Rightarrow$ high
approximate {\it o-}transitivity.

{\rm (2)} Suppose $\pair{N}{G}$ has nonidentity elements $p$ and
$q$ such that $p \in \fs{Lft}(G)$ and $q \in \fs{Rt}(G)$.  Then
$\pair{N}{G}$ has a nonidentity bounded element.
 \end{cor}

\noindent
{\bf Proof }
(1) When $p \in \fs{Rt}(G)$, the proof of (1) is identical to the proof
of (1a) in Theorem~\ref{3H}; the other case is the dual.

(2) Again it suffices to treat the linear case.  We seek $f \in G$ such
that for $q_2 \eqdf q^f$ we have $pq_2 \neq q_2p$.  Then for $q_3
\eqdf q_2^p$ we shall have $q_3 \equiv q_2$ outside some bounded
interval (since $q_2 \in \fs{Rt}(G)$ and $p \in \fs{Lft}(G)$), but $q_3
\neq q_2$.  Therefore $q_2q_3^{-1}$ will be the required nonidentity
bounded element.

To construct such an $f$, proceed as follows.  Replacing $p$ by
$p^{-1}$ if necessary, pick a bounded $L$-interval $J_1 \subseteq
\fs{supp}(p)$ such that $J_1 < p(J_1)$.  Let $a_1 =
\fs{inf}(\fs{supp}(q))$.  Pick $a_2 < a_3$ such that $a_1 < a_2 < a_3$
and such that $a_1$ and $a_2$ lie within the same
interval $I$ of $q$.  Use (1a) to pick $f \in G$ such that $f(a_1) \in
J_1$, and such that $f(a_2) < p(J_1) < f(a_3)$ so that $f(I) \supseteq
p(J_1)$.  Then $f(a_1)$ is fixed by $q^f$ (since $a_1$ is fixed by $q$)
but moved by $p$ (since $f(a_1) \in J_1$), and $p(f(a_1))$ is moved
by $q^f$ (since $f(I) \supseteq p(J_1))$.  Hence $f$ is as required.
\hfill\qed

\bigskip

We turn now to proving Theorem~\ref{CONSEQ}.  Given {\em high}
interval transitivity, we can be a bit more relaxed in our arguments.

\begin{lemma}
\label{}
Suppose $\pair{N}{G}$ is highly interval-transitive.
Then in parts (1c) and (2c) of Definition~\ref{DEFMULTTRANS},
provided that $I_i$'s are bounded and have no common endpoints in
$\barN$ and similarly for the $J_i$'s, it can be arranged in addition
that
\begin{eqnarray*}
g(I_i) \subseteq J_i \mbox{ for each } i \in \CS, \\*
g(I_i) \supseteq J_i \mbox{ for each } i\not\in \CS, \,
\end{eqnarray*}
where $\CS$ is any subset of $\{1,\ldots,n\}$.
\end{lemma}

\noindent
{\bf Proof }
For $\CS = \{1,\ldots,n\}$, the proof is similar to that of
Proposition~\ref{HAT} (with $I_i$ playing the role of $a_i$), and the
modifications for arbitrary $\CS$ are easy.
\hfill\qed
\medskip

For (1a) of \ref{CONSEQ} it suffices to show that for
$\{\fs{Id}\} \neq H \unlhd G$, $\pair{N}{H}$ is highly
interval-transitive.

We deal first with linear permutation groups.  Let $I_1 < \ldots <
I_n$ and $J_1 < \ldots < J_n$ be $L$-intervals.   Replacing these by
subintervals of themselves, we may assume they satisfy the extra
conditions of the lemma.  We will construct an $h \in H$ such that
$h(I_i) \subseteq J_i$ for $i = 1,\ldots,n$.

We first show that such an $h$ exists in the special case in which
$I_n < J_1$.  Pick $h_1 \in H$ such that $K < h_1(K)$ for some
$L$-interval $K$.  Within $K$ pick $L$-intervals
$I_1' < \ldots < I_n'$ satisfying the additional conditions of the
lemma.  Pick $f \in G$ such that for $i = 1,\ldots,n$, $f(I_i) \subseteq
I_i'$ and $f(J_i) \supseteq h_1(I_i')$.  Then $h \eqdf h_1^{f\inverse}$
is as required.

In general, let $K_1 < \ldots < K_n$ be $L$-intervals such that $I_n
\cup J_n < K_1$.  Apply the special case to obtain $h_2 \in H$ such
that $h_2(I_i) \subseteq K_i$ for $i = 1,\ldots,n$.  Apply a similar
special case to obtain $h_3 \in H$ such that $h_3(K_i) \subseteq J_i$
for $i = 1,\ldots,n$.  Then $h \eqdf h_3 h_2$ is as required.  This
concludes the linear case, and the monotonic case is similar.

The circular/monocircular case follows the same outline.  In the
special case it it assumed that there exist disjoint $C$-intervals $I$
and $J$ such that $\fs{Crs}(\fs{cw}(I),I_1,\ldots,I_n,\fs{ccw}(I))$ and
similarly for the $J_i$'s.  We introduce {\em two} sequences
$K_1^I,\ldots,K_n^I$ and
$K_1^J,\ldots,K_n^J$ which are appropriately disjoint relative to the
$I_i$'s (resp., the $J_i$'s) and to each other, and $h$ is the product of
{\em three} elements $h_2,h_3,h_4$.

(1b) follows from Proposition~\ref{BDD} since with
inclusion-transitivity, the existence of a nonidentity bounded
element implies
local movability.

\medskip

Sometimes we will need to ``adjust'' a bounded element $h$ at the
``ends'' while leaving it unchanged in the ``middle'':

\begin{lemma}
\label{ADJUST}
Let $\pair{N}{G}$ be a nearly ordered permutation group which is
highly interval-transitive and has a nonidentity
bounded element.  Let $h \in G$ be bounded, and let
$(a,b)$ and $(c,d)$ be $\barN$-intervals such that
$(a,b) \cup h(a,b) \subset ^S (c,d)$ (see
Definition~\ref{DEFNEST}).
Then there exists a bounded $k \in G$ such that
$g \eqdf h^k$ agrees with $h$ on $(a,b)$ and such that $\fs{supp}(g)
\subseteq (c,d)$.
\end{lemma}

\noindent
{\bf Proof }
Recall that $\pair{N}{G}$ is highly approximately {\it o-}transitive by
Proposition \ref{HAT}.
We deal now with the linear case.  First we ``adjust''
$\fs{inf}(\fs{supp}(h))$, that is, we construct $k_1 \in \fs{Bdd}(G)$
such that $g_1 \eqdf h^{k_1}$ agrees with $h$ on $(a,b)$ and such
that $c < \fs{supp}(g_1)$.

We may assume that $\fs{inf}(\fs{supp}(h)) < c$ (else take $k_1 =
\fs{Id}$), making
\[
\fs{inf}(\fs{supp}(h)) < c < a_1 \eqdf \fs{min}\{a,h(a)\}.
\]
Pick $p \in \fs{Bdd}(G)$ and $\barL$-intervals $I_1$ and $I_2$ such
that
\[
I_1 < I_2 < p(I_1) < \fs{sup}(\fs{supp}(p)).
\]
Now pick $f \in G$ such that $f(\fs{inf}(\fs{supp}(h)) \in I_1$, $f(c)
\in I_2$, and $f(a_1) > \fs{sup}(\fs{supp}(p))$.  Then $k_1 \eqdf
p^{f\inverse}$ is as required since $k_1(\fs{inf}(\fs{supp}(h))) > c$
and $\fs{supp}(k_1) < a_1$.

We can adjust $\fs{sup}(h)$ by a dual construction, involving
counterparts $p_2$ and $f_2$ of $p$ and $f$.  To construct a single
$k$ as required, we combine these two constructions by making
$\fs{sup}(\fs{supp}(p)) < \fs{inf}(\fs{supp}(p_2))$ and then
combining the stipulations about $f$ and $f_2$ to produce a single
new $f$.

The proofs for the other three types of nearly ordered permutation
groups are similar.
\hfill\qed
\medskip

For (2) of \ref{CONSEQ}, we apply (1a) to $H=\fs{Bdd}(G)$ to conclude
that there exists $h \in \fs{Bdd}(G)$ such that $h(a_i) \in \barJ_i$ for
$i = 1,\ldots,n$.  Then we obtain $g$ by adjusting $h$, applying the
lemma with $(a,b) = (a_1,a_n)$ and $(c,d) = K$.

For (3a) the proof entails revisiting the circular argument in the
proof of (1a), with $H = \fs{Bdd}(G)$.  Let $K'$ be a bounded
$\barC$-interval such that $K \subset^S K'$, and pick just one
sequence
$K_1,\ldots,K_n$ such that
$$
\fs{Crs}(\fs{ccw}(K),K_1,\ldots,K_n,\fs{ccw}(K')).
$$
In each instance of the special case, it can be arranged that the $h$
being produced as a conjugate of the (bounded) element $h_1$ have
$\fs{supp}(h) \subseteq K'$; this is accomplished by imposing two
obvious additional stipulations on the conjugator $f$.  Then the final
$h$ produced will also have $\fs{supp}(h) \subseteq K'$, and can be
adjusted as in the proof of (2).

Finally, (3b) follows from the proof of Theorem~\ref{3H}(2a) since in
the induction step the formation of $g$ from $h$ adds just one
bounded factor.

This concludes the proof of Theorem~\ref{CONSEQ}.
\hfill\qed

\bigskip

\noindent
{\bf Proof of Theorem~\ref{EXACT}}
The proof of each analogue is similar to (and easier than) that of its
counterpart, and is left to the reader, except for (1b) and (2b) of
Theorem~\ref{3H}, which will be proved after Theorem \ref{UPWARD}.
\hfill\qed
\medskip

\noindent
{\bf Proof of Theorem~\ref{HOT}}
(1) We treat the linear case, the other cases being similar.  First we
show that for {\em normal} subgroups $\{\fs{Id}\} \neq H \unlhd\, G$,
$\pair{\barN}{H}$ is highly {\it o-}transitive.

Note first that $\pair{N}{\fs{Bdd}(H)}$ is transitive.  For, given $a < b
\in N$, pick $h \in \fs{Bdd}(H)- \{\fs{Id}\}$ by Theorem \ref{CONSEQ},
pick $c \in N$ such that $c < h(c)$ (replacing $h$ by $h\inverse$ if
necessary), and pick $g \in G$ such that $g(c) = a$ and $g(h(c)) = b$.
Then $h^g(a) = b$ and $h^g \in \fs{Bdd}(H)$.

Now let $a_1 < \cdots < a_n$ and $b_1 < \cdots < b_n$ in $N$, and suppose
by induction that there exists $h\in H$ such that $h(a_i) = b_i$ for $i =
1,\dots,n-1$.  Use the transitivity of $\pair{N}{\fs{Bdd}(H)}$ to pick
$h_1 \in \fs{Bdd}(H)$ such that $h_1(h(a_n)) = b_n$.
Then adjust $h_1$ by picking $g \in G$ such that $h_1^g(h(a_n)) =
h_1(h(a_n)) = b_n$ and $b_{n-1} < \fs{supp}(h_1^g)$; this is done by using
Lemma \ref{ADJUST} with $c = b_{n-1}$ and $d > \fs{max} \{h(a_n),b_n \}$.
Now
$(h_1^gh)(a_i) = b_i$ for $i = 1,\dots,n$, and $h_1^gh \in H$.

This establishes (1) for $H \unlhd\, G$, and in view of Theorem
\ref{CONSEQ}(1b), (1) for
$H \unlhd\!\unlhd\, G$ follows by induction on the defect of $H$ (the
minimum length of a subnormal series from $G$ down to $H$).

(2) Let $a_1 < \cdots < a_n$ and $b_1 < \cdots < b_n$ in $\barN$ with each
$b_i$ in the same $\fs{Bdd}(G)$-orbit as $a_i$. An $h_2 \in \fs{Bdd}(G)$
for which $h_2(a_i) = b_i$ for $i = 1,\dots,n$ can be produced as in (1),
taking $H = \fs{Bdd}(G)$, with no need for the 2-{\it o-}transitivity
hypothesis since it is given in advance that $b_n$ and $a_n$ (and thus
also $h(a_n))$ all lie in the same $\fs{Bdd}(G)$-orbit.
\hfill\qed
\medskip

Now we turn to the proofs of parts (1b) and (2b) of
Theorem~\ref{3H}, in which $N$ is Dedekind complete.  Actually, in
each case the transitivity hypotheses can be weakened a bit.  Even
without Dedekind completeness, condition (1a) of the following
definition is obviously a consequence of
approximate 2-{\it o-}transitivity, (1b) of 2-{\it o-}transitivity, and
(2) of 3-{\it o-}transitivity.  (2) is self-dual, and (1a) and (1b) are
equivalent to their duals under the hypotheses of Theorem~\ref{3H}
(because of Theorem~\ref{UPWARD} below).

\begin{defn} \label{}
\begin{rm}
Let $\pair{N}{G}$ be a nearly ordered permutation group.

(1) Let $N=\pair{L}{<}$ or $N=\fs{ED}(L)$.

$\,\,\,\,\,\,\,\,\,\,$(a) We call $\pair{N}{G}$ {\em (exactly)
2-upward-transitive} if $\pair{N}{\fs{Opp}(G)}$ is 1-transitive, and if
for some (hence for any) $a \in L$ and for any $b,c \in L$ such that
$a< b < c$, there exists $g \in \fs{Opp}(G)$ such that
$g(a) = a$ and $g(b) > c$.
This weakening of 2-{\it o-}transitivity was
known in \cite{McPRIM} as {\em o-2-semitransitivity} and in
\cite{GL} as {\em weak double transitivity}.

$\,\,\,\,\,\,\,\,\,\,$(b) We call $\pair{N}{G}$ {\em approximately
2-upward-transitive} if for any $a < b < c \in L$ and any $L$-interval
$I$ such that $a \in I$, there exists $g \in \fs{Opp}(G)$ such that
$g(a) \in I$ and $g(b) > c$.  The existence of such a $g$ then holds
even if $a \not\in I$.  (If $G(a) \cap I = \emptyset$, let $I_1$ be the
largest convex subset of $L$ which contains I and is disjoint from
$G(a)$.  $I_1$ is a proper {\it o-}block of $\pair{N}{G}$, which by the
definition of approximate 2-upward-transitivity cannot have such
blocks.)  This weakening of approximate 2-{\it o-}transitivity is
strict; see Example~\ref{WRS2}.

(2) Let $N = C = \fs{CR}(L)$ or $N = \fs{EO}(C)$.
We call $\pair{N}{G}$ {\em approximately 3-outward-transitive}
if for any $a \in C$, any $b,c \in C$ such that
$\fs{Cr}(b,a,c)$, and any bounded $C$-intervals $I,J$ such that $a \in
I \subseteq J$, there exists $g \in G$ such that $g(a) \in I$ and $g(b),
g(c) \not\in J$.  As before, the existence of such a $g$ then holds
even if $a \not\in I$ (since approximate 3-{\it o-}transitivity implies
inclusion transitivity).

(3) We call $\pair{N}{G}$ {\em locally sweeping} if for all bounded
$N$-intervals I and all $a \in \barI$, there exists $g \in G$ such that
$\fs{supp}(g)\subseteq I$ and $g(a) \neq a$.  (In the
linear/monotonic case this implies that $\setm{g(a)}{\fs{supp}(g) \in
I}$ is coterminal in $\barI$---else apply the definition to $a_1 \eqdf
\fs{sup}\setm{g(a)}{\fs{supp}(g) \subseteq I}$ or dually.  An
analogous statement holds in the circular/monocircular case.)
\end{rm}
\end{defn}

\smallskip

Parts (1b) and (2b) of Theorem~\ref{3H} follow immediately from

\begin{theorem}
\label{UPWARD}
Let $\pair{N}{G}$ be a nearly ordered permutation group having a
nonidentity bounded element.  Suppose that $N$ is Dedekind
complete.  The following are equivalent:

{\rm (1)} $\pair{N}{G}$ is highly approximately {\it o-}transitive.

{\rm (2)} $\pair{N}{G}$ is approximately 2-upward-transitive (in the
linear/monotonic case); or approximately 3-outward-transitive (in
the
circular/monocircular case).

{\rm (3)} $\pair{N}{G}$ is locally sweeping.
\end{theorem}

\noindent
{\bf Proof }
(1) $\Rightarrow$ (2)  This is obvious.

(2) $\Rightarrow$ (3) We deal first with $N=L$ and $N=\fs{ED}(L)$.
Let $a < b < c \in L = \barL$.  We seek a $g \in G$ such that $a <
\fs{supp}(g) < c$ and $g(b) \neq b$.

To obtain $h \in \fs{Bdd}(G)$ such that $\fs{supp}(h) < c$ and $h(b)
\neq b$, pick $p \in \fs{Bdd}(G) - \fs{Id}$, pick an $L$-interval $I
\subseteq \fs{supp}(p)$, and use (2) to pick $f\in G$ such that $f(b)
\in I$ and $f(c) > \fs{supp}(p)$.  Then $h \eqdf p^{f\inverse}$ is as
required.

Hence there exists $a_1 \in L$ such that there exists $g \in G$ with
$a_1 < \fs{supp}(g) < c$ and $g(b) \neq b$, and we let $d$ be the
supremum of such $a_1$'s.  Suppose by way of contradiction that $d
\leq a$.  Apply the previous paragraph to obtain $h_1 \in G$ such
that $\fs{supp}(h_1) < b$ and $h_1(d) > d$.  Pick $a_1 < d$ such that
$h_1(a_1) > d$, and $g_1 \in G$ with $a_1 < \fs{supp}(g_1) < c$ and
$g_1(b) > b$.  Then $g_2 \eqdf g_1^{h_1}$ yields a contradiction since
\begin{eqnarray*}
& d < h_1(a_1) < h_1(\fs{supp}(g_1)) = \fs{supp}(g_1^{h_1}) < c
\mbox{ and } \\*
& g_1^{h_1}(b) = h_1(g_1(b) > h_1(b) = b.
\end{eqnarray*}

Now we deal with $N=C=\fs{CR}(L)$ and $N=\fs{EO}(C)$.  Let
$\fs{Cr}(b,a,c)$ in $\barC = C$.  We seek a $g \in G$ such that
$\fs{Crs}(b,\fs{supp}(g),c)$ and $g(a) \neq a$.

Since the hypotheses force $\pair{N}{G}$ to be locally moving, we
may pick $g_1 \in G - \{\fs{Id}\}$ such that
$\fs{Crs}(b,\fs{supp}(g_1),c)$, pick a $C$-interval $I \subseteq
\fs{supp}(g_1)$, and use (2) to pick $h \in G$ such that $h(a) \in I$
and $h(b),h(c) \not \in (b,c)$.
Then $g \eqdf g_1^{h\inverse}$ is as required.

(3) $\Rightarrow$ (1)
We deal first with $N=L$, the other cases being similar.

We claim that for any bounded $L$-intervals $J \subseteq I$ and
any $a \in I$, there exists $g \in G$ such that $\fs{supp}(g)
\subseteq I$ and $g(a) \in J$.  Suppose we have a counterexample to
the claim.  We may assume that $a < J$.  Then
\[
b \eqdf \fs{sup}\setm{g(a)}{g \in G, \fs{supp}(g) \subseteq I, \mbox
{ and } g(a) < J} < J.
\]
Pick an $L$-interval $K$ such that $\fs{inf}(K) = \fs{inf}(I)$ and
$\fs{sup}(K) \in J$.  Use (3) to pick $h \in G$ such that $\fs{supp}(h)
\subseteq K$ and $h(b) > b$,
pick $c$ such that $c < b < h(c)$, and pick $g_1 \in G$ such that
$\fs{supp}(g_1) \subseteq I$ and $c < g_1(a) < J$.  Let $g = hg_1$.
Then $\fs{supp}(g) \subseteq I$ and $b < g(a) < \fs{sup}(K) <
\fs{sup}(J)$, contradicting the choice of $b$ and proving the claim.

Now let $a_1 < \ldots < a_n \in L$ and let $J_1 < \ldots < J_n$ be
bounded $L$-intervals.  We need to show there exists $g \in G$ such
that $g(a_i) \in J_i$ for $i = 1,\ldots,n$.  The above claim guarantees
this for $n=1$.  Suppose by induction that there exists $h \in G$ such
that $h(a_i) \in J_i$ for $i = 1,\ldots,n-1$.  Pick a bounded
$L$-interval $I$ containing $h(a_n)$ and $J_n$ but not $h(a_{n-1})$.
Use
the claim again to pick $g_1 \in G$ such that $\fs{supp}(g_1)
\subseteq I$ and $g_1(h(a_n)) \in J_n$.  Then $g \eqdf g_1h$ is as
required.
\hfill\qed
\medskip

We remark that under the hypotheses of Theorem~\ref{UPWARD}, an
approximately 2-upward-transitive $\pair{N}{\fs{Opp}(G)}$ must be
(exactly)
2-upward-transitive on each of its orbits.  (This is a consequence of
(1), courtesy of Theorem~\ref{CONSEQ}(2).)  An analogous remark
applies to 3-outward-transitivity.

\medskip
Now we complete the proof of Theorem~\ref{EXACT} by treating
the analogues of (1b) and (2b) of Theorem~\ref{3H}.

For (1b), when $N$ is Dedekind complete, 2-{\it o-}transitivity
$\Rightarrow$ approximate 2-{\it o-}transitivity  $\Rightarrow$ high
approximate {\it o-}transitivity by the original
Theorem~\ref{3H}(1b), and Theorem~\ref{HOT} yields high
{\it o-}transitivity.  For (2b) the argument is analogous.
\hfill\qed

\medskip

The next proposition sharpens Proposition~\ref{NOHELP} for the
special case of full automorphism groups $G = \fs{Aut}(N)$.

\begin{prop}
\label{AUTLC}\

{\rm (1)} Let $L$ be a chain.  If
$L$ is 2-interval-homogeneous, then:

$\,\,\,\,\,\,\,\,\,\,${\rm (a)} $\fs{Aut}(L)$ and even $\fs{Bdd}(L)$ are
highly {\it o-}transitive on each (necessarily dense) orbit $M$ of
$\pair{\barL}{\fs{Aut}(L)}$.  In particular,
2-interval-homogeneity plus transitivity implies high
{\it o-}homogeneity.

$\,\,\,\,\,\,\,\,\,\,${\rm (b)} Each $L$-interval supports a positive $g
\in
\fs{Aut}(L)$.  $\pair{L}{\fs{Aut}(L)} \in K^{\fss{LNn}}$ for all
$n$, and similarly for $\pair{L}{\fs{Bdd}(L)}$.

{\rm (2)} Part (1) holds verbatim for circles $C$ (with ``positive''
replaced by ``signed'', and $K^{\fss{LNn}}$ by $K^{\fss{CRn}}$), the
hypothesis still being only 2-interval-homogeneity.  Here
$\fs{Bdd}(C) = \fs{Aut}(C)$.
\end{prop}

{\bf Note } This proposition is valid even when
2-interval-transitivity is weakened to inclusion-transitivity.

\medskip

\noindent
{\bf Proof }
(1a) for $\fs{Aut}(L)$ is a special case of Proposition~\ref{NOHELP}.
New we prove (1a) for $\fs{Bdd}(L)$.  Let $x_1 < \dots < x_n$ and
$y_1 < \dots
< y_n$ in $M$.  Pick $u < \fs{min} \{x_1,y_1 \}$ and $v > \fs{max}
\{x_n,y_n \}$.  Use (1a) to pick $f \in \fs{Aut}(L)$ such that $f(x_i) =
y_i$ for $i=1,\dots,n$, and such that $f(u) = u$ and $f(v) = v$.  Then
modify $f$ outside $[u,v]$ so that $f \equiv \fs{Id}$ there.  This
makes $f \in \fs{Bdd}(L)$, and $f(x_i) = y_i$ for $i=1,\dots,n$.

(1b) is now obvious.

(2) First we prove (2) for $\fs{Aut}(C)$.  It suffices to show that any
two bounded $C$-intervals with endpoints in $M$ are isomorphic
(when totally ordered in the counterclockwise direction), for then
high {\it o-}transitivity can be achieved by splicing appropriate
isomorphisms together.  For this, it suffices to show that for
$\fs{Cr}(x_1,x_2,y_2)$ in $M$, there exists $f \in \fs{Aut}(C)$ such
that $f(x_1) = x_1$ and $f(x_2) = y_2$.

We claim that the stabilizer subgroup $G_{y_2}$ cannot fix any other
point $z \in \barC$.  For pick bounded $\barC$-intervals $I$ and $J$
such that
$$
\fs{Crs}(\fs{cw}(I),z,J,y_2,\fs{ccw}(I)).
$$
Then use inclusion-transitivity to pick $h \in \fs{Aut}(C)$ such that
$h(I) \subseteq J$.  Then $h^{-1}(J) \supseteq I$, ensuring that
$$
\fs{Cr}(h^{-1}(z),z,y_2,h^{-1}(y_2)).
$$
Now $h^{-1}$ must fix some $u,v \in \barC$ for which
$$
\fs{Cr}(u,h^{-1}(z),z,v,y_2,h^{-1}(y_2)).
$$
Form $k$ by modifying $h^{-1}$
outside $(u,v)$ so that $k \equiv \fs{Id}$ there.  Then $k \in
G_{y_2}$ and $k(z) \neq z$, proving the claim.

Now, back to the construction of the required $f$.  Pick $g \in
\fs{Aut}(C)$ such that $g(x_2) = y_2$.  Since $G_{y_2}$ fixes no $z \in
\barC$, it can be arranged by following $g$ by some $g_1 \in
\fs{Aut}(C)$ that $\fs{Cr}(g(x_1),x_1,x_2,y_2)$ as well as $g(x_2) =
y_2$.  As in the proof of the claim, $g$ must fix some $u$ and $v$ for
which
$$
\fs{Cr}(u,g(x_1),x_1,v,x_2,y_2),
$$
enabling the formation of $f \in \fs{Aut}(C)$ for which $f(x_1) = x_1$
and $f(x_1) = y_2$ as required.

$\fs{Bdd}(C) = \fs{Aut}(C)$ because the above ``modification
techniques'' can be used to write an arbitrary $g\in \fs{Aut}(C)$ as a
product of two elements of $\fs{Bdd}(C)$.
\hfill\qed
\medskip

How can we make worthwhile circles from chains?

\begin{prop}
\label{WRAP}\
Let $L$ be a 2-{\it o-}homogeneous chain.  Suppose that for one
(hence every) $a \in L$, the concatenation $(a,\infty) \cat
(-\infty,a)$ of $L$-rays is order-isomorphic to one (hence every)
bounded $L$-interval.  Then the circular order $C$ based on $L$ is
highly
{\it o-}homogeneous.
\end{prop}

\noindent
{\bf Proof }
The hypotheses guarantee that any bounded $C$-interval $(a,b)$
containing the point $\infty$ at which the ends of $L$ are joined
must be isomorphic to all the other bounded $C$-intervals
(when totally ordered in the counterclockwise direction).  Now high
{\it o-}transitivity can be achieved by splicing appropriate
isomorphisms together.
\hfill\qed
\newpage

%

\section{Counterexamples}
\label{EXAMPLES}

Here we give examples of nearly ordered permutation groups for
which reconstruction does not hold, and examples showing that
certain transitivity properties do not imply others.  After a bit of
preparation, we list the examples and their properties.
Then we turn to constructing them, and that we do in the most
convenient order, which differs from the order in which the
examples are presented.

We begin with linear/monotonic examples.
Most of our linear examples $\pair{L}{G}$ enjoy the following
properties:

(1) $\pair{L}{G}$ is locally moving.

(2) $L$ is a dense subchain of the real line $\R$, and all orbits of
$\pair{\barL}{G}$ ($= \pair{{\R}}{G}$) are dense in $\R$.

(3) $\pair{L}{G}$ is transitive, and even better, $G$ contains a
transitive {\it l-}subgroup of the
lattice-ordered group $\fs{Aut}(L)$.

We refer to the items (1)--(3) as the Linear List.
Just as $n$-interval-transitivity, nest transitivity, and weak
span-transitivity are robust under change of $G$-invariant subchain
and under enlargement of the group (see the remarks preceding
Proposition~\ref{ENLARGEGROUP}), so are the items in the Linear
List, provided in (3) that the new
$G$-invariant subchain is a single orbit of $\pair{\barL}{G}$.
Sometimes we get a sharpening of (2), with $L = \Q$.
It is a consequence of (3) that $\pair{L}{G}$ is coherent (see the remarks
after Definition~\ref{DEFOPRIM}).

Similarly, in most of our monotonic examples $\pair{\fs{ED}(L)}{H}$,
$\pair{L}{G} \eqdf \pair{L}{\fs{Opp}(H)}$ enjoys the properties in
the Linear List (and in particular, $\pair{L}{\fs{Opp}(H)}$ is
transitive, which is not always the case for monotonic permutation
groups).  This we refer to by saying that $\pair{\fs{ED}(L)}{H}$
enjoys the properties in the Monotonic List.

The Linear Reconstruction Theorem~\ref{LRT} says that First
Order Reconstruction holds for 2-interval-transitive locally moving
linear permutation groups.  Here is an example in which
reconstruction does not hold, despite the presence of some milder
transitivity properties.

\begin{example}
\label{WRS2}
A linear permutation group $\pair{L}{G}$ such that:

{\rm (1)} Reconstruction (even second order) does not hold for
$\pair{L}{G}$---there is an automorphism $\alpha$ of $G$
not induced by any monotonic bijection $\tau$ of $\barL$.

{\rm (2)} $\pair{L}{G}$ is inclusion transitive, and is
2-upward-transitive and dually.

{\rm (3)} $\pair{L}{G}$ enjoys all the properties in the Linear List;
also $L= \Q$.
\end{example}

The chain $L$ used in Example \ref{WRS2} will also serve for
Examples \ref{WRSN}, \ref{WR1}, and  \ref{WR2}.

The Monotonic Reconstruction Theorem~\ref{MNRT} says that
First Order Reconstruction holds for 3-interval-transitive locally
moving monotonic permutation groups. In light of the
Linear Reconstruction Theorem, an example showing that 3 cannot be
reduced to 2 must contain an order-reversing permutation:

\begin{example}
\label{WRSN}
A monotonic permutation group $\pair{\fs{ED}(L)}{H}$ containing an
order-reversing permutation, such that:

{\rm (1)} Reconstruction does not hold for
$\pair{\fs{ED}(L)}{H}$---there is an automorphism $\alpha$
of $H$
not induced by any monotonic bijection of $\barL$.

{\rm (2)} $\pair{\fs{ED}(L)}{H}$ is 2-{\it o-}transitive.

{\rm (3)} $\pair{\fs{ED}(L)}{H}$ enjoys all the properties in the
Monotonic List.
\end{example}

Several of the examples occur in pairs related as follows:

\begin{defn}
\label{AFFILIATE}
\begin{rm}
Let $\pair{L}{G}$  be a linear permutation group.

(1) An {\em affiliate} of $\pair{L}{G}$ is a monotonic permutation
group $\pair{\fs{ED}(L_1)}{H}$, with $L \subseteq L_1 \subseteq
\barL$ and $L_1 = H(L)$, for which $\fs{Opp}(H) = G$.  Every
monotonic permutation group $\pair{\fs{ED}(L_1)}{H}$
is an affiliate of some linear permutation group, namely
$\pair{L_1}{\fs{Opp}(H)}$.

(2) An affiliate is {\em proper} if $H$ contains order-reversing
permutations, in which case $H$ is generated by $G$ together with
any one order-reversing $h \in H$.
\end{rm}
\end{defn}

Even for a transitive $\pair{L}{G}$, two distinct orbits of
$\pair{\barL}{G}$ may coalesce into a single orbit of an affiliate
$\pair{\fs{ED}(L)}{H}$; and in particular $L$ may coalesce with some
other orbit $L'$ of $\pair{\barL}{G}$, making  $H(L) = L
\cup L'$.  (See (the linear/monotonic analogue of)
Example~\ref{EXALPHA}, and also Examples \ref{WR1} and
\ref{WR1VARIANT}.)  This
coalescing is limited to two
orbits of $\pair{\barL}{G}$ per orbit of $\pair{\fs{ED}(L_1)}{H}$ since
for $h \in H$ we have $h^2 \in \fs{Opp}(H) = G$.

For $n > 1$, an affiliate $\pair{\fs{ED}(L_1)}{H}$ is
$n$-{\it o-}transitive iff $\pair{L}{G}$ ($= \pair{L_1}{\fs{Opp}(H)}$)
is $n$-{\it o-}transitive (and similarly for approximate and
interval-transitivity), so we omit mention of such multiple
transitivity for the affiliate.

The next example shows for both linear and monotonic permutation
groups that 2-{\it o-}transitivity does not imply 3-{\it o-}transitivity
(or even
3-interval-transitivity) and does not imply nest-transitivity or weak
span-transitivity, and that 2-{\it o-}transitivity and
approximate 2-{\it o-}transitivity are not robust under change of
$G$-invariant subset.
Actually, the $\pair{\fs{ED}(L)}{H}$ of Example~\ref{WR1} will be the
same as that of
Example~\ref{WRSN}, but for the sake of clarity these two examples
are presented separately here.

\begin{example}
\label{WR1}
A linear permutation group $\pair{L}{G}$, together with a proper
monotonic affiliate $\pair{\fs{ED}(L)}{H}$, such that:

{\rm (1)} $\pair{L}{G}$ is 2-{\it o-}transitive but not
3-interval-transitive.

{\rm (2)} $\pair{\barL}{G}$ has precisely three orbits $\neq L$, two
of which coalesce into a single orbit of $\pair{\fs{ED}(L)}{H}$.  On
each of these three orbits, $G$ is 2-interval-transitive but not
approximately 2-{\it o-}transitive or 3-interval-transitive.

{\rm (3)} $\pair{L}{\fs{Lft}(G)}$ is not 2-interval-transitive, and it
has orbits which are not dense in $\barL$.

{\rm (4)} $\pair{L}{G}$ and $\pair{\fs{ED}(L)}{H}$ are neither
nest-transitive nor weakly span-transitive.

{\rm (5)} $\pair{L}{G}$ enjoys all the properties in the Linear List,
and $\pair{\fs{ED}(L)}{H}$ enjoys all those in the Monotonic List.
\end{example}

Almost all the features of the previous example can be accomplished
with $L =
\Q$:

\begin{example}
\label{WR1Q}
A variant of Example \ref{WRSN}/\ref{WR1} in which $L = \Q$, and
all properties
mentioned before still hold (except that besides the previous three
orbits $\neq L$ of $\pair{\barL}{G}$, there are uncountably many
others, on each of
which $G$ is 2-{\it o-}transitive).
\end{example}

2-interval-transitivity does not imply approximate
2-{\it o-}transitivity:

\begin{example}
\label{WR1VARIANT}
A linear permutation group $\pair{L}{G}$, together with a proper
monotonic affiliate $\pair{\fs{ED}(L \cup L')}{H}$, where $L'$
is an orbit $\neq L$ of $\pair{\barL}{G}$, and $L$ and $L'$ coalesce
into a single orbit $L \cup L'$ of $H$, such that:

{\rm (1)} $\pair{L}{G}$ is 2-interval-transitive but not
approximately 2-{\it o-}transitive or
3-interval-transitive.

{\rm (2)} $\pair{L}{G}$ and $\pair{\fs{ED}(L)}{H}$ are neither
nest-transitive nor weakly span-transitive.

{\rm (3)} $\pair{L}{G}$ enjoys all the properties in the Linear List;
and $L \cong \Q$.
\end{example}

\begin{question}
\label{Q2ITIMP2ATORB}
\begin{rm}
Must a locally moving {\it n-}interval-transitive nearly ordered
permutation group $\pair{N}{G}$ be approximately
{\it n-o-}transitive on some orbit in $\pair{\barN}{G}$?
This is of interest for $n=2$ in the linear case, for $n=3$ in the
monotonic and circular cases, and for $n=4$ in the monocircular case.
\end{rm}
\end{question}

Approximate 2-{\it o-}transitivity does not imply exact
2-{\it o-}transitivity, or even the existence of an orbit in
$\pair{\barL}{G}$ on which $G$ is exactly 2-{\it o-}transitive:

\begin{example}
\label{WR2}
A linear permutation group $\pair{L}{G}$, together with a proper
monotonic affiliate $\pair{\fs{ED}(L)}{H}$, such that:

{\rm (1)} $\pair{L}{G}$ is approximately 2-{\it o-}transitive but not
2-{\it o-}transitive or
3-interval-transitive.

{\rm (2)} $\pair{\barL}{G}$ has precisely three orbits $\neq L$, two
of which coalesce into a single orbit of $\pair{\fs{ED}(L)}{H}$.  On
each of these three orbits, $G$ is 2-interval-transitive but not
approximately 2-{\it o-}transitive or 3-interval-transitive.

{\rm (3)} $\pair{L}{G}$ and $\pair{\fs{ED}(L)}{H}$ are neither
nest-transitive nor weakly span-transitive.

{\rm (4)} $\pair{L}{G}$ enjoys all the properties in the Linear List
and $\pair{\fs{ED}(L)}{H}$ all those in the Monotonic List.
\end{example}

2-{\it o-}transitivity is not robust under group enlargement:

\begin{example}
\label{WR2VARIANT}
A variant $\pair{L_1}{G_1}$, with proper affiliate
$\pair{\fs{ED}(L_1)}{H_1}$, of the $\pair{L}{G}$ and
$\pair{\fs{ED}(L)}{H}$ of Example~\ref{WR2}
such that:

(1) $L_1$ is a dense subchain of $L$, $G_1 \leq G$, and
$H_1 \leq H$.

(2) $\pair{L_1}{G_1}$ is 2-{\it o-}transitive but
$\pair{G(L_1)}{G} = \pair{L}{G}$ is not (though it is approximately
2-{\it o-}transitive).

(3) $\pair{L_1}{G_1}$ enjoys all the properties in the Linear List
and $\pair{\fs{ED}(L_1)}{H_1}$ all those in the Monotonic List; and
$L_1 \cong \Q$.
\end{example}

We offer another very different way to achieve
Example~\ref{WR2VARIANT}:

\begin{example}
\label{BACKANDFORTH}
Another $\pair{L_1}{G_1}$, $\pair{\fs{ED}(L_1)}{H_1}$, $\pair{L}{G}$,
$\pair{\fs{ED}(L)}{H}$ enjoying all the properties listed in
Example~\ref{WR2VARIANT} (except that the
$\pair{L}{G}$ and $\pair{\fs{ED}(L)}{H}$ are no longer those of
Example~\ref{WR2}).
Moreover, $\pair{L_1}{G_1}$ is highly {\it o-}transitive but there is
an orbit of $\pair{\barL_1}{G_1}$ on which $G_1$ is not even
2-{\it o-}transitive.
\end{example}

1-transitivity does not suffice for reconstruction even when dealing
with the entire automorphism group:

\begin{example}
\label{AUTRZ}
$\pair{L}{G}$, where $L= \R \times \Z$, ordered lexicographically
from the right, and $G = \fs{Aut}(L)$; together with proper affiliate
$\pair{\fs{ED}(L)}{H}$, where $H= \fs{Aut}(\fs{ED}(L))$.  Here:

{\rm (1)} Reconstruction holds for neither---there is an
automorphism $\alpha$ of $H$ for which
$\alpha(G) = G$ and neither $\alpha$ nor $\alpha \rest G$ is induced
by any monotonic bijection of $\barL$.

{\rm (2)} $\pair{L}{G}$ is a transitive lattice-ordered permutation
group, and
in fact $G = \fs{Aut}(L)$.

{\rm (3)} $\pair{\fs{ED}(L)}{H}$ is a half-lattice-ordered permutation
group (see Definition~\ref{HALFORDERED}), and
in fact $H = \fs{Aut}(\fs{ED}(L))$.

{\rm (4)} $\pair{L}{G}$ and $\pair{\fs{ED}(L)}{H}$ are not
{\it o-}primitive.

{\rm (5)} $\pair{L}{G}$ enjoys all the properties in the Linear List
and $\pair{\fs{ED}(L)}{H}$ all those in the Monotonic List, except that
not all the orbits are dense in $\barL$.
\end{example}

Without local movability, reconstruction need not hold even when
the other hypotheses are very strong:

\begin{example}
\label{FREE}
{\rm (S.H. McCleary \cite[Counterexample 20]{McFREE}) }
The free lattice-ordered group $F_{\eta}$ of any rank $\eta > 1$ has
a faithful representaiton as a highly {\it o-}transitive lattice-ordered
permutation group $\pair{L}{F_{\eta}}$ (not locally moving) for
which there is a {\it l-}group automorphism $\alpha$ not
induced by any monotonic bijection of $\barL$.
\end{example}

We turn now to circular/monocircular examples.
Our circular examples $\pair{C}{G}$ enjoy the following
properties:

(1) $\pair{C}{G}$ is locally moving.

(2) $C$ is a dense subcircle of the real circle $\C$, and all orbits of
$\pair{\barC}{G}$ ($= \pair{{\C}}{G}$) are dense in $\C$.

(3) $\pair{C}{G}$ is transitive.

These items constitute the Circular List.
Robustness of the items parallels the linear/monotonic case.
The Monocircular List is analogous to the Monotonic List.
The definition of {\em affiliate} parallels the
linear/monotonic case, as does the possible coalescing of orbits of
$\pair{\barC}{G}$.  Here it is for $n>2$ that an affiliate
$\pair{\fs{EO}(C_1)}{H}$ is
$n$-{\it o-}transitive iff $\pair{C}{G}$ is $n$-{\it o-}transitive.
Coalescing of orbits can occur as in the linear/monotonic case; see
Example~\ref{EXALPHA}.

The Circular and Monocircular Reconstruction Theorems say that
First Order Reconstruction holds for 3-interval-transitive locally
moving circular permutation groups and for 4-interval-transitive
locally moving monocircular permutation groups.  The following
(rather cheap) example shows only that 1-transitivity does
not suffice, leaving a large gap between theorem and example.

Let ${\C}_{(n)}$ denote the circle obtained by wrapping $n$ copies of
$\R$
in the obvious circular fashion; or equivalently, by deleting $n-1$
points from the real circle $\C$.

\begin{example}
\label{CIRCLE1}
$\pair{C}{G}$, where $C= {\C}_{(n)}$ (for $n \geq 3$) and $G =
\fs{Aut}(C)$;
together with proper affiliate
$\pair{\fs{EO}(C)}{H}$, where $H= \fs{Aut}(\fs{EO}(C))$.  Here:

{\rm (1)} Reconstruction holds for neither---there is an
automorphism $\alpha$ of $H$ for which
$\alpha(G) = G$ and neither $\alpha$ nor $\alpha \rest G$ is induced
by any monocircular bijection of $\barC$.

{\rm (2)} $\pair{C}{G}$ and $\pair{\fs{ED}(C)}{H}$ are not
{\it o-}primitive.

{\rm (3)} $\pair{C}{G}$ enjoys all the properties in the Circular List
and $\pair{\fs{ED}(C)}{H}$ all those in the Monocircular List, except
that not all the orbits are dense in $\barC$.
\end{example}

The last example is more interesting.  It shows that also in the
circular/monocircular case, high {\it o-}transitivity does not imply
2-{\it o-}transitivity on other orbits in $\bar{\C}$, and
2-{\it o-}transitivity is not robust under group enlargement.

\begin{example}
\label{CIRCLE2}
A circular permutation group $\pair{C}{G}$ with proper affiliate
$\pair{\fs{EO}(C)}{H}$, and another $\pair{C_1}{G_1}$ with proper
affiliate $\pair{\fs{EO}(C_1)}{H_1}$,
such that:

{\rm (1)} $C_1$ is a dense subcircle of $C$, $G_1 \leq G$, and
$H_1 \leq H$.

{\rm (2)} $\pair{C_1}{G_1}$ is highly {\it o-}transitive but
$\pair{G(C_1)}{G} = \pair{C}{G}$ is not even 2-{\it o-}transitive.

{\rm (3)} There is an orbit of $\pair{\barC}{G_1}$ on which $G_1$ is
not
2-{\it o-}transitive.

{\rm (4)} $\pair{C_1}{G_1}$ enjoys all the properties in the Circular
List
and $\pair{\fs{EO}(C_1)}{H_1}$ all those in the Monocircular List; and
$C_1 = \C_{\Q}$, the rational circle.
\end{example}

{\bf Construction of the examples }
For any of the linear/monotonic non-reconstruction examples, the
basic
idea is this.  Let $P$ denote the set of elements of $G$ (or of
$\fs{Opp}(H)$ which are positive in the pointwise order, and $N$ the
dual.  Then an automorphism $\alpha$ of $G$ or $H$ which is
conjugation by an order-preserving bijection of $\barL$ must map
$P$ to $P$ and $N$ to $N$,
whereas an $\alpha$ which is conjugation by an order-reversing
bijection of $\barL$ must interchange $P$ and $N$.  The goal is to
find an $\alpha$ doing neither.

We begin with the easiest non-reconstruction example:

\bigskip

{\bf Construction of Example \ref{AUTRZ} }  Clearly each $g \in G =
\fs{Aut}(\R \times \Z)$ permutes the set of local copies of $\R$ by an
integer translation.

Define $\alpha$ to be conjugation by $\psi$, where
$\psi(r,z) = (-r,z)$;
that is, $\alpha(h) = h^{\psi}$.  Then $\alpha \in \fs{Aut}(H)$ and
$\alpha \rest G \in \fs{Aut}(G)$.  Positive elements of $G =
\fs{Opp}(H)$ which fix each local copy of $\R$ setwise are mapped
by $\alpha$ to
negative elements, whereas those which don't are mapped to positive
elements (the element and its image give the same integer
translation), so $\alpha$ is as required.

The rest is clear, using an order-isomorphism to make $L$
literally a subchain of $\R$.
\hfill\qed

\bigskip

{\bf Construction of Example \ref{CIRCLE1} } This is much like the
previous construction.  Here $\psi$ is to fix each local copy of $\R$
setwise, and to negate within each copy.  Then $\alpha$ cannot be
conjugation by an orientation-preserving bijection $\tau$ because of
the negating, nor by an orientation-reversing $\tau$ because of the
setwise effect on the local copies of $\R$.
\hfill\qed

\medskip

Most of the remaining examples are built
from one particular {\em wreath power} of
$\pair{{\Q}}{\fs{Aut}(\Q)}$, which we now define.  This is a special
case of a more general construction well known in the area of
ordered permutation groups (see \cite{McSURVEY} or \cite{GL}), except that
here the group is not equipped with the
pointwise order.

\begin{defn}
\label{}
\begin{rm}
Let $\Gamma = \Z$, which will serve as an index chain and is to be
visualized as running vertically (with larger elements higher).  Let
$K$
be the set of functions $a:\Gamma \raro \Q$ for which $\fs{supp}(a)
\eqdf \setm{\gamma}{a_{\gamma} \neq 0}$ is inversely well
ordered.  For $a \neq b \in K$, the greatest $\gamma$ for which
$a(\gamma) \neq b(\gamma)$ is denoted by $\fs{Val}(a,b)$.
Totally order $K$ lexicographically from the top, that is, by
greatest difference.

For $\gamma \in \Gamma$, let $\CC_{\gamma}$ be the convex
equivalence relation (that is, the equivalence classes are convex)
defined by
$$
 a\CC_{\gamma}b \mbox{\, iff \,} a(\delta) = b(\delta)
\mbox{ for all } \delta \in \Gamma \mbox{ such that }
\delta \geq \gamma.
$$
Let $W$ consist of those $f \in \fs{Aut}(K)$ which respect all these
equivalence relations.  Equivalently, $W$ is the automorphism group
of the structure obtained by adding to $\pair{L}{<}$ the relations
$\setm{\CC_{\gamma}}{\gamma \in \Gamma}$.  The linear
permutation group
$\pair{K}{W}$ is called the {\em wreath power of
$\pair{{\Q}}{\fs{Aut}({\Q})}$ by $\Z$}.
$\pair{K}{W}$ will denote this wreath power throughout the rest of
this paper.
\end{rm}
\end{defn}

We describe $\pair{K}{W}$ in more detail.  By construction, the
$\CC_{\gamma}$'s are convex congruences of $\pair{K}{w}$. Let
$\CC^{\gamma} = \CC_{\gamma+1}$.  Identify each
$a\CC^{\gamma} / \CC_{\gamma}$ with $\Q$ in the obvious fashion.
Ordered by inclusion,
$\setm{a\CC^{\gamma}}{a \in K, \gamma \in \Gamma}$ forms a {\em
root system}, that is,
those elements exceeding any given element form a
tower.

Let $w \in W$.  For each $a \in L$ and $\gamma \in
\Gamma$, $w$ induces an order-isomorphism
$w'_{a\CC^{\gamma}}: a\CC^{\gamma} / \CC_{\gamma} \raro
w(a)\CC^{\gamma} / \CC_{\gamma}$, which in turn induces
$w_{a\CC^{\gamma}} \in \fs{Aut}({\Q})$ via the above identifications.
The $w_{a\CC^{\gamma}}$'s are called the
{\em $\gamma$-components} of $w$.  Thus $w$ has a root system of
components
$w_{a\CC^{\gamma}} \in \fs{Aut}({\Q})$ (where the
$w_{a\CC^{\gamma}}$'s are ordered as are the $a\CC^{\gamma}$'s),
and $w(a)\CC_{\gamma} = (a\CC_{\gamma})w_{a\CC^{\gamma}}$.
We have $(vw)_{a\CC^{\gamma}} =
v_{w(a)\CC^{\gamma}}w_{a\CC^{\gamma}}$.

$\pair{K}{W}$ is locally moving since every $L$-interval contains
some $a\CC^{\gamma}$ and there exists $w \in W$ whose only
nonidentity component is $w_{a\CC^{\gamma}}$, making
$\fs{supp}(w) \subseteq a\CC^{\gamma}$.

Clearly $\pair{K}{W}$ is transitive.
$K$ is an
{\it l-}subgroup of $\fs{Aut}(\pair{K}{<})$ since if $v$ and $w$
respect $\CC^{\gamma}$, so do their pointwise supremum and
infimum.

For any $\gamma
\in \Gamma$ and $a,b \in K$,
$\setm{w_{a\CC^{\gamma}}}{w(a\CC^{\gamma}) = b\CC^{\gamma}}
= \fs{Aut}(\Q)$ under the above identification.
Hence the orbits $\neq K$ of $\pair{\barK}{W}$ are as follows.  For
each $\gamma \in \Gamma$, the suprema of the
$\CC_{\gamma}$-classes form an orbit $\barK_{\gamma}^+$,
the infima of the
$\CC_{\gamma}$-classes form an orbit $\barK_{\gamma}^-$,
and the elements of $\barK$ which correspond to irrational numbers
in the various $a\CC^{\gamma} / \CC_{\gamma}$'s form an orbit
$\barK_{\gamma}^{\irr}$.  There are no other orbits.  For if
$a \in \barK- K$, there must be some $\delta$ for which $a$ does not
lie in (the Dedekind completion of) any $\CC_{\delta}$-class;
otherwise $a \in K$.  Then take $\gamma$ to be the largest such
$\delta$, and $a \in
\barK_{\gamma}^+ \cup \barK_{\gamma}^- \cup
\barK_{\gamma}^{\irr}$.

We mention two subgroups of $W$:

(1) $W_F$ denotes the {\it l-}subgroup of $\fs{Aut}(K)$ consisting of
those $w \in W$ having only finitely many nonidentity components.
$\pair{K}{W_F}$ is not transitive.  The orbit containing any $k \in K$
consists of those $k' \in K$ differing from $k$ at only finitely many
$\gamma$'s, and there are uncountably many such orbits, each
dense in $K$.  The orbits in $\barK - K$ are the same for $W_F$ as
for $W$.

(2) $W_T$ denotes the transitive {\it l-}subgroup of $\fs{Aut}(K)$
consisting of those $w \in W$ such that
for some $\gamma$:
If $\delta \leq \gamma$ all $\delta$-components
of $w$ are the same rational translation
$w^{\delta}$, and if $\delta > \gamma$ all $\delta$-components
are the identity.

We also define two permutations not in $W$:

(1) Let $s \in \fs{Aut}(K)$ be ``shift down one'', that is,
$(s(a))(\gamma) = a(\gamma +1)$.  Clearly $s$ normalizes $W$.

(2) Let $\fs{neg} \in \fs{Aut}(\fs{ED}(K))$ be ``negation'', that is,
$(\fs{neg}(a))(\gamma) = -a(\gamma)$.  Then \fs{neg} also
normalizes $W$ and $\fs{neg}^2 = \fs{Id}$.

\medskip

Construction of affiliates is facilitated by the following easy result:

\begin{prop}
\label{}
{\rm (Scott \cite{Sc2})}
Let $\pair{L}{G}$ be a linear permutation group.  Let $h$ be an
order-reversing permutation of $L$, and let $H$ be the subgroup of
$\fs{Aut}(\fs{ED}(L))$ generated by $G$ and $h$. Then
$\pair{\fs{ED}(L)}{H}$ is a (proper) affiliate of $\pair{L}{G}$ iff $G^h =
G$ and $h^2 \in G$.
\end{prop}

Recall that the following two examples are actually the same, but
for clarity were presented separately.  Many
of the other examples are variations of this one.

\bigskip

{\bf Construction of Example \ref{WRSN}/\ref{WR1} } Take $L =
K$, take $G$ to be the subgroup of $\fs{Aut}(L)$ generated by $W$
and $s$, and take $H$ to be the subgroup of $\fs{Aut}(\fs{ED}(L))$
generated by $G$ and \fs{neg}.  The elements of $G$ are of the form
$ws^z$, where $w \in W$ and $z \in \Z$.  Since \fs{neg} commutes
with $s$, $\pair{\fs{ED}(L)}{H}$ is an affiliate of $\pair{L}{G}$.

Define $\alpha$ by setting $\alpha(h) = h^{\psi}$, where $\psi$ is the
permutation of $L$ given by
\[
(\psi(a))(\gamma) = \left\{
\begin{array}{ll}
-a(\gamma) & \mbox{if $\gamma$ is odd,} \\*
a(\gamma)  & \mbox{if $\gamma$ is even.}
\end{array}
\right.
\]
Then $\alpha(H) \subseteq H$ since $\psi$ normalizes $W$,
$s^{\psi} = \fs{neg} \circ s$, and $\fs{neg}^{\psi} = \fs{neg}$; and
$\alpha \in \fs{Aut}(H)$ since $\psi^{-1} = \psi$.
Define ${\bf 0} \in L$ by setting ${\bf 0}(\gamma) = 0$ for all
$\gamma \in \Gamma$.
Pick an odd $\gamma_1$ and an even
$\gamma_2$.
For $i = 1,2$, pick $h_i \in \fs{Opp}(H) = G$ such that
$h_{{\bf 0}\CC^{\gamma_i}}$ is its only nonidentity component and
such that
$h_{{\bf 0}\CC^{\gamma_i}}$ is
positive.  Then $h_1,h_2 \in P$, and $\alpha(h_1) \in N$
whereas $\alpha(h_2)
\in P$, so $\alpha$ is as required.

Now we show that $\pair{L}{G}$ (and thus also
$\pair{\fs{ED}(L)}{H}$) are 2-{\it o-}transitive.  Let $a_1 < a_2$ and
$b_1 < b_2$ in $L$.  Let $z = \fs{Val}(a_1,a_2) - \fs{Val}(b_1,b_2)$.
Then $\fs{Val}(s^z(a_1),s^z(a_2)) = \fs{Val}(b_1,b_2)$, and we denote
it by $\gamma$.  We pick $w \in W$ such that
$w(s^z(a_1))\CC_{\gamma} = b_1\CC_{\gamma}$ and
$w(s^z(a_2))\CC_{\gamma} = b_2\CC_{\gamma}$.
For any one $\delta < \gamma$, the $\delta$-components
$w_{s^z(a_1)\CC^{\delta}}$ and
$w_{s^z(a_2)\CC^{\delta}}$ can be picked independently of each
other.  It follows that $\pair{L}{G}$ is 2-{\it o-}transitive.

On the other hand, $\pair{L}{G}$ is not 3-interval-transitive.  In
some $\CC^\gamma$-class $M$, pick $\barL$-intervals $I_1 < I_2 <
I_3$ contained in distinct $\CC_{\gamma}$-classes within $M$.  Take
$J_1 = I_1$ and $J_2 = I_2$, and pick an $\barL$-interval $J_3 > M$.
An element $g \in G$ for which
$g(I_1)$ meets $J_1$ and $g(I_2)$ meets $J_2$ cannot involve a
shift, so $g(M) = M$, and then $g(I_3)$
does not meet $J_3$.

Because of the shift $s$, the $W$-orbits
$\setm{\barK_{\gamma}^+}{\gamma \in \Gamma}$ form a single
orbit
$\barK^+$ of $\pair{\barL}{G}$, and similarly for $\barK^-$ and
$\barK^{\irr}$.  $\barK^+$, $\barK^-$, and $\barK^{\irr}$ are the three
orbits $\neq L$ of $\pair{\barL}{G}$.  Each of these orbits is dense in
$\barL$.    Thanks to \fs{neg},
$\barK^+$ and $\barK^-$ coalesce into a single
orbit of $\pair{\fs{ED}(L)}{G}$.

$G$ is not approximately 2-{\it o-}transitive on $\barK^+$.  For
pick $\CC_{\gamma}$-classes $I_1 < I_2$ within the same
$\CC^{\gamma}$-class $M$, and let $a_1 = \fs{sup}(I_1)$ and
$a_2 = \fs{sup}(I_2)$.  Any $g = ws^z \in G$ for which $g(a_1) \in
I_1$ must have $z>0$, forcing $g(M) \subseteq I_1$ and precluding
$g(a_2) \in I_2$.  Similar arguments apply for
$\barK^-$ and $\barK^{\irr}$.

Since $G$ is 2-{\it o-}transitive and thus 2-interval-transitive on $L$,
$G$ must also be 2-interval-transitive on each of the other three
orbits, and it cannot be 3-interval-transitive on any of them because
then it would be 3-interval-transitive on $L$ as well.

Condition (3) of \ref{WR1} holds because no element of $\fs{Lft}(G)$
involves a shift.
$\pair{L}{G}$ and $\pair{\fs{ED}(L)}{G}$ are neither nest-transitive
nor weakly span-transitive
because the $\CC_{\gamma}$-classes are convex semi-blocks.  ($I$ is
a {\em semi-block} if for all $g \in G$ such that $g(I) \cap I \neq
\emptyset$, either $g(I) \subseteq I$ or $g(I) \supseteq I$.)
$\pair{L}{W}$ is transitive and $W$ is an {\it l-}subgroup of
$\fs{Aut}(L)$.
$\barL$ has a countable dense subchain, namely $\barK^+$, so we
may use an order-isomorphism to make $L$ literally a dense
subchain of $\R$.
The rest is clear.
\hfill\qed

\bigskip

{\bf Construction of Example \ref{WR1Q} }
Here we take $L$ to be the dense subchain of $K$ consisting of those
points $a \in K$ for which $\fs{supp}(a)$ is finite, and take $G$ to be
the subgroup generated by $W_F$ and $s$.  $K$ decomposes into
uncountable many dense orbits, one of which is $L$, and the orbits
within
$\barK - K$ remain as in Example~\ref{WR1}.  The rest of the
argument is as before,
except that since $L$ is a countable dense chain we may take it to be
$\Q$.
\hfill\qed

\bigskip

{\bf Construction of Example \ref{WR1VARIANT} }
In Example \ref{WR1} (or \ref{WR1Q}) change the point of view,
taking for $L$ the orbit $\barK^+$.
\hfill\qed

\bigskip

{\bf Construction of Example \ref{WRS2} }
Take $L= K$ as in Example~\ref{WRSN}, and take $G$ to be the
subgroup of $\fs{Aut}(L)$ generated by $W$ and $s^2$.  The
elements of $G$ are of the form $w(s^2)^z$, where $w \in W$ and $z
\in \Z$.  Take $\psi$ to be as in Example~\ref{WRSN}, and take
$\alpha$ to be the restriction to $G$ of the $\alpha$ of
Example~\ref{WRSN}.  Then $\alpha(G) \subseteq G$ since $\psi$
normalizes $W$ and $(s^2)^{\psi} = s^2$ (which is why we use $s^2$
here rather than $s$); and $\alpha \in \fs{Aut}(G)$ since $\psi^{-1} =
\psi$.  Pick $h_1$ and $h_2$ as in Example~\ref{WRSN}.  As before,
$\alpha$ is as required.

Applying appropriate powers of $s$ reveals that $\pair{L}{G}$ is
inclusion-transitive and
2-upward-transitive (taking $a = {\bf 0}$ in the definition of the
latter).  The rest is
clear, the requirement that $L = \Q$ being arranged as
in Example~\ref{WR1Q}.
\break
\rule{1pt}{0pt}\hfill\qed

\medskip

{\bf Construction of Example \ref{WR2} }
Take $L = K$ as in Example~\ref{WR1}, take $G$ to be
the subgroup of $\fs{Aut}(L)$ generated by $W_F$, $W_T$, and $s$,
and take $H$ to be the subgroup of $\fs{Aut}(\fs{ED}(L))$ generated
by $G$ and \fs{neg}.  $G$ is a
subgroup of the group of Example~\ref{WR1}, and the
elements of $G$ are of the form $fts^z$, where $f \in W_F$, $t \in
W_T$, and $z \in \Z$.

$\pair{L}{G}$ is transitive because $\pair{L}{W_T}$ is transitive, and
then the argument for 2-{\it o-}transitivity in Example~\ref{WR1}
gives approximate 2-{\it o-}transitivity.  (Actually, it gives a bit
more:  Given $a_1 < a_2$ and $b_1 < I_2$  in $L$, there exists $g \in
G$ such that $g(a_1) = b_1$ and $g(a_2) \in I_2$.)  However,
$\pair{L}{G}$ is not 2-{\it o-}transitive, for no $g$ fixing {\bf 0} can
map a point having finite support to a point having infinite support.

The rest of the argument parallels that for Example~\ref{WR1}.
\hfill\qed

\bigskip

{\bf Construction of Example \ref{WR2VARIANT} }
Modify Example~\ref{WR2} as Example~\ref{WR1} was modified to
get Example ~\ref{WR1Q}.
\hfill\qed

\bigskip

{\bf Construction of Example \ref{BACKANDFORTH} }
Form the direct sum $M_1 \eqdf \Q \oplus \Q\sqrt 2 \leq
\trpl{\R}{<}{+}$.  Take $L_1$ to be any coset $d+M_1$ where $d$ is
algebraic and $d \not\in M_1$.  For each $q \in \Q$,
let $M_1^q  = {\Q} + q\sqrt 2$.  $M_1$, $L_1$, and the $M_1^q$'s are
countable dense subchains of $\R$.

Take $G_1$ to be the set of all $g  \in\fs{Aut}(L_1)$ such that $g$
permutes $\setm{M_1^q}{q \in \Q}$.  Using the group $T$ of
translations by elements
of $M_1$, we see that $\pair{L_1}{G_1}$ is transitive and that $M_1$
is an orbit of $\pair{\barL_1}{G_1}$.

We claim that $\pair{L_1}{G_1}$ is highly {\it o-}transitive.  For each
pair $a_1 < \dots < a_n$ and $b_1 < \dots < b_n$ in $L_1$, use a
standard back-and-forth argument to extend the map $a_i \mapsto b_i$ to
some $g
\in \fs{Aut}(L_1)$ such that

(1) $g(a_i) = b_i$ for $i=1,\dots,n$,

(2) $g(M_1^q) = M_1^q$ for each $q \in \Q$, and

(3) $g$ has bounded support.
\newline
Then $g \in G_1$, establishing the claim.  However, on the orbit
$M_1$ of
$\pair{\barL_1}{G_1}$, $G_1$ is not 2-{\it o-}transitive because the
$M_1^q$'s are (non-convex) blocks.  By (3), $\pair{L_1}{G_1}$ is
locally moving.

Define an order-reversing permutation $h_1$ of $L_1 = d + M_1$ by
setting $h_1(d+m_1) = d-m_1$.  Take $H_1$ to be the subgroup of
$\fs{Aut}(\R)$ generated by $G_1$ and $h_1$;
$\pair{\fs{ED}(L_1)}{H_1}$ is an affiliate of $\pair{L_1}{G_1}$.

We obtain the larger $\pair{L}{G}$ as follows.  Take $L = L_1 \cup
M_1$.  Use a back-and-forth argument to produce $f \in \fs{Aut}(L)$
interchanging $L_1$ and $M_1$.  Take $G$ to be the subgroup of
$\fs{Aut}(L)$ generated by $G_1$ and $f$.  Then $\pair{G(L_1)}{G} =
\pair{L}{G}$ is not
2-{\it o-}transitive since $L_1$ and $M_1$ are blocks.
Take $H$ to be the subgroup of $\fs{Aut}(\R)$ generated by $G$ and
$h_1$; $\pair{\fs{ED}(L)}{H}$ is an affiliate of $\pair{L}{G}$.

For $\pair{L_1}{G_1}$, the translations by elements of $M_1$ form a
transitive {\it l-}subgroup of $G_1$; even better, this subgroup is
{\em totally} ordered.  The rest is clear.
\hfill\qed

\medskip

This last construction continues the previous one, and deals with the
circular/monocircular case.

\bigskip

{\bf Construction of Example \ref{CIRCLE2} }
View the real circle $\C$ of circumference $\pi$ as ${\R}/\pi$, the
set of cosets of ${\Z}\pi$.  Denote $b+{\Z}\pi$ by $\hat{b}$.
$\hat{L}_1$ and $\hat{M}_1$ are countable dense subsets of $\C$.
Since $d$ is algebraic, distinct points in $L_1 \cup M_1$ yield
distinct points in $\C$.

Each $g \in \fs{Aut}({\R})$ having period $\pi$ (that is,
$g(a+\pi)=g(a)+\pi$ for all $a \in \R$) can be transferred to $\C$,
obtaining $\hat{g} \in \fs{Aut}({\C})$ defined by
$\hat{g}(\hat{b})=\widehat{g(b)}$.  Let
\[
\hat{G}_1 = \setm{\hat{g}}{g \in G_1 \mbox{ and } g \mbox{ has
period } \pi} \leq \fs{Aut}({\C}).
\]
Since $\hat{T} \leq \hat{G}_1$, $\pair{\hat{L}_1}{\hat{G}_1}$ is
transitive and $\hat{M}_1$ is an orbit of
$\pair{\bar{\hat{L}}_1}{\hat{G}_1}$

We claim that $\pair{\hat{L}_1}{\hat{G}_1}$ is highly
{\it o-}transitive.  For $\fs{Cr}(a_1,\dots,a_n)$ and
$\fs{Cr}(b_1,\dots,b_n)$ in $\C$, a back-and-forth argument similar
to that in the previous construction produces $g \in \fs{Aut}({\C})$
such that $g(a_i)=b_i$ for each $i$ and such that $g$ fixes each
$\widehat{M_1^q}$, and when $n=1$ it can be arranged that $g$ has
bounded support.

The $h_1$ used for Example~\ref{BACKANDFORTH} has a different
kind of periodicity that also permits transfer to $\CC$:
$h_1(a+\pi)=h_1(a)-\pi$.  Take
$\hat{H}_1$ to be the subgroup of $\fs{Aut}(\fs{EO}({\C}))$
generated
by $\hat{G}_1$ and $\hat{h}_1$; $\pair{\hat{L}_1}{\hat{H}_1}$ is an
affiliate of $\pair{\hat{L}_1}{\hat{G}_1}$.

The rest parallels the previous construction.
\hfill\qed
\newpage

%

\font\goth=eufm10 scaled 1300
\newfont{\sgoth}{eufm10}

\end{document}